\documentclass[preprint,12pt]{elsarticle}

\usepackage{amssymb}
\usepackage{relsize}
\usepackage{enumerate}
\usepackage{enumitem}
\usepackage{bm}
\usepackage{amsmath}
\usepackage{color}
\usepackage{verbatim}

\newtheorem{rmk}{Remark}

\usepackage[section]{placeins}

\usepackage{xcolor}

\usepackage{natbib}

\usepackage{booktabs} % For better table rules
\usepackage{multirow}

\journal{Journal of Computational Physics}

\begin{document}

\begin{frontmatter}



\title{A higher-order three-scale computational method for efficient nonlinear thermo-mechanical coupling simulation of heterogeneous structures with multiple spatial scales}


\author[label1]{Hao Dong\corref{cor1}}\ead{donghao@mail.nwpu.edu.cn}
\author[label1]{Yanqi Wang}
\author[label1]{Jiale Linghu}
\author[label2]{Qiang Ma\corref{cor1}}\ead{maqiang809@scu.edu.cn}
\cortext[cor1]{Corresponding author.}

\address[label1]{School of Mathematics and Statistics, Xidian University, Xi'an 710071, PR China}
\address[label2]{School of Mathematics, Sichuan University, Chengdu 610041, PR China}

%%\affiliation{organization={},%Department and Organization
            %%addressline={},
            %%city={},
            %%postcode={},
            %%state={},
            %%country={}}

%% Abstract
\begin{abstract}
Classical multi-scale methods involving two spatial scales face significant challenges when simulating heterogeneous structures with complicated three-scale spatial configurations. This study proposes an innovative higher-order three-scale (HOTS) computational method, aimed at accurately and efficiently computing the transient nonlinear thermo-mechanical coupling problems of heterogeneous structures with multiple spatial scales. In these heterogeneous structures, temperature-dependent material properties have an important impact on the thermo-mechanical coupling responses, which is the particular interest in this work. At first, the detailed macro-meso-micro correlative model with higher-order correction terms is established by recursively two-scale analysis between macro-meso and meso-micro scales, which enables high-accuracy analysis of temperature-dependent nonlinear thermo-mechanical behaviors of heterogeneous structures with complicated three-scale configurations. The local error analysis mathematically illustrates the well-balanced property of HOTS computational model, endowing it with high computational accuracy. In addition, a two-stage numerical algorithm with off-line and on-line stages is proposed in order to efficiently simulate the nonlinear thermo-mechanical responses of heterogeneous structures with three-level spatial scales and accurately capture their highly oscillatory information at micro-scale. Finally, the high computational efficiency, high numerical accuracy and low computational cost of the presented higher-order three-scale computational approach are substantiated via representative numerical experiments. It can be summarized that this scalable and robust HOTS computational approach offers a reliably numerical tool for nonlinear multiphysics simulation of large-scale heterogeneous structures in real-world applications.
\end{abstract}

%%\item Research highlight 1
%%\item Research highlight 2

%% Keywords

\end{frontmatter}


\section{Introduction}
Owing to the superior physical and mechanical properties compared to conventional homogeneous materials, composite materials with hierarchical configurations enable their extensive applications in diverse fields, including aviation, aerospace, electromechanical engineering, and civil engineering, etc \cite{R1,R2,R3,R4}. Meanwhile, with the monumental progress in materials science and technology, heterogeneous structures with complicated multi-scale features have been designed and fabricated, which display complicated multi-scale configurations \cite{R5,R6,R7,R8,R9}, as displayed in Fig.\hspace{1mm}1. Particularly in the aerospace field, it is common for inhomogeneous structures to be subjected to complicated thermal and mechanical coupling environments, wherein the physical and mechanical properties of these composites exhibit significant changes with variations in temperature, resulting in prominent nonlinearities in material performances \cite{R1,R2,R3,R4}. Moreover, in superconducting field, significant temperature changes during magnetization will bring about drastic variations in specific heat and thermal conductivity of $\mathrm{MgB_2}$ superconducting bulk \cite{R10}. In such cases, employing invariable material properties of the investigated heterogeneous structures would introduce substantial inaccuracies in multi-scale simulations. Hence, it is of significant values to develop effective numerical approaches for nonlinear multi-scale computation of composite structures with complicated three-scale configurations.
\begin{figure}[!htb]
	\centering
	\begin{minipage}[c]{0.48\textwidth}
		\centering
		\includegraphics[width=45mm]{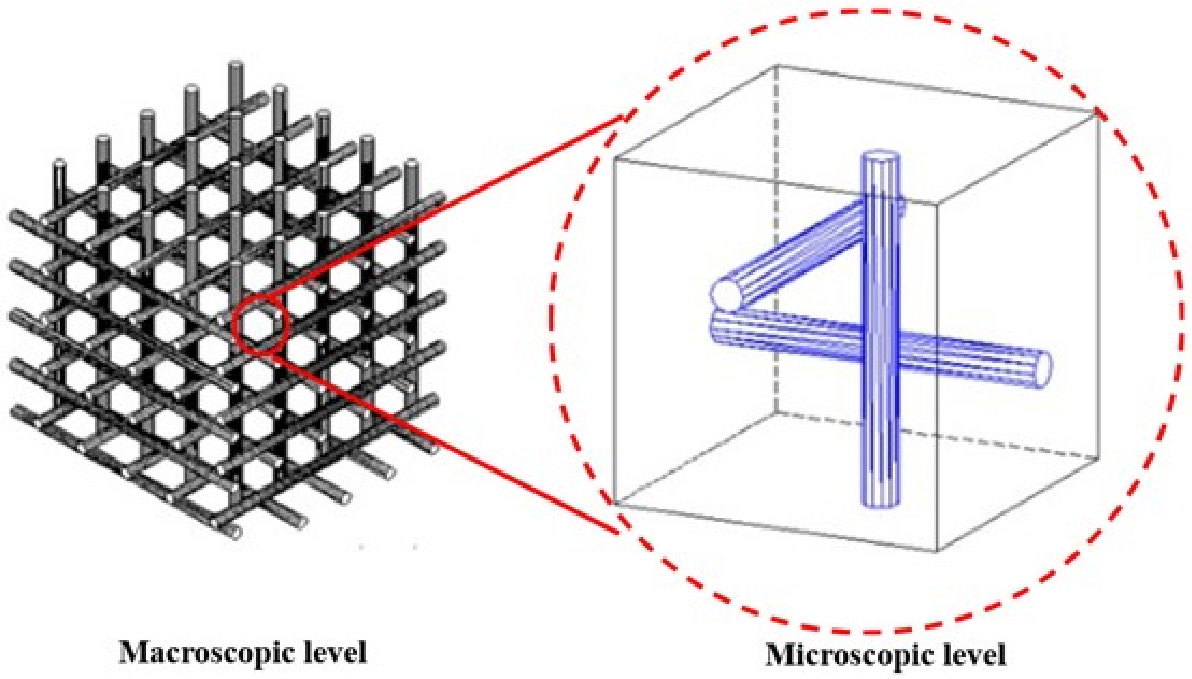} \\
		(a)
	\end{minipage}
	\begin{minipage}[c]{0.48\textwidth}
		\centering
		\includegraphics[width=80mm]{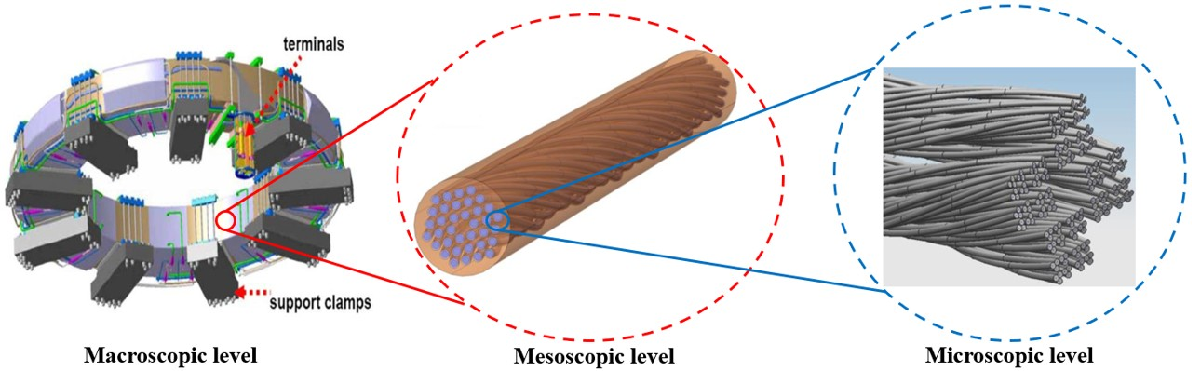} \\
		(b)
	\end{minipage}
	\caption {The illustration of inhomogeneous structures: (a) braided structure with two-scale spatial hierarchy; (b) superconducting bulk with three-scale spatial hierarchy \cite{R11}.}\label{f1}
\end{figure}

It is widely known that accurate analysis and evaluation of thermo-mechanical coupling behaviors in heterogeneous structures necessitates solving initial-boundary value problems governed by PDE with highly oscillatory and sharply discontinuous coefficients. However, conventional numerical approaches for addressing these problems face two fundamental challenges: (1) They require prohibitively fine mesh to resolve the sufficiently oscillatory information of heterogeneous structures at the smallest scale, resulting in excessive computational demands. (2) The material coefficients of composite structures with complex multi-scale features often vary significantly at different scales, posing tough challenges to the convergence and stability of classical numerical methods. Moreover, some numerical techniques, including Laplace transform, for linear system are inapplicable for theoretical analysis and numerical simulation of the nonlinear thermo-mechanical coupling system due to its temperature-dependent material properties (nonlinear equation coefficients). Thus, establishing an efficient computational approach for simulating the transient nonlinear thermo-mechanical problems of composite structures with spatially three-scale hierarchy has enormous engineering value and scientific significance. To overcome the challenging issues inherent in the multi-scale nature of composite structures, scientists and engineers have presented lots of multi-scale methods to preserve the balance between efficiency and accuracy . The asymptotic homogenization method (AHM) \cite{R12} is initially developed by Bensoussan and Oleinik et al. based on rigorous mathematical theories. The essential idea of AHM is to mathematically replace an equation with oscillatory coefficients, particularly in space, with one having homogeneous (uniform) coefficients, while physically establishing equivalence between periodically heterogeneous materials and corresponding homogeneous materials, which can effectively simplify intricate multi-scale problems and offer reliable macroscopic models. It should be mentioned that the effective properties (homogeneous coefficients) of composite materials are computed by solving a fine-scale problem over the representative volume element (RVE) or unit cell (UC) based on the homogenization strategy. Building upon the foregoing theoretical principle of AHM, researchers have developed a variety of multi-scale methods to conduct effective computation and evaluation of the multi-scale responses of composite materials with hierarchical configurations, such as the numerical homogenization method (NHM) \cite{R13}, the multi-scale finite element method (MsFEM) \cite{R14,R15,R16}, the heterogeneous multi-scale method (HMM) \cite{R17}, the variational asymptotic multi-scale method (VAMM) \cite{R18}, the local orthogonal decomposition method (LOD) \cite{R19}, the multi-scale eigenelement method (MEM) \cite{R20}, the computational homogenization method (CHM) \cite{R21}, asymptotic-preserving (AP) schemes \cite{R22} and multi-scale reduced basis method (MRBM) \cite{R23}, etc. To improve the calculation precision of multi-scale approach for engineering applications, Cui et al. established a series of second-order two-scale (SOTS) and higher-order three-scale (HOTS) methodologies. By virtue of applying the higher-order cell functions and corresponding correction terms, these approaches can precisely capture the locally physical and mechanical behaviors inside constituent materials, providing a feasible multi-scale computational framework for inhomogeneous structures \cite{R24,R25,R26,R27}. As far as we know, multi-scale approaches have been employed in limited studies addressing nonlinear multi-scale problems in composite structures. In references \cite{R28,R29}, Gaka and Chung et al. studied the stationary heat equation of periodically composite materials characterized by temperature-dependent properties and particularly emphasized on evaluating the effective thermal conductivity of layered composites characterized by temperature-dependent behaviors. Abdulle and colleagues developed an innovative reduced basis finite element heterogeneous multi-scale method (RB-FE-HMM) for a serious of nonlinear homogenization elliptic problems, which can drastically lower the computational demands of the finite element heterogeneous multi-scale method (FE-HMM) in \cite{R30,R31}. Sengupta et al. proposes an advanced two-scale formulation for fully coupled continuum thermo-mechanics employing the finite element method at macro- and micro- scales, which has been successfully used to deal with a strongly coupled problem in \cite{R32}. In reference \cite{R33}, researcher presented a multi-scale approach for strongly nonlinear monotone equations stemming from localized orthogonal decomposition. The proposed new approach gives optimal a priori error estimates up to linearization errors. In reference \cite{R34}, researchers established an efficient domain decomposition-based approach for nonlinear multi-scale problems,which is derived from manifold learning techniques and exploits the tangent spaces. Furthermore, Dong et al. proposed a higher-order two-scale computational model and associated efficient computational algorithm for nonlinear problems of composite materials, whose calculation precision is higher than the homogenization method and lower-order two-scale approach \cite{R25,R35,R36}. But aforementioned studies have predominantly focused on nonlinear thermo-mechanical problems in heterogeneous structures with two-scale spatial configurations, while remaining inapplicable to three-scale hierarchical systems. The critical limitations of two-scale methodologies manifest in two key aspects: First, direct extension of conventional two-scale approaches to composite structures with three-level spatial configuration necessitates prohibitively expensive computational efforts for solving auxiliary cell problems at meso-scale. Second, existing two-scale methods fail to adequately resolve the micro-scale oscillatory information of composite structures with three-scale spatial configurations. Currently, there remains a significant research gap in the investigation of temperature-dependent nonlinear thermo-mechanical coupling problems in composite structures exhibiting more than two-scale spatial configurations.

For addressing the deficiencies of the two-scale methods, there have been several three-scale computational approaches to efficiently simulate the material behaviors of heterogeneous structures with three-scale spatial configurations. In references \cite{R12,R37,R38,R39,R40}, the reiterated homogenization method was presented, which was used to evaluate the thermal transfer performance of composite materials with three levels of spatial heterogeneity, early-age concrete with randomly distributed inclusions, heterogeneous materials with perfect thermal contact and interfacial thermal resistance, and three-scale ordered arrays separately. Although these three-scale methods are capable of predicting the macro-scale thermal conductivity parameters of heterogeneous structures exhibiting a three-level spatial scales, they cannot fully decouple the scale coupling between the micro-scale and the meso-scale due to their improper assumptions (mesoscopic coordinate $\bm{y} = {{\bm{x}} \mathord{\left/{\vphantom {x {{\varepsilon}}}} \right.\kern-\nulldelimiterspace} {{\varepsilon}}}$ and microscopic coordinate $\bm{z} =  {{\bm{y}} \mathord{\left/{\vphantom {x {{\varepsilon}}}} \right.\kern-\nulldelimiterspace} {{\varepsilon}}}= {{\bm{x}} \mathord{\left/{\vphantom {x {{\varepsilon^2}}}} \right.\kern-\nulldelimiterspace} {{\varepsilon^2}}}$). Thus, Cui and his team pioneered a higher-order three-scale (HOTS) computational framework and developed associated multi-scale computational algorithms based on innovative assumptions (mesoscopic coordinate $\bm{y} = {{\bm{x}} \mathord{\left/{\vphantom {x {{\varepsilon _1}}}} \right.\kern-\nulldelimiterspace} {{\varepsilon _1}}}$ and microscopic coordinate $\bm{z} = {{{\varepsilon _1}{\bm{y}}} \mathord{\left/{\vphantom {{{\varepsilon_1}y} {{\varepsilon_2}}}} \right.\kern-\nulldelimiterspace} {{\varepsilon _2}}} = {{\bm{x}} \mathord{\left/{\vphantom {x {{\varepsilon_2}}}} \right.\kern-\nulldelimiterspace} {{\varepsilon _2}}}$) \cite{R41}. Building upon this novel three-scale computational framework, different higher-order three-scale approaches are developed for high-accuracy and efficient three-scale simulation of mechanical problem \cite{R42}, coupled conduction-radiation problem \cite{R43}, thermo-mechanical coupling problem \cite{R44}, hydro-mechanical coupling problem \cite{R45} and eigenvalue problem \cite{R46}, etc. The obtained results clearly indicated that the proposed higher-order three-scale methods accurately acquire the essential oscillatory behaviors of heterogeneous structures at the smallest scale, thereby satisfying the high-precision numerical solutions required for practical engineering applications. Notably, the existing higher-order three-scale methods exhibit two significant limitations: (1) They fail to account for temperature-dependent material properties in nonlinear thermo-mechanical coupling simulation of heterogeneous structures. (2) Existing methods lack rigorous theoretical analyses, including error estimation and convergence proof, when applying to nonlinear thermo-mechanical problems in composite structures with three-level spatial scales. In summary, while these nonlinear multi-scale problems have garnered increasing research attention over the past two decades owing to their profound theoretical implication and practical engineering significance, investigations specifically targeting the effective simulation and analysis of nonlinear coupled thermo-mechanical behaviors in three-scale heterogeneous structures with temperature-dependent properties remain notably limited.

This paper is dedicated to developing an efficient and high-precision three-scale computational method for investigating the time-dependent thermo-mechanical problems of composite structures with spatial three-level configurations and temperature-dependent properties. Firstly, by adopting SOTS model between the macro-scale and the meso-scale, macroscopic and mesoscopic nonlinear thermo-mechanical responses of composite structures with three-scale spatial configurations are precisely acquired. Next, SOTS model is further employed to mesoscopic cell functions for deriving lower-order and higher-order microscopic cell functions and corresponding microscopic correction terms. Through the top-down successive SOTS modeling, an innovative higher-order macro-meso-micro correlative model is designed for the nonlinear coupled thermo-mechanical computation of heterogeneous structures. Subsequently, to evaluate the numerical precision of the newly developed HOTS model, local error analysis is derived from the point-wise sense. Finally, an effectively two-stage numerical algorithm is presented to reduce the significant amount of computational resources, since the direct use of classical numerical methods for three-scale simulation needs an exceedingly fine mesh to catch the oscillating behaviors at micro-scale. In this two-stage algorithm, the mesoscopic cell functions and the microscopic cell functions are calculated via the finite element method (FEM) in the off-line phase. And then, in the on-line phase, the homogenization equations at macro-scale are computed using the hybrid finite element-finite difference method (FEM-FDM) and decoupling method whilst the higher-order three-scale solutions are computed through classical interpolation methods.

The outline of this study is introduced as below. Section 2 establishes the governing equations for nonlinear temperature-dependent thermo-mechanical coupling problems in three-scale heterogeneous structures with temperature-dependent properties. Next, we systematically establish an innovative three-scale asymptotic model capable of precisely capturing microscopic oscillatory behaviors in nonlinear multi-scale composites through successive higher-order two-scale analysis. Moreover, the local error analysis in the pointwise sense is given, mathematically illustrating the well-balanced property of HOTS computational model. Section 3 develops an efficient two-stage numerical algorithm for simulating nonlinear temperature-dependent thermo-mechanical coupling problems in three-scale heterogeneous structures. The algorithm incorporates the off-line computation of auxiliary cell functions at meso- and micro-scales, and the on-line simulation of macroscopic homogenized problems and HOTS solutions. Section 4 numerically validates both the accuracy and computational efficiency of the presented higher-order three-scale model and two-stage algorithm through several numerical experiments. The paper concludes with discussion and summary remarks in Section 5.

Throughout this study,  the Einstein summation convention is adopted to streamline expressions with repetitive indices.
\section{Higher-order three-scale asymptotic model for nonlinear thermo-mechanical problems}
\subsection{Problem setting and governing equations}
Inspired by multi-scale theory and thermo-mechanical coupling theory, the governing equations for the nonlinear thermo-mechanical problems of heterogeneous structures characterized by three-level spatial hierarchy are described as below, which possess the temperature-dependent material coefficients \cite{R3,R12,R35}.
\begin{equation}
	\left\{
	\begin{aligned}
		&{\rho ^{{\zeta _1}{\zeta _2}}}(\bm{x},{\theta ^{{\zeta _1}{\zeta _2}}}){c^{{\zeta _1}{\zeta _2}}}(\bm{x},\!{\theta ^{{\zeta _1}{\zeta  _2}}})\frac{{\partial {\theta ^{{\zeta  _1}{\zeta  _2}}}(\bm{x},\!t)}}{{\partial t}} \!-\! \frac{\partial }{{\partial {x_i}}}(k_{ij}^{{\zeta  _1}{\zeta  _2}}(\bm{x},{\theta ^{{\zeta  _1}{\zeta  _2}}})\frac{{\partial {\theta ^{{\zeta  _1}{\zeta  _2}}}(\bm{x},\!t)}}{{\partial {x_j}}})\\
		&\quad \quad \quad \quad   + \vartheta _{ij}^{{\zeta  _1}{\zeta  _2}}(\bm{x},\!{\theta ^{{\zeta  _1}{\zeta  _2}}})\frac{\partial }{{\partial t}}(\frac{{\partial u_i^{{\zeta  _1}{\zeta  _2}}(\bm{x},t)}}{{\partial {x_j}}}) = h(\bm{x},t),\;\;{\rm{in}}\;\;\Omega  \times (0,T),\\
		&{\rho ^{{\zeta  _1}{\zeta  _2}}}(\bm{x},{\theta ^{{\zeta  _1}{\zeta  _2}}})\frac{{{\partial ^2}u_i^{{\zeta  _1}{\zeta  _2}}(\bm{x},\!t)}}{{\partial {t^2}}} \!-\! \frac{\partial }{{\partial {x_j}}}[C_{ijkl}^{{\zeta  _1}{\zeta  _2}}(\bm{x},{\theta ^{{\zeta  _1}{\zeta  _2}}})\frac{{\partial u_k^{{\zeta  _1}{\zeta  _2}}(\bm{x},\!t)}}{{\partial {x_l}}}\\
		&\quad \quad \quad \quad   - \beta _{ij}^{{\zeta  _1}{\zeta  _2}}(\bm{x},{\theta ^{{\zeta  _1}{\zeta  _2}}})({\theta ^{{\zeta  _1}{\zeta  _2}}}(\bm{x},t) - \tilde \theta )] = {f_i}(\bm{x},t),\;\;{\rm{in}}\;\;\Omega  \times (0,T),\\
		&{\theta ^{{\zeta  _1}{\zeta  _2}}}(\bm{x},t) = \widehat \theta (\bm{x},t),\;\;{\rm{on}}\;\;\partial {\Omega _\theta } \times (0,T),\\
		&{{\bm {u}}^{{\zeta  _1}{\zeta  _2}}}(\bm{x},t) = \widehat {\bm{ u}}(\bm{x},t),\;\;{\rm{on}}\;\;\partial {\Omega _u} \times (0,T),\\
		&k_{ij}^{{\zeta  _1}{\zeta  _2}}(\bm{x},{\theta ^{{\zeta  _1}{\zeta  _2}}})\frac{{\partial {\theta ^{{\zeta  _1}{\zeta  _2}}}(\bm{x},t)}}{{\partial {x_j}}}{n_i} = \bar q(\bm{x},t),\;\;{\rm{on}}\;\;\partial {\Omega _q} \times (0,T),\\
		&\big[C_{ijkl}^{{\zeta  _1}{\zeta  _2}}(\bm{x},{\theta ^{{\zeta  _1}{\zeta  _2}}})\frac{{\partial u_k^{{\zeta  _1}{\zeta  _2}}(\bm{x},t)}}{{\partial {x_l}}} - \beta _{ij}^{{\zeta  _1}{\zeta  _2}}(\bm{x},{\theta ^{{\zeta  _1}{\zeta  _2}}})({\theta ^{{\zeta  _1}{\zeta  _2}}}(\bm{x},t) - \tilde \theta )\big]{n_j} \\
		&\quad \quad \quad \quad \quad \quad \quad \quad \quad \quad \quad \quad = {{\bar \sigma }_i}(\bm{x},t),\;\;{\rm{on}}\;\;\partial {\Omega _\sigma } \times (0,T),\\
		&{\theta ^{{\zeta  _1}{\zeta  _2}}}(\bm{x},0) = \tilde \theta,\;\;{{\bm{u}}^{{\zeta  _1}{\zeta  _2}}}(\bm{x},0) = {{\bm{u}}^0},\;\;\frac{{\partial {{\bm{u}}^{{\zeta  _1}{\zeta  _2}}}(\bm{x},t)}}{{\partial t}}{|_{t = 0}} = {{\bm{u}}^1}(\bm{x}),\;\;{\rm{in}}\;\;\Omega ,
	\end{aligned} \right.
\end{equation}
where $\Omega$ is a bounded convex domain in $\mathbb{R}^\mathcal{N}(\mathcal{N}=2,3)$ with a boundary $\partial\Omega = \partial {\Omega _\theta } \cup\partial {\Omega _u} \cup \partial {\Omega _q}\cup \partial {\Omega _\sigma}$. In this work, the composite structure $\Omega$ is made up of periodic mesoscopic UC ${Y}$ with characteristic size $\zeta _1$ in the meso-scale. Furthermore, mesoscopic UC $Y$ can be composed of periodic microscopic UC $Z$ with characteristic size $\zeta _2$ in the micro-scale. By introducing the scaling parameters $\zeta _1$ and $\zeta _2$, the macroscopic coordinate $\bm{x}$, mesoscopic coordinate $\bm{y}$ and microscopic coordinate $\bm{z}$ are given by $\bm{y} = {{\bm{x}} \mathord{\left/{\vphantom {x {{\zeta  _1}}}} \right.\kern-\nulldelimiterspace} {{\zeta  _1}}}$ and $\bm{z} = {{{\zeta  _1}{\bm{y}}} \mathord{\left/{\vphantom {{{\zeta _1}y} {{\zeta _2}}}} \right.\kern-\nulldelimiterspace} {{\zeta  _2}}} = {{\bm{x}} \mathord{\left/{\vphantom {x {{\zeta _2}}}} \right.\kern-\nulldelimiterspace} {{\zeta  _2}}}$ for the investigated composite structures characterized by three well-separated spatial scales, as exhibited in Fig.\hspace{1mm}1(b).

Moreover, ${\rho ^{{\zeta  _1}{\zeta  _2}}}$ is the mass density; ${c^{{\zeta  _1}{\zeta  _2}}}$ is the specific heat capacity; $\{k_{ij}^{{\zeta  _1}{\zeta  _2}}\}$ is the second-order thermal conductivity tensor; $\{\vartheta _{ij}^{{\zeta  _1}{\zeta  _2}}\}$ is the nonlinear two-way thermo-mechanical coupled effect tensor; $\{C_{ijkl}^{{\zeta  _1}{\zeta  _2}}\}$ is the fourth-order elastic tensor; $ \{\beta _{ij}^{{\zeta  _1}{\zeta  _2}}\} $ is the second-order thermal modulus tensor; $h(\bm{x},t)$ and $f_i(\bm{x},t)$ are the internal heat source and body forces, respectively; $\widehat\theta(\bm{x},t)$ is the prescribed temperature on the boundary $\partial\Omega_{\theta}$, $\widehat {\bm{u}}(\bm{x},t)$ is the prescribed displacements on the boundary $\partial {\Omega _u}$, $\bar{q}(\bm{x},t)$ is the prescribed heat flux normal to the boundary $\partial\Omega_{q}$ with the normal vector $n_i$ and $ {{\bar \sigma }_i}(\bm{x},t)$ is the prescribed traction on the boundary $\partial {\Omega _\sigma }$ with the normal vector ${n_j}$; $\tilde \theta$, ${{\bm{u}}^0}$ and ${{\bm{u}}^1}(\bm{x})$ represent the initial temperature, initial displacement and velocity conditions of composite structure $\Omega$, respectively. The research target is to solve the temperature field ${\theta ^{{\zeta  _1}{\zeta  _2}}}(\bm{x},t)$ and displacement field ${u_i^{{\zeta  _1}{\zeta  _2}}}(\bm{x},t)$ of the investigated composite structures.

Furthermore, the following assumptions for the governing equations (1) are presented according to the existing works \cite{R26,R27}.
\begin{enumerate}[label=(\textbf{\Alph*}), leftmargin=*, align=left]
\item[(A)]
Functions $k_{ij}^{{\zeta  _1}{\zeta  _2}}( {\bm{x},{\theta ^{{\zeta  _1}{\zeta  _2}}}} ),\vartheta _{ij}^{{\zeta  _1}{\zeta  _2}}(\bm{x},{\theta ^{{\zeta  _1}{\zeta  _2}}}),C_{ijkl}^{{\zeta  _1}{\zeta  _2}}(\bm{x},{\theta ^{{\zeta  _1}{\zeta  _2}}})$ and $\beta _{ij}^{{\zeta  _1}{\zeta  _2}}(\bm{x},{\theta ^{{\zeta  _1}{\zeta  _2}}})$ are all symmetric and uniformly elliptic such that
\begin{displaymath}
	\begin{aligned}
			&k_{ij}^{{\zeta  _1}{\zeta  _2}}=k_{ji}^{{\zeta  _1}{\zeta  _2}},\;\underline{\gamma}|\bm{\xi}|^2\leq {k_{ij}^{{\zeta  _1}{\zeta  _2}}}( {{\bm{x}},{\theta ^{{\zeta  _1}{\zeta  _2}}}} )\xi_i\xi_j  \le\overline{\gamma}|\bm{\xi}|^2,\\
			&\vartheta_{ij}^{{\zeta  _1}{\zeta  _2}}=\vartheta_{ji}^{{\zeta  _1}{\zeta  _2}},\;\underline{\gamma}|\bm{\xi}|^2\leq {\vartheta_{ij}^{{\zeta  _1}{\zeta  _2}}}( {\bm{x},{\theta ^{{\zeta  _1}{\zeta  _2}}}} )\xi_i\xi_j  \le\overline{\gamma}|\bm{\xi}|^2,\\
			&C_{ijkl}^{{\zeta  _1}{\zeta  _2}}=C_{ijlk}^{{\zeta  _1}{\zeta  _2}}=C_{klij}^{{\zeta  _1}{\zeta  _2}},\;\underline{\gamma}{\lambda _{ij}}{\lambda _{ij}}\leq {C_{ijkl}^{{\zeta  _1}{\zeta  _2}}}( {\bm{x},{\theta ^{{\zeta  _1}{\zeta  _2}}}} ){\lambda _{ij}}{\lambda _{kl}}  \le\overline{\gamma}{\lambda _{ij}}{\lambda _{ij}},\\
			&\beta_{ij}^{{\zeta  _1}{\zeta  _2}}=\beta_{ji}^{{\zeta  _1}{\zeta  _2}},\;\underline{\gamma}|\bm{\xi}|^2\leq {\beta_{ij}^{{\zeta  _1}{\zeta  _2}}}( {\bm{x},{\theta ^{{\zeta  _1}{\zeta  _2}}}} )\xi_i\xi_j  \le\overline{\gamma}|\bm{\xi}|^2,
	\end{aligned}
\end{displaymath}
for positive constants $\underline{\gamma}$ and $\overline{\gamma}$ independent of $\zeta _1$ and $\zeta _2$, arbitrary vector $\bm{\xi}=(\xi_1,\cdots,\xi_\mathcal{N}) \in \mathbb{R}^\mathcal{N}$, arbitrary symmetric matrix $\{\lambda_{ij}\} \in \mathbb{R}^{\mathcal{N}\times\mathcal{N}}$ and arbitrary point $({\bm{x},{\theta ^{{\zeta  _1}{\zeta  _2}}}})$  $\rm{in}$  $\Omega  \times [{\theta _{min}},{\theta _{max}} + {C_*}]$.
\item[(B)]
Material property parameters ${\rho ^{{\zeta  _1}{\zeta  _2}}}(\bm{x},{\theta ^{{\zeta  _1}{\zeta  _2}}})$, ${c^{{\zeta  _1}{\zeta  _2}}}(\bm{x},{\theta ^{{\zeta  _1}{\zeta  _2}}})$, ${k_{ij}^{{\zeta  _1}{\zeta  _2}}}( {{\bm{x}},{\theta ^{{\zeta  _1}{\zeta  _2}}}} ) $, $\vartheta _{ij}^{{\zeta  _1}{\zeta  _2}}(\bm{x},{\theta ^{{\zeta  _1}{\zeta  _2}}})$, $C_{ijkl}^{{\zeta  _1}{\zeta  _2}}(\bm{x},{\theta ^{{\zeta  _1}{\zeta  _2}}})$ and $\beta _{ij}^{{\zeta  _1}{\zeta  _2}}(\bm{x},{\theta ^{{\zeta  _1}{\zeta  _2}}})$ belong to $L^\infty (\Omega)$, and there exist four positive constants $\underline{\rho}$, $\underline{c}$, $\overline{\rho}$ and $\overline{c}$ such that
\begin{displaymath}
\underline{\rho}\leq {\rho ^{{\zeta  _1}{\zeta  _2}}}({\bm{x},{\theta ^{{\zeta  _1}{\zeta  _2}}}}) \leq \overline{\rho},\;\;\underline{c}\leq {c^{{\zeta  _1}{\zeta  _2}}}( {\bm{x},{\theta ^{{\zeta  _1}{\zeta  _2}}}}) \leq \overline{c}.
\end{displaymath}
\item[(C)]
$h(\bm{x},t)\in L^2(\Omega\times(0,T))$, $f_i(\bm{x},t)\in L^2(\Omega\times(0,T))$, $\widehat{\theta}(\bm{x},t)\in L^2(0,T;H^1(\Omega))$, $\widehat{{\bm{u}}}(\bm{x}, t)\in L^2(0,T;(H^1(\Omega))^{\mathcal{N}})$, $\bar{q}(\bm{x},t) \in L^2(\Omega\times(0,T))$, $\bar{\sigma}_i(\bm{x}, t) \in L^2(\Omega\times(0,T))$, ${\bm{u}}^1(\bm{x}) \in (L^2(\Omega))^{\mathcal{N}}$.
\end{enumerate}
\subsection{The second-order two-scale asymptotic model of governing equations}
To implement three-scale computable modeling, based on these mathematical settings, the chain rules for arbitrary function ${\Psi^{{\zeta  _1}{\zeta  _2}}}(\bm{x},t) = {\Psi ^{{\zeta  _2}}}(\bm{x},\bm{y},t) = \Psi (\bm{x},\bm{y},\bm{z},t)$ achieve as follows.
\begin{equation}
\frac{{\partial {\Psi ^{{\zeta  _1}{\zeta  _2}}}(\bm{x},t)}}{{\partial {x_i}}} = \frac{{\partial {\Psi ^{{\zeta  _2}}}(\bm{x},\bm{y},t)}}{{\partial {x_i}}} + \zeta  _1^{ - 1}\frac{{\partial {\Psi ^{{\zeta  _2}}}(\bm{x},\bm{y},t)}}{{\partial {y_i}}},
\end{equation}
\begin{equation}
\frac{{\partial {\Psi ^{{\zeta  _2}}}(\bm{x},\bm{y},t)}}{{\partial {y_i}}} = \frac{{\partial \Psi (\bm{x},\bm{y},\bm{z},t)}}{{\partial {y_i}}} +{\zeta  _1}\zeta  _2^{ - 1}\frac{{\partial \Psi (\bm{x},\bm{y},\bm{z},t)}}{{\partial {z_i}}},
\end{equation}
\begin{equation}
\begin{aligned}
\frac{{\partial {\Psi ^{{\zeta  _1}{\zeta  _2}}}(\bm{x},t)}}{{\partial {x_i}}} &= \frac{{\partial \Psi (\bm{x},\bm{y},\bm{z},t)}}{{\partial {x_i}}}\!+ \zeta  _1^{ - 1}\frac{{\partial \Psi (\bm{x},\bm{y},\bm{z},t)}}{{\partial {y_i}}} \!+\zeta  _2^{ - 1}\frac{{\partial \Psi (\bm{x},\bm{y},\bm{z},t)}}{{\partial {z_i}}},
\end{aligned}
\end{equation}
which will be widely utilized in the sequel.

Motivated by the pervious works in \cite{R46,R47}, ${\theta ^{{\zeta  _1}{\zeta  _2}}}(\bm{x},t)$ and ${{u_i}^{{\zeta  _1}{\zeta  _2}}}(\bm{x},t)$ are sought in the forms of the asymptotic expansions in macro-meso scales.
\begin{equation}
	\left\{
	\begin{aligned}
		&{\theta ^{{\zeta  _1}{\zeta  _2}}}({\bm{x}},t) = {\theta _0}({\bm{x}},{\bm{y}},t) + {\zeta  _1}{\theta _1}({\bm{x}},{\bm{y}},t) + \zeta  _1^2{\theta _2}({\bm{x}},{\bm{y}},t) + {\rm{O}}(\zeta  _1^3),\\
		&u_i^{{\zeta  _1}{\zeta  _2}}({\bm{x}},t) = {u_{i0}}({\bm{x}},{\bm{y}},t) + {\zeta  _1}{u_{i1}}({\bm{x}},{\bm{y}},t) + \zeta  _1^2{u_{i2}}({\bm{x}},{\bm{y}},t) + {\rm{O}}(\zeta  _1^3),
	\end{aligned}\right.
\end{equation}
where the terms with subscript 0 are zeroth-order expansion terms, the terms with subscript 1 are first-order
asymptotic terms (lower-order asymptotic terms), and the terms with subscript 2 are second-order asymptotic terms (higher-order asymptotic terms).

After that, employing the Taylor's formula, and multi-index
notations $f_y(x,y)={\bf{D}}^{(0,1)}f(x,y)$ and $f_{yy}(x,y)={\bf{D}}^{(0,2)}f(x,y)$ in \cite{R25,R35,R47}, temperature-dependent three-scale function ${\Phi^{{\zeta  _1}{\zeta  _2}}}(\bm{x},{\theta ^{{\zeta  _1}{\zeta  _2}}})$ has the asymptotic expansion formulation $\displaystyle{\Phi^{{\zeta  _1}{\zeta  _2}}}(\bm{x},{\theta ^{{\zeta  _1}{\zeta  _2}}}) = {\Phi^{{\zeta  _2}}}(\bm{y},{\theta ^{{\zeta  _1}{\zeta  _2}}}) = {\Phi^{{\zeta  _2}}}(\bm{y},{\theta _0} + {\zeta  _1}{\theta _1} + \zeta  _1^2{\theta _2} + {\rm{O}}(\zeta  _1^3))= {\Phi^{{\zeta  _2}}}(\bm{y},{\theta _0}) + {\zeta  _1}{\theta _1}{{\bf{D}}^{(0,1)}}{\Phi^{{\zeta  _2}}}(\bm{y},{\theta _0})+ \zeta  _1^2[{\theta _2}{{\bf{D}}^{(0,1)}}{\Phi^{{\zeta  _2}}}(\bm{y},{\theta _0}) + \frac{1}{2}\theta _1^2{{\bf{D}}^{(0,2)}}{\Phi ^{{\zeta  _2}}}(\bm{y},{\theta _0})] + {\rm{O}}(\zeta  _1^3)= {\Phi^{(0)}}(\bm{y},{\theta _0}) + {\zeta  _1}{\Phi^{(1)}}(\bm{x},\bm{y},{\theta _0}) + \zeta  _1^2{\Phi^{(2)}}(\bm{x},\bm{y},{\theta _0}) + {\rm{O}}(\zeta  _1^3)$. Then, material property functions ${\rho ^{{\zeta_1}{\zeta_2}}}$, ${c^{{\zeta_1}{\zeta_2}}}$, ${k_{ij}^{{\zeta  _1}{\zeta _2}}}$, $\vartheta _{ij}^{{\zeta_1}{\zeta_2}}$, $C_{ijkl}^{{\zeta_1}{\zeta_2}}$ and $\beta _{ij}^{{\zeta_1}{\zeta_2}}$ shall be expanded as
\begin{equation}
	\begin{aligned}
		{\rho ^{{\zeta  _1}{\zeta  _2}}}(\bm{x},{\theta ^{{\zeta  _1}{\zeta  _2}}}) = {\rho ^{(0)}}(\bm{y},{\theta _0}) + {\zeta  _1}{\rho ^{(1)}}(\bm{x},\bm{y},{\theta _0}) + \zeta  _1^2{\rho ^{(2)}}(\bm{x},\bm{y},{\theta _0}) + {\rm{O}}(\zeta  _1^3).
	\end{aligned}
\end{equation}
\begin{equation}
	\begin{aligned}
		{c^{{\zeta  _1}{\zeta  _2}}}(\bm{x},{\theta ^{{\zeta  _1}{\zeta  _2}}}) = {c^{(0)}}(\bm{y},{\theta _0}) + {\zeta  _1}{c^{(1)}}(\bm{x},\bm{y},{\theta _0}) + \zeta  _1^2{c^{(2)}}(\bm{x},\bm{y},{\theta _0}) + {\rm{O}}(\zeta  _1^3).
	\end{aligned}
\end{equation}
\begin{equation}
	\begin{aligned}
		k_{ij}^{{\zeta  _1}{\zeta  _2}}(\bm{x},{\theta ^{{\zeta  _1}{\zeta  _2}}}) = k_{ij}^{(0)}(\bm{y},{\theta _0}) + {\zeta  _1}k_{ij}^{(1)}(\bm{x},\bm{y},{\theta _0}) + \zeta  _1^2k_{ij}^{(2)}(\bm{x},\bm{y},{\theta _0}) + {\rm{O}}(\zeta  _1^3).
	\end{aligned}
\end{equation}
\begin{equation}
	\begin{aligned}
		\vartheta _{ij}^{{\zeta  _1}{\zeta  _2}}({\bm{x}},{\theta ^{{\zeta  _1}{\zeta  _2}}})= \vartheta _{ij}^{(0)}({\bm{y}},{\theta _0}) + {\zeta  _1}\vartheta _{ij}^{(1)}({\bm{x}},{\bm{y}},{\theta _0}) + \zeta  _1^2\vartheta _{ij}^{(2)}({\bm{x}},{\bm{y}},{\theta _0}) + {\rm{O}}(\zeta  _1^3).
	\end{aligned}
\end{equation}
\begin{equation}
	\begin{aligned}
		C_{ijkl}^{{\zeta  _1}{\zeta  _2}}({\bm{x}},{\theta ^{{\zeta  _1}{\zeta  _2}}})= C_{ijkl}^{(0)}({\bm{y}},\!{\theta _0})\! + \!{\zeta  _1}C_{ijkl}^{(1)}({\bm{x}},\!{\bm{y}},\!{\theta _0}) \!+\!\zeta  _1^2C_{ijkl}^{(2)}({\bm{x}},\!{\bm{y}},\!{\theta _0})\! + \!{\rm{O}}(\zeta  _1^3).
	\end{aligned}
\end{equation}
\begin{equation}
	\begin{aligned}
		\beta _{ij}^{{\zeta  _1}{\zeta  _2}}({\bm{x}},{\theta ^{{\zeta  _1}{\zeta  _2}}}) = \beta _{ij}^{(0)}({\bm{y}},{\theta _0})\! + \!{\zeta  _1}\beta _{ij}^{(1)}({\bm{x}},{\bm{y}},{\theta _0}) \!+\! \zeta  _1^2\beta _{ij}^{(2)}({\bm{x}},{\bm{y}},{\theta _0}) \!+ \!{\rm{O}}(\zeta  _1^3).
	\end{aligned}
\end{equation}
Then expanding all spatial derivatives by the macro-meso chain rule (2) and substituting (5)-(11) into multi-scale nonlinear thermo-mechanical equations (1), a series of equations can be systematically derived via collecting the power-like terms of $\zeta_1$
\begin{equation}
	{\rm O}({\zeta _1^{-2}}): \left\{
	\begin{aligned}
		&\frac{\partial }{{\partial {y_i}}}\Bigl(k_{ij}^{(0)}\frac{{\partial {\theta _0}}}{{\partial {y_j}}}\Bigl) = 0,\\
		&\frac{\partial }{{\partial {y_j}}}\Bigl(C_{ijkl}^{(0)}\frac{{\partial {u_{k0}}}}{{\partial {y_l}}}\Bigl)= 0.
	\end{aligned}\right.
\end{equation}
\begin{equation}
	{\rm{O}}(\zeta _1^{-1}):
	\left\{
	\begin{aligned}
		&\frac{\partial }{{\partial {x_i}}}\Bigl(k_{ij}^{(0)}\frac{{\partial {\theta _0}}}{{\partial {y_j}}}\Bigl) + \frac{\partial }{{\partial {y_i}}}\Bigl(k_{ij}^{(0)}(\frac{{\partial {\theta _0}}}{{\partial {x_j}}} + \frac{{\partial {\theta _1}}}{{\partial {y_j}}})\Bigl)\\
		&+ \frac{\partial }{{\partial {y_i}}}\Bigl(k_{ij}^{(1)}\frac{{\partial {\theta _0}}}{{\partial {y_j}}}\Bigl) - \vartheta _{ij}^{(0)}\frac{\partial }{{\partial t}}(\frac{{\partial {u_{i0}}}}{{\partial {y_j}}}) = 0,\\
		&\frac{\partial }{{\partial {x_j}}}\Bigl(C_{ijkl}^{(0)}\frac{{\partial {u_{k0}}}}{{\partial {y_l}}}\Bigl) + \frac{\partial }{{\partial {y_j}}}\Bigl(C_{ijkl}^{(0)}(\frac{{\partial {u_{k0}}}}{{\partial {x_l}}} + \frac{{\partial {u_{k1}}}}{{\partial {y_l}}})\Bigl)\\
		&+ \frac{\partial }{{\partial {y_j}}}\Bigl(C_{ijkl}^{(1)}\frac{{\partial {u_{k0}}}}{{\partial {y_l}}}\Bigl) - \frac{\partial }{{\partial {y_j}}}\Bigl(\beta _{ij}^{(0)}({\theta _0} - \tilde \theta )\Bigl) = 0.
	\end{aligned} \right.
\end{equation}
\begin{equation}
	{\rm{O}}(\zeta _1^{0}):\left\{
	\begin{aligned}
		&{\rho ^{(0)}}{c^{(0)}}\frac{{\partial {\theta _0}}}{{\partial t}} = \frac{\partial }{{\partial {x_i}}}\Bigl(k_{ij}^{(1)}\frac{{\partial {\theta _0}}}{{\partial {y_j}}}\Bigl) + \frac{\partial }{{\partial {y_i}}}\Bigl(k_{ij}^{(2)}\frac{{\partial {\theta _0}}}{{\partial {y_j}}}\Bigl)\\
		&+ \frac{\partial }{{\partial {x_i}}}\Bigl(k_{ij}^{(0)}(\frac{{\partial {\theta _0}}}{{\partial {x_j}}} + \frac{{\partial {\theta _1}}}{{\partial {y_j}}})\Bigl) + \frac{\partial }{{\partial {y_i}}}\Bigl(k_{ij}^{(0)}(\frac{{\partial {\theta _1}}}{{\partial {x_j}}} + \frac{{\partial {\theta _2}}}{{\partial {y_j}}})\Bigl) \\
		&+ \frac{\partial }{{\partial {y_i}}}\Bigl(k_{ij}^{(1)}(\frac{{\partial {\theta _0}}}{{\partial {x_j}}} + \frac{{\partial {\theta _1}}}{{\partial {y_j}}})\Bigl) - \vartheta _{ij}^{(1)}\frac{\partial }{{\partial t}}(\frac{{\partial {u_{i0}}}}{{\partial {y_j}}})\\
		&- \vartheta _{ij}^{(0)}\frac{\partial }{{\partial t}}(\frac{{\partial {u_{i0}}}}{{\partial {x_j}}}) - \vartheta _{ij}^{(0)}\frac{\partial }{{\partial t}}(\frac{{\partial {u_{i1}}}}{{\partial {y_j}}}) + h,\\
		&{\rho ^{(0)}}\frac{{{\partial ^2}{u_{i0}}}}{{\partial {t^2}}} = \frac{\partial }{{\partial {y_j}}}\Bigl(C_{ijkl}^{(2)}\frac{{\partial {u_{k0}}}}{{\partial {y_l}}}\Bigl) + \frac{\partial }{{\partial {x_j}}}\Bigl(C_{ijkl}^{(1)}\frac{{\partial {u_{k0}}}}{{\partial {y_l}}}\Bigl) \\
		&+ \frac{\partial }{{\partial {x_j}}}\Bigl(C_{ijkl}^{(0)}(\frac{{\partial {u_{k0}}}}{{\partial {x_l}}} + \frac{{\partial {u_{k1}}}}{{\partial {y_l}}})\Bigl) + \frac{\partial }{{\partial {y_j}}}\Bigl(C_{ijkl}^{(0)}(\frac{{\partial {u_{k1}}}}{{\partial {x_l}}} + \frac{{\partial {u_{k2}}}}{{\partial {y_l}}})\Bigl) \\
		& + \frac{\partial }{{\partial {y_j}}}\Bigl(C_{ijkl}^{(1)}(\frac{{\partial {u_{k0}}}}{{\partial {x_l}}} + \frac{{\partial {u_{k1}}}}{{\partial {y_l}}})\Bigl)- \frac{\partial }{{\partial {x_j}}}\Bigl(\beta _{ij}^{(0)}({\theta _0} - \tilde \theta )\Bigl) \\
		& - \frac{\partial }{{\partial {y_j}}}(\beta _{ij}^{(0)}{\theta _1}) - \frac{\partial }{{\partial {y_j}}}\Bigl(\beta _{ij}^{(1)}({\theta _0} - \tilde \theta )\Bigl) + {f_i}.
	\end{aligned} \right.
\end{equation}

Analogous to existing developments in second-order two-scale nonlinear models \cite{R25,R35}, from ${\rm O}({\zeta _1^{-2}})$-order equations (12), we can derive that
\begin{equation}
	{\theta _0}({\bm{x}},{\bm{y}},t) = {\theta _0}({\bm{x}},t),{u_{i0}}({\bm{x}},{\bm{y}},t) = {u_{i0}}({\bm{x}},t).
\end{equation}
Then using the above conclusion (15) and ${\rm{O}}(\zeta _1^{-1})$-order equations (13), the macro-meso decoupling expressions for $\theta _1$ and $u_{i1}$ are presented as below
\begin{equation}
	\left\{
	\begin{aligned}
		&{\theta _1}({\bm{x}},{\bm{y}},t) = {M_{{\alpha _1}}}({\bm{y}},{\theta _0})\frac{{\partial {\theta _0}}}{{\partial {x_{{\alpha _1}}}}},\\
		&{u_{i1}}({\bm{x}},{\bm{y}},t) = N_{im}^{{\alpha _1}}({\bm{y}},{\theta _0})\frac{{\partial {u_{m0}}}}{{\partial {x_{{\alpha _1}}}}} - {P_i}({\bm{y}},{\theta _0})({\theta _0} - \tilde \theta ),\;\;m,n = 1,2,3,
	\end{aligned} \right.
\end{equation}
where ${M_{\alpha_1}}$, $N_{im}^{{\alpha _1}}$ and $P_i$ are the first-order auxiliary cell functions defined on mesoscopic unit cell $Y^{\zeta_2}$. By the way, the mesoscopic unit cell $Y^{\zeta_2}$ differs significantly from the aforementioned mesoscopic unit cell $Y$ due to the fact that mesoscopic unit cell $Y^{\zeta_2}$ has complicatedly configurations at micro-scale. Furthermore, from ${\rm{O}}(\zeta _1^{-1})$-order equations (13), the governing equations for first-order auxiliary cell functions ${M_{\alpha_1}}$, $N_{im}^{{\alpha _1}}$ and $P_i$ are established via imposing homogeneous Dirichlet boundary condition, described next.
\begin{equation}
	\left\{
	\begin{aligned}
		&\frac{\partial }{{\partial {y_i}}}[k_{ij}^{(0)}({\bm{y},{\theta _0}})\frac{{\partial {M_{{\alpha _1}}}}}{{\partial {y_j}}}] = - \frac{{\partial k_{i{\alpha _1}}^{(0)}( {\bm{y},{\theta _0}})}}{{\partial {y_i}}},&\;\;\;&\bm{y}\in Y^{\zeta  _2}, \\
		&{M_{{\alpha _1}}}(\bm{y},\theta _0) = 0,&\;\;\;&\bm{y}\in\partial Y^{\zeta  _2}.
	\end{aligned} \right.
\end{equation}
\begin{equation}
	\left\{
	\begin{aligned}
		&\frac{\partial }{{\partial {y_j}}}[C_{ijkl}^{(0)}({\bm{y}},{\theta _0})\frac{{\partial N_{km}^{{\alpha _1}}}}{{\partial {y_l}}}] =  - \frac{{\partial C_{ijm{\alpha _1}}^{(0)}({\bm{y}},{\theta _0})}}{{\partial {y_j}}},&\;\;\;&\bm{y}\in Y^{\zeta  _2}, \\
		&N_{km}^{{\alpha _1}}({\bm{y}},{\theta _0}) = 0,&\;\;\;&\bm{y}\in\partial Y^{\zeta  _2}.
	\end{aligned} \right.
\end{equation}
\begin{equation}
	\left\{
	\begin{aligned}
		&\frac{\partial }{{\partial {y_j}}}[C_{ijkl}^{(0)}({\bm{y}},{\theta _0})\frac{{\partial {P_k}}}{{\partial {y_l}}}] =  - \frac{{\partial \beta _{ij}^{(0)}({\bm{y}},{\theta _0})}}{{\partial {y_j}}},&\;\;\;&\bm{y}\in Y^{\zeta  _2}, \\
		&{P_k}({\bm{y}},{\theta _0}) = 0,&\;\;\;&\bm{y}\in\partial Y^{\zeta  _2}.
	\end{aligned} \right.
\end{equation}
\begin{rmk}
Unlike those in linear periodic composites, above-mentioned first-order mesoscopic cell functions with quasi-periodic features exhibit a significant dependence on the macroscopic temperature $\theta _0$.
\end{rmk}

Next, integrating both sides of ${\rm O}({\zeta _1^{0}})$-order equations (14) on mesoscopic unit cell $Y^{\zeta  _2}$ and removing the $\bm{y}$-variable, we can derive the following macroscopic homogenized equations corresponding to nonlinear multi-scale problem (1)
\begin{equation}
	\left\{
	\begin{aligned}
		&\hat S({\theta _0})\frac{{\partial {\theta _0}}}{{\partial t}} - \frac{\partial }{{\partial {x_i}}}({{\hat k}_{ij}}({\theta _0})\frac{{\partial {\theta _0}}}{{\partial {x_j}}}) + {{\hat \vartheta }_{ij}}({\theta _0})\frac{\partial }{{\partial t}}(\frac{{\partial {u_{i0}}}}{{\partial {x_j}}}) = h,\;{\rm{in}}\;\Omega  \times (0,T),\\
		&\hat \rho ({\theta _0})\frac{{{\partial ^2}{u_{i0}}}}{{\partial {t^2}}}\! - \!\frac{\partial }{{\partial {x_j}}}\big[{{\hat C}_{ijkl}}({\theta _0})\frac{{\partial {u_{k0}}}}{{\partial {x_l}}}\! -\! {{\hat \beta }_{ij}}({\theta _0})({\theta _0}\! - \!\tilde \theta )\big]\! = \!{f_i},\;{\rm{in}}\;\Omega  \times (0,T),\\
		&{\theta _0}({\bm{x}},t) = \hat \theta ({\bm{x}},t),\;\;{\rm{on}}\;\;\partial {\Omega _\theta } \times (0,T),\\
		&{{\bm{u}}_{0}}({\bm{x}},t) = \hat {\bm{u}}({\bm{x}},t),\;\;{\rm{on}}\;\;\partial {\Omega _u} \times (0,T),\\
		&{{\hat k}_{ij}}({\theta _0})\frac{{\partial {\theta _0}({\bm{x}},t)}}{{\partial {x_j}}}{n_i} = \bar q({\bm{x}},t),\;\;{\rm{on}}\;\;\partial {\Omega _q} \times (0,T),\\
		&\big[{{\hat C}_{ijkl}}({\theta _0})\frac{{\partial {u_{k0}}({\bm{x}},t)}}{{\partial {x_l}}} - \hat \beta_{ij} ({\theta _0})({\theta _0}({\bm{x}},t) - \tilde \theta )\big]{n_j}\\
		&\quad \quad \quad \quad \quad \quad \quad \quad \quad \quad \quad \quad  = {{\bar \sigma }_i}(\bm{x},t),\;\;{\rm{on}}\;\;\partial {\Omega _\sigma } \times (0,T),\\
		&{\theta _{0}}(\bm{x},0) = \tilde \theta,\;\;{{\bm{u}}_{0}}(\bm{x},0) = {{\bm{u}}^0},\;\;\frac{{\partial {{\bm{u}}_{0}}(\bm{x},t)}}{{\partial t}}{|_{t = 0}} = {{\bm{u}}^1}(\bm{x}),\;\;{\rm{in}}\;\;\Omega,
	\end{aligned}\right.
\end{equation}
in which macroscopic homogenized material parameters can be derived respectively via the following formulas
\begin{equation}
	\begin{aligned}
		&\hat S({\theta _0}) = \frac{1}{{|{Y^{{\zeta  _2}}}|}}\int_{{Y^{{\zeta  _2}}}} {\Bigl({\rho ^{(0)}}({\bm{y}},{\theta _0}){c^{(0)}}({\bm{y}},{\theta _0}) - \vartheta _{ij}^{(0)}({\bm{y}},{\theta _0})\frac{{\partial {P_i}}}{{\partial {y_j}}}\Bigl)d{Y^{{\zeta  _2}}}} ,\\
		&{{\hat k}_{ij}}({\theta _0}) = \frac{1}{{|{Y^{{\zeta  _2}}}|}}\int_{{Y^{{\zeta  _2}}}} {\Bigl(k_{ij}^{(0)}({\bm{y}},{\theta _0}) + k_{ik}^{(0)}({\bm{y}},{\theta _0})\frac{{\partial {M_j}}}{{\partial {y_k}}}\Bigl)d{Y^{{\zeta  _2}}}} ,\\
		&{{\hat \vartheta }_{ij}}({\theta _0}) = \frac{1}{{|{Y^{{\zeta  _2}}}|}}\int_{{Y^{{\zeta  _2}}}} \Bigl(\vartheta _{ij}^{(0)}(\bm{y},{\theta _0}) + \vartheta _{mn}^{(0)}({\bm{y}},{\theta _0})\frac{{\partial N_{mi}^j}}{{\partial {y_n}}}\Bigl)d{Y^{{\zeta  _2}}},\\
		&\hat \rho ({\theta _0}) = \frac{1}{{|{Y^{{\zeta  _2}}}|}}\int_{{Y^{{\zeta  _2}}}}{{\rho ^{(0)}}} ({\bm{y}},{\theta _0})d{Y^{{\zeta  _2}}},\\
		&{{\hat C}_{ijkl}}({\theta _0}) = \frac{1}{{|{Y^{{\zeta  _2}}}|}}\int_{{Y^{{\zeta  _2}}}} {\Bigl(C_{ijkl}^{(0)}({\bm{y}},{\theta _0}) + C_{ijmn}^{(0)}({\bm{y}},{\theta _0})\frac{{\partial N_{mk}^l}}{{\partial {y_n}}}\Bigl)d{Y^{{\zeta  _2}}}} ,\\
		&{{\hat \beta }_{ij}}({\theta _0}) = \frac{1}{{|{Y^{{\zeta  _2}}}|}}\int_{{Y^{{\zeta  _2}}}} {\Bigl(\beta _{ij}^{(0)}({\bm{y}},{\theta _0}) + C_{ijkl}^{(0)}({\bm{y}},{\theta _0})\frac{{\partial {P_k}}}{{\partial {y_l}}}\Bigl)d{Y^{{\zeta  _2}}}} .
	\end{aligned}
\end{equation}
\begin{rmk}
Owing to the quasi-periodic features of the first-order mesoscopic cell functions, all macroscopic homogenized material parameters exhibit a dependence on the macroscopic temperature solution $\theta_0$.
\end{rmk}

By eliminating the terms $h$ and $f_i$ involving in ${\rm O}({\zeta _1^{0}})$-order equations (14) from the macroscopic homogenized equations (20), we can establish the concrete decoupling expressions for $\theta_2$ and $u_{i2}$ as below
\begin{equation}
	\left\{
	\begin{aligned}
		{\theta _2}({\bm{x}},{\bm{y}},t) &= A({\bm{y}},{\theta _0})\frac{{\partial {\theta _0}}}{{\partial t}} \!+\! {M_{{\alpha _1}{\alpha _2}}}({\bm{y}},{\theta _0})\frac{{{\partial ^2}{\theta _0}}}{{\partial {x_{{\alpha _1}}}\partial {x_{{\alpha _2}}}}}\! +\! {C_{{\alpha _1}}}({\bm{y}},{\theta _0})\frac{{\partial {\theta _0}}}{{\partial {x_{{\alpha _1}}}}}\\
		&- {B_{{\alpha _1}{\alpha _2}}}({\bm{y}},{\theta _0})\frac{{\partial {\theta _0}}}{{\partial {x_{{\alpha _1}}}}}\frac{{\partial {\theta _0}}}{{\partial {x_{{\alpha _2}}}}} - {E_{{\alpha _1}{\alpha _2}}}({\bm{y}},{\theta _0})\frac{{{\partial ^2}{u_{{\alpha _1}0}}}}{{\partial {x_{{\alpha _2}}}\partial t}},\\
		{u_{i2}}({\bm{x}},{\bm{y}},t) &= F_i^{{\alpha _1}}({\bm{y}},{\theta _0})\frac{{{\partial ^2}{u_{{\alpha _1}0}}}}{{\partial {t^2}}} + N_{im}^{{\alpha _1}{\alpha _2}}({\bm{y}},{\theta _0})\frac{{{\partial ^2}{u_{m0}}}}{{\partial {x_{{\alpha _1}}}\partial {x_{{\alpha _2}}}}} \\
		&+ Z_{im}^{{\alpha _1}}({\bm{y}},{\theta _0})\frac{{\partial {u_{m0}}}}{{\partial {x_{{\alpha _1}}}}} - {Q_i}({\bm{y}},{\theta _0})({\theta _0} - \tilde \theta ) - H_i^{{\alpha _1}}({\bm{y}},{\theta _0})\frac{{\partial {\theta _0}}}{{\partial {x_{{\alpha _1}}}}} \\
		&+ W_i^{{\alpha _1}}({\bm{y}},{\theta _0})\frac{{\partial {\theta _0}}}{{\partial {x_{{\alpha _1}}}}}({\theta _0} - \tilde \theta ) - J_{im}^{{\alpha _1}{\alpha _2}}({\bm{y}},{\theta _0})\frac{{\partial {\theta _0}}}{{\partial {x_{{\alpha _1}}}}}\frac{{\partial {u_{m0}}}}{{\partial {x_{{\alpha _2}}}}},
	\end{aligned}\right.
\end{equation}
where $A$, $M_{\alpha_1\alpha_2}$, $C_{\alpha_1}$, ${B_{\alpha_1\alpha_2}}$, ${E_{\alpha_1\alpha_2}}$, $F_i^{\alpha _1}$, $N_{im}^{{\alpha _1}{\alpha _2}}$, $ Z_{im}^{{\alpha _1}}$, ${Q_i}$, $H_i^{{\alpha _1}}$, $W_i^{{\alpha _1}}$ and $J_{im}^{{\alpha _1}{\alpha _2}}$ are the second-order auxiliary cell functions defined on mesoscopic unit cell $Y^{\zeta  _2}$. After that, imposing homogeneous Dirichlet boundary conditions, we can derive the following class of equations for the second-order auxiliary cell functions.
\begin{equation}
	\left\{
	\begin{aligned}
		&\frac{\partial}{\partial y_i}[ { k_{ij}^{(0)}({\bm{y}},\theta _0){\frac{\partial A}{\partial y_j}}} ] =  {\rho ^{(0)}}{c^{(0)}} - \hat S - \vartheta _{ij}^{(0)}\frac{{\partial {P_i}}}{{\partial {y_j}}},&\;\;\;&\bm{y}\in Y^{\zeta  _2},\\
		&A({\bm{y}},\theta _0)=0,&\;\;\;&\bm{y}\in\partial Y^{\zeta  _2}.
	\end{aligned} \right.	
\end{equation}
\begin{equation}
	\left\{
	\begin{aligned}
		&\frac{\partial }{{\partial {y_i}}}[k_{ij}^{(0)}({\bm{y}},{\theta _0})\frac{{\partial {M_{{\alpha _1}{\alpha _2}}}}}{{\partial {y_j}}}] = {{\hat k}_{{\alpha _1}{\alpha _2}}} - k_{{\alpha _1}{\alpha _2}}^{(0)} - \frac{\partial }{{\partial {y_i}}}(k_{i{\alpha _1}}^{(0)}{M_{{\alpha _2}}}) - k_{{\alpha _1}j}^{(0)}\frac{{\partial {M_{{\alpha _2}}}}}{{\partial {y_j}}},&\;\;\;&\bm{y}\in Y^{\zeta  _2},\\
		&{M_{{\alpha _1}{\alpha _2}}}(\bm{y},\theta _0) = 0,&\;\;\;&\bm{y}\in\partial Y^{\zeta  _2}.
	\end{aligned} \right.
\end{equation}
\begin{equation}
	\left\{
	\begin{aligned}
		&\frac{\partial }{{\partial {y_i}}}[{k_{ij}^{(0)}(\bm{y},\theta _0)\frac{{\partial {C_{\alpha _1}}}}{{\partial {y_j}}}}] = {\frac{{\partial {{\hat k}_{i{\alpha _1}}}}}{{\partial {x_i}}} - \frac{{\partial k_{i{\alpha _1}}^{(0)}}}{{\partial {x_i}}} - \frac{\partial }{{\partial {y_i}}}(k_{ij}^{(0)}\frac{{\partial {M_{{\alpha _1}}}}}{{\partial {x_j}}}) - \frac{\partial }{{\partial {x_i}}}(k_{ij}^{(0)}\frac{{\partial {M_{{\alpha _1}}}}}{{\partial {y_j}}})},&\;\;\;&\bm{y}\in Y^{\zeta  _2},\\
		&{C_{{\alpha _1}}}(\bm{y},\theta _0) = 0,&\;\;\;&\bm{y}\in\partial Y^{\zeta  _2}.
	\end{aligned} \right.
\end{equation}
\begin{equation}
	\left\{
	\begin{aligned}
		&\frac{\partial }{{\partial {y_i}}}[k_{ij}^{(0)}(\bm{y},\theta _0)\frac{{\partial {B_{{\alpha _1}{\alpha _2}}}}}{{\partial {y_j}}}]=  \frac{\partial }{{\partial {y_i}}}\Big({M_{{\alpha _1}}}{\mathbf{D}^{({0,1})}}k_{i{\alpha _2}}^{(0)} + {M_{{\alpha _1}}}{\mathbf{D}^{({0,1})}}k_{ij}^{(0)}\frac{{\partial {M_{{\alpha _2}}}}}{{\partial {y_j}}}\Big),&\;\;\;&\bm{y}\in Y^{\zeta  _2},\\
		&{B_{{\alpha _1}{\alpha _2}}}(\bm{y},\theta _0) = 0,&\;\;\;&\bm{y}\in\partial Y^{\zeta  _2}.
	\end{aligned} \right.
\end{equation}
\begin{equation}
	\left\{
	\begin{aligned}
		&\frac{\partial }{{\partial {y_i}}}[k_{ij}^{(0)}({\bm{y}},{\theta _0})\frac{{\partial {E_{{\alpha _1}{\alpha _2}}}}}{{\partial {y_j}}}] = {{\hat \vartheta }_{{\alpha _1}{\alpha _2}}} - \vartheta _{{\alpha _1}{\alpha _2}}^{(0)} - \vartheta _{ij}^{(0)}\frac{{\partial N_{i{\alpha _1}}^{{\alpha _2}}}}{{\partial {y_j}}},&\;\;\;&\bm{y}\in Y^{\zeta  _2},\\
		&{E_{{\alpha _1}{\alpha _2}}}({\bm{y}},{\theta _0}) = 0,&\;\;\;&\bm{y}\in\partial Y^{\zeta  _2}.
	\end{aligned} \right.
\end{equation}
\begin{equation}
	\left\{
	\begin{aligned}
		&\frac{\partial }{{\partial {y_j}}}[C_{ijkl}^{(0)}({\bm{y}},{\theta _0})\frac{{\partial F_k^{{\alpha _1}}}}{{\partial {y_l}}}] = {\rho ^{(0)}} - \hat \rho ,&\;\;\;&\bm{y}\in Y^{\zeta  _2},\\
		&F_k^{{\alpha _1}}({\bm{y}},{\theta _0}) = 0,&\;\;\;&\bm{y}\in\partial Y^{\zeta  _2}.
	\end{aligned} \right.
\end{equation}
\begin{equation}
	\left\{
	\begin{aligned}
		&\frac{\partial }{{\partial {y_j}}}[C_{ijkl}^{(0)}({\bm{y}},{\theta _0})\frac{{\partial N_{km}^{{\alpha _1}{\alpha _2}}}}{{\partial {y_l}}}] = {{\hat C}_{i{\alpha _1}m{\alpha _2}}} - C_{i{\alpha _1}m{\alpha _2}}^{(0)} \\
        &\quad\quad\quad\quad\quad- \frac{\partial }{{\partial {y_j}}}(C_{ijk{\alpha _1}}^{(0)}N_{km}^{{\alpha _2}}) - C_{i{\alpha _1}kj}^{(0)}\frac{{\partial N_{km}^{{\alpha _2}}}}{{\partial {y_j}}},&\;\;\;&\bm{y}\in Y^{\zeta  _2},\\
		&N_{km}^{{\alpha _1}{\alpha _2}}({\bm{y}},{\theta _0}) = 0,&\;\;\;&\bm{y}\in\partial Y^{\zeta  _2}.
	\end{aligned}\right.
\end{equation}
\begin{equation}
	\left\{
	\begin{aligned}
		&\frac{\partial }{{\partial {y_j}}}[C_{ijkl}^{(0)}({\bm{y}},{\theta _0})\frac{{\partial Z_{km}^{{\alpha _1}}}}{{\partial {y_l}}}] = \frac{{\partial {{\hat C}_{ijm{\alpha _1}}}}}{{\partial {x_j}}} - \frac{{\partial C_{ijm{\alpha _1}}^{(0)}}}{{\partial {x_j}}} \\
        &\quad\quad\quad- \frac{\partial }{{\partial {y_j}}}(C_{ijkl}^{(0)}\frac{{\partial N_{km}^{{\alpha _1}}}}{{\partial {x_l}}}) - \frac{\partial }{{\partial {x_j}}}(C_{ijkl}^{(0)}\frac{{\partial N_{km}^{{\alpha _1}}}}{{\partial {y_l}}}),&\;\;\;&\bm{y}\in Y^{\zeta  _2},\\
		&Z_{km}^{{\alpha _1}}({\bm{y}},{\theta _0}) = 0,&\;\;\;&\bm{y}\in\partial Y^{\zeta  _2}.
	\end{aligned} \right.
\end{equation}
\begin{equation}
	\left\{
	\begin{aligned}
		&\frac{\partial }{{\partial {y_j}}}[C_{ijkl}^{(0)}({\bm{y}},{\theta _0})\frac{{\partial {Q_k}}}{{\partial {y_l}}}] = \frac{{\partial {{\hat \beta }_{ij}}}}{{\partial {x_j}}} - \frac{{\partial \beta _{ij}^{(0)}}}{{\partial {x_j}}}\\
        &\quad\quad\quad - \frac{\partial }{{\partial {y_j}}}(C_{ijkl}^{(0)}\frac{{\partial {P_k}}}{{\partial {x_l}}}) - \frac{\partial }{{\partial {x_j}}}(C_{ijkl}^{(0)}\frac{{\partial {P_k}}}{{\partial {y_l}}}),&\;\;\;&\bm{y}\in Y^{\zeta  _2},\\
		&{Q_k}({\bm{y}},{\theta _0}) = 0,&\;\;\;&\bm{y}\in\partial Y^{\zeta  _2}.
	\end{aligned}\right.
\end{equation}
\begin{equation}
	\left\{
	\begin{aligned}
		&\frac{\partial }{{\partial {y_j}}}[C_{ijkl}^{(0)}({\bm{y}},{\theta _0})\frac{{\partial H_k^{{\alpha _1}}}}{{\partial {y_l}}}] = {{\hat \beta }_{i{\alpha _1}}} - \beta _{i{\alpha _1}}^{(0)}- \frac{\partial }{{\partial {y_j}}}(C_{ijk{\alpha _1}}^{(0)}{P_k}) \\
        &\quad\quad\quad\quad- C_{i{\alpha _1}kl}^{(0)}\frac{{\partial {P_k}}}{{\partial {y_l}}} - \frac{\partial }{{\partial {y_j}}}(\beta _{ij}^{(0)}{M_{{\alpha _1}}}),&\;&\bm{y}\in Y^{\zeta  _2},\\
		&H_k^{{\alpha _1}}({\bm{y}},{\theta _0}) = 0,&\;&\bm{y}\in\partial Y^{\zeta  _2}.
	\end{aligned} \right.
\end{equation}
\begin{equation}
	\left\{
	\begin{aligned}
		&\frac{\partial }{{\partial {y_j}}}[C_{ijkl}^{(0)}({\bm{y}},{\theta _0})\frac{{\partial W_k^{{\alpha _1}}}}{{\partial {y_l}}}] = \frac{\partial }{{\partial {y_j}}}\Bigl({M_{{\alpha _1}}}{{\bf{D}}^{(0,1)}}\beta _{ij}^{(0)} + {M_{{\alpha _1}}}{{\bf{D}}^{(0,1)}}C_{ijkl}^{(0)}\frac{{\partial {P_k}}}{{\partial {y_l}}}\Bigl),&\;\;\;&\bm{y}\in Y^{\zeta  _2},\\
		&W_k^{{\alpha _1}}({\bm{y}},{\theta _0}) = 0,&\;\;\;&\bm{y}\in\partial Y^{\zeta  _2}.
	\end{aligned} \right.
\end{equation}
\begin{equation}
	\left\{
	\begin{aligned}
		&\frac{\partial }{{\partial {y_j}}}[C_{ijkl}^{(0)}({\bm{y}},{\theta _0})\frac{{\partial J_{km}^{{\alpha _1}{\alpha _2}}}}{{\partial {y_l}}}] = \frac{\partial }{{\partial {y_j}}}\Bigl({M_{{\alpha _1}}}{{\bf{D}}^{(0,1)}}C_{ijm{\alpha _2}}^{(0)} + {M_{{\alpha _1}}}{{\bf{D}}^{(0,1)}}C_{ijkl}^{(0)}\frac{{\partial N_{km}^{{\alpha _2}}}}{{\partial {y_l}}}\Bigl),&\;\;\;&\bm{y}\in Y^{\zeta  _2},\\
		&J_{km}^{{\alpha _1}{\alpha _2}}({\bm{y}},{\theta _0}) = 0,&\;\;\;&\bm{y}\in\partial Y^{\zeta  _2}.
	\end{aligned}  \right.
\end{equation}
To wrap up, we have successfully established the macro-meso correlative second-order two-scale asymptotic solutions for the multi-scale nonlinear thermo-mechanical problems (1), which are constructed without meso-micro decoupling.
\begin{equation}
	\begin{aligned}
		{\theta ^{{\zeta  _1}{\zeta  _2}}}({\bm{x}},t) &= {\theta _0}({\bm{x}},t) + {\zeta  _1}{M_{{\alpha _1}}}({\bm{y}},{\theta _0})\frac{{\partial {\theta _0}}}{{\partial {x_{{\alpha _1}}}}} + \zeta  _1^2[A({\bm{y}},{\theta _0})\frac{{\partial {\theta _0}}}{{\partial t}} \\
		&+ {M_{{\alpha _1}{\alpha _2}}}({\bm{y}},{\theta _0})\frac{{{\partial ^2}{\theta _0}}}{{\partial {x_{{\alpha _1}}}\partial {x_{{\alpha _2}}}}} + {C_{{\alpha _1}}}({\bm{y}},{\theta _0})\frac{{\partial {\theta _0}}}{{\partial {x_{{\alpha _1}}}}} \\
		&- {B_{{\alpha _1}{\alpha _2}}}({\bm{y}},{\theta _0})\frac{{\partial {\theta _0}}}{{\partial {x_{{\alpha _1}}}}}\frac{{\partial {\theta _0}}}{{\partial {x_{{\alpha _2}}}}} - {E_{{\alpha _1}{\alpha _2}}}({\bm{y}},{\theta _0})\frac{{{\partial ^2}{u_{{\alpha _1}0}}}}{{\partial {x_{{\alpha _2}}}\partial t}}]+{\rm O}(\zeta  _1^3).
	\end{aligned}
\end{equation}
\begin{equation}
	\begin{aligned}
		u_i^{{\zeta  _1}{\zeta  _2}}({\bm{x}},t) &= {u_{i0}}({\bm{x}},t) + {\zeta  _1}[N_{im}^{{\alpha _1}}({\bm{y}},{\theta _0})\frac{{\partial {u_{m0}}}}{{\partial {x_{{\alpha _1}}}}} - {P_i}({\bm{y}},{\theta _0})({\theta _0} - \tilde \theta )]\\
		&+ \zeta  _1^2[F_i^{{\alpha _1}}({\bm{y}},{\theta _0})\frac{{{\partial ^2}{u_{{\alpha _1}0}}}}{{\partial {t^2}}} + N_{im}^{{\alpha _1}{\alpha _2}}({\bm{y}},{\theta _0})\frac{{{\partial ^2}{u_{m0}}}}{{\partial {x_{{\alpha _1}}}\partial {x_{{\alpha _2}}}}} \\
		&+ Z_{im}^{{\alpha _1}}({\bm{y}},{\theta _0})\frac{{\partial {u_{m0}}}}{{\partial {x_{{\alpha _1}}}}}- {Q_i}({\bm{y}},{\theta _0})({\theta _0} - \tilde \theta ) - H_i^{{\alpha _1}}({\bm{y}},{\theta _0})\frac{{\partial {\theta _0}}}{{\partial {x_{{\alpha _1}}}}}\\
		&+\! W_i^{{\alpha _1}}({\bm{y}},{\theta _0})\frac{{\partial {\theta _0}}}{{\partial {x_{{\alpha _1}}}}}({\theta _0} - \tilde \theta ) \!-\! J_{im}^{{\alpha _1}{\alpha _2}}({\bm{y}},{\theta _0})\frac{{\partial {\theta _0}}}{{\partial {x_{{\alpha _1}}}}}\frac{{\partial {u_{m0}}}}{{\partial {x_{{\alpha _2}}}}}]\!+\!{\rm O}(\zeta  _1^3) .
	\end{aligned}
\end{equation}
However, there is a key limitation that aforementioned second-order two-scale model without meso-micro decoupling is inadequate for analyzing nonlinear thermo-mechanical problems in three-scale composite structures, primarily because the mesoscopic unit cells in these structures still possess intricate microscopic configurations and numerically solving the corresponding mesoscopic cell problems, which involve rapidly oscillatory thermal conductivity $k_{ij}^{(0)}(\bm{y},\theta_{0})$ and stiffness coefficient $C_{ijkl}^{(0)}(\bm{y},\theta_{0})$, still demands prohibitive computational resources. Furthermore, the pivotal idea of macro-meso two-scale modeling is illustrated in Fig.\hspace{1mm}2.

\begin{figure}[!htb]
	\centering
	\begin{minipage}[c]{1.0\textwidth}
		\centering
		\includegraphics[
		width=0.8\linewidth
		]{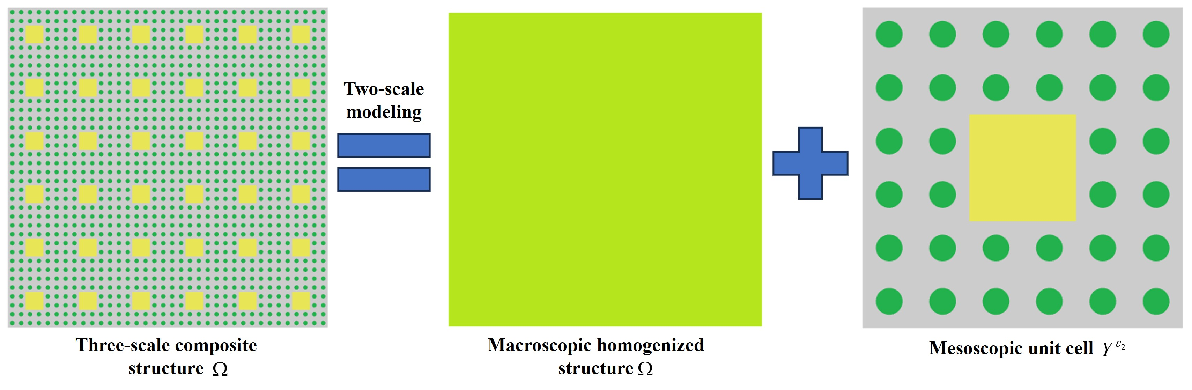}
	\end{minipage}
	\caption {The illustration of macro-meso two-scale modeling.}\label{f1}
\end{figure}

\subsection{The second-order two-scale asymptotic model of mesoscopic cell problems defined on $Y^{\zeta_2}$}
Given the multiple spatial scales of the composite structure, its mesoscopic unit cell $Y^{\zeta_2}$ has underlying microscopic configurations. Hence, we continuously conduct second-order two-scale analysis on the mesoscopic cell problems (17)-(19) and (23)-(34) in the subsequent subsection. Moreover, the core idea of meso-micro two-scale modeling is shown in Fig.\hspace{1mm}3.
\begin{figure}[!htb]
	\centering
	\begin{minipage}[c]{1.0\textwidth}
		\centering
		\includegraphics[
		width=0.8\linewidth
		]{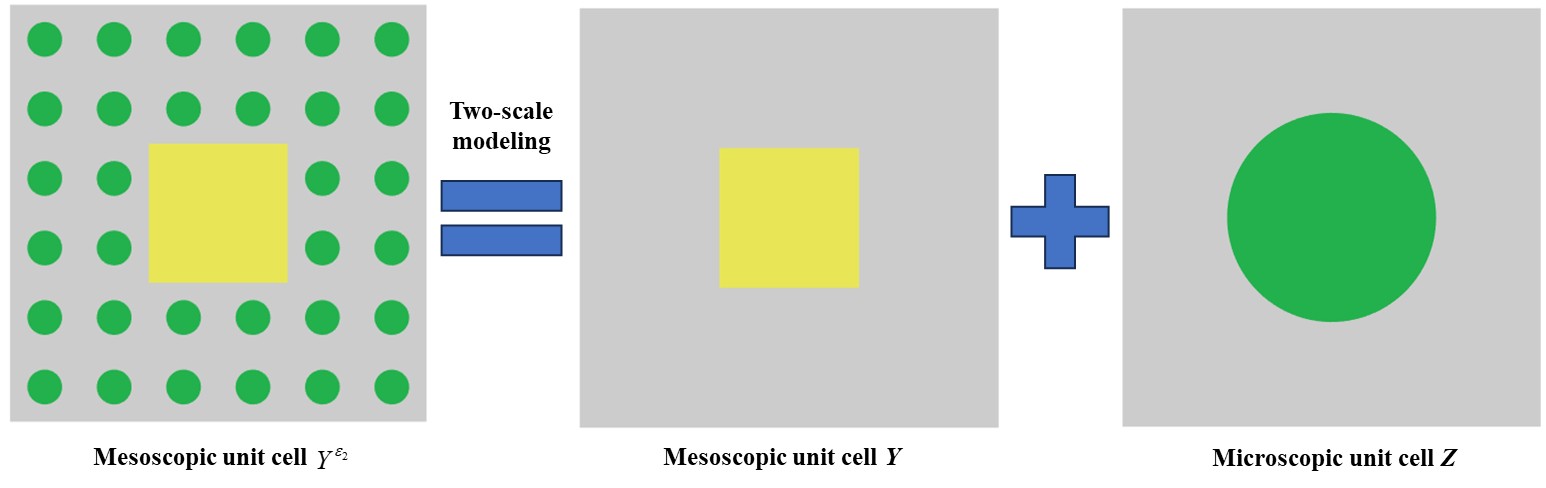}
	\end{minipage}
	\caption {The illustration of meso-micro two-scale modeling.}\label{f1}
\end{figure}

Using mesoscopic cell problem (18) as a representative case, we assume that the succeeding asymptotic expansion form are derived for first-order mesoscopic cell function $N_{im}^{{\alpha _1}}({\bm{y}},{\theta _0})$ in meso-micro scales.
\begin{equation}
	\begin{aligned}
	N_{im}^{{\alpha _1}}({\bm{y}},{\theta _0}) = N_{im}^{0{\alpha _1}}({\bm{y}},{\bm{z}},{\theta _0}) + \zeta  _1^{ - 1}{\zeta  _2}N_{im}^{1{\alpha _1}}({\bm{y}},{\bm{z}},{\theta _0}) + \zeta  _1^{ - 2}\zeta  _2^2N_{im}^{2{\alpha _1}}({\bm{y}},{\bm{z}},{\theta _0}) + {\rm{O}}(\zeta  _1^{ - 3}\zeta  _2^3).
	\end{aligned}
\end{equation}
Before proceeding, it is remarked that material parameters have the equivalent forms as below
\begin{equation}
	\left\{
	\begin{aligned}
&{k_{ij}^{{\zeta  _2}}}(\bm{y},{\theta _0})={k_{ij}^{(0)}}(\bm{y},{\theta _0})=k_{ij}(\bm{y},\bm{z},{\theta _0}),\\
&{C_{ijkl}^{{\zeta  _2}}}(\bm{y},{\theta _0})={C_{ijkl}^{(0)}}(\bm{y},{\theta _0})=C_{ijkl}(\bm{y},\bm{z},{\theta _0}).
	\end{aligned} \right.
\end{equation}
Next, expanding all spatial derivatives through the macro-meso chain rule (3) in mesoscopic cell problem (18), we hence derive a sequence of equations via balancing the power-like terms of $\displaystyle{\zeta  _1^{ - 1}{\zeta  _2}}$
\begin{equation}
	{\rm{O}}({(\zeta  _1^{ - 1}{\zeta  _2})^{ - 2}}):\frac{\partial }{{\partial {z_j}}}\Bigl(C_{ijkl}(\bm{y},\bm{z},{\theta _0})\frac{{\partial N_{km}^{0{\alpha _1}}}}{{\partial {z_l}}}\Bigl) = 0.
\end{equation}
\begin{equation}
	{\rm{O}}({(\zeta  _1^{ - 1}{\zeta  _2})^{ - 1}}):
	\begin{aligned}		
		&\frac{\partial }{{\partial {y_j}}}\Bigl(C_{ijkl}(\bm{y},\bm{z},{\theta _0})\frac{{\partial N_{km}^{0{\alpha _1}}}}{{\partial {z_l}}}\Bigl) \!+\! \frac{\partial }{{\partial {z_j}}}\Bigl(C_{ijkl}(\bm{y},\bm{z},{\theta _0})\frac{{\partial N_{km}^{0{\alpha _1}}}}{{\partial {y_l}}}\Bigl) \\
		&+ \frac{\partial }{{\partial {z_j}}}\Bigl(C_{ijkl}(\bm{y},\bm{z},{\theta _0})\frac{{\partial N_{km}^{1{\alpha _1}}}}{{\partial {z_l}}}\Bigl) =  - \frac{{\partial C_{ijm{\alpha _1}}(\bm{y},\bm{z},{\theta _0})}}{{\partial {z_j}}}.
	\end{aligned}
\end{equation}
\begin{equation}
	{\rm{O}}({(\zeta  _1^{ - 1}{\zeta  _2})^0}):
	\begin{aligned}
		&\frac{\partial }{{\partial {y_j}}}\Bigl(C_{ijkl}(\bm{y},\bm{z},{\theta _0})\frac{{\partial N_{km}^{0{\alpha _1}}}}{{\partial {y_l}}}\Bigl) + \frac{\partial }{{\partial {y_j}}}\Bigl(C_{ijkl}(\bm{y},\bm{z},{\theta _0})\frac{{\partial N_{km}^{1{\alpha _1}}}}{{\partial {z_l}}}\Bigl)\\
		&+\! \frac{\partial }{{\partial {z_j}}}\Bigl(C_{ijkl}(\bm{y},\!\bm{z},\!{\theta _0})\frac{{\partial N_{km}^{1{\alpha _1}}}}{{\partial {y_l}}}\Bigl)\!+ \frac{\partial }{{\partial {z_j}}}\Bigl(C_{ijkl}(\bm{y},\!\bm{z},\!{\theta _0})\frac{{\partial N_{km}^{2{\alpha _1}}}}{{\partial {z_l}}}\Bigl) \\
		&=  - \frac{{\partial C_{ijm{\alpha _1}}(\bm{y},\bm{z},{\theta _0})}}{{\partial {y_j}}}.
	\end{aligned}
\end{equation}

From ${\rm{O}}\displaystyle({(\zeta  _1^{ - 1}{\zeta  _2})^{ - 2}})$-order equation (39), one can determine that
\begin{equation}
N_{im}^{0{\alpha _1}}({\bm{y}},{\bm{z}},{\theta _0}) = N_{im}^{0{\alpha _1}}({\bm{y}},{\theta _0}).
\end{equation}
After that, inserting (42) into (40), equation (40) can be further simplified as follows
\begin{equation}
	\begin{aligned}
\frac{\partial }{{\partial {z_j}}}\Bigl(C_{ijkl}(\bm{y},\bm{z},{\theta _0})\frac{{\partial N_{km}^{1{\alpha _1}}}}{{\partial {z_l}}}\Bigl) =  - \frac{\partial }{{\partial {z_j}}}\Bigl(C_{ijkl}(\bm{y},\bm{z},{\theta _0})\frac{{\partial N_{km}^{0{\alpha _1}}}}{{\partial {y_l}}}\Bigl) - \frac{{\partial C_{ijm{\alpha _1}}(\bm{y},\bm{z},{\theta _0})}}{{\partial {z_j}}}.
	\end{aligned}
\end{equation}
According to (43), we can construct the following expression for mesoscopic cell function ${N_{im}^{1{\alpha _1}}}$ as below
\begin{equation}
	\begin{aligned}
N_{im}^{1{\alpha _1}}({\bm{y}},{\bm{z}},{\theta _0}) = T_{ip}^n({\bm{y}},{\bm{z}},{\theta _0})\frac{{\partial N_{pm}^{0{\alpha _1}}}}{{\partial {y_n}}} + T_{im}^{{\alpha _1}}({\bm{y}},{\bm{z}},{\theta _0}),
	\end{aligned}
\end{equation}
where $T_{im}^{{\alpha _1}}$ is denoted as first-order auxiliary cell function defined in microscale. Then, putting (44) into (43), the governing equation for $T_{im}^{{\alpha _1}}$ can be established by means of imposing homogeneous Dirichlet boundary condition as below
\begin{equation}
	\left\{
	\begin{aligned}
		&\frac{\partial }{{\partial {z_j}}}[C_{ijkl}({\bm{y}},{\bm{z}},{\theta _0})\frac{{\partial T_{km}^{{\alpha _1}}}}{{\partial {z_l}}}] =  - \frac{{\partial C_{ijm{\alpha _1}}({\bm{y}},{\bm{z}},{\theta _0})}}{{\partial {z_j}}},&\;\;\;&\bm{z}\in Z,\\
		&T_{km}^{{\alpha _1}}({\bm{y}},{\bm{z}},{\theta _0}) = 0,&\;\;\;&\bm{z}\in\partial Z.
	\end{aligned} \right.
\end{equation}
Then integrating both sides of ${\rm{O}}({(\zeta  _1^{ - 1}{\zeta  _2})^0})$-order equation (41) on microscopic unit cell $Z$ and removing the $\bm{z}$-variable, we can obtain the following mesoscopic homogenized equations corresponding to mesoscopic cell problem (18)
\begin{equation}
	\left\{
	\begin{aligned}
		&\frac{\partial }{{\partial {y_j}}}\Bigl({{\bar C}_{ijkl}}({\bm{y}},{\theta _0})\frac{{\partial N_{km}^{0{\alpha _1}}}}{{\partial {y_l}}}\Bigl) =  - \frac{{\partial {{\bar C}_{ijm{\alpha _1}}}({\bm{y}},{\theta _0})}}{{\partial {y_j}}},&\;\;\;&\bm{y}\in Y,\\
		&N_{km}^{0{\alpha _1}}({\bm{y}},{\theta _0}) = 0,&\;\;\;&\bm{y}\in\partial Y,
	\end{aligned} \right.
\end{equation}
where mesoscopic homogenized material parameter ${\bar C_{ijkl}}({\bm{y}},{\theta _0})$ is defined as follows
\begin{equation}
	\begin{aligned}
		{{\bar C}_{ijkl}}({\bm{y}},{\theta _0}) = \frac{1}{{|Z|}}\int_Z \Bigl( C_{ijkl}^{(0)}({\bm{y}},{\bm{z}},{\theta _0}) + C_{ijpn}^{(0)}({\bm{y}},{\bm{z}},{\theta _0})\frac{{\partial T_{pk}^l}}{{\partial {z_n}}}\Bigl)dZ.
	\end{aligned}
\end{equation}
It should be noted that we introduce a new set of macroscopic homogenized material parameters $\hat S^*({\theta _0})$, ${{\hat k}_{ij}}^*({\theta _0})$, ${{\hat \vartheta }_{ij}}^*({\theta _0})$, $\hat \rho ^* ({\theta _0})$, ${{\hat C}_{ijkl}}^*({\theta _0})$ and ${{\hat \beta }_{ij}}^*({\theta _0}) $ via the up-scaling procedure. These novel parameters are distinct from those in (21) and their derivation is thoroughly analyzed in Appendix A.

After that, we can obtain the decoupling expression of essential second-order correction term $N_{im}^{2{\alpha _1}}$. By inserting (42) and (44) into (41), and subtracting (41) from (46), the equation is derived equivalently
\begin{equation}
	\begin{aligned}
		\frac{\partial }{{\partial {z_j}}}\Bigl(C_{ijkl}\frac{{\partial N_{km}^{2{\alpha _1}}}}{{\partial {z_l}}}\Bigl) &= \Big[{{\bar C}_{iqpn}} - C_{iqpn} - C_{iqkj}\frac{{\partial T_{kp}^n}}{{\partial {z_j}}} - \frac{\partial }{{\partial {z_j}}}(C_{ijkq}T_{kp}^n)\Big]\frac{{{\partial ^2}N_{pm}^{0{\alpha _1}}}}{{\partial {y_q}\partial {y_n}}}\\
		&+ \Big[\frac{{\partial {{\bar C}_{ijpn}}}}{{\partial {y_j}}} - \frac{{\partial C_{ijpn}}}{{\partial {y_j}}} - \frac{\partial }{{\partial {y_j}}}(C_{ijkl}\frac{{\partial T_{kp}^n}}{{\partial {z_l}}}) - \frac{\partial }{{\partial {z_j}}}(C_{ijkl}\frac{{\partial T_{kp}^n}}{{\partial {y_l}}})\Big]\frac{{\partial N_{pm}^{0{\alpha _1}}}}{{\partial {y_n}}}\\
		&+ \Big[\frac{{\partial {{\bar C}_{ijm{\alpha _1}}}}}{{\partial {y_j}}} - \frac{{\partial C_{ijm{\alpha _1}}}}{{\partial {y_j}}} - \frac{\partial }{{\partial {y_j}}}(C_{ijkl}\frac{{\partial T_{km}^{{\alpha _1}}}}{{\partial {z_l}}}) - \frac{\partial }{{\partial {z_j}}}(C_{ijkl}\frac{{\partial T_{km}^{{\alpha _1}}}}{{\partial {y_l}}})\Big].
	\end{aligned}
\end{equation}
On account of (48), one can conclude that $N_{im}^{2{\alpha _1}}$ has the detailed decoupling expression as below
\begin{equation}
	\begin{aligned}
	N_{im}^{2{\alpha _1}}({\bm{y}},{\bm{z}},{\theta _0})\!=\!T_{ip}^{qn}({\bm{y}},\!{\bm{z}},\!{\theta _0})\frac{{{\partial ^2}N_{pm}^{0{\alpha _1}}}}{{\partial {y_q}\partial {y_n}}} \!+ \! L_{ip}^n({\bm{y}},\!{\bm{z}},\!{\theta _0})\frac{{\partial N_{pm}^{0{\alpha _1}}}}{{\partial {y_n}}}\! + \! L_{im}^{{\alpha _1}}({\bm{y}},\!{\bm{z}},\!{\theta _0}),
	\end{aligned}
\end{equation}
where $T_{im}^{{\alpha _1}{\alpha _2}}$ and $L_{im}^{{\alpha _1}}$ are second-order auxiliary cell functions defined on microscopic UC $Z$. Then substituting (49) into (48), a serious of equations for second-order auxiliary cell functions $T_{im}^{{\alpha _1}{\alpha _2}}$ and $L_{im}^{{\alpha _1}}$ can be established by attaching homogeneous Dirichlet boundary condition as below
\begin{equation}
	\left\{
	\begin{aligned}
		&\frac{\partial }{{\partial {z_j}}}[C_{ijkl}({\bm{y}},{\bm{z}},{\theta _0})\frac{{\partial T_{km}^{{\alpha _1}{\alpha _2}}}}{{\partial {z_l}}}] = {{\bar C}_{i{\alpha _1}m{\alpha _2}}} - C_{i{\alpha _1}m{\alpha _2}} - C_{i{\alpha _1}kj}\frac{{\partial T_{km}^{{\alpha _2}}}}{{\partial {z_j}}} - \frac{\partial }{{\partial {z_j}}}(C_{ijk{\alpha _1}}T_{km}^{{\alpha _2}}),&\;\;\;&\bm{z}\in Z,\\
		&T_{km}^{{\alpha _1}{\alpha _2}}({\bm{y}},{\bm{z}},{\theta _0}) = 0,&\;\;\;&\bm{z}\in\partial Z.
	\end{aligned} \right.
\end{equation}
\begin{equation}
	\left\{
	\begin{aligned}
		&\frac{\partial }{{\partial {z_j}}}[C_{ijkl}({\bm{y}},{\bm{z}},{\theta _0})\frac{{\partial L_{km}^{{\alpha _1}}}}{{\partial {z_l}}}] = \frac{{\partial {{\bar C}_{ijm{\alpha _1}}}}}{{\partial {y_j}}} - \frac{{\partial C_{ijm{\alpha _1}}}}{{\partial {y_j}}} - \frac{\partial }{{\partial {y_j}}}(C_{ijkl}\frac{{\partial T_{km}^{{\alpha _1}}}}{{\partial {z_l}}})\\
 &\quad\quad\quad\quad- \frac{\partial }{{\partial {z_j}}}(C_{ijkl}\frac{{\partial T_{km}^{{\alpha _1}}}}{{\partial {y_l}}}),&\;\;\;&\bm{z}\in Z,\\
		&L_{km}^{{\alpha _1}}({\bm{y}},{\bm{z}},{\theta _0}) = 0,&\;\;\;&\bm{z}\in\partial Z.
	\end{aligned}  \right.
\end{equation}
In conclusion, we can finally obtain the meso-micro decoupling expressions for $N_{im}^{{\alpha _1}}$ as follows
\begin{equation}
	\begin{aligned}
		N_{im}^{{\alpha _1}}({\bm{y}},{\theta _0}) &= N_{im}^{0{\alpha _1}}({\bm{y}},{\theta _0}) + \zeta  _1^{ - 1}{\zeta  _2}[T_{ip}^n({\bm{y}},{\bm{z}},{\theta _0})\frac{{\partial N_{pm}^{0{\alpha _1}}}}{{\partial {y_n}}} + T_{im}^{{\alpha _1}}({\bm{y}},{\bm{z}},{\theta _0})]\\
		&+ \zeta  _1^{ - 2}\zeta  _2^2[T_{ip}^{qn}({\bm{y}},{\bm{z}},{\theta _0})\frac{{{\partial ^2}N_{pm}^{0{\alpha _1}}}}{{\partial {y_q}\partial {y_n}}} + L_{ip}^n({\bm{y}},{\bm{z}},{\theta _0})\frac{{\partial N_{pm}^{0{\alpha _1}}}}{{\partial {y_n}}} + L_{im}^{{\alpha _1}}({\bm{y}},{\bm{z}},{\theta _0})] + {\rm{O}}(\zeta  _1^{ - 3}\zeta  _2^3).
	\end{aligned}
\end{equation}

By means of the preceding second-order two-scale analysis for $N_{im}^{{\alpha _1}}$, we can derive the following meso-micro decoupling forms for the remaining mesoscopic cell functions ${M_{{\alpha _1}}}$, ${P_i}$, $A$, $M_{\alpha_1\alpha_2}$, $C_{\alpha_1}$, ${B_{\alpha_1\alpha_2}}$, ${E_{\alpha_1\alpha_2}}$, $F_i^{\alpha _1}$, $N_{im}^{{\alpha _1}{\alpha _2}}$, $ Z_{im}^{{\alpha _1}}$, ${Q_i}$, $H_i^{{\alpha _1}}$, $W_i^{{\alpha _1}}$ and $J_{im}^{{\alpha _1}{\alpha _2}}$.
\begin{equation}
	\begin{aligned}
		{M_{{\alpha _1}}}({\bm{y}},{\theta _0}) &= M_{{\alpha _1}}^0({\bm{y}},{\theta _0}) + {\zeta  _1^{ - 1}{\zeta  _2}}[{R_k}({\bm{y}},{\bm{z}},{\theta _0})\frac{{\partial M_{{\alpha _1}}^0}}{{\partial {y_k}}} + {R_{{\alpha _1}}}({\bm{y}},{\bm{z}},{\theta _0})]\\
		& + \zeta  _1^{ - 2}\zeta  _2^2\Big[{R_{ij}}({\bm{y}},{\bm{z}},{\theta _0})\frac{{{\partial ^2}M_{{\alpha _1}}^0}}{{\partial {y_i}\partial {y_j}}} + {D_k}({\bm{y}},{\bm{z}},{\theta _0})\frac{{\partial M_{{\alpha _1}}^0}}{{\partial {y_k}}} + {D_{{\alpha _1}}}({\bm{y}},{\bm{z}},{\theta _0})\Big] + {\rm{O}}(\zeta  _1^{ - 3}\zeta  _2^3).		
	\end{aligned}
\end{equation}
\begin{equation}
	\begin{aligned}
		{P_i}({\bm{y}},{\theta _0}) &= P_i^0({\bm{y}},{\theta _0}) + \zeta  _1^{ - 1}{\zeta  _2}[T_{im}^{{\alpha _1}}({\bm{y}},{\bm{z}},{\theta _0})\frac{{\partial P_m^0}}{{\partial {y_{{\alpha _1}}}}} + {O_i}({\bm{y}},{\bm{z}},{\theta _0})]\\
		&+ \zeta  _1^{ - 2}\zeta  _2^2[T_{im}^{{\alpha _1}{\alpha _2}}({\bm{y}},{\bm{z}},{\theta _0})\frac{{{\partial ^2}P_m^0}}{{\partial {y_{{\alpha _1}}}\partial {y_{{\alpha _2}}}}} + L_{im}^{{\alpha _1}}({\bm{y}},{\bm{z}},{\theta _0})\frac{{\partial P_m^0}}{{\partial {y_{{\alpha _1}}}}} + P_i^ * ({\bm{y}},{\bm{z}},{\theta _0})] + {\rm{O}}(\zeta  _1^{ - 3}\zeta  _2^3).
	\end{aligned}
\end{equation}
\begin{equation}
	\begin{aligned}
		A({\bm{y}},{\theta _0}) &= {A^0}({\bm{y}},{\theta _0}) + \zeta  _1^{ - 1}{\zeta  _2}{R_k}({\bm{y}},{\bm{z}},{\theta _0})\frac{{\partial {A^0}}}{{\partial {y_k}}} \\
		&+ \zeta  _1^{ - 2}\zeta  _2^2[{R_{ij}}({\bm{y}},{\bm{z}},{\theta _0})\frac{{{\partial ^2}{A^0}}}{{\partial {y_i}\partial {y_j}}} + {D_k}({\bm{y}},{\bm{z}},{\theta _0})\frac{{\partial {A^0}}}{{\partial {y_k}}} \\
		&+ {A^ * }({\bm{y}},{\bm{z}},{\theta _0})] + {\rm{O}}(\zeta  _1^{ - 3}\zeta  _2^3).
	\end{aligned}
\end{equation}
\begin{equation}
	\begin{aligned}
		{M_{{\alpha _1}{\alpha _2}}}({\bm{y}},{\theta _0}) &= \!M_{{\alpha _1}{\alpha _2}}^0({\bm{y}},{\theta _0})\!+\! \zeta  _1^{ - 1}{\zeta  _2}[{R_k}({\bm{y}},{\bm{z}},{\theta _0})\frac{{\partial M_{{\alpha _1}{\alpha _2}}^0}}{{\partial {y_k}}}\! + \!{R_{{\alpha _1}}}({\bm{y}},{\bm{z}},{\theta _0})M_{{\alpha _2}}^0]\\
		&+ \zeta  _1^{ - 2}\zeta  _2^2[{R_{ij}}({\bm{y}},{\bm{z}},{\theta _0})\frac{{{\partial ^2}M_{{\alpha _1}{\alpha _2}}^0}}{{\partial {y_i}\partial {y_j}}} + {D_k}({\bm{y}},{\bm{z}},{\theta _0})\frac{{\partial M_{{\alpha _1}{\alpha _2}}^0}}{{\partial {y_k}}} \\
		&+ {R_{{\alpha _1}{\alpha _2}}}({\bm{y}},{\bm{z}},{\theta _0}) + {R_{{\alpha _1}j}}({\bm{y}},{\bm{z}},{\theta _0})\frac{{\partial M_{{\alpha _2}}^0}}{{\partial {y_j}}} \!+ \!{R_{j{\alpha _1}}}({\bm{y}},{\bm{z}},{\theta _0})\frac{{\partial M_{{\alpha _2}}^0}}{{\partial {y_j}}}\\
		&  + {D_{{\alpha _1}}}({\bm{y}},{\bm{z}},{\theta _0})M_{{\alpha _2}}^0({\bm{y}},{\theta _0})] + {\rm{O}}(\zeta  _1^{ - 3}\zeta  _2^3).
	\end{aligned}
\end{equation}
\begin{equation}
	\begin{aligned}
		{C_{{\alpha _1}}}({\bm{y}},{\theta _0}) &= C_{{\alpha _1}}^0({\bm{y}},{\theta _0}) + \zeta  _1^{ - 1}{\zeta  _2}[{R_k}({\bm{y}},{\bm{z}},{\theta _0})\frac{{\partial C_{{\alpha _1}}^0}}{{\partial {y_k}}}\! +\! {R_k}({\bm{y}},{\bm{z}},{\theta _0})\frac{{\partial M_{{\alpha _1}}^0}}{{\partial {x_k}}}]\\
		&+ \zeta  _1^{ - 2}\zeta  _2^2[{R_{ij}}({\bm{y}},{\bm{z}},{\theta _0})\frac{{{\partial ^2}C_{{\alpha _1}}^0}}{{\partial {y_i}\partial {y_j}}} + {D_k}({\bm{y}},{\bm{z}},{\theta _0})\frac{{\partial C_{{\alpha _1}}^0}}{{\partial {y_k}}} \\
		& \!+\! D_{{\alpha _1}}^*({\bm{y}},{\bm{z}},{\theta _0}) \!+\! {R_{ij}}({\bm{y}},{\bm{z}},{\theta _0})\frac{{\partial M_{{\alpha _1}}^0}}{{\partial {y_i}\partial {x_j}}}\!+ \!{R_{ij}}({\bm{y}},{\bm{z}},{\theta _0})\frac{{\partial M_{{\alpha _1}}^0}}{{\partial {x_i}\partial {y_j}}} \\
		& + {D_k}({\bm{y}},{\bm{z}},{\theta _0})\frac{{\partial M_{{\alpha _1}}^0}}{{\partial {x_k}}} + D_k^ * ({\bm{y}},{\bm{z}},{\theta _0})\frac{{\partial M_{{\alpha _1}}^0}}{{\partial {y_k}}}] + {\rm{O}}(\zeta  _1^{ - 3}\zeta  _2^3).
	\end{aligned}
\end{equation}
\begin{equation}
	\begin{aligned}
		{B_{{\alpha _1}{\alpha _2}}}({\bm{y}},{\theta _0}) &= B_{{\alpha _1}{\alpha _2}}^0({\bm{y}},{\theta _0}) + \zeta  _1^{ - 1}{\zeta  _2}[{R_k}({\bm{y}},{\bm{z}},{\theta _0})\frac{{\partial B_{{\alpha _1}{\alpha _2}}^0}}{{\partial {y_k}}} \\
		&+ M_{{\alpha _1}}^0({\bm{y}},{\theta _0})R_{{\alpha _2}}^ * ({\bm{y}},{\bm{z}},{\theta _0}) + M_{{\alpha _1}}^0({\bm{y}},{\theta _0})R_k^ * ({\bm{y}},{\bm{z}},{\theta _0})\frac{{\partial M_{{\alpha _2}}^0}}{{\partial {y_k}}}]\\
		&+ \zeta  _1^{ - 2}\zeta  _2^2[{R_{ij}}({\bm{y}},{\bm{z}},{\theta _0})\frac{{{\partial ^2}B_{{\alpha _1}{\alpha _2}}^0}}{{\partial {y_i}\partial {y_j}}} + {D_k}({\bm{y}},{\bm{z}},{\theta _0})\frac{{\partial B_{{\alpha _1}{\alpha _2}}^0}}{{\partial {y_k}}} \\
		&+ {G_{{\alpha _1}{\alpha _2}}}({\bm{y}},{\bm{z}})] + {\rm{O}}(\zeta  _1^{ - 3}\zeta  _2^3).
	\end{aligned}
\end{equation}
\begin{equation}
	\begin{aligned}
		{E_{{\alpha _1}{\alpha _2}}}({\bm{y}},{\theta _0}) &= E_{{\alpha _1}{\alpha _2}}^0({\bm{y}},{\theta _0}) + \zeta  _1^{ - 1}{\zeta  _2}{R_k}({\bm{y}},{\bm{z}},{\theta _0})\frac{{\partial E_{{\alpha _1}{\alpha _2}}^0}}{{\partial {y_k}}} \\
		&+ \zeta  _1^{ - 2}\zeta  _2^2[{R_{ij}}({\bm{y}},{\bm{z}},{\theta _0})\frac{{{\partial ^2}E_{{\alpha _1}{\alpha _2}}^0}}{{\partial {y_i}\partial {y_j}}} + {D_k}({\bm{y}},{\bm{z}},{\theta _0})\frac{{\partial E_{{\alpha _1}{\alpha _2}}^0}}{{\partial {y_k}}}  \\
		&+ E_{{\alpha _1}{\alpha _2}}^ * ({\bm{y}},{\bm{z}},{\theta _0}) + E_{pn}^ * ({\bm{y}},{\bm{z}},{\theta _0})\frac{{\partial N_{p{\alpha _1}}^{0{\alpha _2}}}}{{\partial {y_n}}}] + {\rm{O}}(\zeta  _1^{ - 3}\zeta  _2^3).
	\end{aligned}
\end{equation}
\begin{equation}
	\begin{aligned}
		F_i^{{\alpha _1}}({\bm{y}},{\theta _0}) &= F_{i}^{{0\alpha _1}}({\bm{y}},{\theta _0}) + \zeta  _1^{ - 1}{\zeta  _2}T_{im}^n({\bm{y}},{\bm{z}},{\theta _0})\frac{{\partial F_{m}^{{0\alpha _1}}}}{{\partial {y_n}}}\\
		&+ \zeta  _1^{ - 2}\zeta  _2^2[T_{im}^{qn}({\bm{y}},{\bm{z}},{\theta _0})\frac{{{\partial ^2}F_{m}^{{0\alpha _1}}}}{{\partial {y_q}\partial {y_n}}} + L_{im}^n({\bm{y}},{\bm{z}},{\theta _0})\frac{{\partial F_{m}^{{0\alpha _1}}}}{{\partial {y_n}}} \\
		&+ {F_i^ * }({\bm{y}},{\bm{z}},{\theta _0})] + {\rm{O}}(\zeta  _1^{ - 3}\zeta  _2^3).
	\end{aligned}
\end{equation}
\begin{equation}
	\begin{aligned}
	N_{im}^{{\alpha _1}{\alpha _2}}({\bm{y}},{\theta _0}) &= N_{im}^{0{\alpha _1}{\alpha _2}}({\bm{y}},{\theta _0}) + \zeta  _1^{ - 1}{\zeta  _2}[T_{ip}^n({\bm{y}},{\bm{z}},{\theta _0})\frac{{\partial N_{pm}^{0{\alpha _1}{\alpha _2}}}}{{\partial {y_n}}} + T_{ip}^{{\alpha _1}}({\bm{y}},{\bm{z}},{\theta _0})N_{pm}^{0{\alpha _2}}]\\
	&+ \zeta  _1^{ - 2}\zeta  _2^2[ T_{ip}^{qn}({\bm{y}},{\bm{z}},{\theta _0})\frac{{{\partial ^2}N_{pm}^{0{\alpha _1}{\alpha _2}}}}{{\partial {y_q}\partial {y_n}}}+ T_{ip}^{{\alpha _1}n}({\bm{y}},{\bm{z}},{\theta _0})\frac{{\partial N_{pm}^{0{\alpha _2}}}}{{\partial {y_n}}}\\
	& + L_{ip}^n({\bm{y}},{\bm{z}},{\theta _0})\frac{{\partial N_{pm}^{0{\alpha _1}{\alpha _2}}}}{{\partial {y_n}}} + L_{ip}^{{\alpha _1}}({\bm{y}},{\bm{z}},{\theta _0})N_{pm}^{0{\alpha _2}}({\bm{y}},{\theta _0})\\
	&+T_{im}^{{\alpha _1}{\alpha _2}}({\bm{y}},{\bm{z}},{\theta _0})+ T_{ip}^{n{\alpha _1}}({\bm{y}},{\bm{z}},{\theta _0})\frac{{\partial N_{pm}^{0{\alpha _2}}}}{{\partial {y_n}}}] + {\rm{O}}(\zeta  _1^{ - 3}\zeta  _2^3).
	\end{aligned}
\end{equation}
\begin{equation}
	\begin{aligned}
		Z_{im}^{{\alpha _1}}({\bm{y}},{\theta _0}) &= Z_{im}^{0{\alpha _1}}({\bm{y}},{\theta _0}) + \zeta  _1^{ - 1}{\zeta  _2}[T_{ip}^n({\bm{y}},{\bm{z}},{\theta _0})\frac{{\partial Z_{pm}^{0{\alpha _1}}}}{{\partial {y_n}}} + T_{ip}^n({\bm{y}},{\bm{z}},{\theta _0})\frac{{\partial N_{pm}^{0{\alpha _1}}}}{{\partial {x_n}}}]\\
		&+ \zeta  _1^{ - 2}\zeta  _2^2[U_{im}^{{\alpha _1}}({\bm{y}},{\bm{z}},{\theta _0}) + T_{ip}^{qn}({\bm{y}},{\bm{z}},{\theta _0})\frac{{{\partial ^2}Z_{pm}^{0{\alpha _1}}}}{{\partial {y_q}\partial {y_n}}} + L_{ip}^n({\bm{y}},{\bm{z}},{\theta _0})\frac{{\partial Z_{pm}^{0{\alpha _1}}}}{{\partial {y_n}}}\\
		&+ T_{ip}^{qn}({\bm{y}},{\bm{z}},{\theta _0})\frac{{{\partial ^2}N_{pm}^{0{\alpha _1}}}}{{\partial {y_q}\partial {x_n}}} + L_{ip}^n({\bm{y}},{\bm{z}},{\theta _0})\frac{{\partial N_{pm}^{0{\alpha _1}}}}{{\partial {x_n}}} + T_{ip}^{qn}({\bm{y}},{\bm{z}},{\theta _0})\frac{{{\partial ^2}N_{pm}^{0{\alpha _1}}}}{{\partial {x_q}\partial {y_n}}}\\
		&+ U_{ip}^n({\bm{y}},{\bm{z}},{\theta _0})\frac{{\partial N_{pm}^{0{\alpha _1}}}}{{\partial {y_n}}}] + {\rm{O}}(\zeta  _1^{ - 3}\zeta  _2^3).
	\end{aligned}
\end{equation}
\begin{equation}
\begin{aligned}
	{Q_i}({\bm{y}},{\theta _0}) &= {Q_{i}^0}({\bm{y}},{\theta _0}) + \zeta  _1^{ - 1}{\zeta  _2}[T_{im}^n({\bm{y}},{\bm{z}},{\theta _0})\frac{{\partial {Q_{m}^0}}}{{\partial {y_n}}} + T_{im}^n({\bm{y}},{\bm{z}},{\theta _0})\frac{{\partial P_m^0}}{{\partial {x_n}}}]\\
	&+ \zeta  _1^{ - 2}\zeta  _2^2[Q_k^ * ({\bm{y}},{\bm{z}},{\theta _0}) + T_{im}^{qn}({\bm{y}},{\bm{z}},{\theta _0})\frac{{{\partial ^2}{Q_{m}^0}}}{{\partial {y_q}\partial {y_n}}} + L_{im}^n({\bm{y}},{\bm{z}},{\theta _0})\frac{{\partial {Q_{m}^0}}}{{\partial {y_n}}}\\
	& + T_{im}^{qn}({\bm{y}},{\bm{z}},{\theta _0})\frac{{{\partial ^2}P_m^0}}{{\partial {x_q}\partial {y_n}}} + U_{im}^n({\bm{y}},{\bm{z}},{\theta _0})\frac{{\partial P_m^0}}{{\partial {y_n}}} + T_{im}^{qn}({\bm{y}},{\bm{z}},{\theta _0})\frac{{{\partial ^2}P_m^0}}{{\partial {y_q}\partial {x_n}}} \\
	&+ L_{im}^n({\bm{y}},{\bm{z}},{\theta _0})\frac{{\partial P_m^0}}{{\partial {x_n}}}] +{\rm{O}}(\zeta  _1^{ - 3}\zeta  _2^3).
\end{aligned}
\end{equation}
\begin{equation}
	\begin{aligned}
		H_i^{{\alpha _1}}({\bm{y}},{\theta _0}) &= H_{i}^{{0\alpha _1}}({\bm{y}},{\theta _0})\! +\! \zeta  _1^{ - 1}{\zeta  _2}[T_{im}^n({\bm{y}},{\bm{z}},{\theta _0})\frac{{\partial H_{m}^{{0\alpha _1}}}}{{\partial {y_n}}}\!+\! T_{im}^{{\alpha _1}}({\bm{y}},{\bm{z}},{\theta _0})P_m^0 \\
		&+ {O_i}({\bm{y}},{\bm{z}},{\theta _0})M_{{\alpha _1}}^0]+ \zeta  _1^{ - 2}\zeta  _2^2[V_i^{{\alpha _1}}({\bm{y}},{\bm{z}},{\theta _0}) \\
		&+ T_{im}^{qn}({\bm{y}},{\bm{z}},{\theta _0})\frac{{{\partial ^2}H_{m}^{{0\alpha _1}}}}{{\partial {y_q}\partial {y_n}}} + L_{im}^n({\bm{y}},{\bm{z}},{\theta _0})\frac{{\partial H_{m}^{{0\alpha _1}}}}{{\partial {y_n}}}\\
		&\!+ \!T_{im}^{{\alpha _1}n}({\bm{y}},{\bm{z}},{\theta _0})\frac{{\partial P_m^0}}{{\partial {y_n}}} \!+ \!T_{im}^{n{\alpha _1}}({\bm{y}},{\bm{z}},{\theta _0})\frac{{\partial P_m^0}}{{\partial {y_n}}} \!+ \! L_{im}^{{\alpha _1}}({\bm{y}},{\bm{z}},{\theta _0})P_m^0\\
		&+ V_i^n({\bm{y}},{\bm{z}},{\theta _0})\frac{{\partial M_{{\alpha _1}}^0}}{{\partial {y_n}}} + P_i^ * ({\bm{y}},{\bm{z}},{\theta _0})M_{{\alpha _1}}^0] + {\rm{O}}(\zeta  _1^{ - 3}\zeta  _2^3).
	\end{aligned}
\end{equation}
\begin{equation}
	\begin{aligned}
		W_i^{{\alpha _1}}({\bm{y}},{\theta _0}) &= W_{i}^{{0\alpha _1}}({\bm{y}},{\theta _0}) + \zeta  _1^{ - 1}{\zeta  _2}[T_{im}^n({\bm{y}},{\bm{z}},{\theta _0})\frac{{\partial W_{m}^{{0\alpha _1}}}}{{\partial {y_n}}}\\
		&+ M_{{\alpha _1}}^0({\bm{y}},{\theta _0})O_i^ * ({\bm{y}},{\bm{z}},{\theta _0}) + M_{{\alpha _1}}^0({\bm{y}},{\theta _0})T_{im}^{ * n}({\bm{y}},{\bm{z}},{\theta _0})\frac{{\partial P_m^0}}{{\partial {y_n}}}]\\
		&+ \zeta  _1^{ - 2}\zeta  _2^2[T_{im}^{qn}({\bm{y}},{\bm{z}},{\theta _0})\frac{{{\partial ^2}W_{m}^{{0\alpha _1}}}}{{\partial {y_q}\partial {y_n}}} + L_{im}^n({\bm{y}},{\bm{z}},{\theta _0})\frac{{\partial W_{m}^{{0\alpha _1}}}}{{\partial {y_n}}} \\
		&+ X_i^{{\alpha _1}}({\bm{y}},{\bm{z}})] + {\rm{O}}(\zeta  _1^{ - 3}\zeta  _2^3).
	\end{aligned}
\end{equation}
\begin{equation}
	\begin{aligned}
		J_{im}^{{\alpha _1}{\alpha _2}}({\bm{y}},{\bm{z}},{\theta _0}) &= J_{im}^{0{\alpha _1}{\alpha _2}}({\bm{y}},{\theta _0}) + \zeta  _1^{ - 1}{\zeta  _2}[T_{kp}^n({\bm{y}},{\bm{z}},{\theta _0})\frac{{\partial J_{pm}^{0{\alpha _1}{\alpha _2}}}}{{\partial {y_n}}}\\
		&+ M_{{\alpha _1}}^0({\bm{y}},{\theta _0})T_{im}^{ * {\alpha _2}}({\bm{y}},{\bm{z}},{\theta _0}) + M_{{\alpha _1}}^0({\bm{y}},{\theta _0})T_{ip}^{ * n}({\bm{y}},{\bm{z}},{\theta _0})\frac{{\partial N_{pm}^{0{\alpha _2}}}}{{\partial {y_n}}}]\\
		&+ \zeta  _1^{ - 2}\zeta  _2^2[T_{ip}^{qn}({\bm{y}},{\bm{z}},{\theta _0})\frac{{{\partial ^2}J_{pm}^{0{\alpha _1}{\alpha _2}}}}{{\partial {y_q}\partial {y_n}}} + L_{ip}^n({\bm{y}},{\bm{z}},{\theta _0})\frac{{\partial J_{pm}^{0{\alpha _1}{\alpha _2}}}}{{\partial {y_n}}}\\
		& + J_{im}^{ * {\alpha _1}{\alpha _2}}({\bm{y}},{\bm{z}})] + {\rm{O}}(\zeta  _1^{ - 3}\zeta  _2^3).
	\end{aligned}
\end{equation}
Moreover, Appendix B provides the detailed formulations of the auxiliary cell functions $M_{{\alpha _1}}^0$, $P_{i}^0$, ${A^0}$, ${A^*}$, $M_{{\alpha _1}{\alpha _2}}^0$, $C_{{\alpha _1}}^0$, $D_{{\alpha _1}}^*$, $B_{{\alpha _1}{\alpha _2}}^0$, $R_{{\alpha _1}}^*$, $E_{{\alpha _1}{\alpha _2}}^0$, $E_{{\alpha _1}{\alpha _2}}^*$, $F_{i}^{{0\alpha _1}}$, ${F_i^*}$, $N_{im}^{0{\alpha _1}{\alpha _2}}$, $Z_{im}^{0{\alpha _1}}$, $U_{im}^{{\alpha _1}}$, ${Q_{i}^0}$, $Q_i^*$, $H_{i}^{{0\alpha _1}}$, $V_i^{{\alpha _1}}$, $P_i^*$, $W_{i}^{{0\alpha _1}}$, $O_i^*$, $T_{im}^{*n}$ and $J_{im}^{0{\alpha _1}{\alpha _2}}$ in (53)-(66).

\subsection{The higher-order three-scale asymptotic model of nonlinear thermo-mechanical coupling problems}
By the macro-meso correlative SOTS modeling and meso-micro correlative SOTS modeling, the core idea of macro-meso-micro three-scale modeling is shown in Fig.\hspace{1mm}4.
\begin{figure}[!htb]
	\centering
	\begin{minipage}[c]{1.0\textwidth}
		\centering
		\includegraphics[
		width=0.8\linewidth
		]{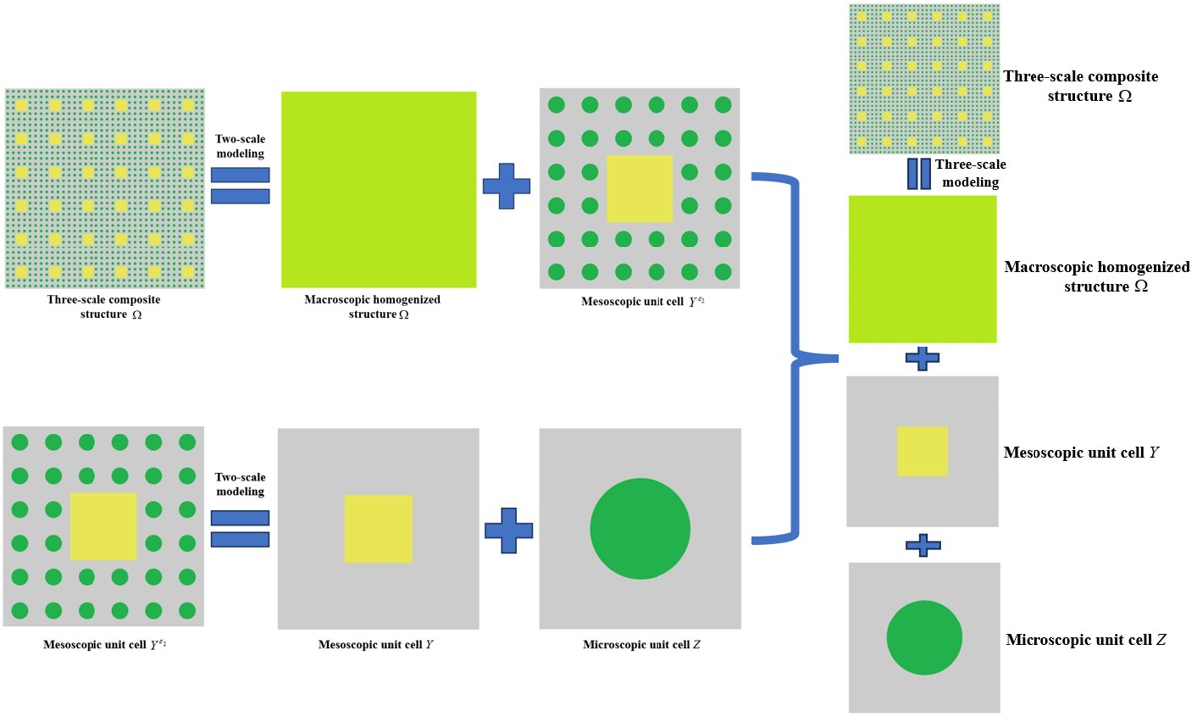}
	\end{minipage}
	\caption {The illustration of macro-meso-micro three-scale modeling.}\label{f1}
\end{figure}

Then, the macro-meso-micro coupled HOTS solutions for nonlinear thermo-mechanical problem (1) of heterogeneous structures with multiple spatial scales are established as follow.
\begin{equation}
	\begin{aligned}
		{\theta ^{{\zeta  _1}{\zeta  _2}}}({\bm{x}},t) &= {\theta ^{(0)}}({\bm{x}},t) + {\zeta  _1}{\theta ^{(1)}}({\bm{x}},{\bm{y}},{t}) + {\zeta  _2}{\theta ^{(2)}}({\bm{x}},{\bm{y}},{\bm{z}},{t})\\
		& + \zeta  _1^{ - 1}\zeta  _2^2{\theta ^{(3)}}({\bm{x}},{\bm{y}},{\bm{z}},{t}) + \zeta  _1^2{\theta ^{(4)}}({\bm{x}},{\bm{y}},{t}) + {\zeta  _1}{\zeta  _2}{\theta ^{(5)}}({\bm{x}},{\bm{y}},{\bm{z}},{t}) \\
		& + \zeta  _2^2{\theta ^{(6)}}({\bm{x}},{\bm{y}},{\bm{z}},{t}) + {\rm{O}}(\zeta  _1^3) + {\rm{O}}(\zeta  _1^{ - 2}\zeta  _2^3) + {\rm{O}}(\zeta  _1^{ - 1}\zeta  _2^3)\\
		&=: \theta _T^{{\zeta  _1}{\zeta  _2}}(\bm{x},t) + {\rm{O}}(\zeta  _1^{ - 2}\zeta  _2^3) + {\rm{O}}(\zeta  _1^{ - 1}\zeta  _2^3) + {\rm{O}}(\zeta  _1^3),
	\end{aligned}
\end{equation}
\begin{equation}
	\begin{aligned}
		u_i^{{\zeta  _1}{\zeta  _2}}({\bm{x}},t) &= {u_{i}^{(0)}}({\bm{x}},t) + {\zeta  _1}u_i^{(1)}({\bm{x}},{\bm{y}},{t}) + {\zeta  _2}u_i^{(2)}({\bm{x}},{\bm{y}},{\bm{z}},{t})\\
		& + \zeta  _1^{ - 1}\zeta  _2^2u_i^{(3)}({\bm{x}},{\bm{y}},{\bm{z}},{t})+ \zeta  _1^2u_i^{(4)}({\bm{x}},{\bm{y}},{t}) + {\zeta  _1}{\zeta  _2}u_i^{(5)}({\bm{x}},{\bm{y}},{\bm{z}},{t})\\
		& + \zeta  _2^2u_i^{(6)}({\bm{x}},{\bm{y}},{\bm{z}},{t})+ {\rm{O}}(\zeta  _1^3) + {\rm{O}}(\zeta  _1^{ - 2}\zeta  _2^3) + {\rm{O}}(\zeta  _1^{ - 1}\zeta  _2^3)\\
		&=: {u}_{iT}^{{\zeta  _1}{\zeta  _2}}(\bm{x},t) + {\rm{O}}(\zeta  _1^{ - 2}\zeta  _2^3) + {\rm{O}}(\zeta  _1^{ - 1}\zeta  _2^3) + {\rm{O}}(\zeta  _1^3),
	\end{aligned}
\end{equation}
where ${\theta ^{(0)}}$, ${\theta ^{(1)}}$, ${\theta ^{(2)}}$, ${\theta ^{(3)}}$, ${\theta ^{(4)}}$, ${\theta ^{(5)}}$ and ${\theta ^{(6)}} $  in (67) and $u_i^{(0)}$, $u_i^{(1)}$, $u_i^{(2)}$, $u_i^{(3)}$, $u_i^{(4)}$, $u_i^{(5)}$ and $u_i^{(6)}$ in (68) are defined as follows
\begin{small}
\begin{equation}
	{\theta ^{(0)}}({\bm{x}},t) = {\theta _0}({\bm{x}},t),
\end{equation}
\begin{equation}
	{\theta ^{(1)}}({\bm{x}},{\bm{y}},{t}) = M_{{\alpha _1}}^0\frac{{\partial {\theta _0}}}{{\partial {x_{{\alpha _1}}}}},
\end{equation}
\begin{equation}
	{\theta ^{(2)}}({\bm{x}},{\bm{y}},{\bm{z}},{t}) = [{R_k}\frac{{\partial M_{{\alpha _1}}^0}}{{\partial {y_k}}} + {R_{{\alpha _1}}}]\frac{{\partial {\theta _0}}}{{\partial {x_{{\alpha _1}}}}},
\end{equation}
\begin{equation}
	{\theta ^{(3)}}({\bm{x}},{\bm{y}},{\bm{z}},{t}) = [{R_{ij}}\frac{{{\partial ^2}M_{{\alpha _1}}^0}}{{\partial {y_i}\partial {y_j}}} + {D_k}\frac{{\partial M_{{\alpha _1}}^0}}{{\partial {y_k}}} + {D_{{\alpha _1}}}]\frac{{\partial {\theta _0}}}{{\partial {x_{{\alpha _1}}}}},
\end{equation}
\begin{equation}
	\begin{aligned}
		{\theta ^{(4)}}({\bm{x}},{\bm{y}},{t}) &= {A^0}\frac{{\partial {\theta _0}}}{{\partial t}} + M_{{\alpha _1}{\alpha _2}}^0\frac{{{\partial ^2}{\theta _0}}}{{\partial {x_{{\alpha _1}}}\partial {x_{{\alpha _2}}}}} + C_{{\alpha _1}}^0\frac{{\partial {\theta _0}}}{{\partial {x_{{\alpha _1}}}}} \\
		& - B_{{\alpha _1}{\alpha _2}}^0\frac{{\partial {\theta _0}}}{{\partial {x_{{\alpha _1}}}}}\frac{{\partial {\theta _0}}}{{\partial {x_{{\alpha _2}}}}} - E_{{\alpha _1}{\alpha _2}}^0\frac{{{\partial ^2}{u_{{\alpha _1}0}}}}{{\partial {x_{{\alpha _2}}}\partial t}},
	\end{aligned}
\end{equation}
\begin{equation}
	\begin{aligned}
		{\theta ^{(5)}}({\bm{x}},{\bm{y}},{\bm{z}},{t})& = {R_k}\frac{{\partial {A^0}}}{{\partial {y_k}}}\frac{{\partial {\theta _0}}}{{\partial t}} - {R_k}\frac{{\partial E_{{\alpha _1}{\alpha _2}}^0}}{{\partial {y_k}}}\frac{{{\partial ^2}{u_{{\alpha _1}0}}}}{{\partial {x_{{\alpha _2}}}\partial t}} + [{R_k}\frac{{\partial M_{{\alpha _1}{\alpha _2}}^0}}{{\partial {y_k}}}  \\
		&+ {R_{{\alpha _1}}}M_{{\alpha _2}}^0]\frac{{{\partial ^2}{\theta _0}}}{{\partial {x_{{\alpha _1}}}\partial {x_{{\alpha _2}}}}}+ [{R_k}\frac{{\partial C_{{\alpha _1}}^0}}{{\partial {y_k}}} + {R_k}\frac{{\partial M_{{\alpha _1}}^0}}{{\partial {x_k}}}]\frac{{\partial {\theta _0}}}{{\partial {x_{{\alpha _1}}}}}\\
		&- [{R_k}\frac{{\partial B_{{\alpha _1}{\alpha _2}}^0}}{{\partial {y_k}}} + M_{{\alpha _1}}^0R_{{\alpha _2}}^ *+ M_{{\alpha _1}}^0R_k^ *\frac{{\partial M_{{\alpha _2}}^0}}{{\partial {y_k}}}]\frac{{\partial {\theta _0}}}{{\partial {x_{{\alpha _1}}}}}\frac{{\partial {\theta _0}}}{{\partial {x_{{\alpha _2}}}}},
	\end{aligned}
\end{equation}
\begin{equation}
	\begin{aligned}
		{\theta ^{(6)}}({\bm{x}},{\bm{y}},{\bm{z}},{t}) &= [{R_{ij}}\frac{{{\partial ^2}{A^0}}}{{\partial {y_i}\partial {y_j}}} + {D_k}\frac{{\partial {A^0}}}{{\partial {y_k}}} + {A^ * }]\frac{{\partial {\theta _0}}}{{\partial t}}+ [{R_{ij}}\frac{{{\partial ^2}M_{{\alpha _1}{\alpha _2}}^0}}{{\partial {y_i}\partial {y_j}}}\\
		&  + {D_k}\frac{{\partial M_{{\alpha _1}{\alpha _2}}^0}}{{\partial {y_k}}}+ {R_{{\alpha _1}{\alpha _2}}}+ {R_{{\alpha _1}j}}\frac{{\partial M_{{\alpha _2}}^0}}{{\partial {y_j}}} + {R_{j{\alpha _1}}}\frac{{\partial M_{{\alpha _2}}^0}}{{\partial {y_j}}} \\
		&+ {D_{{\alpha _1}}}M_{{\alpha _2}}^0]\frac{{{\partial ^2}{\theta _0}}}{{\partial {x_{{\alpha _1}}}\partial {x_{{\alpha _2}}}}}+ [{R_{ij}}\frac{{{\partial ^2}C_{{\alpha _1}}^0}}{{\partial {y_i}\partial {y_j}}} + {D_k}\frac{{\partial C_{{\alpha _1}}^0}}{{\partial {y_k}}} + {R_{ij}}\frac{{\partial M_{{\alpha _1}}^0}}{{\partial {y_i}\partial {x_j}}}\\
		&+ {R_{ij}}\frac{{\partial M_{{\alpha _1}}^0}}{{\partial {x_i}\partial {y_j}}} + {D_k}\frac{{\partial M_{{\alpha _1}}^0}}{{\partial {x_k}}}+ D_{{\alpha _1}}^ *+ D_k^ *\frac{{\partial M_{{\alpha _1}}^0}}{{\partial {y_k}}}]\frac{{\partial {\theta _0}}}{{\partial {x_{{\alpha _1}}}}}\\
		&- [{R_{ij}}\frac{{{\partial ^2}B_{{\alpha _1}{\alpha _2}}^0}}{{\partial {y_i}\partial {y_j}}} + {D_k}\frac{{\partial B_{{\alpha _1}{\alpha _2}}^0}}{{\partial {y_k}}} + {G_{{\alpha _1}{\alpha _2}}}]\frac{{\partial {\theta _0}}}{{\partial {x_{{\alpha _1}}}}}\frac{{\partial {\theta _0}}}{{\partial {x_{{\alpha _2}}}}}\\
		&- [{R_{ij}}\frac{{{\partial ^2}E_{{\alpha _1}{\alpha _2}}^0}}{{\partial {y_i}\partial {y_j}}} + {D_k}\frac{{\partial E_{{\alpha _1}{\alpha _2}}^0}}{{\partial {y_k}}} + E_{{\alpha _1}{\alpha _2}}^ *+ E_{pn}^*\frac{{\partial N_{p{\alpha _1}}^{0{\alpha _2}}}}{{\partial {y_n}}}]\frac{{{\partial ^2}{u_{{\alpha _1}0}}}}{{\partial {x_{{\alpha _2}}}\partial t}},
	\end{aligned}
\end{equation}
\begin{equation}
	u_i^{(0)}({\bm{x}},t) = u_{i0}({\bm{x}},t),
\end{equation}
\begin{equation}
	u_i^{(1)}({\bm{x}},{\bm{y}},{t}) = N_{im}^{0{\alpha _1}}\frac{{\partial {u_{m0}}}}{{\partial {x_{{\alpha _1}}}}} - P_i^0({\theta _0} - \tilde \theta ),
\end{equation}
\begin{equation}
	\begin{aligned}
		u_i^{(2)}({\bm{x}},{\bm{y}},{\bm{z}},{t}) &= [T_{ip}^n\frac{{\partial N_{pm}^{0{\alpha _1}}}}{{\partial {y_n}}} + T_{im}^{{\alpha _1}}]\frac{{\partial {u_{m0}}}}{{\partial {x_{{\alpha _1}}}}} + [T_{im}^{{\alpha _1}}\frac{{\partial P_m^0}}{{\partial {y_{{\alpha _1}}}}} + {O_i}]({\theta _0} - \tilde \theta ),
	\end{aligned}
\end{equation}
\begin{equation}
	\begin{aligned}
	u_i^{(3)}({\bm{x}},{\bm{y}},{\bm{z}},{t}) &= [T_{ip}^{qn}\frac{{{\partial ^2}N_{pm}^{0{\alpha _1}}}}{{\partial {y_q}\partial {y_n}}} + L_{ip}^n\frac{{\partial N_{pm}^{0{\alpha _1}}}}{{\partial {y_n}}} + L_{im}^{{\alpha _1}}]\frac{{\partial {u_{m0}}}}{{\partial {x_{{\alpha _1}}}}}\\
	&- [T_{im}^{{\alpha _1}{\alpha _2}}\frac{{{\partial ^2}P_m^0}}{{\partial {y_{{\alpha _1}}}\partial {y_{{\alpha _2}}}}} + L_{im}^{{\alpha _1}}\frac{{\partial P_m^0}}{{\partial {y_{{\alpha _1}}}}} + P_i^ * ]({\theta _0} - \tilde \theta ),	
	\end{aligned}
\end{equation}
\begin{equation}
	\begin{aligned}
		u_i^{(4)}({\bm{x}},{\bm{y}},{t}) &= F_{i}^{{0\alpha _1}}\frac{{{\partial ^2}{u_{{\alpha _1}0}}}}{{\partial {t^2}}} + N_{im}^{0{\alpha _1}{\alpha _2}}\frac{{{\partial ^2}{u_{m0}}}}{{\partial {x_{{\alpha _1}}}\partial {x_{{\alpha _2}}}}} + Z_{im}^{0{\alpha _1}}\frac{{\partial {u_{m0}}}}{{\partial {x_{{\alpha _1}}}}}- {Q_{i}^0}({\theta _0} - \tilde \theta )\\
		& - H_{i}^{{0\alpha _1}}\frac{{\partial {\theta _0}}}{{\partial {x_{{\alpha _1}}}}}+ W_{i}^{{0\alpha _1}}\frac{{\partial {\theta _0}}}{{\partial {x_{{\alpha _1}}}}}({\theta _0} - \tilde \theta ) - J_{im}^{0{\alpha _1}{\alpha _2}}\frac{{\partial {\theta _0}}}{{\partial {x_{{\alpha _1}}}}}\frac{{\partial {u_{m0}}}}{{\partial {x_{{\alpha _2}}}}},
	\end{aligned}
\end{equation}
\begin{equation}
	\begin{aligned}
		u_i^{(5)}({\bm{x}},{\bm{y}},{\bm{z}},{t})&= T_{im}^n\frac{{\partial F_{m}^{{0\alpha _1}}}}{{\partial {y_n}}}\frac{{{\partial ^2}{u_{{\alpha _1}0}}}}{{\partial {t^2}}}+ [T_{ip}^n\frac{{\partial N_{pm}^{0{\alpha _1}{\alpha _2}}}}{{\partial {y_n}}} + T_{ip}^{{\alpha _1}}N_{pm}^{0{\alpha _2}}]\frac{{{\partial ^2}{u_{m0}}}}{{\partial {x_{{\alpha _1}}}\partial {x_{{\alpha _2}}}}}\\
		&+ [T_{ip}^n\frac{{\partial Z_{pm}^{0{\alpha _1}}}}{{\partial {y_n}}} + T_{ip}^n\frac{{\partial N_{pm}^{0{\alpha _1}}}}{{\partial {x_n}}}]\frac{{\partial {u_{m0}}}}{{\partial {x_{{\alpha _1}}}}}- [T_{im}^n\frac{{\partial {Q_{m}^{0}}}}{{\partial {y_n}}} \\
		& + T_{im}^n\frac{{\partial P_m^0}}{{\partial {x_n}}}]({\theta _0} - \tilde \theta )- [T_{im}^n\frac{{\partial H_{m}^{{0\alpha _1}}}}{{\partial {y_n}}} + T_{im}^{{\alpha _1}}P_m^0+ {O_i}M_{{\alpha _1}}^0]\frac{{\partial {\theta _0}}}{{\partial {x_{{\alpha _1}}}}}\\
		&+ [T_{im}^n\frac{{\partial W_{m}^{{0\alpha _1}}}}{{\partial {y_n}}} + M_{{\alpha _1}}^0O_i^ * + M_{{\alpha _1}}^0T_{im}^{ * n}\frac{{\partial P_m^0}}{{\partial {y_n}}}]\frac{{\partial {\theta _0}}}{{\partial {x_{{\alpha _1}}}}}({\theta _0} - \tilde \theta )\\
		&- [T_{kp}^n\frac{{\partial J_{pm}^{0{\alpha _1}{\alpha _2}}}}{{\partial {y_n}}} + M_{{\alpha _1}}^0T_{im}^{ * {\alpha _2}}+ M_{{\alpha _1}}^0T_{ip}^{ * n}\frac{{\partial N_{pm}^{0{\alpha _2}}}}{{\partial {y_n}}}]\frac{{\partial {\theta _0}}}{{\partial {x_{{\alpha _1}}}}}\frac{{\partial {u_{m0}}}}{{\partial {x_{{\alpha _2}}}}},
	\end{aligned}
\end{equation}
\begin{equation}
	\begin{aligned}
		u_i^{(6)}({\bm{x}},{\bm{y}},{\bm{z}},{t}) &= [T_{im}^{qn}\frac{{{\partial ^2}F_{m}^{{0\alpha _1}}}}{{\partial {y_q}\partial {y_n}}} + L_{im}^n\frac{{\partial F_{m}^{{0\alpha _1}}}}{{\partial {y_n}}} + {F_i^ * }]\frac{{{\partial ^2}{u_{{\alpha _1}0}}}}{{\partial {t^2}}}+ [T_{im}^{{\alpha _1}{\alpha _2}}\\
		&+ T_{ip}^{qn}\frac{{{\partial ^2}N_{pm}^{0{\alpha _1}{\alpha _2}}}}{{\partial {y_q}\partial {y_n}}} + T_{ip}^{{\alpha _1}n}\frac{{\partial N_{pm}^{0{\alpha _2}}}}{{\partial {y_n}}}+ L_{ip}^n\frac{{\partial N_{pm}^{0{\alpha _1}{\alpha _2}}}}{{\partial {y_n}}} + L_{ip}^{{\alpha _1}}N_{pm}^{0{\alpha _2}}\\
		&+ T_{ip}^{n{\alpha _1}}\frac{{\partial N_{pm}^{0{\alpha _2}}}}{{\partial {y_n}}}]\frac{{{\partial ^2}{u_{m0}}}}{{\partial {x_{{\alpha _1}}}\partial {x_{{\alpha _2}}}}}+ [U_{im}^{{\alpha _1}}+ T_{ip}^{qn}\frac{{{\partial ^2}Z_{pm}^{0{\alpha _1}}}}{{\partial {y_q}\partial {y_n}}} + L_{ip}^n\frac{{\partial Z_{pm}^{0{\alpha _1}}}}{{\partial {y_n}}} \\
		&+ T_{ip}^{qn}\frac{{{\partial ^2}N_{pm}^{0{\alpha _1}}}}{{\partial {y_q}\partial {x_n}}}+ L_{ip}^n\frac{{\partial N_{pm}^{0{\alpha _1}}}}{{\partial {x_n}}} + T_{ip}^{qn}\frac{{{\partial ^2}N_{pm}^{0{\alpha _1}}}}{{\partial {x_q}\partial {y_n}}}+ U_{ip}^n\frac{{\partial N_{pm}^{0{\alpha _1}}}}{{\partial {y_n}}}]\frac{{\partial {u_{m0}}}}{{\partial {x_{{\alpha _1}}}}}\\
        &- [Q_k^ *+ T_{im}^{qn}\frac{{{\partial ^2}{Q_{m}^{0}}}}{{\partial {y_q}\partial {y_n}}} + L_{im}^n\frac{{\partial {Q_{m}^{0}}}}{{\partial {y_n}}}+ T_{im}^{qn}\frac{{{\partial ^2}P_m^0}}{{\partial {x_q}\partial {y_n}}}\\
		&+ U_{im}^n\frac{{\partial P_m^0}}{{\partial {y_n}}} + T_{im}^{qn}\frac{{{\partial ^2}P_m^0}}{{\partial {y_q}\partial {x_n}}}+ L_{im}^n\frac{{\partial P_m^0}}{{\partial {x_n}}}]({\theta _0} - \tilde \theta )\\
		& - [V_i^{{\alpha _1}}+ T_{im}^{qn}\frac{{{\partial ^2}H_{m}^{{0\alpha _1}}}}{{\partial {y_q}\partial {y_n}}}+ L_{im}^n\frac{{\partial H_{m}^{{0\alpha _1}}}}{{\partial {y_n}}}+ T_{im}^{{\alpha _1}n}\frac{{\partial P_m^0}}{{\partial {y_n}}} \\
		& + T_{im}^{n{\alpha _1}}\frac{{\partial P_m^0}}{{\partial {y_n}}} + L_{im}^{{\alpha _1}}P_m^0+ V_i^n\frac{{\partial M_{{\alpha _1}}^0}}{{\partial {y_n}}}+P_i^ *M_{{\alpha _1}}^0]\frac{{\partial {\theta _0}}}{{\partial {x_{{\alpha _1}}}}}\\
		&+ [T_{im}^{qn}\frac{{{\partial ^2}W_{m}^{{0\alpha _1}}}}{{\partial {y_q}\partial {y_n}}} +  L_{im}^n\frac{{\partial W_{m}^{{0\alpha _1}}}}{{\partial {y_n}}} + X_i^{{\alpha _1}}]\frac{{\partial {\theta _0}}}{{\partial {x_{{\alpha _1}}}}}({\theta _0} - \tilde \theta )\\
		&- [T_{ip}^{qn}\frac{{{\partial ^2}J_{pm}^{0{\alpha _1}{\alpha _2}}}}{{\partial {y_q}\partial {y_n}}} + L_{ip}^n\frac{{\partial J_{pm}^{0{\alpha _1}{\alpha _2}}}}{{\partial {y_n}}} + J_{im}^{ * {\alpha _1}{\alpha _2}}]\frac{{\partial {\theta _0}}}{{\partial {x_{{\alpha _1}}}}}\frac{{\partial {u_{m0}}}}{{\partial {x_{{\alpha _2}}}}}.
	\end{aligned}
\end{equation}
\end{small}

%\section{Main convergence result and its proof}
\subsection{The local error analysis of multi-scale asymptotic solutions}
Before carrying out the error estimate, it is necessary to define different kinds of multi-scale solutions for nonlinear thermo-mechanical coupling simulation of heterogeneous structures with three-scale spatial hierarchy. At first, the two-scale asymptotic solutions $\theta_D^{\zeta _1}(\bm{x},t)$ and $u_{iD}^{{\zeta  _1}}(\bm{x},t)$ with meso-micro decoupling are given as follows
\begin{equation}
	\begin{aligned}
		\theta_D^{\zeta _1}(\bm{x},t)={\theta _0}({\bm{x}},t) + {\zeta  _1}{\theta ^{(1)}}({\bm{x}},{\bm{y}},{t}) + {\zeta  _1^2}{\theta ^{(4)}}({\bm{x}},{\bm{y}},{t}),
	\end{aligned}
\end{equation}
\begin{equation}
	\begin{aligned}
		u_{iD}^{{\zeta  _1}}(\bm{x},t)= {u_{i}^{(0)}}({\bm{x}},t) + {\zeta  _1}u_i^{(1)}({\bm{x}},{\bm{y}},{t}) + {\zeta  _1^2}u_i^{(4)}({\bm{x}},{\bm{y}},{t}).
	\end{aligned}
\end{equation}
Moreover, the lower-order three-scale solutions $\theta_t^{\zeta _1\zeta _2}(\bm{x},t)$ and $u_{it}^{{\zeta  _1}{\zeta  _2}}(\bm{x},t)$ can be presented in the following
\begin{equation}
	\begin{aligned}
		\theta_t^{\zeta _1\zeta _2}(\bm{x},t)={\theta _0}({\bm{x}},t) + {\zeta  _1}{\theta ^{(1)}}({\bm{x}},{\bm{y}},{t}) + {\zeta  _2}{\theta ^{(2)}}({\bm{x}},{\bm{y}},{\bm{z}},{t}),
	\end{aligned}
\end{equation}
\begin{equation}
	\begin{aligned}
		u_{it}^{{\zeta  _1}{\zeta  _2}}(\bm{x},t)= {u_{i}^{(0)}}({\bm{x}},t) + {\zeta  _1}u_i^{(1)}({\bm{x}},{\bm{y}},{t}) + {\zeta  _2}u_i^{(2)}({\bm{x}},{\bm{y}},{\bm{z}},{t}).
	\end{aligned}
\end{equation}
Additionally, the novel higher-order three-scale solutions $\theta_T^{\zeta _1\zeta _2}(\bm{x},t)$ and $	u_{iT}^{{\zeta  _1}{\zeta  _2}}(\bm{x},t)$ have been proposed as the formulas (67) and (68).

Furthermore, the residual functions of different multi-scale solutions are defined in the following way
\begin{equation}
	\begin{aligned}
		&\theta _\Delta ^{1{\zeta  _1}{\zeta  _2}}(\bm{x},t) = {\theta ^{{\zeta  _1}{\zeta  _2}}}(\bm{x},t)\! -\! \theta _D^{{\zeta  _1}}(\bm{x},t),u_{i\Delta }^{1{\zeta  _1}{\zeta  _2}}(\bm{x},t) = u_i^{{\zeta  _1}{\zeta  _2}}(\bm{x},t)\! - \!u_{iD}^{{\zeta  _1}}(\bm{x},t),\\
		&\theta _\Delta ^{2{\zeta  _1}{\zeta  _2}}(\bm{x},t)\! = \!{\theta ^{{\zeta  _1}{\zeta  _2}}}(\bm{x},t) \!- \!\theta _t^{{\zeta  _1}{\zeta  _2}}(\bm{x},t),u_{i\Delta }^{2{\zeta  _1}{\zeta  _2}}(\bm{x},t) \!= \!u_i^{{\zeta  _1}{\zeta  _2}}(\bm{x},\!t)\!-\! u_{it}^{{\zeta  _1}{\zeta  _2}}(\bm{x},\!t),\\
		&\theta _\Delta ^{3{\zeta  _1}{\zeta  _2}}(\bm{x},t)\! =\! {\theta ^{{\zeta  _1}{\zeta  _2}}}(\bm{x},t) \!-\! \theta _T^{{\zeta  _1}{\zeta  _2}}(\bm{x},t),u_{i\Delta }^{3{\zeta  _1}{\zeta  _2}}(\bm{x},t)\! =\! u_i^{{\zeta  _1}{\zeta  _2}}(\bm{x},\!t)\! -\! u_{iT}^{{\zeta  _1}{\zeta  _2}}(\bm{x},\!t).
	\end{aligned}
\end{equation}
Next, the following residual equations of different multi-scale solutions are derived by individually substituting $\theta _\Delta ^{1{\zeta  _1}{\zeta  _2}}$, $\theta _\Delta ^{2{\zeta  _1}{\zeta  _2}}$, $\theta _\Delta ^{3{\zeta  _1}{\zeta  _2}}$ for $ {\theta ^{{\zeta  _1}{\zeta  _2}}}$ and $u_{i\Delta }^{1{\zeta  _1}{\zeta  _2}}$,  $u_{i\Delta }^{2{\zeta  _1}{\zeta  _2}}$, $u_{i\Delta }^{3{\zeta  _1}{\zeta  _2}}$ for  $u_i^{{\zeta  _1}{\zeta  _2}}$ in the three-scale nonlinear thermo-mechanical problem (1), which hold in the distribution sense
\begin{equation}
	\left\{
	\begin{aligned}
		&{\rho ^{{\zeta  _1}{\zeta  _2}}}{c^{{\zeta  _1}{\zeta  _2}}}\frac{{\partial \theta _\Delta ^{1{\zeta  _1}{\zeta  _2}}}}{{\partial t}} - \frac{\partial }{{\partial {x_i}}}(k_{ij}^{{\zeta  _1}{\zeta  _2}}\frac{{\partial \theta _\Delta ^{1{\zeta  _1}{\zeta  _2}}}}{{\partial {x_j}}}) + \vartheta _{ij}^{{\zeta  _1}{\zeta  _2}}\frac{\partial }{{\partial t}}(\frac{{\partial u_{i\Delta }^{1{\zeta  _1}{\zeta  _2}}}}{{\partial {x_j}}}) = {\rm{O}}(\zeta  _2^{ - 1}),\;\;{\rm{in}}\;\;\Omega  \times (0,T),\\
		&{\rho ^{{\zeta  _1}{\zeta  _2}}}\frac{{{\partial ^2}u_{i\Delta }^{1{\zeta  _1}{\zeta  _2}}}}{{\partial {t^2}}} - \frac{\partial }{{\partial {x_j}}}\Big[C_{ijkl}^{{\zeta  _1}{\zeta  _2}}\frac{{\partial u_{k\Delta }^{1{\zeta  _1}{\zeta  _2}}}}{{\partial {x_l}}} - \beta _{ij}^{{\zeta  _1}{\zeta  _2}}(\theta _\Delta ^{1{\zeta  _1}{\zeta  _2}} - \tilde \theta )\Big] = {\rm{O}}(\zeta  _2^{ - 1}),\;\;{\rm{in}}\;\;\Omega  \times (0,T).
	\end{aligned}  \right.
\end{equation}
\begin{equation}
	\left\{
	\begin{aligned}
		&{\rho ^{{\zeta  _1}{\zeta  _2}}}{c^{{\zeta  _1}{\zeta  _2}}}\frac{{\partial \theta _\Delta ^{2{\zeta  _1}{\zeta  _2}}}}{{\partial t}} - \frac{\partial }{{\partial {x_i}}}(k_{ij}^{{\zeta  _1}{\zeta  _2}}\frac{{\partial \theta _\Delta ^{2{\zeta  _1}{\zeta  _2}}}}{{\partial {x_j}}}) + \vartheta _{ij}^{{\zeta  _1}{\zeta  _2}}\frac{\partial }{{\partial t}}(\frac{{\partial u_{i\Delta }^{2{\zeta  _1}{\zeta  _2}}}}{{\partial {x_j}}}) = {\rm{O}}(\zeta  _1^2\zeta  _2^{ - 2}),\;\;{\rm{in}}\;\;\Omega  \times (0,T),\\
		&{\rho ^{{\zeta  _1}{\zeta  _2}}}\frac{{{\partial ^2}u_{i\Delta }^{2{\zeta  _1}{\zeta  _2}}}}{{\partial {t^2}}} - \frac{\partial }{{\partial {x_j}}}\Big[C_{ijkl}^{{\zeta  _1}{\zeta  _2}}\frac{{\partial u_{k\Delta }^{2{\zeta  _1}{\zeta  _2}}}}{{\partial {x_l}}} - \beta _{ij}^{{\zeta  _1}{\zeta  _2}}(\theta _\Delta ^{2{\zeta  _1}{\zeta  _2}} - \tilde \theta )\Big] = {\rm{O}}(\zeta  _1^2\zeta  _2^{ - 2}),\;\;{\rm{in}}\;\;\Omega  \times (0,T).
	\end{aligned}  \right.
\end{equation}
\begin{equation}
	\left\{
	\begin{aligned}
		&{\rho ^{{\zeta  _1}{\zeta  _2}}}{c^{{\zeta  _1}{\zeta  _2}}}\frac{{\partial \theta _\Delta ^{3{\zeta  _1}{\zeta  _2}}}}{{\partial t}} - \frac{\partial }{{\partial {x_i}}}(k_{ij}^{{\zeta  _1}{\zeta  _2}}\frac{{\partial \theta _\Delta ^{3{\zeta  _1}{\zeta  _2}}}}{{\partial {x_j}}})+ \vartheta _{ij}^{{\zeta  _1}{\zeta  _2}}\frac{\partial }{{\partial t}}(\frac{{\partial u_{i\Delta }^{3{\zeta  _1}{\zeta  _2}}}}{{\partial {x_j}}})\\
		&\quad\quad\quad\quad  = {\rm{O}}({\zeta  _1}) + {\rm{O}}(\zeta  _1^{ - 2}{\zeta  _2}) + {\rm{O}}(\zeta  _1^{ - 1}{\zeta  _2}),\;\;{\rm{in}}\;\;\Omega  \times (0,T),\\
		&{\rho ^{{\zeta  _1}{\zeta  _2}}}\frac{{{\partial ^2}u_{i\Delta }^{3{\zeta  _1}{\zeta  _2}}}}{{\partial {t^2}}} - \frac{\partial }{{\partial {x_j}}}\Big[C_{ijkl}^{{\zeta  _1}{\zeta  _2}}\frac{{\partial u_{k\Delta }^{3{\zeta  _1}{\zeta  _2}}}}{{\partial {x_l}}}- \beta _{ij}^{{\zeta  _1}{\zeta  _2}}(\theta _\Delta ^{3{\zeta  _1}{\zeta  _2}} - \tilde \theta )\Big] \\
		&\quad\quad\quad\quad = {\rm{O}}({\zeta  _1}) + {\rm{O}}(\zeta  _1^{ - 2}{\zeta  _2}) + {\rm{O}}(\zeta  _1^{ - 1}{\zeta  _2}),\;\;{\rm{in}}\;\;\Omega  \times (0,T).
	\end{aligned} \right.
\end{equation}

To wrap up, owing to the magnitude relationship ${\zeta  _2} \ll \zeta  _1^2 \ll {\zeta  _1} \ll 1$, it is clear that the residual terms ${\rm{O}}(\zeta  _2^{ - 1})$ and ${\rm{O}}(\zeta  _1^2\zeta  _2^{ - 2})$ in right sides of equations (88) and (89) are significantly greater than constant 1,rendering both the two-scale solutions and the lower-order three-scale solutions incapable of maintaining the local balance for precise simulation of three-scale thermo-mechanical problems. In addition, the residual terms of higher-order three-scale solutions in equation (90) are of order ${\rm{O}}({\zeta  _1}) + {\rm{O}}(\zeta  _1^{ - 2}{\zeta  _2}) + {\rm{O}}(\zeta  _1^{ - 1}{\zeta  _2})$, which remain smaller than constant 1 and converge to zero when the characteristic parameters $\zeta_1$ and $\zeta_2$ approach zero. It demonstrates that novel HOTS solutions can physically guarantee the well-balanced property of local stress and heat flux, thereby effectively capturing the oscillatory information of composite structures with multiple spatial scales at the smallest scale. This is the primary reason and driven to develop the HOTS computational model.

\section{Two-stage numerical algorithm}
In this section, an efficient two-stage numerical algorithm is developed for effectively simulating the time-dependent nonlinear temperature-dependent thermo-mechanical coupling behaviors in three-scale heterogeneous structures. The detailed algorithm procedures with off-line and on-line stages are presented as below.
\subsection{Off-line computation for microscopic and mesoscopic UC problems}
\begin{enumerate}[label=(\arabic*), leftmargin=*, align=left]
	\item Determine the geometric structure of microscopic UC $Z={(0,1)^{\mathcal{N}}}$, and partition $Z$ into finite element meshes. In addition, denote the linear conforming finite element space ${V_{{h_1}}}(Z) = \{ \nu  \in {C^0}(\bar Z):{ \nu |_{\partial Z}} = 0,{ \nu |_K} \in {P_1}(K)\}  \subset H_0^1(Z)$ for solving microscopic UC problems.
	\item Verify the material parameters of composite materials in microscopic UC $Z$, and then simulate the lower-order microscopic UC problems defined by (45), (B.1) and (B.2) on $V_{{h_1}}(Z)$ which correspond to distinct representative macroscopic temperature ${\bar \theta _s} \in [\theta_{min },\theta_{max }]$. Subsequently, evaluate the mesoscopic homogenized material parameters $\bar S({\bm{y}},{\theta _0})$, $\bar{k}_{ij}(\bm{y},{\theta_0})$, $\bar \vartheta _{ij}({\bm{y}},{\theta _0})$, $\bar \rho({\bm{y}},{\theta _0})$, ${\bar C_{ijkl}}({\bm{y}},{\theta _0})$ and ${\bar \beta _{ij}}({\bm{y}},{\theta _0})$ by making integral of (47), (B.3)-(B.7) which correspond to distinct temperature ${\bar \theta _s}$ in the macroscale.
	\item Simulate the higher-order microscopic UC problems defined by (50)-(51) and (B.8)-(B.20) associated with distinct macroscopic temperature ${\bar \theta _s} $ on $V_{{h_1}}(Z)$ by means of the identical mesh used for the lower-order microscopic UC problems.
	\item Determine the geometric structure of mesoscopic UC $Y={(0,1)^{\mathcal{N}}}$, and partition $Y$ into finite element meshes. In addition, denote the linear conforming finite element space ${V_{{h_2}}}(Y) = \{ \nu  \in {C^0}(\bar Y):{ \nu |_{\partial Y}} = 0,{ \nu |_K} \in {P_1}(K)\}  \subset H_0^1(Y)$ for solving mesoscopic UC problems.
	\item Verify the material parameters of composite materials in mesoscopic UC $Y$, and then simulate lower-order mesoscopic UC problems defined by (46), (B.21) and (B.22) on $V_{{h_2}}(Y)$ corresponding to distinct representative macroscale temperature ${\bar \theta _s} \in [\theta_{min },\theta_{max }]$. Then, simulate the macroscopic homogenized material parameters $\hat S({\theta _0})$, ${{\hat k}_{ij}}({\theta _0}) $, ${{\hat \vartheta }_{ij}}({\theta _0}) $, $\hat \rho ({\theta _0}) $, ${{\hat C}_{ijkl}}({\theta _0}) $ and ${{\hat \beta }_{ij}}({\theta _0})$ by making integral of (A.1)-(A.6) associated with distinct temperature ${\bar \theta _s} $ in the macroscale.
	\item Evaluate the higher-order mesoscopic UC problems defined by (B.23)-(B.32) associated with distinct macroscopic temperature ${\bar \theta _s} $ on $V_{{h_2}}(Y)$ via the identical mesh used for the lower-order mesoscopic UC problems.
\end{enumerate}
\subsection{On-line computation for macroscopic homogenized problem and HOTS solutions}
\begin{enumerate}[label=(\arabic*), leftmargin=*, align=left]
	\item Determine the geometric structure of macroscopic homogenized region $\Omega$ in $\mathbb{R}^\mathcal{N}$, and generate the triangular finite element mesh in $\mathbb{R}^{2}$ or tetrahedral mesh in $\mathbb{R}^{3}$. Next, denote the linear conforming finite element spaces ${V_{{h_0}}}(\Omega) = \{ \nu  \in {C^0}(\bar \Omega):{ \nu |_{\partial \Omega_{\theta}}} = 0,{ \nu |_e} \in {P_1}(e)\}  \subset H^1(\Omega)$ and ${V_{{h_0}}^*}(\Omega) = \{ \nu  \in {C^0}(\bar \Omega):{ \nu |_{\partial \Omega_{u}}} = 0,{ \nu |_e} \in {P_1}(e)\}  \subset H^1(\Omega)$.
	\item  By means of the mixed FDM-FEM method, we can solve macroscopic homogenized problem (20) in a coarse mesh and with a large time step on the domain $\Omega \times (0,T)$. Employing the equidistant time step $\Delta t= T/P$ to discretize the time-domain $(0,T)$ as $0=t_0 < t_1 <\cdots< t_P =T$ and $t_N=N\Delta t (N=0,\cdots,P)$, then we define $h^N=h(\bm{x},t_N)$ and $f_i^N=f_i(\bm{x},t_N)$. To ensure the numerical stability of the algorithm, the implicit FDM scheme and the FEM scheme are applied in temporal domain and spatial domain, respectively. Then, we here present the concrete FDM-FEM scheme for macroscopic homogenized problem (20) as below
\begin{equation}
	\left\{
	\begin{aligned}
		\int_{\Omega} \hat{S}&(\theta_0^{N+1}) \frac{\theta_0^{N+1} - \theta_0^N}{\Delta t} {\widetilde \varphi}^{h_0} d\Omega \\
		& + \int_{\Omega} \hat{k}_{ij}(\theta_0^{N+1}) \frac{\partial \theta_0^{N+1}}{\partial x_j} \frac{\partial {\widetilde \varphi}^{h_0}}{\partial x_i} d\Omega \\
		& + \int_{\Omega} \hat{\vartheta}_{ij}(\theta_0^{N+1}) \frac{1}{\Delta t} \left( \frac{\partial u_{i0}^{N+1}}{\partial x_j} - \frac{\partial u_{i0}^N}{\partial x_j} \right) {\widetilde \varphi}^{h_0} d\Omega \\
		& = \int_{\Omega} h^{N+1}{\widetilde \varphi}^{h_0} d\Omega + \int_{\partial \Omega_q} \bar{q}^{N+1} {\widetilde \varphi}^{h_0} ds, \quad \forall {\widetilde \varphi}^{h_0} \in V_{h_0}(\Omega),\\
      	\int_{\Omega} \hat{\rho}&(\theta_0^{N+1}) \frac{u_{i0}^{N+1} - 2u_{i0}^N + u_{i0}^{N-1}}{(\Delta t)^2} v_i^{h_0} d\Omega \\
    	& + \int_{\Omega} \hat{C}_{ijkl}(\theta_0^{N+1}) \frac{\partial u_{k0}^{N+1}}{\partial x_l} \frac{\partial v_i^{h_0}}{\partial x_j} d\Omega \\
   	  	& - \int_{\Omega} \hat{\beta}_{ij}(\theta_0^{N+1})(\theta_0^{N+1} - \widetilde{\theta}) \frac{\partial v_i^{h_0}}{\partial x_j} d\Omega \\
     	& = \int_{\Omega} f_i^{N+1} v_i^{h_0} d\Omega + \int_{\partial \Omega_\sigma} \bar{\sigma}_i^{N+1} v_i^{h_0} ds, \quad \forall {\bf v }^{h_0} \in (V_{h_0}^*(\Omega))^{\mathcal{N}},\\
     	\theta_0^N &= \widehat{\theta}(\bm{x}, t_N), \quad \text{on } \partial \Omega_\theta,\\
     	{\bm u}_{0}^N &= \widehat{\bm u}(\bm{x}, t_N), \quad \text{on } \partial \Omega_u.
	\end{aligned} \right.
\end{equation}
Since the nonlinear coupled system (91) failed to simulate directly, we can employ the decoupling approach to decompose the system into two sub-problems, which can maintain the accuracy and stability \cite{R33}. First we set that
\begin{equation}
	u_{i0}^{N+1}=u_{i0}^{N}+\varpi (u_{i0}^{N}-u_{i0}^{N-1}),
\end{equation}
and $\varpi$ is a correction coefficient. Then we can derive two sub-systems by decomposing the above system (91).
\begin{equation}
	\left\{
	\begin{aligned}
		\int_{\Omega} \hat{S}&(\theta_0^{N+1}) \frac{\theta_0^{N+1} - \theta_0^N}{\Delta t} {\widetilde \varphi}^{h_0} d\Omega \\
		& + \int_{\Omega} \hat{k}_{ij}(\theta_0^{N+1}) \frac{\partial \theta_0^{N+1}}{\partial x_j} \frac{\partial {\widetilde \varphi}^{h_0}}{\partial x_i} d\Omega \\
		& + \int_{\Omega} \varpi \hat{\vartheta}_{ij}(\theta_0^{N+1}) \frac{1}{\Delta t} \left( \frac{\partial u_{i0}^{N}}{\partial x_j} - \frac{\partial u_{i0}^{N-1}}{\partial x_j} \right) {\widetilde \varphi}^{h_0} d\Omega \\
		& = \int_{\Omega} h^{N+1} {\widetilde \varphi}^{h_0} d\Omega + \int_{\partial \Omega_q} \bar{q}^{N+1} {\widetilde \varphi}^{h_0} ds, \quad \forall {\widetilde \varphi}^{h_0} \in V_{h_0}(\Omega),\\
		\theta_0^N &= \widehat{\theta}(\bm{x}, t_N), \quad \text{on } \partial \Omega_\theta.
	\end{aligned} \right.
\end{equation}
\begin{equation}
	\left\{
	\begin{aligned}
		\int_{\Omega} \hat{\rho}&(\theta_0^{N+1}) \frac{u_{i0}^{N+1} - 2u_{i0}^N + u_{i0}^{N-1}}{(\Delta t)^2} v_i^{h_0} d\Omega \\
		& + \int_{\Omega} \hat{C}_{ijkl}(\theta_0^{N+1}) \frac{\partial u_{k0}^{N+1}}{\partial x_l} \frac{\partial v_i^{h_0}}{\partial x_j} d\Omega \\
		& - \int_{\Omega} \hat{\beta}_{ij}(\theta_0^{N+1})(\theta_0^{N+1} - \widetilde{\theta}) \frac{\partial v_i^{h_0}}{\partial x_j} d\Omega \\
		& = \int_{\Omega} f_i^{N+1} v_i^{h_0} d\Omega + \int_{\partial \Omega_\sigma} \bar{\sigma}_i^{N+1} v_i^{h_0} ds, \quad \forall {\bf v}^{h_0} \in (V_{h_0}^*(\Omega))^{\mathcal{N}},\\
		{\bm u}_{0}^N &= \widehat{\bm u}(\bm{x}, t_N), \quad \text{on } \partial \Omega_u.
	\end{aligned} \right.
\end{equation}
	\item Next, the  fixed-point iterative method is introduced for evaluating the sub-systems (93) and (94).
	\begin{enumerate}[label=\textbf{Step \arabic*:}, leftmargin=*, align=left]
		\item Denote $\tilde{\theta}_\eta (\bm{x})$ and $\tilde{\bm u}_\eta(\bm{x})$ as the temperature and displacement fields at $\eta$-th iterative step, and $\tilde{\theta}_0(\bm{x})$ as the initial function. Let the iteration threshold be $E_{tol}$ and $E_{tol}^*$, respectively, and start to iterate.
		\item Employ $\tilde{\theta}_{\eta - 1}(\bm{x})$ to linearize the sub-systems (93) and (94) at the $\eta$-th iteration step.
\begin{equation}
	\left\{
	\begin{aligned}
		\int_{\Omega} \hat{S}&(\tilde{\theta}_{\eta - 1}) \frac{\tilde{\theta}_\eta - \theta_0^N}{\Delta t} {\widetilde \varphi}^{h_0} d\Omega \\
		& + \int_{\Omega} \hat{k}_{ij}(\tilde{\theta}_{\eta - 1}) \frac{\partial \tilde{\theta}_\eta}{\partial x_j} \frac{\partial {\widetilde \varphi}^{h_0}}{\partial x_i} d\Omega \\
		& + \int_{\Omega} \varpi \hat{\vartheta}_{ij}(\tilde{\theta}_{\eta - 1}) \frac{1}{\Delta t} \left( \frac{\partial u_{i0}^{N}}{\partial x_j} - \frac{\partial u_{i0}^{N-1}}{\partial x_j} \right) {\widetilde \varphi}^{h_0} d\Omega \\
		& = \int_{\Omega} h^{N+1} {\widetilde \varphi}^{h_0} d\Omega + \int_{\partial \Omega_q} \bar{q}^{N+1} {\widetilde \varphi}^{h_0} ds, \quad \forall {\widetilde \varphi}^{h_0} \in V_{h_0}(\Omega),\\
		\theta_0^N &= \widehat{\theta}(\bm{x}, t_N), \quad \text{on } \partial \Omega_\theta.
	\end{aligned} \right.
\end{equation}
\begin{equation}
	\left\{
	\begin{aligned}
		\int_{\Omega} \hat{\rho}&(\tilde{\theta}_{\eta - 1}) \frac{\tilde{u}_{\eta,i} - 2u_{i0}^N + u_{i0}^{N-1}}{(\Delta t)^2} v_i^{h_0} d\Omega \\
		& + \int_{\Omega} \hat{C}_{ijkl}(\tilde{\theta}_{\eta - 1}) \frac{\partial \tilde{u}_{\eta,k}}{\partial x_l} \frac{\partial v_i^{h_0}}{\partial x_j} d\Omega \\
		& - \int_{\Omega} \hat{\beta}_{ij}(\tilde{\theta}_{\eta - 1})(\tilde{\theta}_{\eta} - \widetilde{\theta}) \frac{\partial v_i^{h_0}}{\partial x_j} d\Omega \\
		& = \int_{\Omega} f_i^{N+1} v_i^{h_0} d\Omega + \int_{\partial \Omega_\sigma} \bar{\sigma}_i^{N+1} v_i^{h_0} ds, \quad \forall {\bf v}^{h_0} \in (V_{h_0}^*(\Omega))^N,\\
		{\bm u}_{0}^N &= \widehat{\bm u}(\bm{x}, t_N), \quad \text{on } \partial \Omega_u.
	\end{aligned} \right.
\end{equation}
	\item  If  $||\tilde{\theta}_\eta(\bm{x}) - \tilde{\theta}_{\eta-1}(\bm{x})||_{L^\infty(\Omega)} \leq E_{tol}$ and $||\tilde{\bm u}_\eta(\bm{x}) - \tilde{\bm u}_{\eta-1}(\bm{x})||_{L^\infty(\Omega)} \leq E_{tol}^*$, stop; otherwise $ \eta = \eta + 1$,  return to Step 2.
	\item Let  $\theta_0^{N+1} = \theta_{sat}$ and ${\bm u}_0^{N+1} = {\bm u}_{sat}$,  where  $\theta_{sat}$ and ${\bm{u}}_{sat}$ are the solutions of the above systems (95) and (96) that meet the iteration threshold as $ E_{tol}$ and $E_{tol}^*$, respectively.
	\end{enumerate}
	\item For any point $(\bm{x},t) \in \Omega \times (0,T)$, we employed the interpolation technique for computing the values of mesoscopic cell functions, microscopic cell functions alongside macroscopic homogenized solutions.
	\item In formulas (67) and (68), the spatial derivatives $ \displaystyle\frac{\partial \theta_0}{\partial x_{\alpha_1}} $, $ \displaystyle \frac{\partial^2 \theta_0}{\partial x_{\alpha_1} \partial x_{\alpha_2}}$, $ \displaystyle\frac{\partial u_{m0}}{\partial x_{\alpha_1}} $ and $ \displaystyle \frac{\partial^2 u_{m0}}{\partial x_{\alpha_1} \partial x_{\alpha_2}}$ are simulated via the average technique on interdependent elements \cite{R48,R49,R50}, and the temporal derivatives $ \displaystyle\frac{\partial \theta_0}{\partial t} $, $\displaystyle \frac{{{\partial ^2}u_{i0}}}{{\partial {t^2}}}$ and $\displaystyle \frac{{{\partial ^2}{u_{{\alpha _1}0}}}}{{\partial {x_{{\alpha _2}}}\partial t}}$ are simulated via the difference schemes (93) and (94) at every time step.
	\item Finally, we can evaluate the temperature field $\theta^{{\zeta  _1}{\zeta  _2}}(\bm{x},t)$ and displacement field ${\bm u}^{{\zeta  _1}{\zeta  _2}}(\bm{x},t)$ according to the formulas (67) and (68).
\end{enumerate}

\section{Numerical experiments}
To validate the computational cost, computational efficiency and numerical accuracy of the HOTS approach, this section presents four numerical experiments. All numerical examples are conducted using the FreeFem++ software on a HP desktop workstation equipped with 16.0 GB RAM and a 2.20 GHz Intel Core i7-8750 processor. Given that it is nearly impossible to simulate the exact solutions for nonlinear three-scale multi-physics problem (1), the precise FEM solutions in a highly refined mesh, which we denote as $\theta_{Ref}^{\zeta _1\zeta _2}$ and ${\bm u}_{Ref}^{{\zeta  _1}{\zeta  _2}}$, served as the reference solutions. Moreover, some error notations without confusion are expressed as below.
\begin{equation}
	\begin{aligned}
		&{\rm{T_{{L^2}}^0}}(t) = \frac{{||\theta _{Ref}^{{\zeta  _1}{\zeta  _2}} - {\theta _0}|{|_{{L^2}(\Omega )}}}}{{||\theta _{Ref}^{{\zeta  _1}{\zeta  _2}}|{|_{{L^2}(\Omega )}}}},\;{\rm{T_{{L^2}}^1}}(t) = \frac{{||\theta _{Ref}^{{\zeta  _1}{\zeta  _2}} - \theta _D^{{\zeta  _1}}|{|_{{L^2}(\Omega )}}}}{{||\theta _{Ref}^{{\zeta  _1}{\zeta  _2}}|{|_{{L^2}(\Omega )}}}},\\
		&{\rm{T_{{L^2}}^2}}(t) = \frac{{||\theta _{Ref}^{{\zeta  _1}{\zeta  _2}} - \theta _t^{{\zeta  _1}{\zeta  _2}}|{|_{{L^2}(\Omega )}}}}{{||\theta _{Ref}^{{\zeta  _1}{\zeta  _2}}|{|_{{L^2}(\Omega )}}}},\;{\rm{T_{{L^2}}^3}}(t) = \frac{{||\theta _{Ref}^{{\zeta  _1}{\zeta  _2}} - \theta _T^{{\zeta  _1}{\zeta  _2}}|{|_{{L^2}(\Omega )}}}}{{||\theta _{Ref}^{{\zeta  _1}{\zeta  _2}}|{|_{{L^2}(\Omega )}}}},
	\end{aligned}
\end{equation}
\begin{equation}
\begin{aligned}
		&{\rm{T_{{H^1}}^0}}(t) = \frac{{|\theta _{Ref}^{{\zeta  _1}{\zeta  _2}} - {\theta _0}{|_{{H^1}(\Omega )}}}}{{|\theta _{Ref}^{{\zeta  _1}{\zeta  _2}}{|_{{H^1}(\Omega )}}}},\;{\rm{T_{{H^1}}^1}}(t) = \frac{{|\theta _{Ref}^{{\zeta  _1}{\zeta  _2}} - \theta _D^{{\zeta  _1}}{|_{{H^1}(\Omega )}}}}{{|\theta _{Ref}^{{\zeta  _1}{\zeta  _2}}{|_{{H^1}(\Omega )}}}},\\
		&{\rm{T_{{H^1}}^2}}(t) = \frac{{|\theta _{Ref}^{{\zeta  _1}{\zeta  _2}} - \theta _t^{{\zeta  _1}{\zeta  _2}}{|_{{H^1}(\Omega )}}}}{{|\theta _{Ref}^{{\zeta  _1}{\zeta  _2}}{|_{{H^1}(\Omega )}}}},\;{\rm{T_{{H^1}}^3}}(t) = \frac{{|\theta _{Ref}^{{\zeta  _1}{\zeta  _2}} - \theta _T^{{\zeta  _1}{\zeta  _2}}{|_{{H^1}(\Omega )}}}}{{|\theta _{Ref}^{{\zeta  _1}{\zeta  _2}}{|_{{H^1}(\Omega )}}}},
\end{aligned}
\end{equation}
\begin{equation}
\begin{aligned}
		&{\rm{{D}_{{L^2}}^0}}(t) = \frac{{||{\bm u}_{Ref}^{{\zeta  _1}{\zeta  _2}} - {\bm u_{0}}|{|_{({L^2}(\Omega ))^{\mathcal{N}}}}}}{{||{\bm u}_{Ref}^{{\zeta  _1}{\zeta  _2}}|{|_{({L^2}(\Omega ))^{\mathcal{N}}}}}},\;{\rm{{D}_{{L^2}}^1}}(t) = \frac{{||{\bm u}_{Ref}^{{\zeta  _1}{\zeta  _2}} - {\bm u}_D^{{\zeta  _1}}|{|_{({L^2}(\Omega ))^{\mathcal{N}}}}}}{{||{\bm u}_{Ref}^{{\zeta  _1}{\zeta  _2}}|{|_{({L^2}(\Omega ))^{\mathcal{N}}}}}},\\
		&{\rm{{D}_{{L^2}}^2}}(t) = \frac{{||{\bm u}_{Ref}^{{\zeta  _1}{\zeta  _2}} - {\bm u}_t^{{\zeta  _1}{\zeta  _2}}|{|_{({L^2}(\Omega ))^{\mathcal{N}}}}}}{{||{\bm u}_{Ref}^{{\zeta  _1}{\zeta  _2}}|{|_{({L^2}(\Omega ))^{\mathcal{N}}}}}},\;{\rm{{D}_{{L^2}}^3}}(t) = \frac{{||{\bm u}_{Ref}^{{\zeta  _1}{\zeta  _2}} - {\bm u}_T^{{\zeta  _1}{\zeta  _2}}|{|_{({L^2}(\Omega ))^{\mathcal{N}}}}}}{{||{\bm u}_{Ref}^{{\zeta  _1}{\zeta  _2}}|{|_{({L^2}(\Omega ))^{\mathcal{N}}}}}},
\end{aligned}
\end{equation}
\begin{equation}
\begin{aligned}
	&{\rm{{D}_{{H^1}}^0}}(t) = \frac{{|{\bm u}_{Ref}^{{\zeta  _1}{\zeta  _2}} - {\bm u_{0}}{|_{({H^1}(\Omega ))^{\mathcal{N}}}}}}{{|{\bm u}_{Ref}^{{\zeta  _1}{\zeta  _2}}{|_{({H^1}(\Omega ))^{\mathcal{N}}}}}},\;{\rm{{D}_{{H^1}}^1}}(t) = \frac{{|{\bm u}_{Ref}^{{\zeta  _1}{\zeta  _2}} - {\bm u}_D^{{\zeta  _1}}{|_{({H^1}(\Omega ))^{\mathcal{N}}}}}}{{|{\bm u}_{Ref}^{{\zeta  _1}{\zeta  _2}}{|_{({H^1}(\Omega ))^{\mathcal{N}}}}}},\\
	&{\rm{{D}_{{H^1}}^2}}(t) = \frac{{|{\bm u}_{Ref}^{{\zeta  _1}{\zeta  _2}} - {\bm u}_t^{{\zeta  _1}{\zeta  _2}}{|_{({H^1}(\Omega ))^{\mathcal{N}}}}}}{{|{\bm u}_{Ref}^{{\zeta  _1}{\zeta  _2}}{|_{({H^1}(\Omega ))^{\mathcal{N}}}}}},\;{\rm{{D}_{{H^1}}^3}}(t) = \frac{{|{\bm u}_{Ref}^{{\zeta  _1}{\zeta  _2}} - {\bm u}_T^{{\zeta  _1}{\zeta  _2}}{|_{({H^1}(\Omega ))^{\mathcal{N}}}}}}{{|{\bm u}_{Ref}^{{\zeta  _1}{\zeta  _2}}{|_{({H^1}(\Omega ))^{\mathcal{N}}}}}}.
\end{aligned}
\end{equation}

\subsection{Example 1: 2D composite structure I with multiple spatial scales}
\begin{figure}[!htb]
	\centering
	\begin{minipage}[c]{0.46\textwidth}
		\centering
		\includegraphics[
		width=0.9\linewidth,
		trim=2cm 1.3cm 1.8cm 1.3cm, % 从左、下、右、上各裁剪1厘米
		clip % 应用裁剪
		]{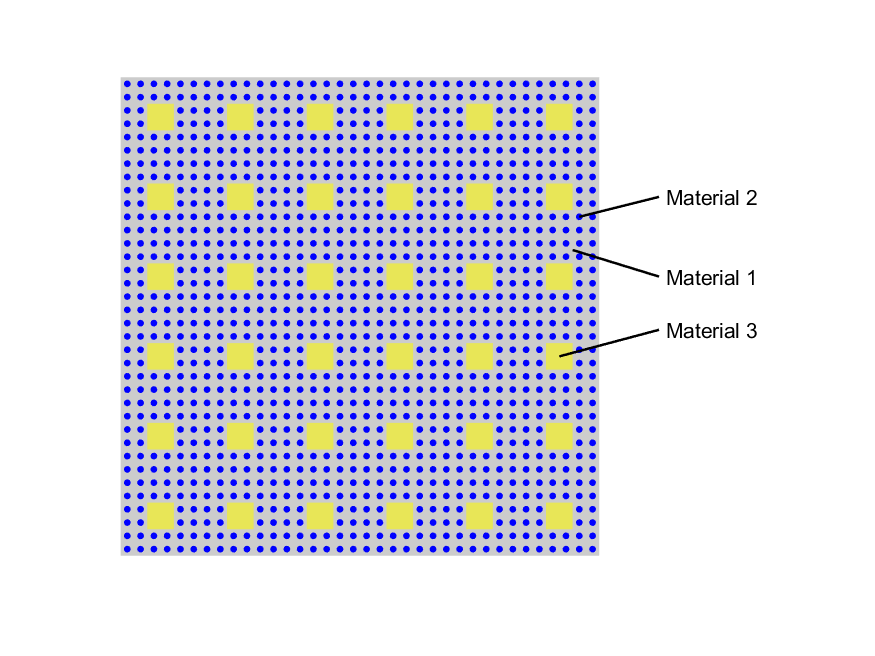} \\
		(a)
	\end{minipage}
	\begin{minipage}[c]{0.46\textwidth}
		\centering
		\includegraphics[
		width=0.9\linewidth,
		trim=2cm 1.1cm 1.5cm 1.5cm, % 从左、下、右、上各裁剪1厘米
		clip % 应用裁剪
		]{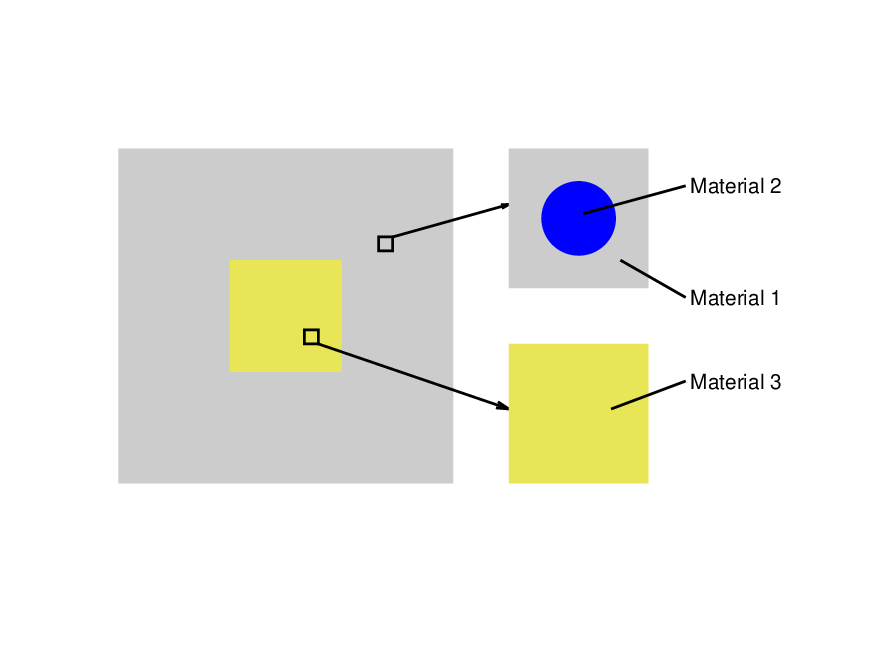} \\
		(b)
	\end{minipage}
	\caption{(a) Macroscopic domain $\Omega$ ; (b) Mesoscopic unit cell $Y$ (left side) and the microscopic unit cell $Z$ (right side).}\label{E1f1}
\end{figure}

The nonlinear thermo-mechanical behaviors of a 2D heterogeneous structure I with three-level structural scales are simulated by the presented HOTS computational approach and FEM. The investigated heterogeneous structure is exhibited in Fig.\hspace{1mm}\ref{E1f1}, where the characteristic periodic parameters $\zeta _1=1/6$ and $\zeta _2=1/36$ and $\Omega=(x_1,x_2)=[0,1]\times[0,1]\mathrm{cm^2}$.

As for this example, the boundary conditions on $\partial\Omega$ are prescribed as $\widehat{\theta}(\bm{x},t)=373.15\mathrm{K}$ and $\widehat{\bm{u}}(\bm{x},t)=0$ separately. Next, the initial conditions are prescribed as $\tilde{\theta}=373.15\mathrm{K}$, ${\bm{u}}^{0}=0$ and ${\bm{u}}^{1}(\bm{x})=0$, and then the heat source and body forces are given by $h(\bm{x},t)=10000.0\mathrm{J/(cm^2\cdot s)}$ and $f_i(\bm{x},t)=(-8000.0\mathrm{N/cm^2},-8000.0\mathrm{N/cm^2})$. In addition, the detailed temperature-dependent material parameters of the heterogeneous structure I are given in Table \ref{E1t1}.
\begin{table}[!htb]
	\centering
	\caption{Dimensionless material parameters ($\theta$ stands for temperature).}\label{E1t1}
	\begin{tabular}{ll}
		\hline
		Material Property     & Material 1       \\ \hline
		Mass density ${\rho ^{{\zeta _1}{\zeta _2}}}$ ($\mathrm{kg/m^3}$)     &4410        \\
		Specific heat ${c^{{\zeta _1}{\zeta _2}}}$ (J/(kg$\cdot$K))      & $808.3+0.081\theta+0.00008\theta^2 $       \\
		Thermal conductivity $k_{ij}^{{\zeta  _1}{\zeta  _2}}$ (W/(m$\cdot$K))   & $1000+0.1\theta+10^{-5}\theta^2$       \\
		Young's modulus $E$ (GPa)      &$3.0 \times {10^7} - 300\theta - 0.03\theta^2$    \\
		Poisson's ratio $v$         & 0.30              \\
		Thermal modulus $\beta _{ij}^{{\zeta  _1}{\zeta  _2}}$ (MPa/K)     &$19-1.9\times10^{-3}\theta-1.9\times10^{-7}\theta^2$     \\ \hline
		Material Property        & Material 2      \\ \hline
		Mass density ${\rho ^{{\zeta _1}{\zeta _2}}}$ ($\mathrm{kg/m^3}$)      &5600                 \\
		Specific heat ${c^{{\zeta _1}{\zeta _2}}}$ (J/(kg$\cdot$K))        & $615.6+0.062\theta+0.00006\theta^2 $       \\
		Thermal conductivity $k_{ij}^{{\zeta  _1}{\zeta  _2}}$  (W/(m$\cdot$K))   &$1+10^{-4}\theta$+$10^{-8}\theta^2$       \\
		Young's modulus $E$ (GPa)    &$6.0 \times 10^{6}-60\theta-0.006\theta^2$       \\
		Poisson's ratio $v$         & 0.20                \\
		Thermal modulus $\beta _{ij}^{{\zeta  _1}{\zeta  _2}}$ (MPa/K)   &$17-1.7\times10^{-3}\theta-1.7\times10^{-7}\theta^2$      \\ \hline
		Material Property          &Material 3    \\ \hline
		Mass density ${\rho ^{{\zeta _1}{\zeta _2}}}$ ($\mathrm{kg/m^3}$)     &5800       \\
		Specific heat ${c^{{\zeta _1}{\zeta _2}}}$ (J/(kg$\cdot$K))      & $590.9+0.059\theta+0.00006\theta^2 $       \\
		Thermal conductivity $k_{ij}^{{\zeta  _1}{\zeta  _2}}$ (W/(m$\cdot$K))   &$200+0.02\theta+2\times10^{-6}\theta^2$      \\
		Young's modulus $E$ (GPa)  &$2.5 \times 10^{7}-250\theta-0.025\theta^2$             \\
		Poisson's ratio $v$         & 0.25            \\
		Thermal modulus $\beta _{ij}^{{\zeta  _1}{\zeta  _2}}$ (MPa/K)  &$18-1.8\times10^{-3}\theta-1.8\times10^{-7}\theta^2$              \\ \hline
	\end{tabular}
\end{table}

By conducting the reference FEM and HOTS method, the computational cost can be presented as shown in Table \ref{E1t2}, including the numbers of nodes and elements, as well as the computational time.
\begin{table}[!htb]
	\centering
	\caption{Summary of computational cost ($\Delta t = 0.001, t \in [0,1] $).}\label{E1t2}
	\begin{tabular}{lcccc}
		\hline
		 & \multirow{2}{*}{Reference FEM} & \multicolumn{3}{c}{HOTS method} \\
		\cmidrule(l){3-5}
		 & & $Z$ & $Y$ & $\Omega$ \\
		\hline
		Number of nodes & 43,633 & 1,179 & 2,056 & 961 \\
		Number of elements & 86,400 & 3,396 & 3,942 & 1,800 \\
		Computational time & 33,339.6s &\multicolumn{3}{c}{22,127.8s} \\
		\hline
	\end{tabular}
\end{table}

\begin{figure}[!htb]
	\centering
	\begin{minipage}[b]{0.35\textwidth}
		\centering
		\includegraphics[
		width=4cm,
		trim=2.5cm 4.5cm 0cm 4.5cm, % 从左、下、右、上各裁剪1厘米
		clip % 应用裁剪
		]{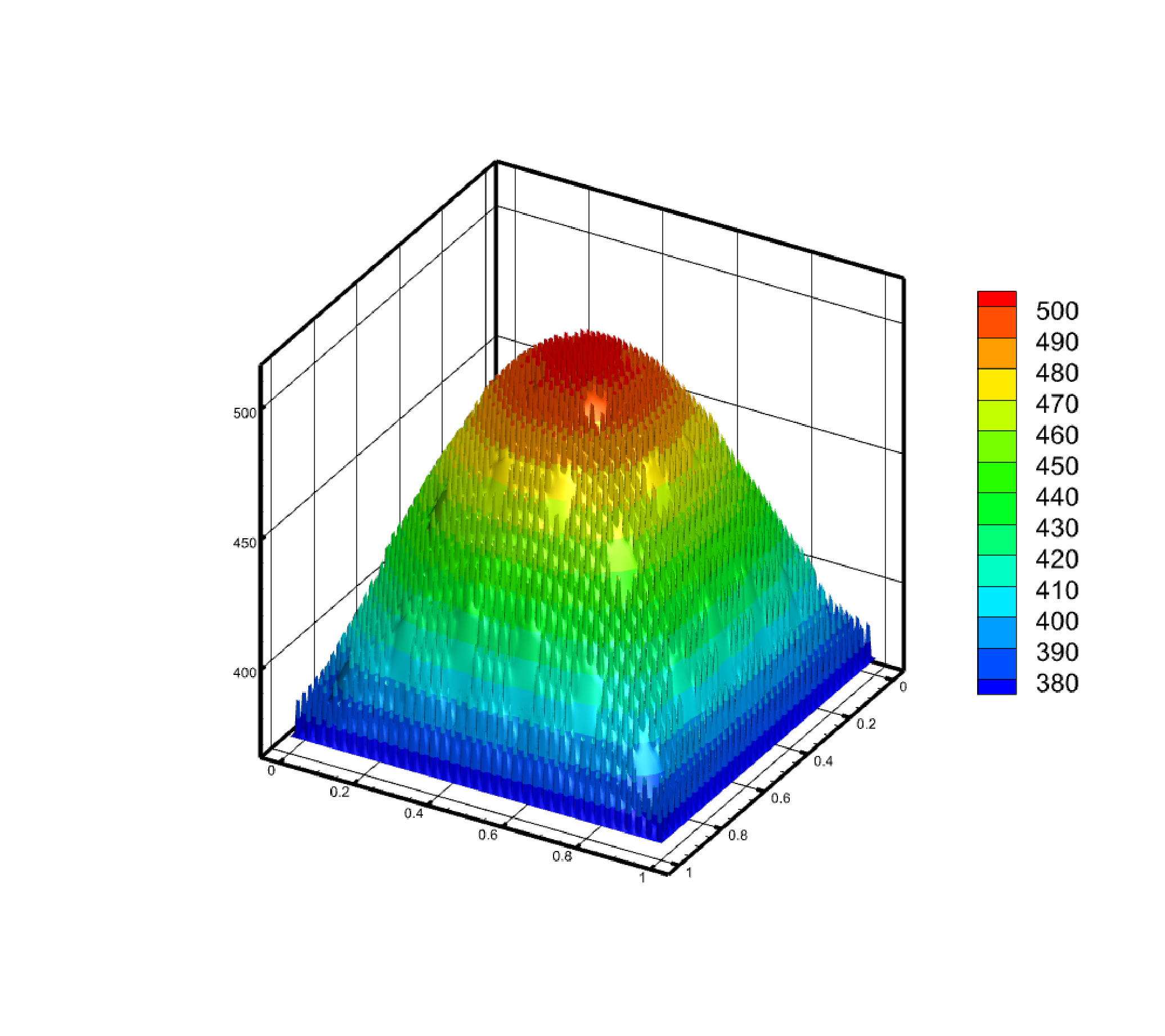} \\
		(a)
	\end{minipage}
	\begin{minipage}[b]{0.35\textwidth}
		\centering
		\includegraphics[
		width=4cm,
		trim=2.5cm 4.5cm 0cm 4.5cm, % 从左、下、右、上各裁剪1厘米
		clip % 应用裁剪
		]{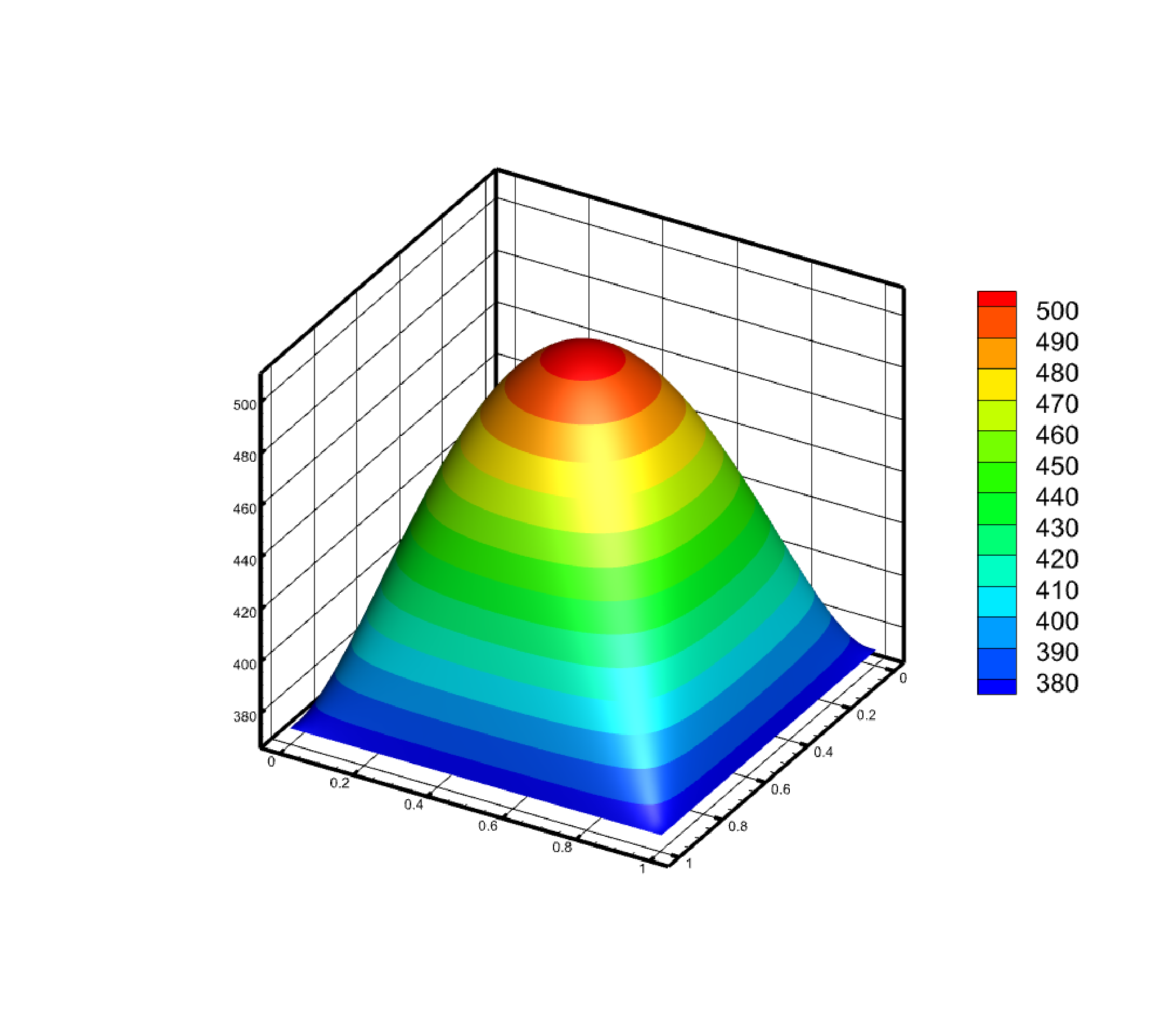} \\
		(b)
	\end{minipage}
	\begin{minipage}[b]{0.32\textwidth}
		\centering
		\includegraphics[
		width=4cm,
		trim=2.5cm 4.5cm 0cm 4.5cm, % 从左、下、右、上各裁剪1厘米
		clip % 应用裁剪
		]{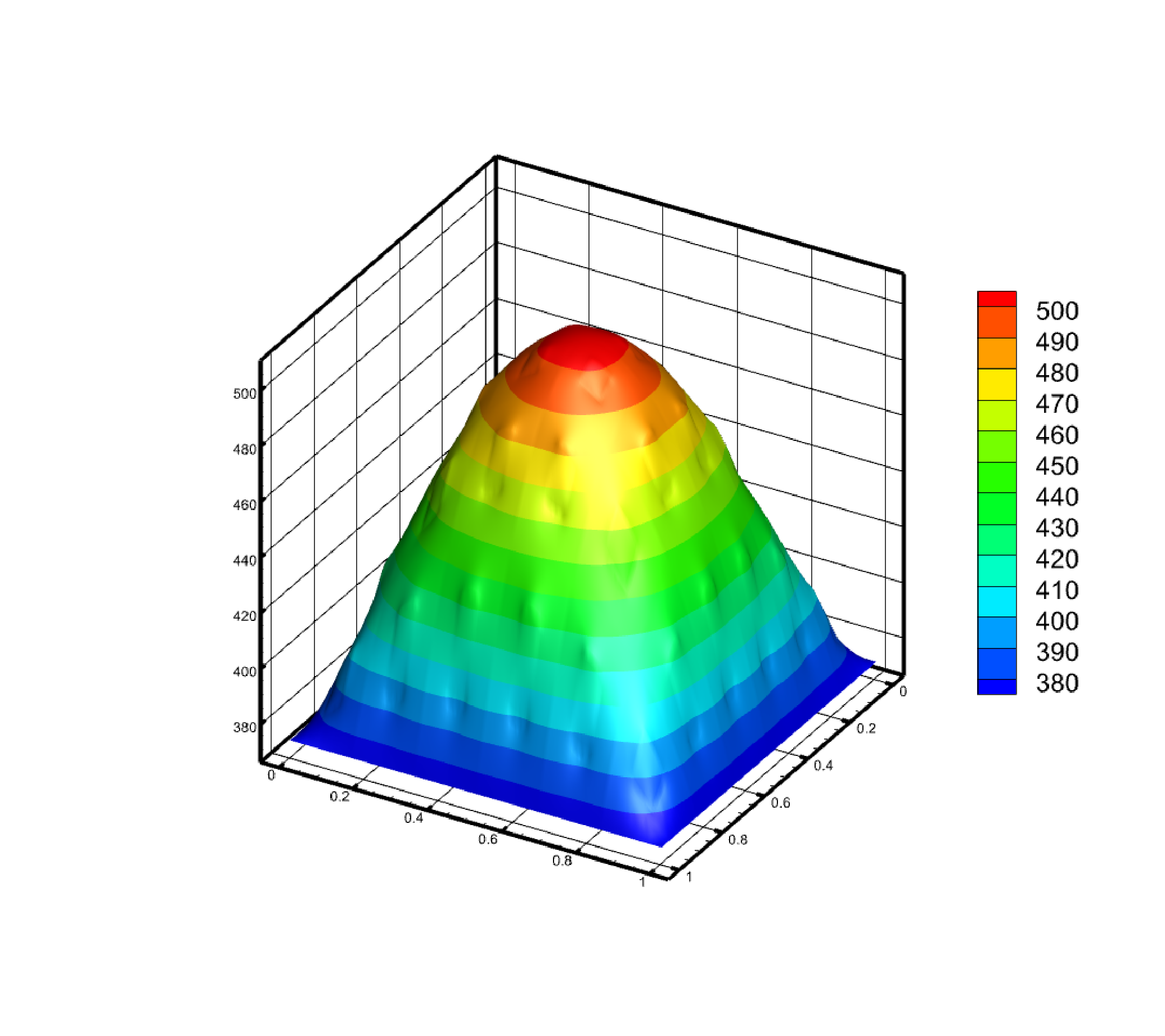} \\
		(c)
	\end{minipage}
	\begin{minipage}[b]{0.32\textwidth}
		\centering
		\includegraphics[
		width=4cm,
		trim=2.5cm 4.5cm 0cm 4.5cm, % 从左、下、右、上各裁剪1厘米
		clip % 应用裁剪
		]{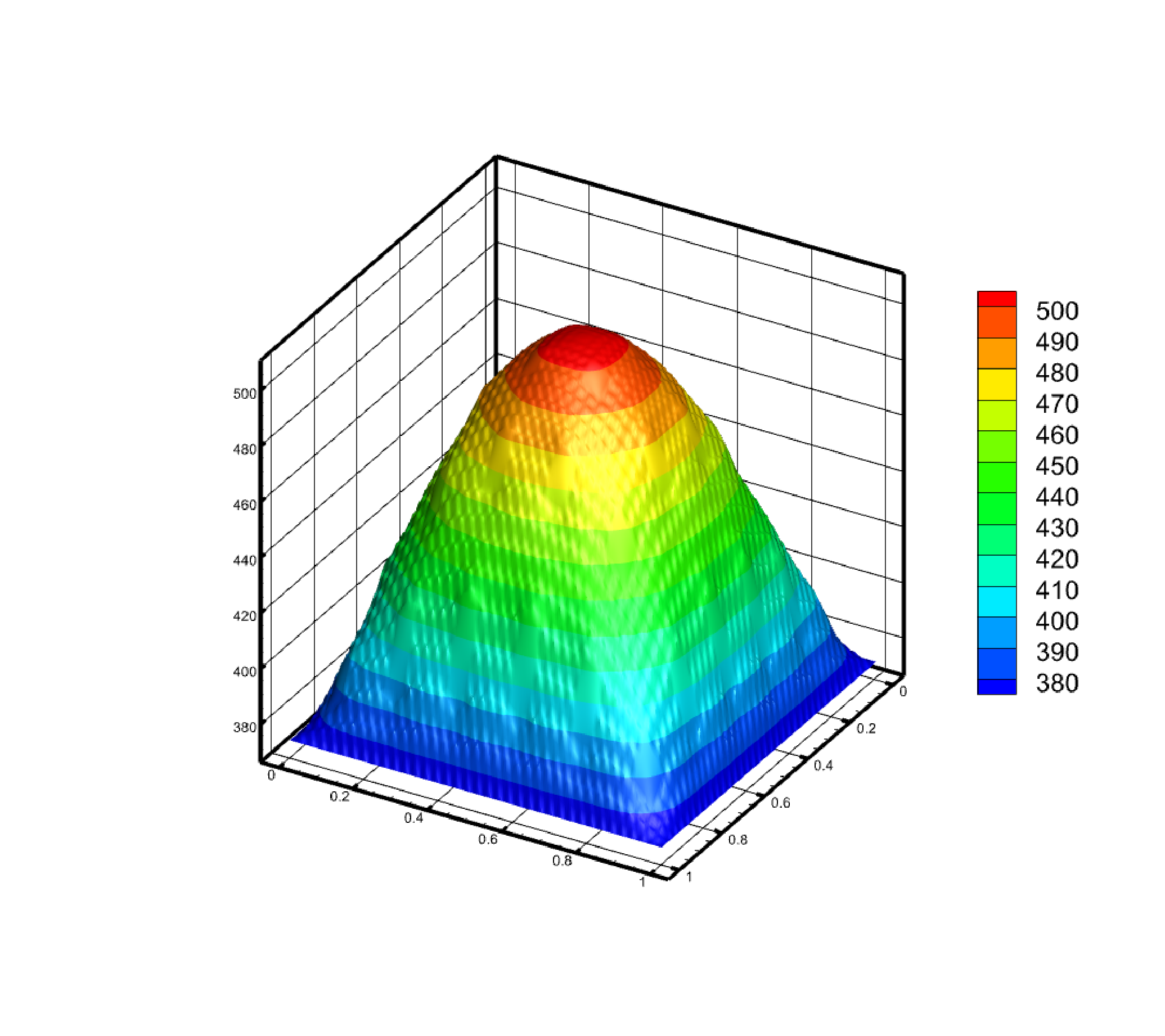} \\
		(d)
	\end{minipage}
	\begin{minipage}[b]{0.32\textwidth}
		\centering
		\includegraphics[
		width=4cm,
		trim=2.5cm 4.5cm 0cm 4.5cm, % 从左、下、右、上各裁剪1厘米
		clip % 应用裁剪
		]{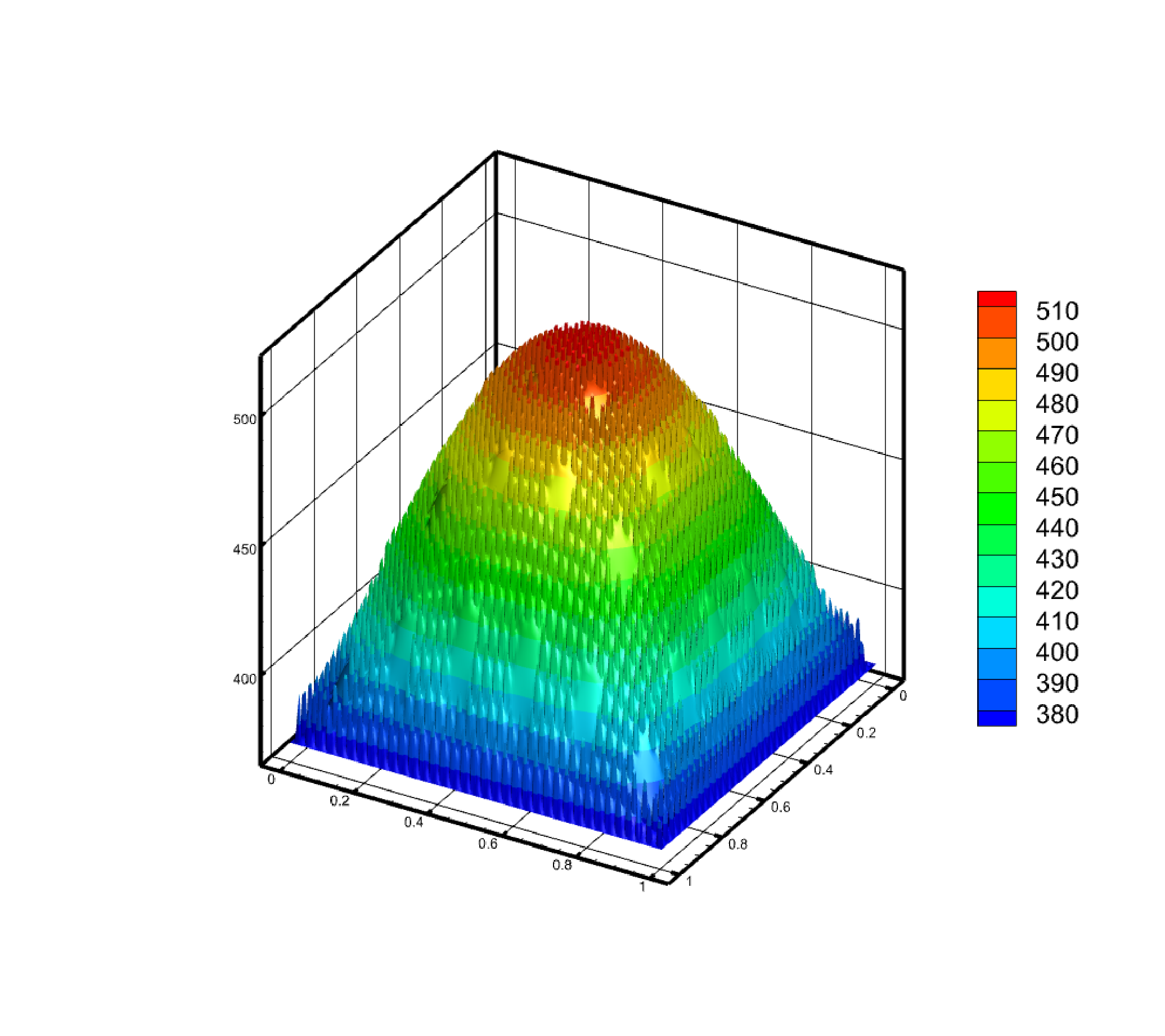} \\
		(e)
	\end{minipage}
	\caption{The temperature field: (a) $\theta_{Ref}^{\zeta _1\zeta _2}$; (b) $\theta_{0}$; (c) ${\theta _D^{{\zeta  _1}}}$; (d) ${\theta _t^{{\zeta  _1}{\zeta  _2}}}$; (e) ${\theta _T^{{\zeta  _1}{\zeta  _2}}}$.}\label{E1f2}
\end{figure}

\begin{figure}[!htb]
	\centering
	\begin{minipage}[b]{0.35\textwidth}
		\centering
		\includegraphics[
		width=4cm,
		trim=2.5cm 6.5cm 0cm 6.5cm, % 从左、下、右、上各裁剪1厘米
		clip % 应用裁剪
		]{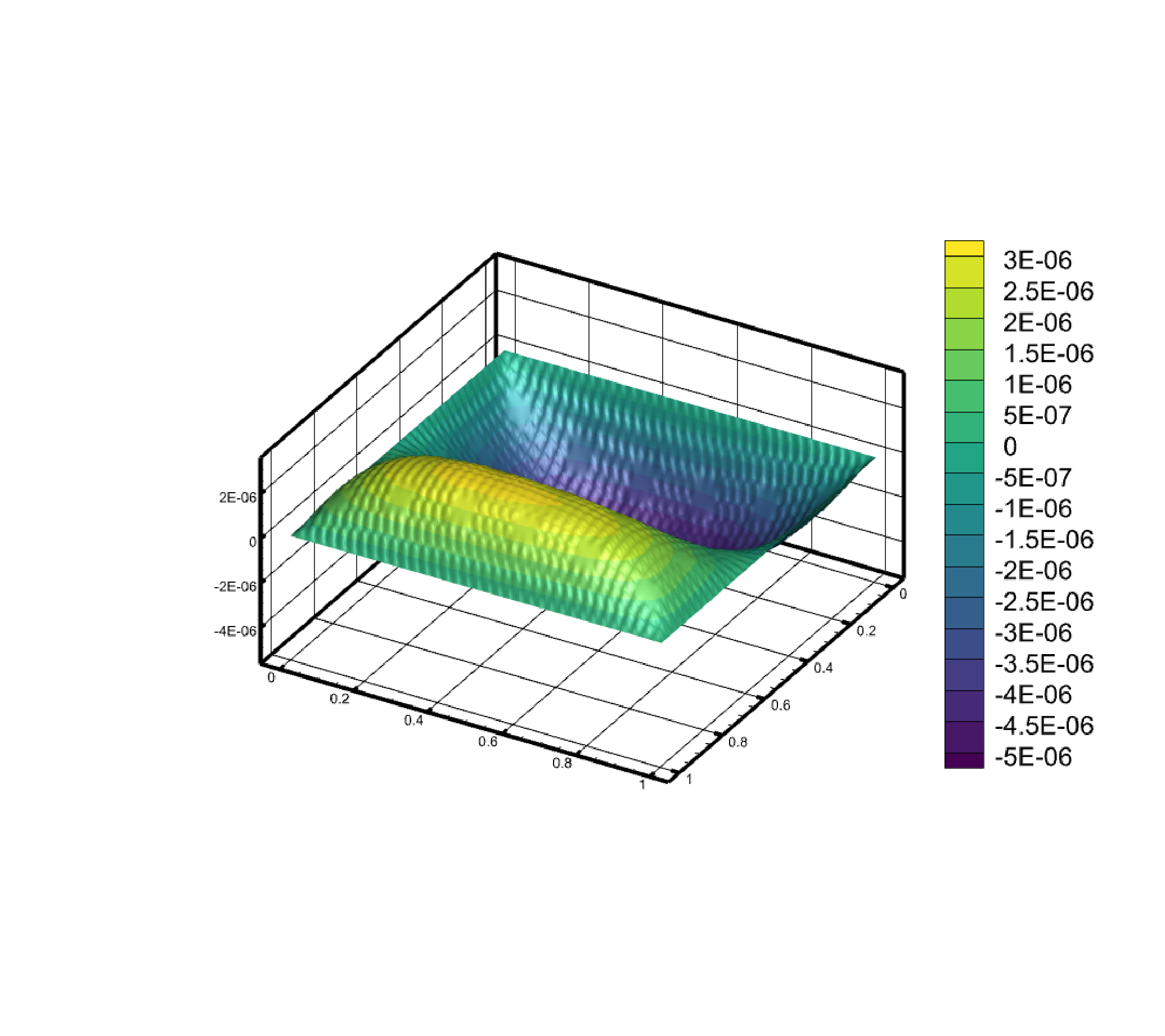} \\
		(a)
	\end{minipage}
	\begin{minipage}[b]{0.35\textwidth}
		\centering
		\includegraphics[
		width=4cm,
		trim=2.5cm 6.5cm 0cm 6.5cm, % 从左、下、右、上各裁剪1厘米
		clip % 应用裁剪
		]{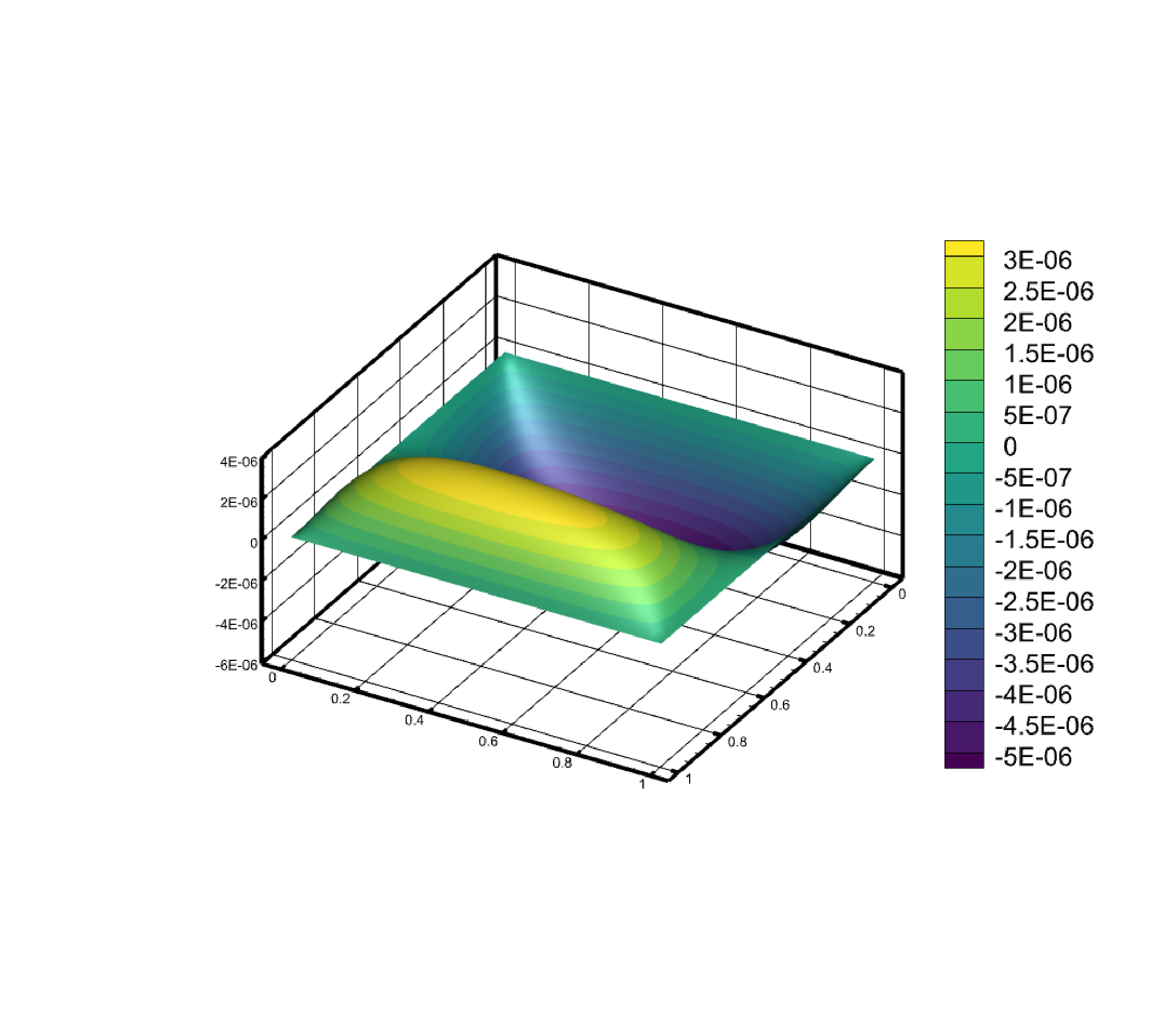} \\
		(b)
	\end{minipage}
	\begin{minipage}[b]{0.32\textwidth}
		\centering
		\includegraphics[
		width=4cm,
		trim=2.5cm 6.5cm 0cm 6.5cm, % 从左、下、右、上各裁剪1厘米
		clip % 应用裁剪
		]{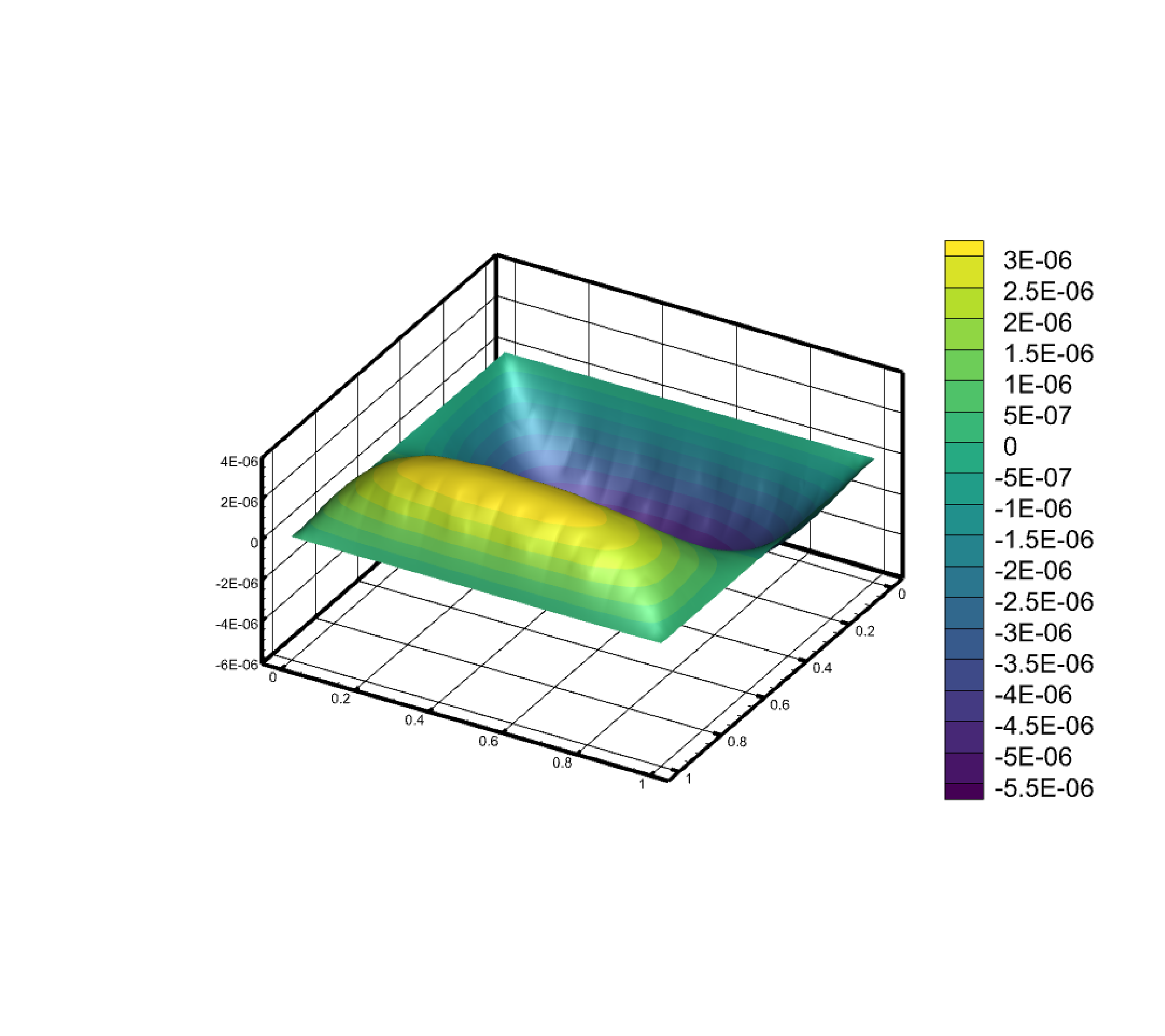} \\
		(c)
	\end{minipage}
	\begin{minipage}[b]{0.32\textwidth}
		\centering
		\includegraphics[
		width=4cm,
		trim=2.5cm 6.5cm 0cm 6.5cm, % 从左、下、右、上各裁剪1厘米
		clip % 应用裁剪
		]{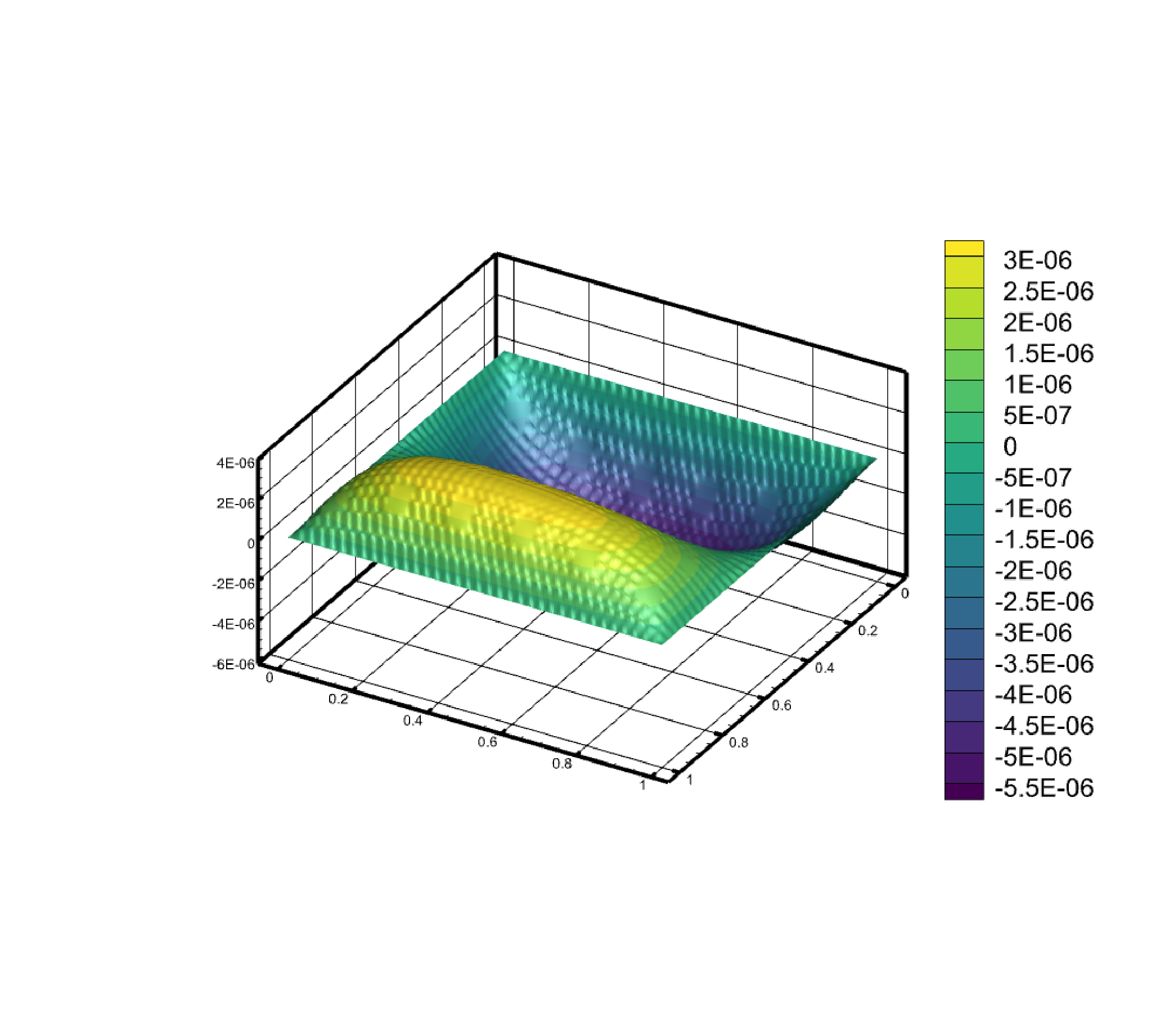} \\
		(d)
	\end{minipage}
	\begin{minipage}[b]{0.32\textwidth}
		\centering
		\includegraphics[
		width=4cm,
		trim=2.5cm 6.5cm 0cm 6.5cm, % 从左、下、右、上各裁剪1厘米
		clip % 应用裁剪
		]{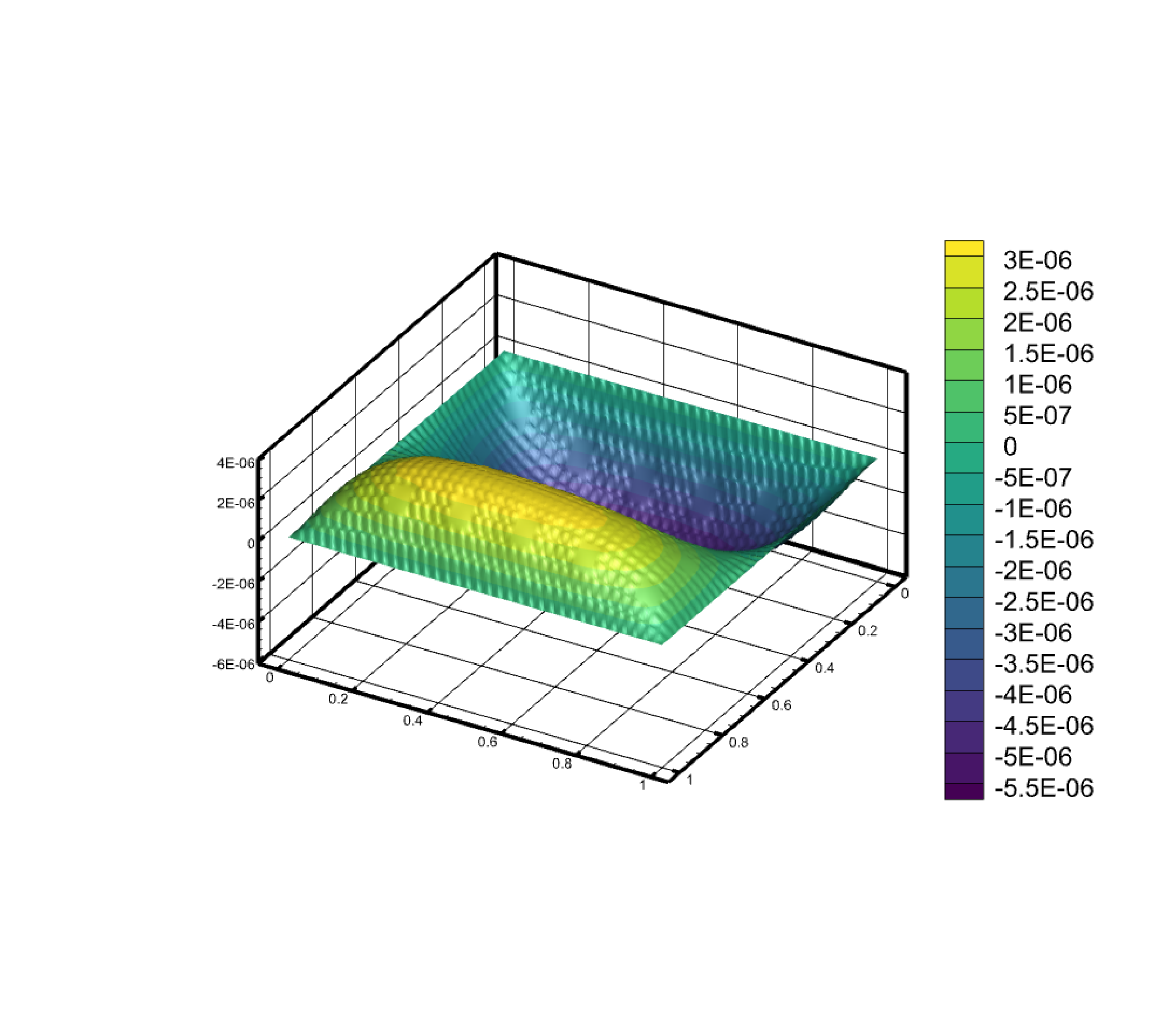} \\
		(e)
	\end{minipage}
	\caption{The first component for the displacement field: (a) ${u _{1Ref}^{{\zeta  _1}{\zeta  _2}}}$; (b) ${{u_{1}^{(0)}}}$;
		(c) ${u_{1D}^{{\zeta  _1}}}$; (d) ${u_{1t}^{{\zeta  _1}{\zeta  _2}}}$; (e) ${u_{1T}^{{\zeta  _1}{\zeta  _2}}}$.}\label{E1f3}
\end{figure}

\begin{figure}[!htb]
	\centering
	\begin{minipage}[b]{0.35\textwidth}
		\centering
		\includegraphics[
		width=4cm,
		trim=2.5cm 6.5cm 0cm 6.5cm, % 从左、下、右、上各裁剪1厘米
		clip % 应用裁剪
		]{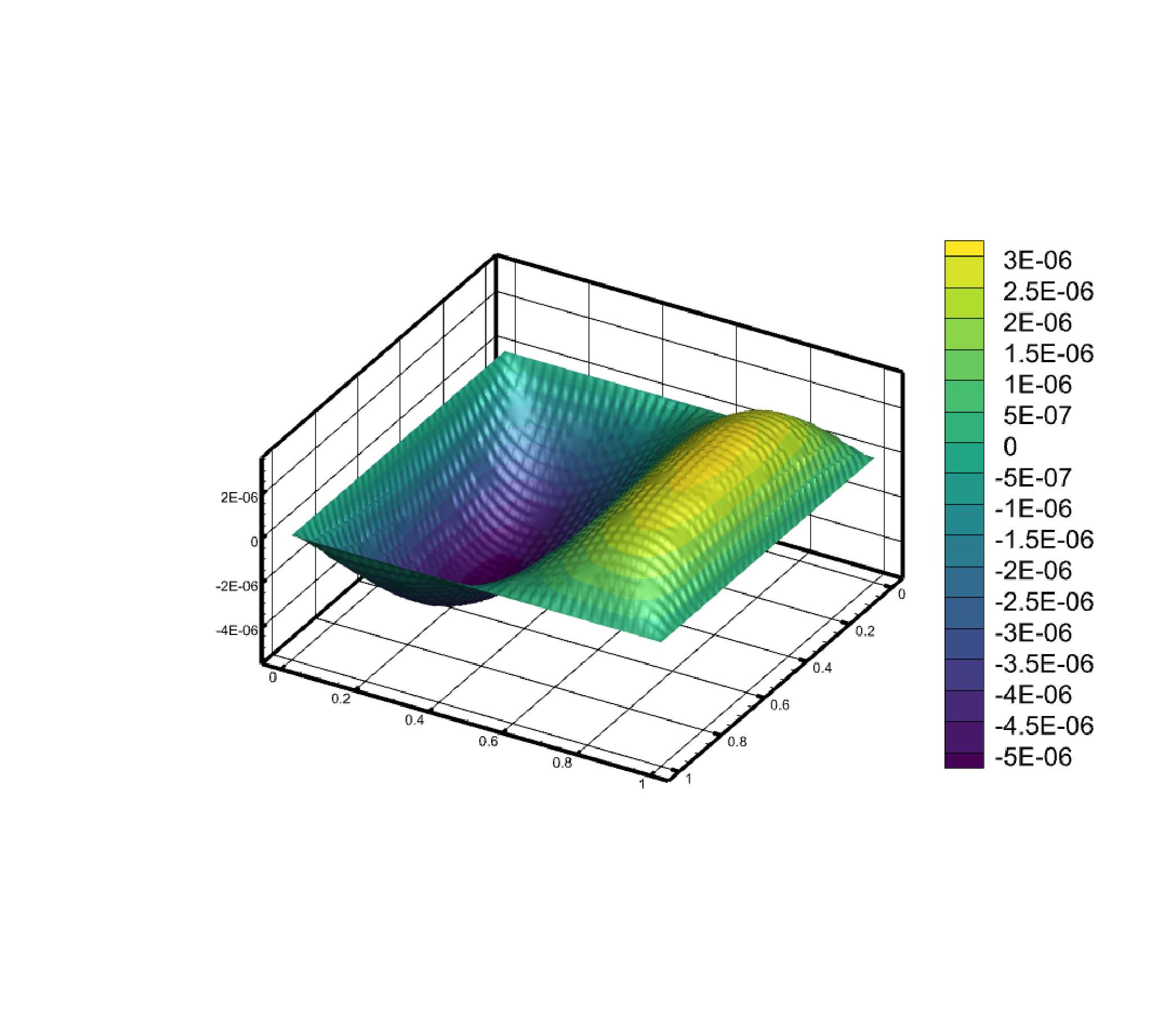} \\
		(a)
	\end{minipage}
	\begin{minipage}[b]{0.35\textwidth}
		\centering
		\includegraphics[
		width=4cm,
		trim=2.5cm 6.5cm 0cm 6.5cm, % 从左、下、右、上各裁剪1厘米
		clip % 应用裁剪
		]{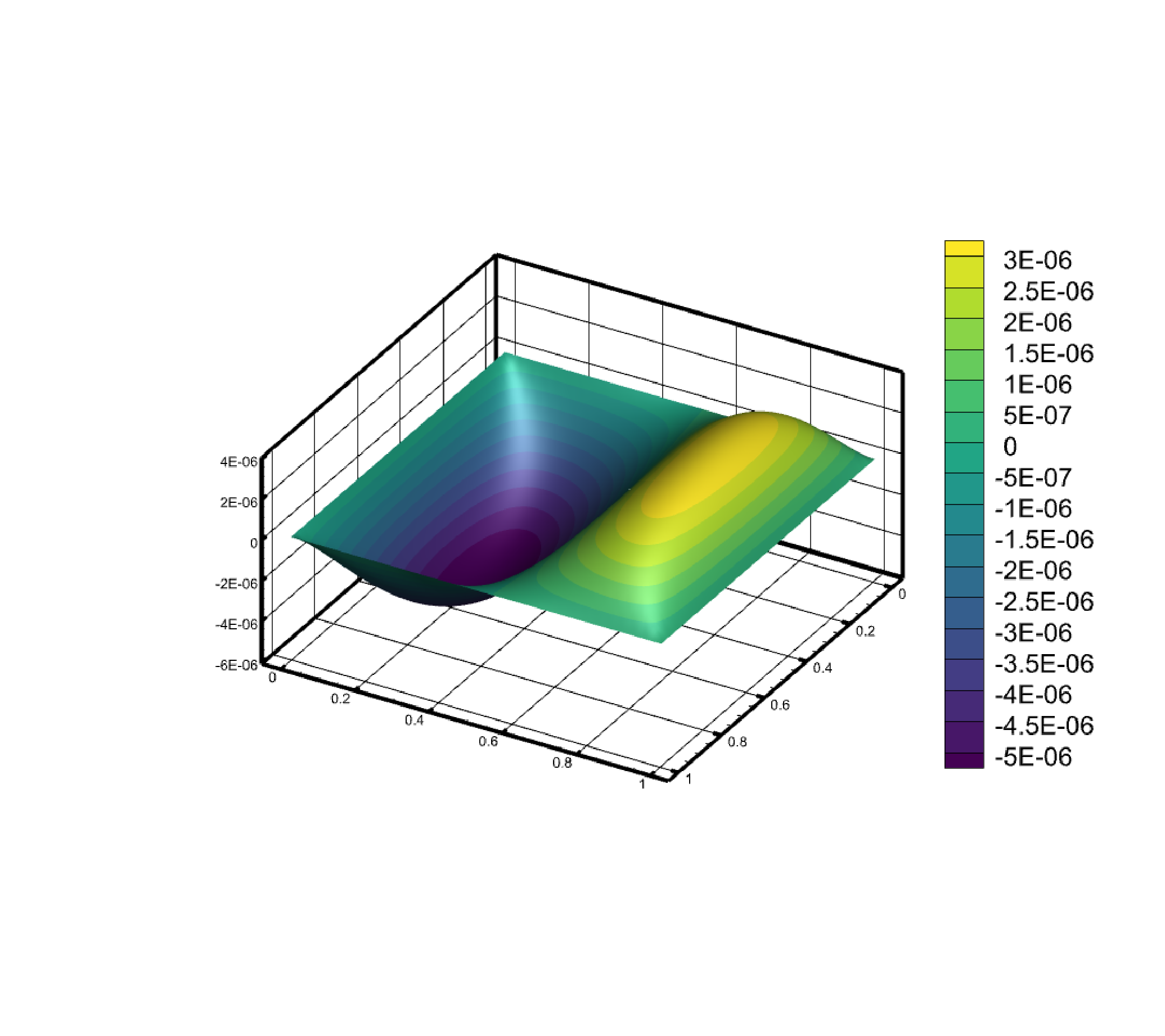} \\
		(b)
	\end{minipage}
	\begin{minipage}[b]{0.32\textwidth}
		\centering
		\includegraphics[
		width=4cm,
		trim=2.5cm 6.5cm 0cm 6.5cm, % 从左、下、右、上各裁剪1厘米
		clip % 应用裁剪
		]{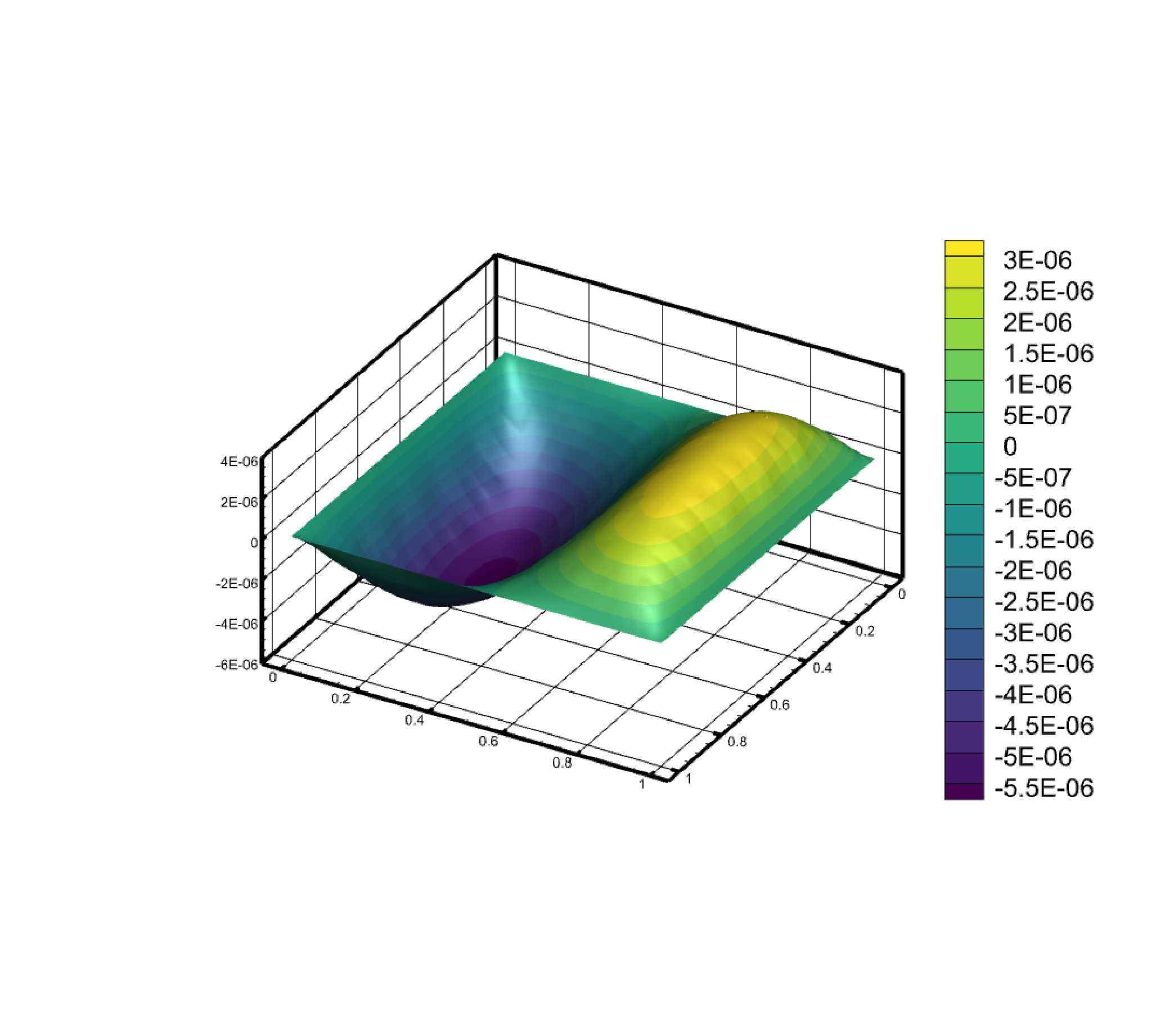} \\
		(c)
	\end{minipage}
	\begin{minipage}[b]{0.32\textwidth}
		\centering
		\includegraphics[
		width=4cm,
		trim=2.5cm 6.5cm 0cm 6.5cm, % 从左、下、右、上各裁剪1厘米
		clip % 应用裁剪
		]{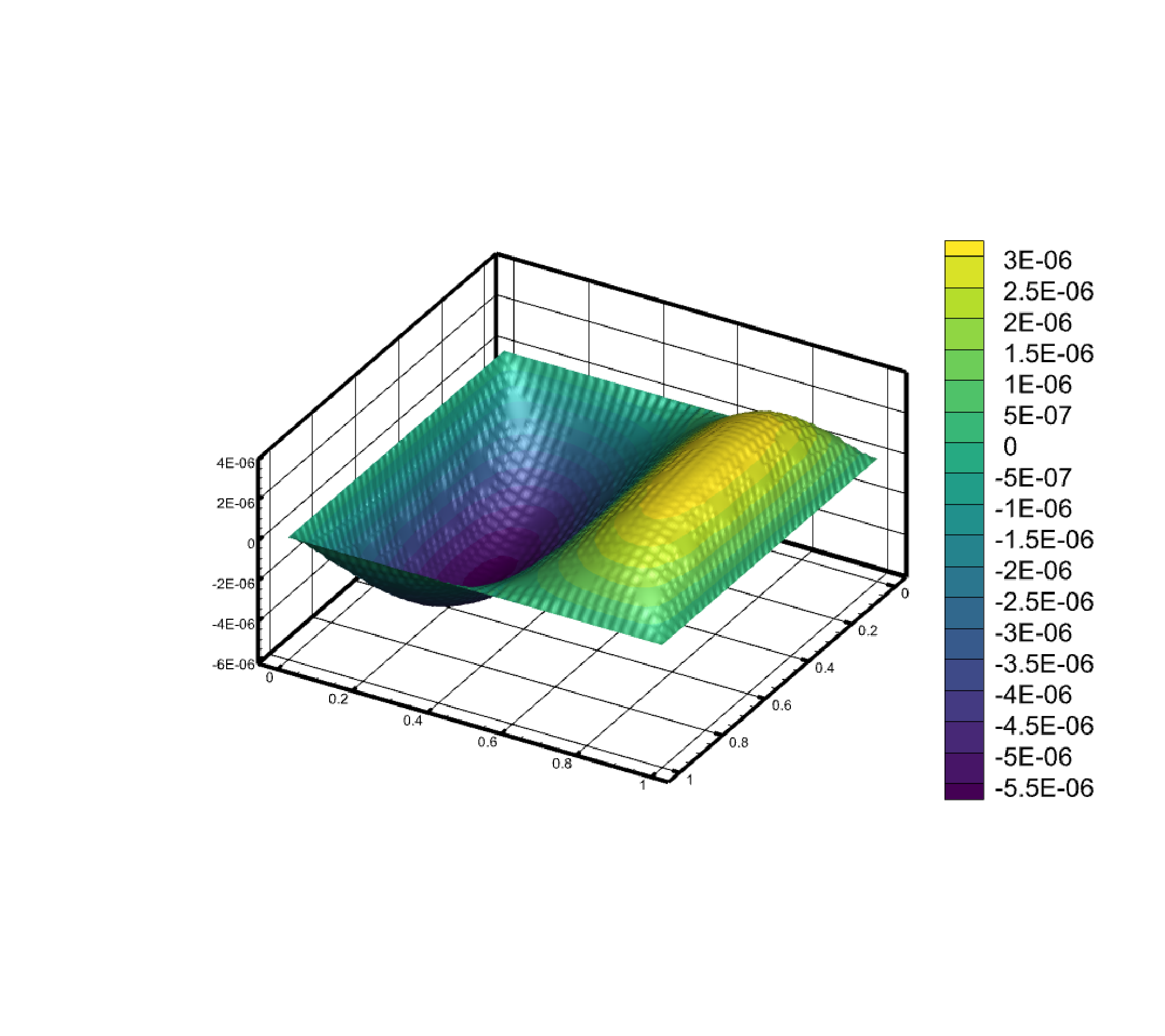} \\
		(d)
	\end{minipage}
	\begin{minipage}[b]{0.32\textwidth}
		\centering
		\includegraphics[
		width=4cm,
		trim=2.5cm 6.5cm 0cm 6.5cm, % 从左、下、右、上各裁剪1厘米
		clip % 应用裁剪
		]{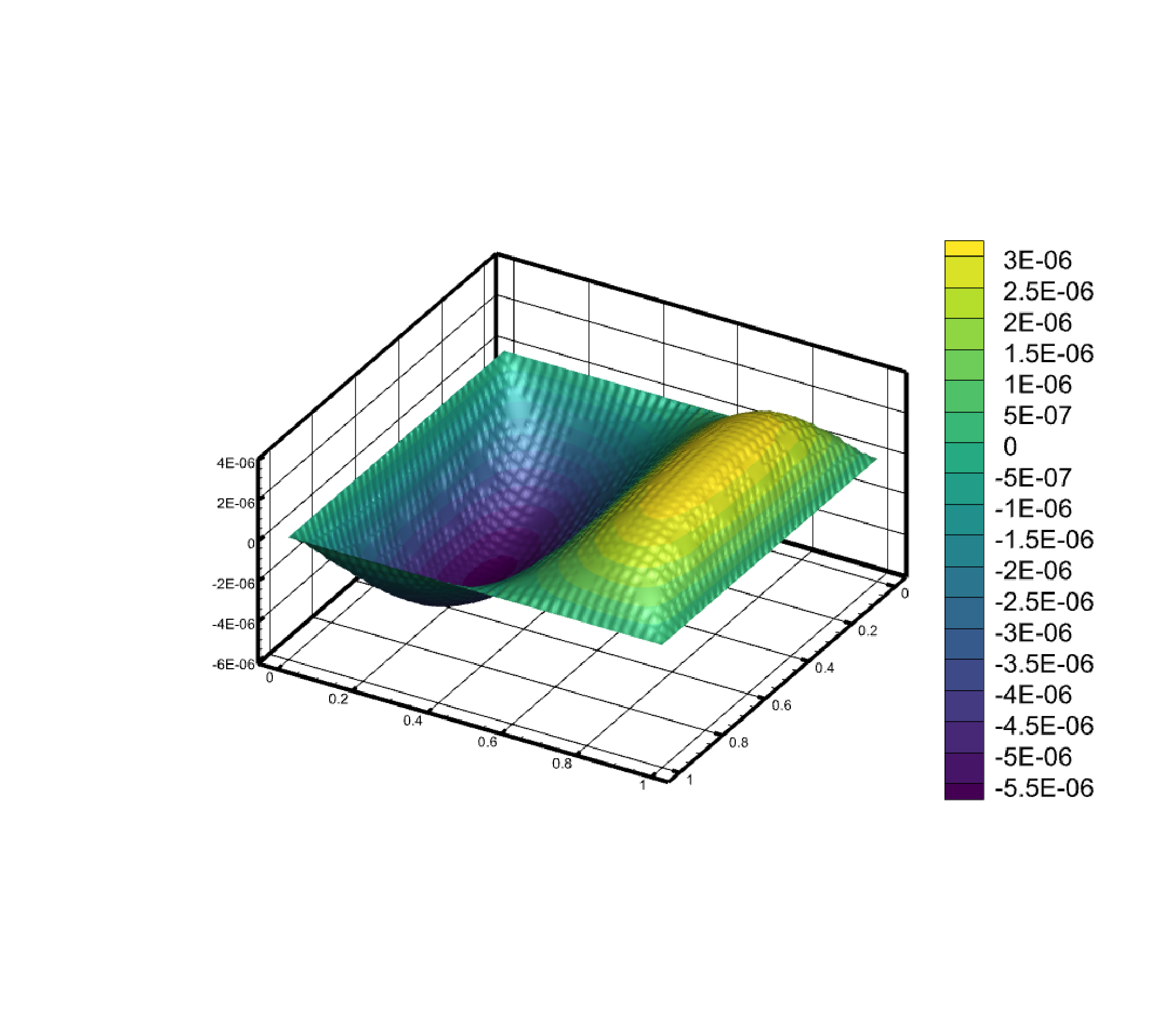} \\
		(e)
	\end{minipage}
	\caption{The second component for the displacement field: (a) ${u _{2Ref}^{{\zeta  _1}{\zeta  _2}}}$; (b) ${{u_{2}^{(0)}}}$;
		(c) ${u_{2D}^{{\zeta  _1}}}$; (d) ${u_{2t}^{{\zeta  _1}{\zeta  _2}}}$; (e) ${u_{2T}^{{\zeta  _1}{\zeta  _2}}}$.}\label{E1f4}
\end{figure}

After simulating for 2D composite structure I, Figs.\hspace{1mm}\ref{E1f2}-\ref{E1f4} display the simulation results at time $t= 1.0 s$ including temperature field and displacement fields, respectively.

\begin{figure}[!htb]
	\centering
	\begin{minipage}[b]{0.4\textwidth}
		\centering
		\includegraphics[width=4cm]{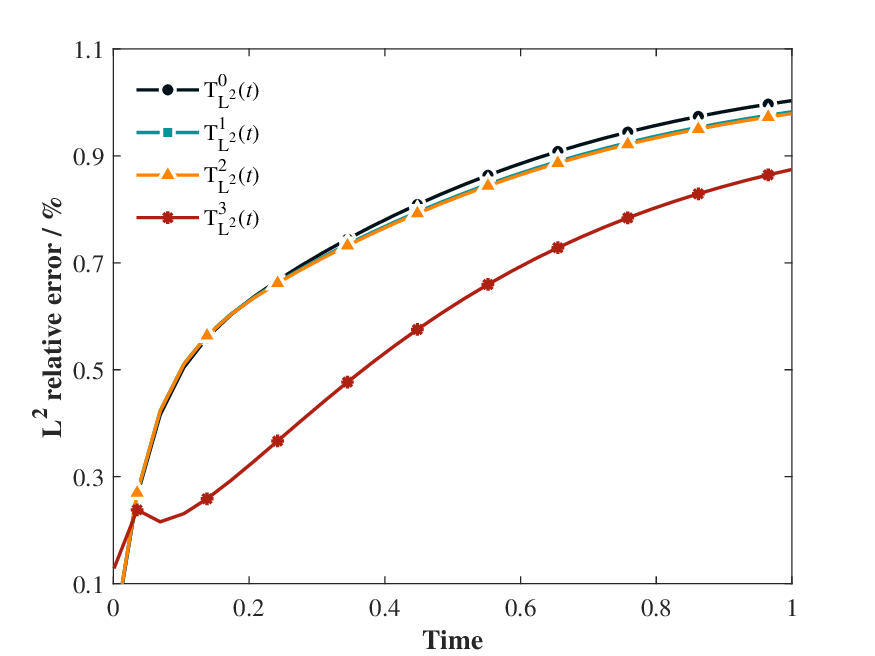} \\
		(a)
	\end{minipage}
	\begin{minipage}[b]{0.4\textwidth}
		\centering
		\includegraphics[width=4cm]{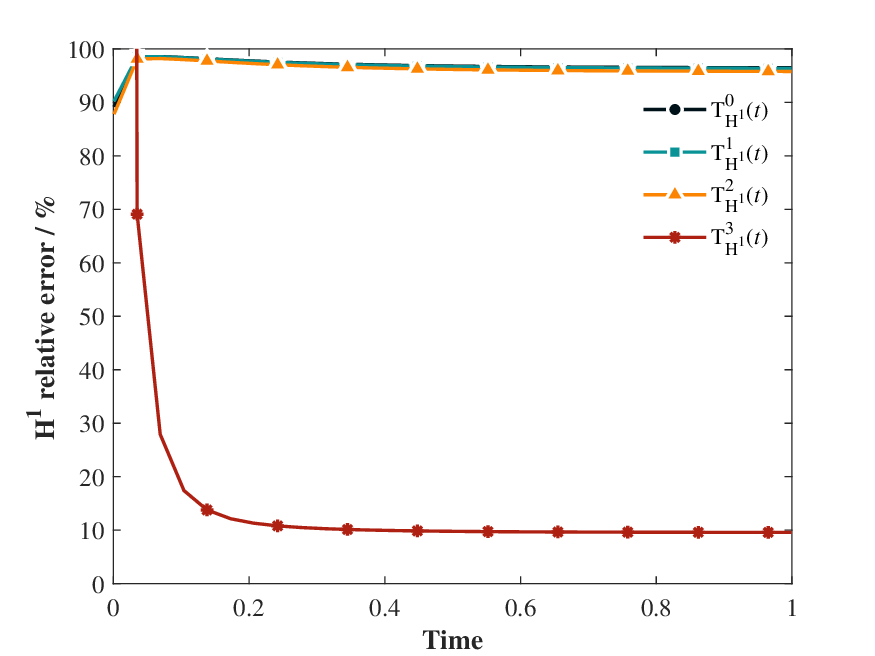} \\
		(b)
	\end{minipage}
	\begin{minipage}[b]{0.4\textwidth}
		\centering
		\includegraphics[width=4cm]{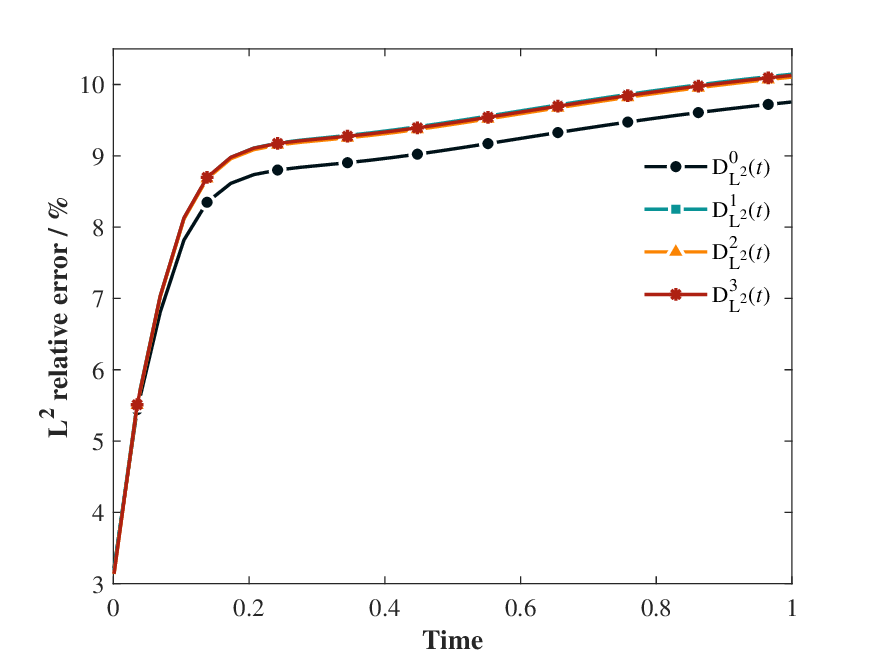} \\
		(c)
	\end{minipage}
	\begin{minipage}[b]{0.4\textwidth}
		\centering
		\includegraphics[width=4cm]{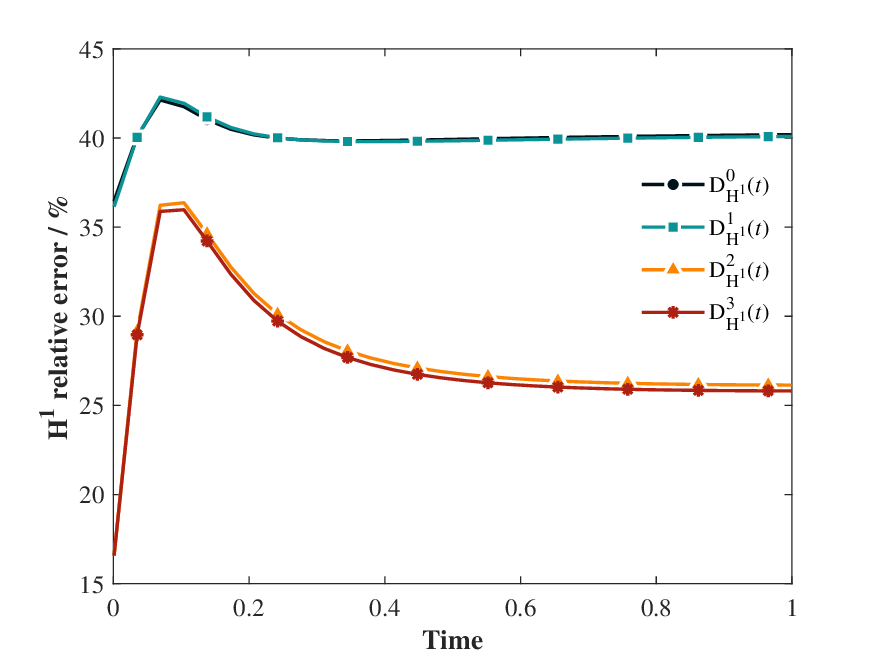} \\
		(d)
	\end{minipage}
	\caption{The evolutive relative errors of temperature and displacement fields: (a) ${\rm{T_{{L^2}}}}(t)$; (b) ${\rm{T_{{H^1}}}}(t)$; (c) ${\rm{D_{{L^2}}}}(t)$; (d) ${\rm{D_{{H^1}}}}(t)$.}\label{E1f5}
\end{figure}

Next, the evolutive relative errors of distinct multi-scale solutions for temperature and displacement field are depicted in Fig.\hspace{1mm}\ref{E1f5}.

To wrap up, as demonstrated in Table \ref{E1t2}, it is apparent that our HOTS methodology can reduce the computational expense and improve the computational efficiency. As shown in Figs.\hspace{1mm}\ref{E1f2}-\ref{E1f4}, it is clearly that only the solutions of the HOTS method demonstrate strong agreement with the reference FEM solutions and accurately capture the highly microscopic oscillating information. Furthermore, other multi-scale solutions failed to capture the essential oscillating information at the smallest scale, yielding significantly lower numerical accuracy compared to that of the higher-order three-scale solutions. From Fig.\hspace{1mm}\ref{E1f5}, it is plainly visible that the accuracy of the HOTS solutions of physical fields are much better than that of the other three solutions. Furthermore, we can also conclude that the proposed two-stage numerical algorithm is robust without blow-up phenomenon even if it runs for a long time.

\subsection{Example 2: 2D composite structure II with multiple spatial scales}
\begin{figure}[!htb]
	\centering
	\begin{minipage}[c]{0.42\textwidth}
		\centering
		\includegraphics[
		width=0.95\linewidth,
		trim=2cm 1.5cm 1.8cm 1.3cm, % 从左、下、右、上各裁剪1厘米
		clip % 应用裁剪
		]{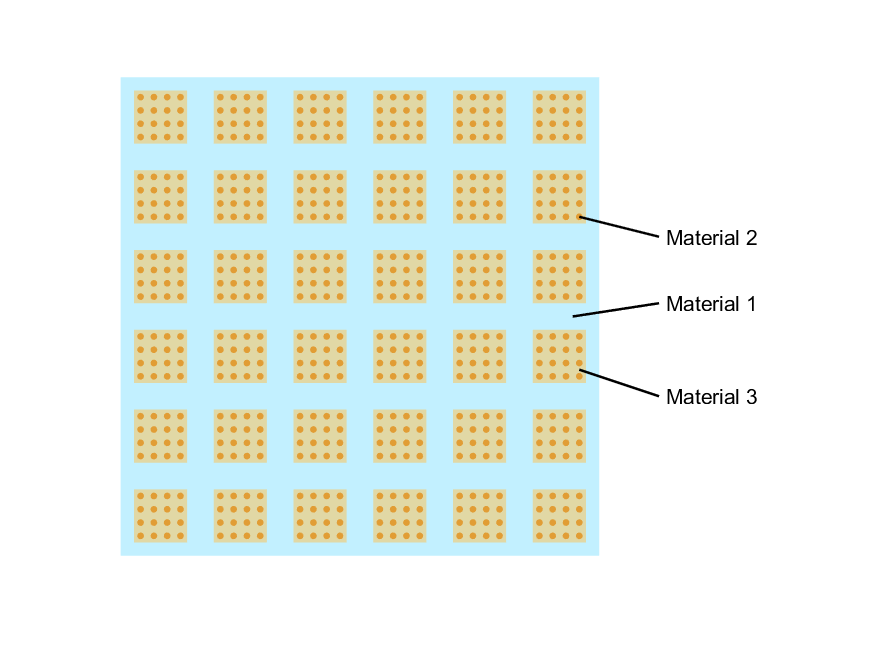} \\
		(a)
	\end{minipage}
	\begin{minipage}[c]{0.5\textwidth}
		\centering
		\includegraphics[
		width=0.95\linewidth,
		trim=2cm 1.5cm 1.5cm 1.5cm, % 从左、下、右、上各裁剪1厘米
		clip % 应用裁剪
		]{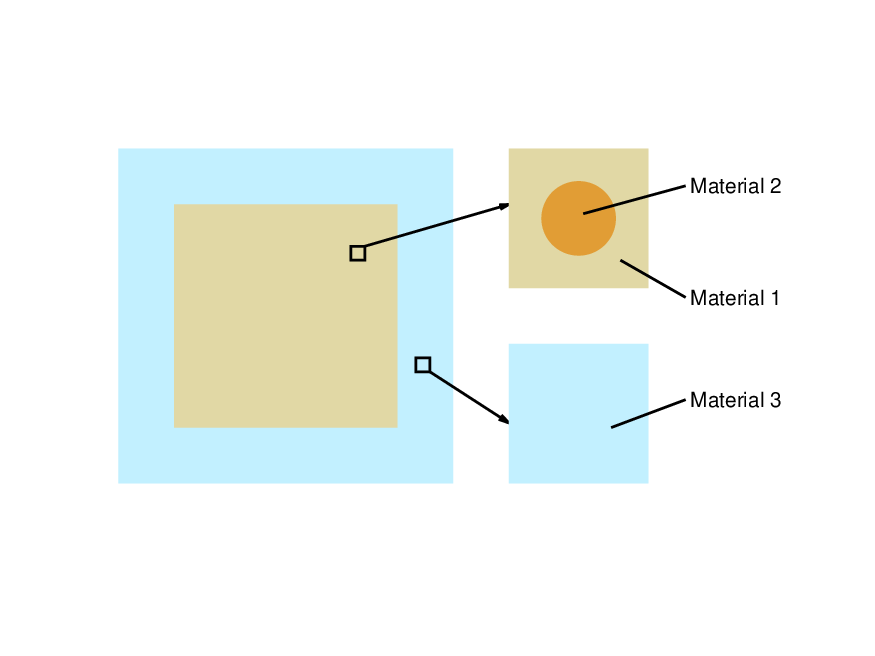} \\
		(b)
	\end{minipage}
	\caption{(a) Macroscopic domain $\Omega$; (b) Mesoscopic unit cell $Y$ (left side) and microscopic unit cell $Z$ (right side).}\label{E2f1}
\end{figure}

In this example, we consider the nonlinear thermo-mechanical behaviors of a novel 2D composite structure II with $\zeta _1=1/6$ and $\zeta _2=1/36$, whose specific three-scale configuration is depicted in Fig.\hspace{1mm}\ref{E2f1}. Additionally, this example employs the same physical configurations as Example 1, including material properties, internal sources, as well as initial-boundary conditions.

After that, we execute the triangular meshing to three-scale nonlinear problem, mesoscopic UC problems, microscopic UC problems and associated macroscopic homogenized problem. Following that, an assessment of the computational costs for the reference FEM and HOTS method is provided in Table \ref{E2t1}.
\begin{table}[!htb]
	\centering
	\caption{Summary of computational cost ($\Delta t = 0.001, t \in [0,1] $).}\label{E2t1}
	\begin{tabular}{lcccc}
		\hline
		& \multirow{2}{*}{Reference FEM} & \multicolumn{3}{c}{HOTS method} \\
		\cmidrule(l){3-5}
		& & $Z$ & $Y$ & $\Omega$ \\
		\hline
		Number of nodes & 105,049 & 1,179 & 1,838 & 10,201 \\
		Number of elements & 208,800 & 3,396 & 3,514 & 20,000 \\
		Computational time & 134,860s &\multicolumn{3}{c}{61,127.9s} \\
		\hline
	\end{tabular}
\end{table}

\begin{figure}[!htb]
	\centering
	\begin{minipage}[b]{0.35\textwidth}
		\centering
		\includegraphics[
		width=4cm,
		trim=2.5cm 4.5cm 0cm 4.5cm, % 从左、下、右、上各裁剪1厘米
		clip % 应用裁剪
		]{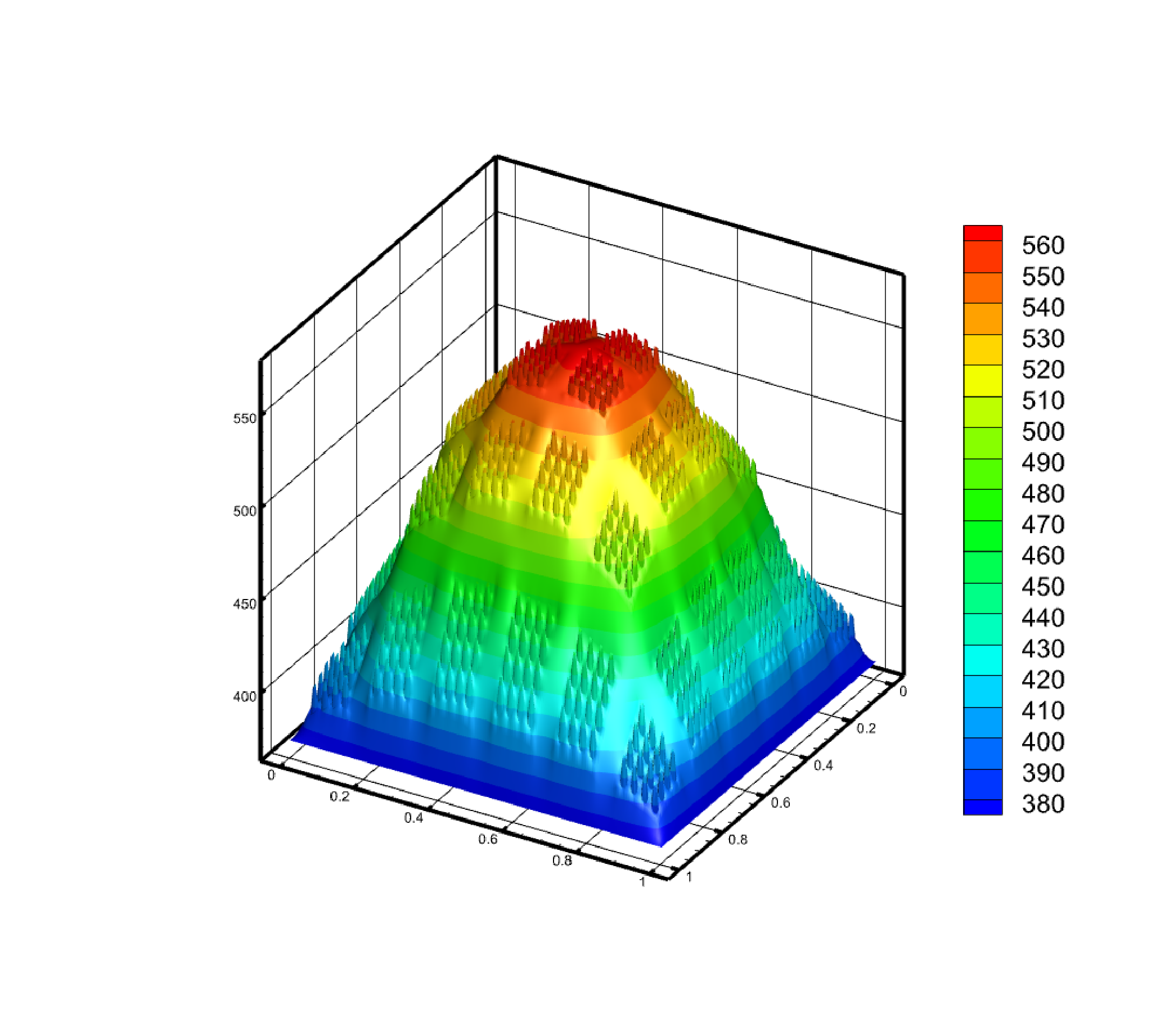} \\
		(a)
	\end{minipage}
	\begin{minipage}[b]{0.35\textwidth}
		\centering
		\includegraphics[
		width=4cm,
		trim=2.5cm 4.5cm 0cm 4.5cm, % 从左、下、右、上各裁剪1厘米
		clip % 应用裁剪
		]{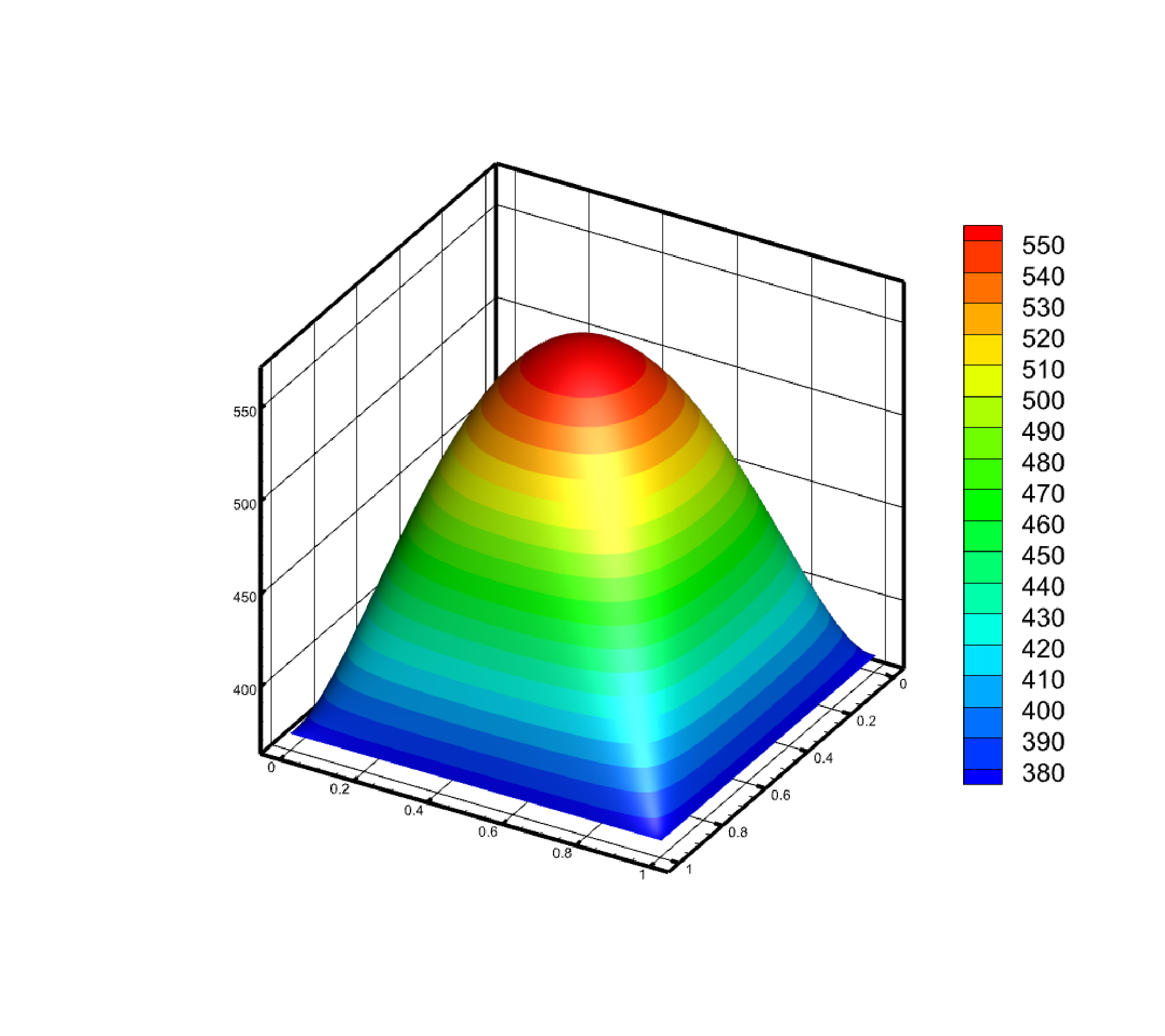} \\
		(b)
	\end{minipage}
	\begin{minipage}[b]{0.32\textwidth}
		\centering
		\includegraphics[
		width=4cm,
		trim=2.5cm 4.5cm 0cm 4.5cm, % 从左、下、右、上各裁剪1厘米
		clip % 应用裁剪
		]{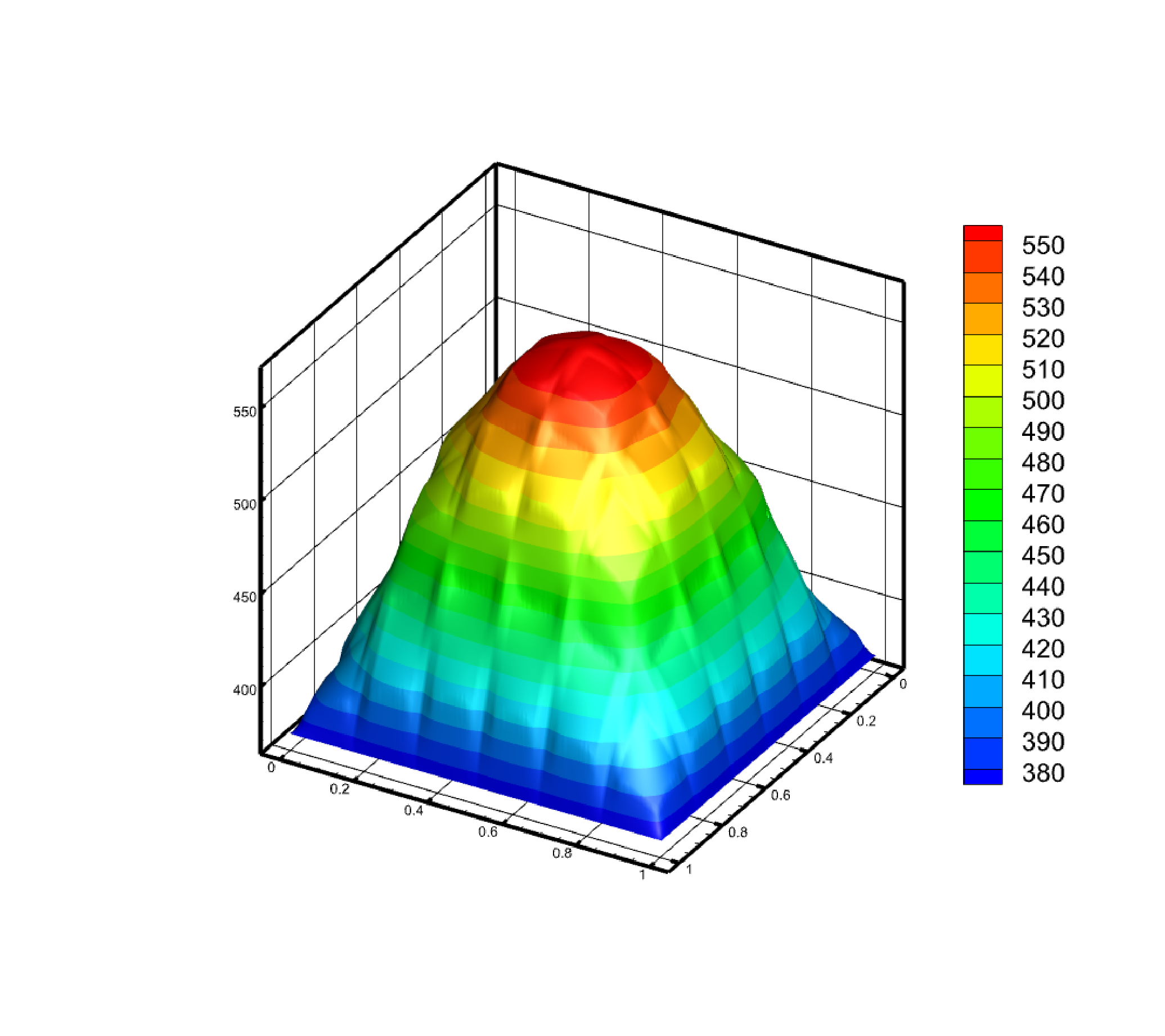} \\
		(c)
	\end{minipage}
	\begin{minipage}[b]{0.32\textwidth}
		\centering
		\includegraphics[
		width=4cm,
		trim=2.5cm 4.5cm 0cm 4.5cm, % 从左、下、右、上各裁剪1厘米
		clip % 应用裁剪
		]{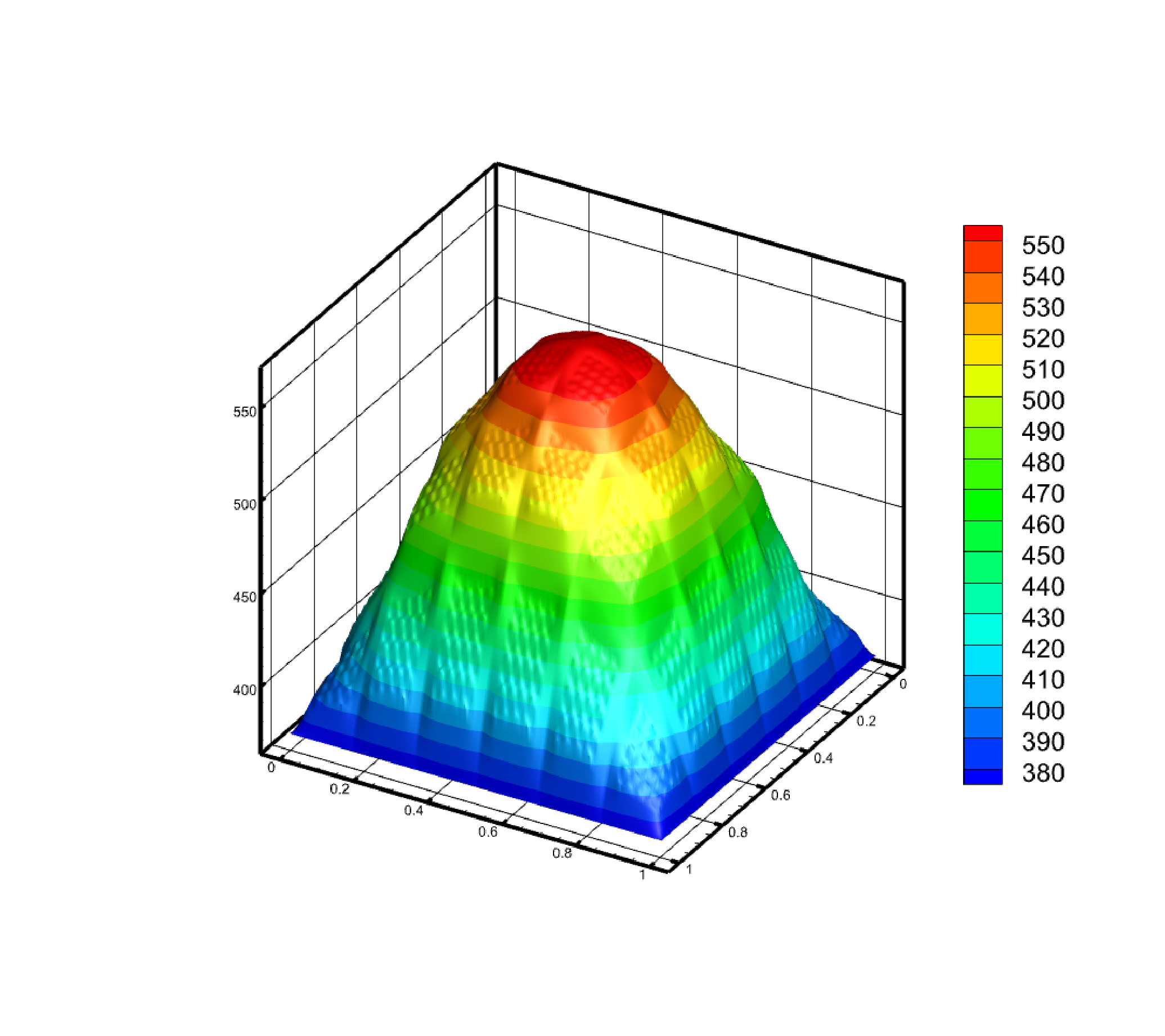} \\
		(d)
	\end{minipage}
	\begin{minipage}[b]{0.32\textwidth}
		\centering
		\includegraphics[
		width=4cm,
		trim=2.5cm 4.5cm 0cm 4.5cm, % 从左、下、右、上各裁剪1厘米
		clip % 应用裁剪
		]{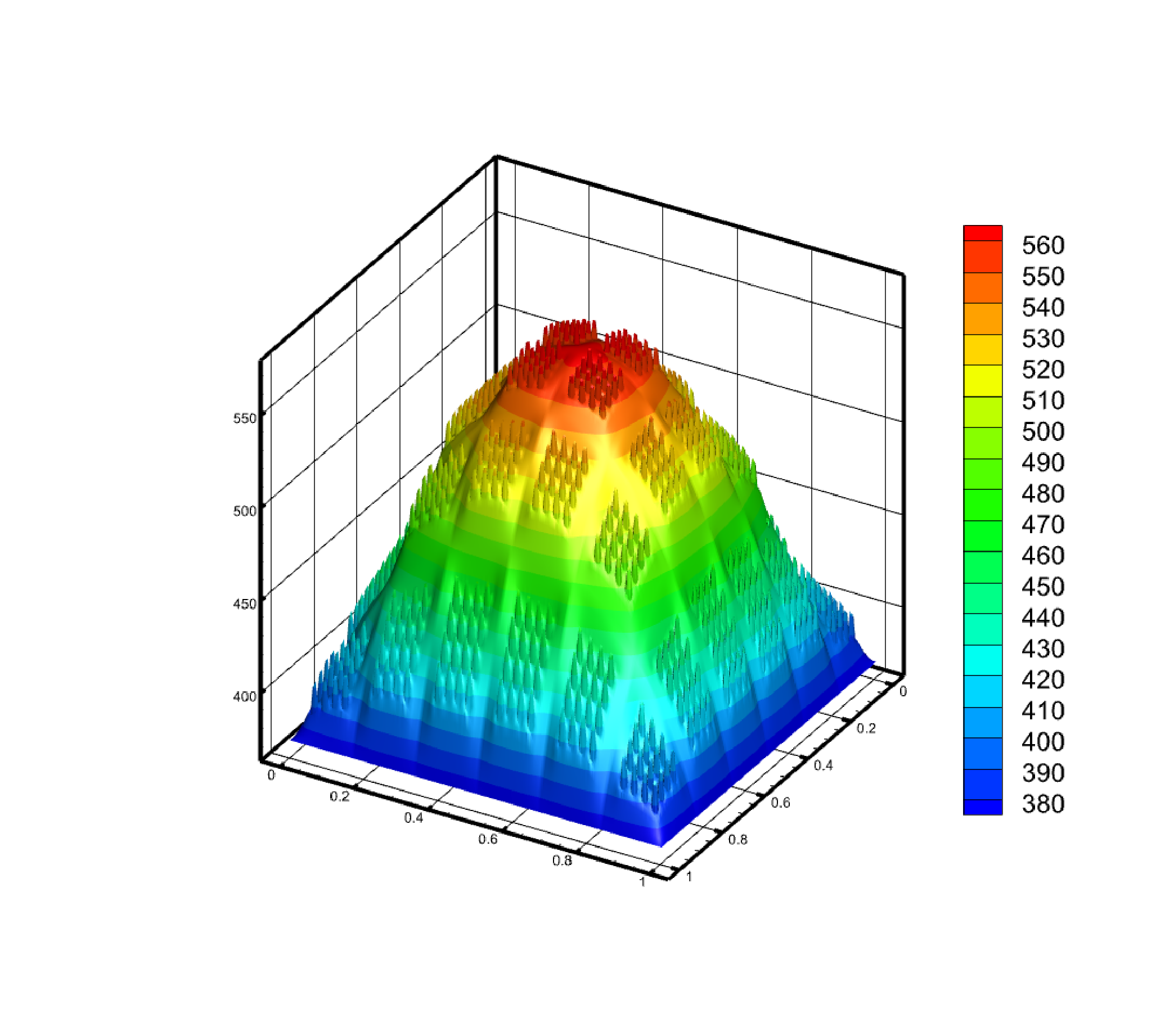} \\
		(e)
	\end{minipage}
	\caption{The temperature field: (a) $\theta_{Ref}^{\zeta _1\zeta _2}$; (b) $\theta_{0}$; (c) ${\theta _D^{{\zeta  _1}}}$; (d) ${\theta _t^{{\zeta  _1}{\zeta  _2}}}$; (e) ${\theta _T^{{\zeta  _1}{\zeta  _2}}}$.}\label{E2f2}
\end{figure}

\begin{figure}[!htb]
	\centering
	\begin{minipage}[b]{0.35\textwidth}
		\centering
		\includegraphics[
		width=4cm,
		trim=2.5cm 6.5cm 0cm 6.5cm, % 从左、下、右、上各裁剪1厘米
		clip % 应用裁剪
		]{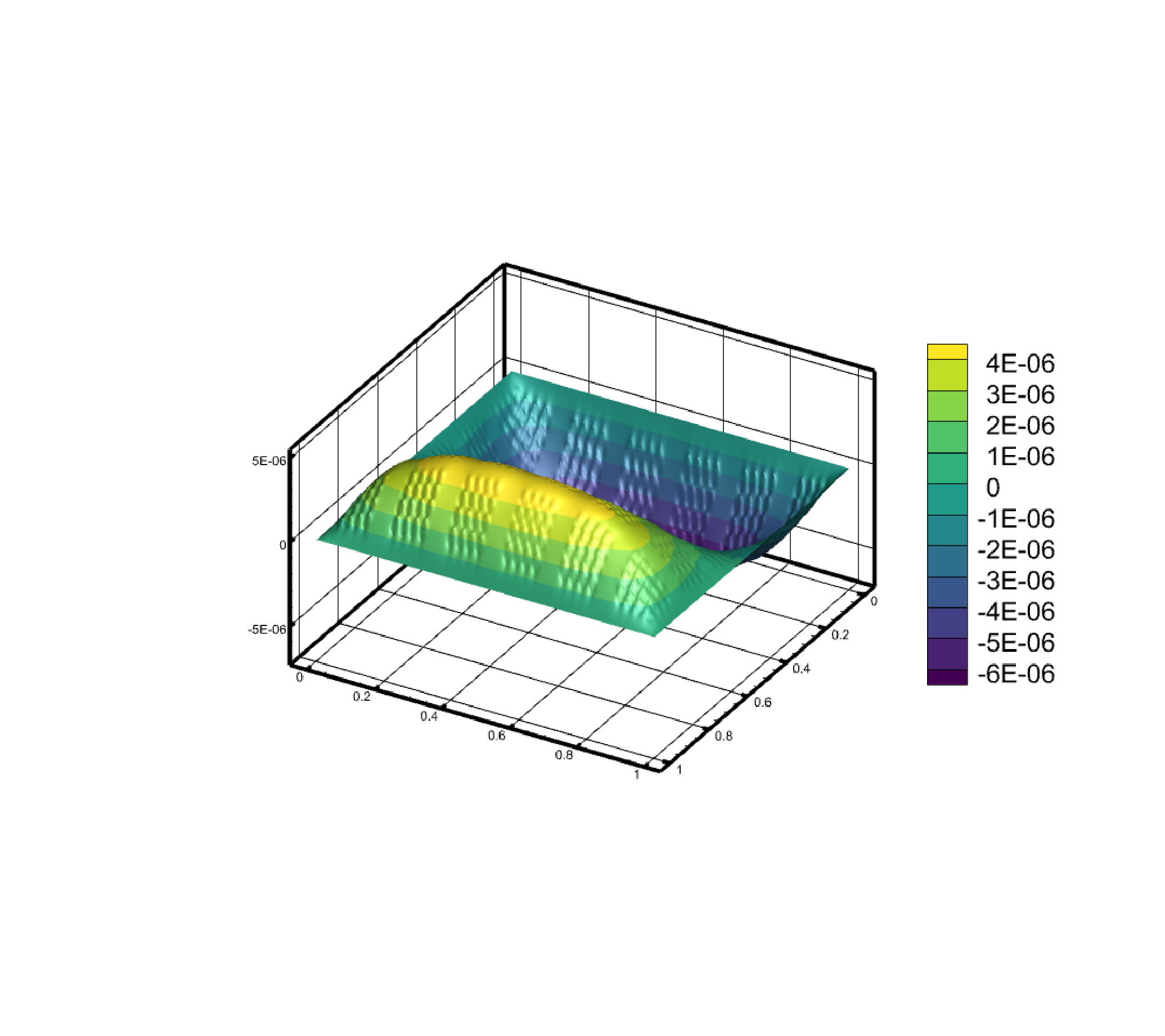} \\
		(a)
	\end{minipage}
	\begin{minipage}[b]{0.35\textwidth}
		\centering
		\includegraphics[
		width=4cm,
		trim=2.5cm 6.5cm 0cm 6.5cm, % 从左、下、右、上各裁剪1厘米
		clip % 应用裁剪
		]{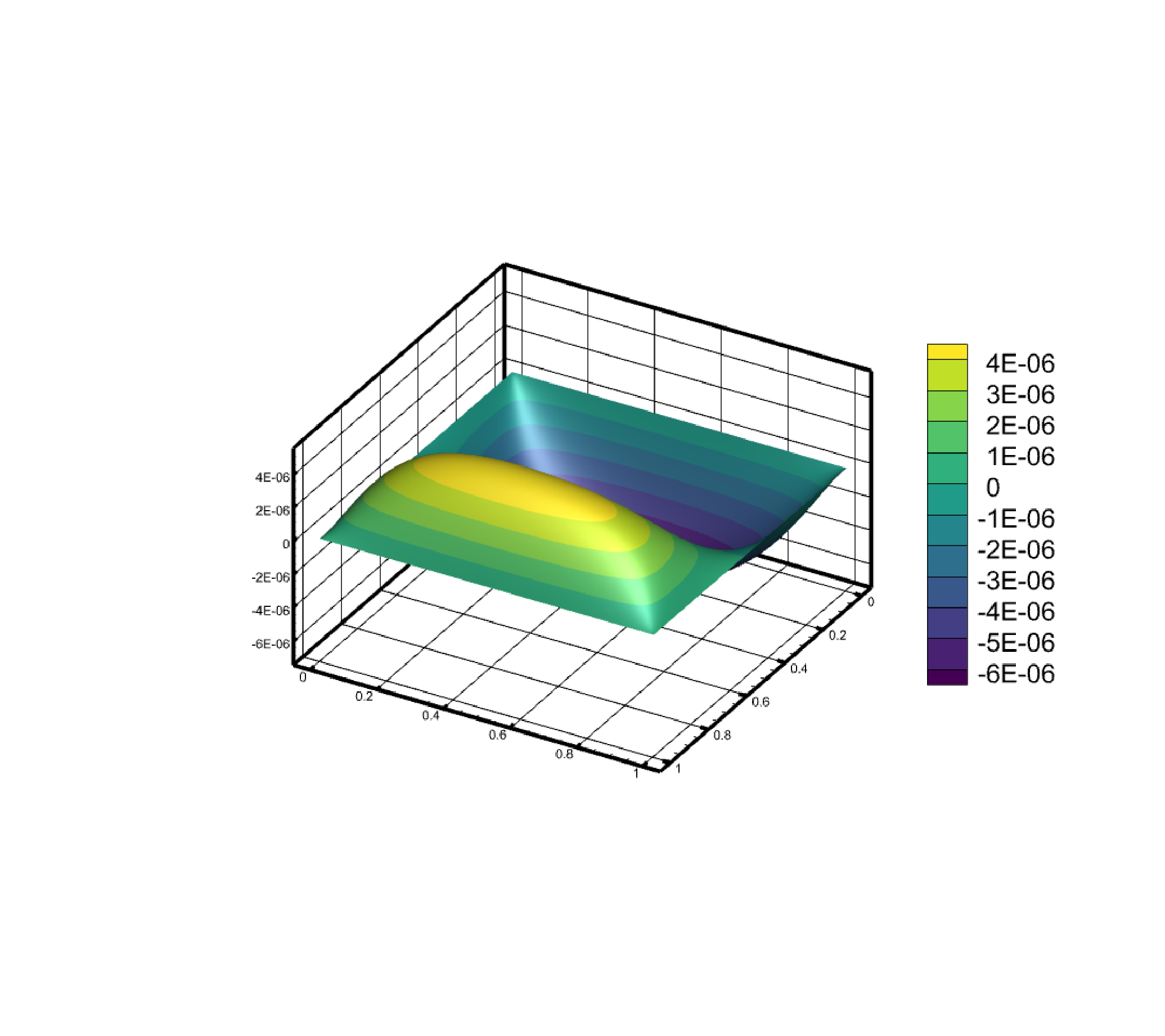} \\
		(b)
	\end{minipage}
	\begin{minipage}[b]{0.32\textwidth}
		\centering
		\includegraphics[
		width=4cm,
		trim=2.5cm 6.5cm 0cm 6.5cm, % 从左、下、右、上各裁剪1厘米
		clip % 应用裁剪
		]{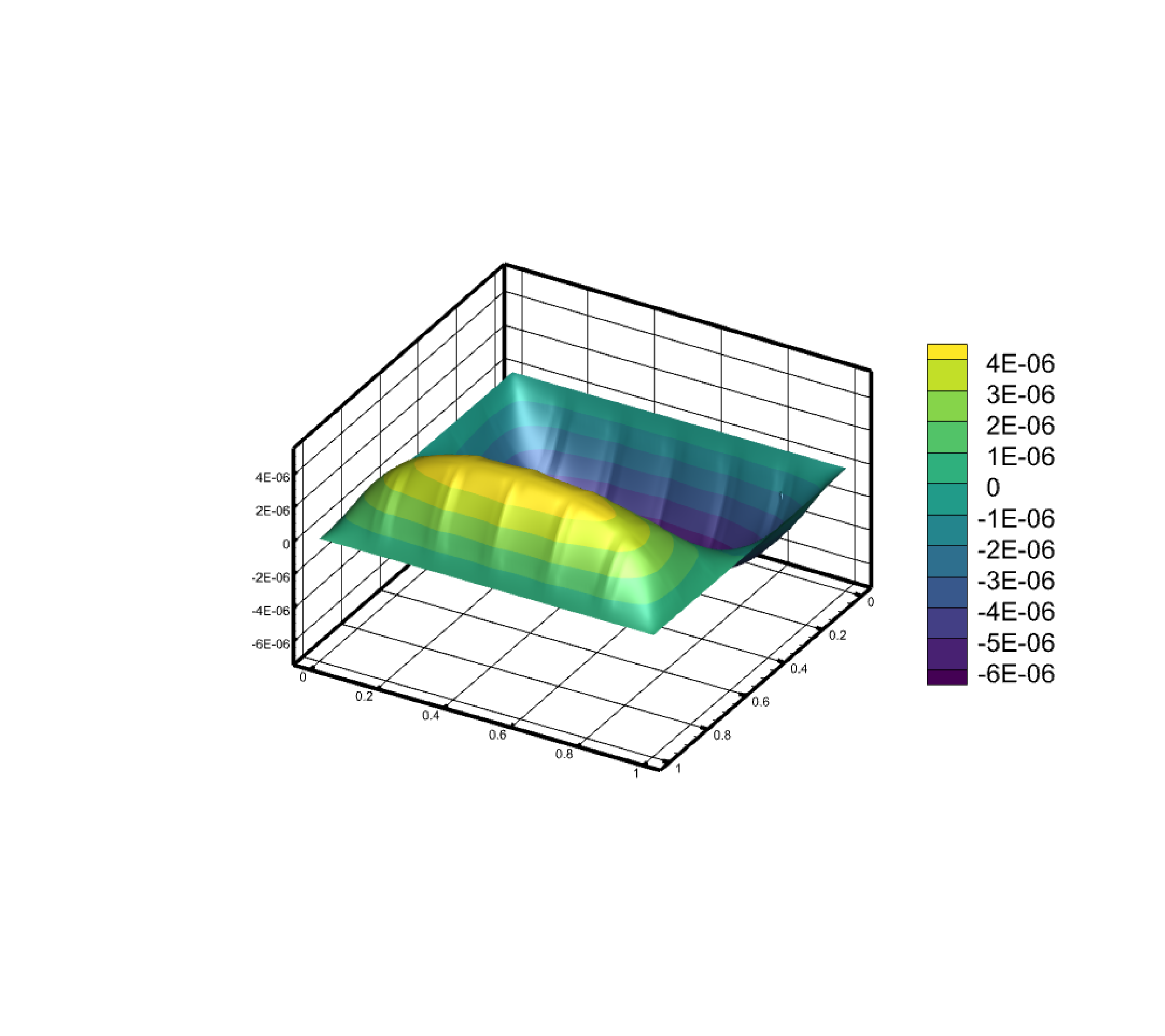} \\
		(c)
	\end{minipage}
	\begin{minipage}[b]{0.32\textwidth}
		\centering
		\includegraphics[
		width=4cm,
		trim=2.5cm 6.5cm 0cm 6.5cm, % 从左、下、右、上各裁剪1厘米
		clip % 应用裁剪
		]{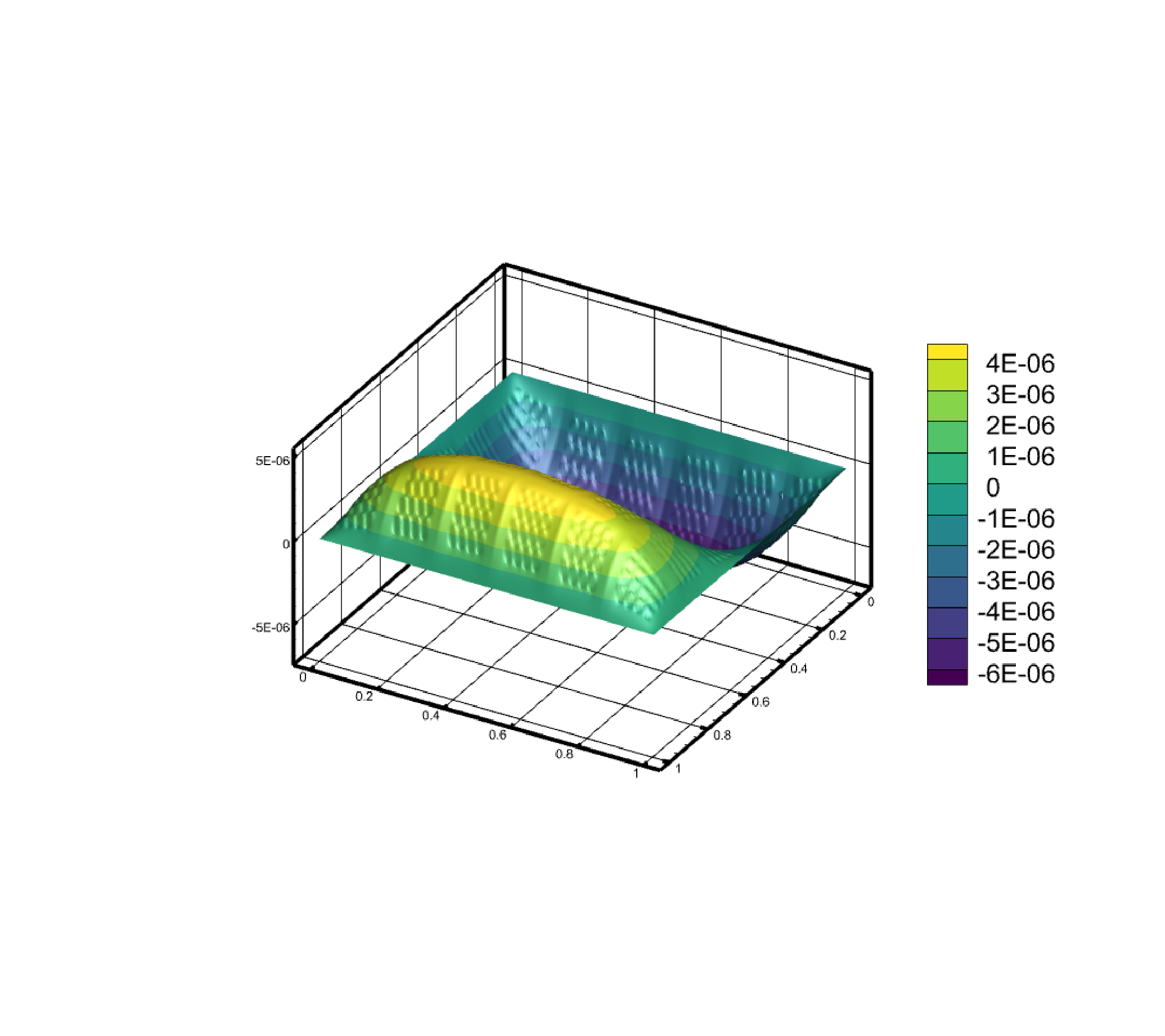} \\
		(d)
	\end{minipage}
	\begin{minipage}[b]{0.32\textwidth}
		\centering
		\includegraphics[
		width=4cm,
		trim=2.5cm 6.5cm 0cm 6.5cm, % 从左、下、右、上各裁剪1厘米
		clip % 应用裁剪
		]{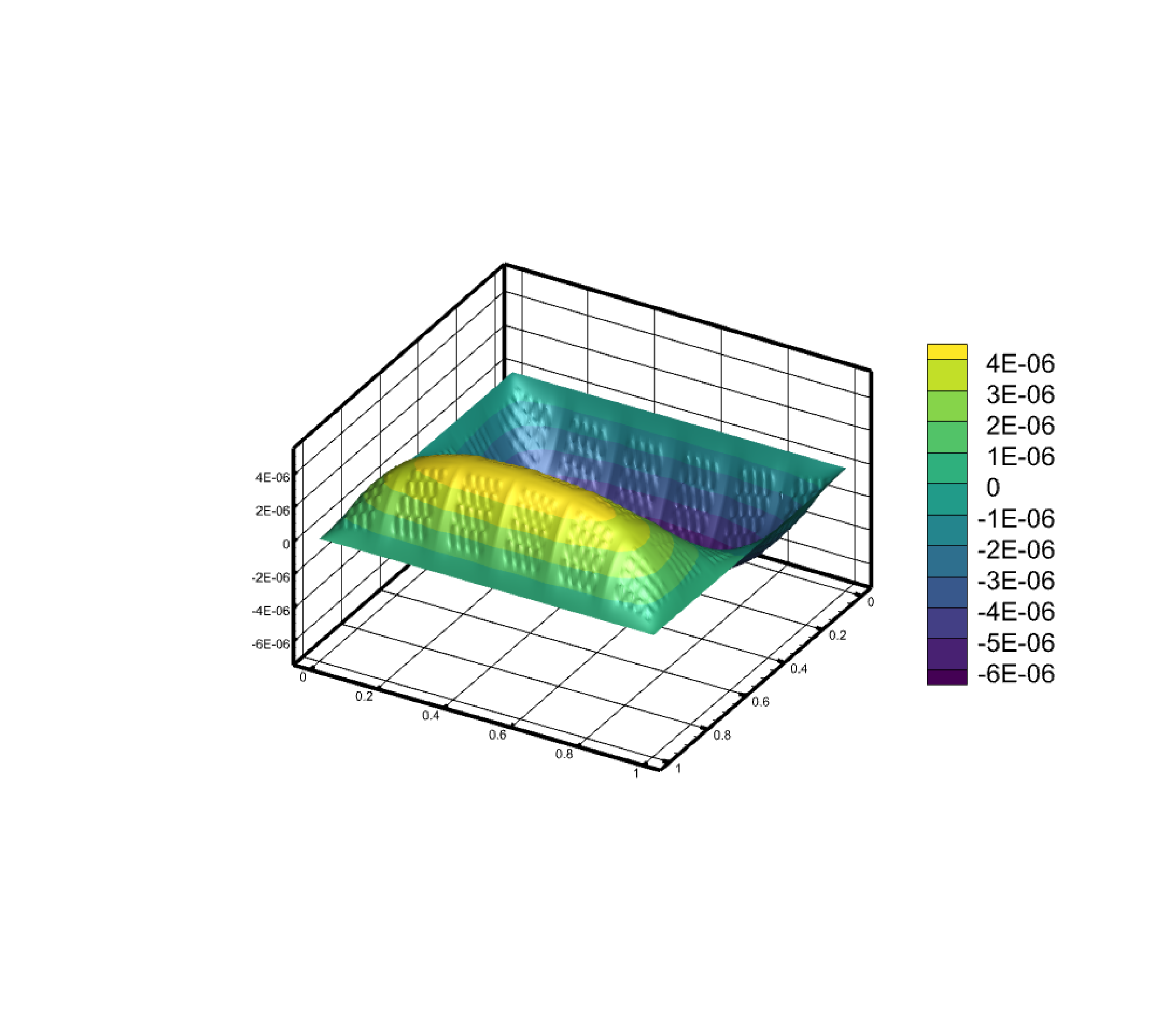} \\
		(e)
	\end{minipage}
	\caption{The first component for the displacement field: (a) ${u _{1Ref}^{{\zeta  _1}{\zeta  _2}}}$; (b) ${{u_{1}^{(0)}}}$;
		(c) ${u_{1D}^{{\zeta  _1}}}$; (d) ${u_{1t}^{{\zeta  _1}{\zeta  _2}}}$; (e) ${u_{1T}^{{\zeta  _1}{\zeta  _2}}}$.}\label{E2f3}
\end{figure}

\begin{figure}[!htb]
	\centering
	\begin{minipage}[b]{0.35\textwidth}
		\centering
		\includegraphics[
		width=4cm,
		trim=2.5cm 6.5cm 0cm 6.5cm, % 从左、下、右、上各裁剪1厘米
		clip % 应用裁剪
		]{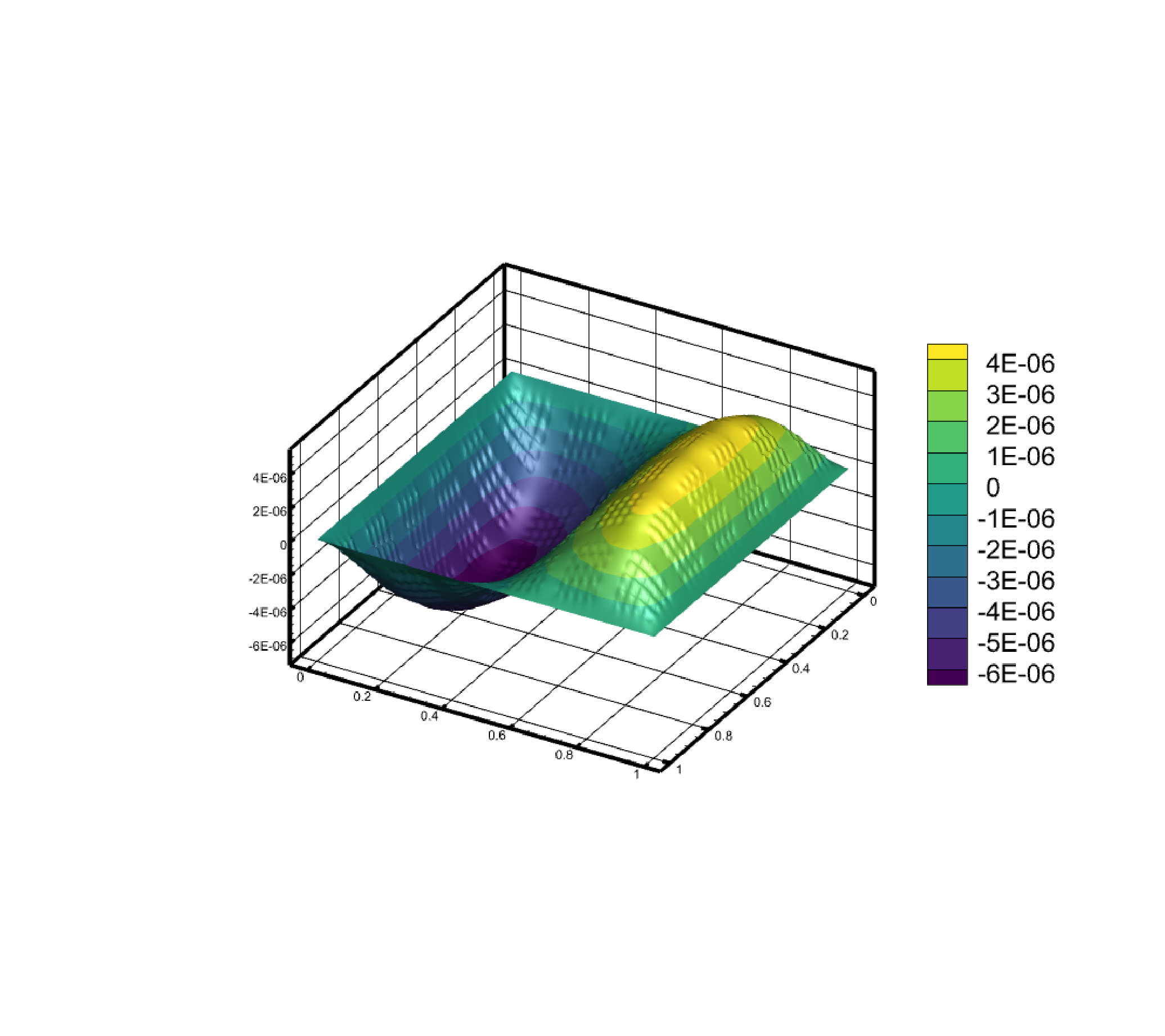} \\
		(a)
	\end{minipage}
	\begin{minipage}[b]{0.35\textwidth}
		\centering
		\includegraphics[
		width=4cm,
		trim=2.5cm 6.5cm 0cm 6.5cm, % 从左、下、右、上各裁剪1厘米
		clip % 应用裁剪
		]{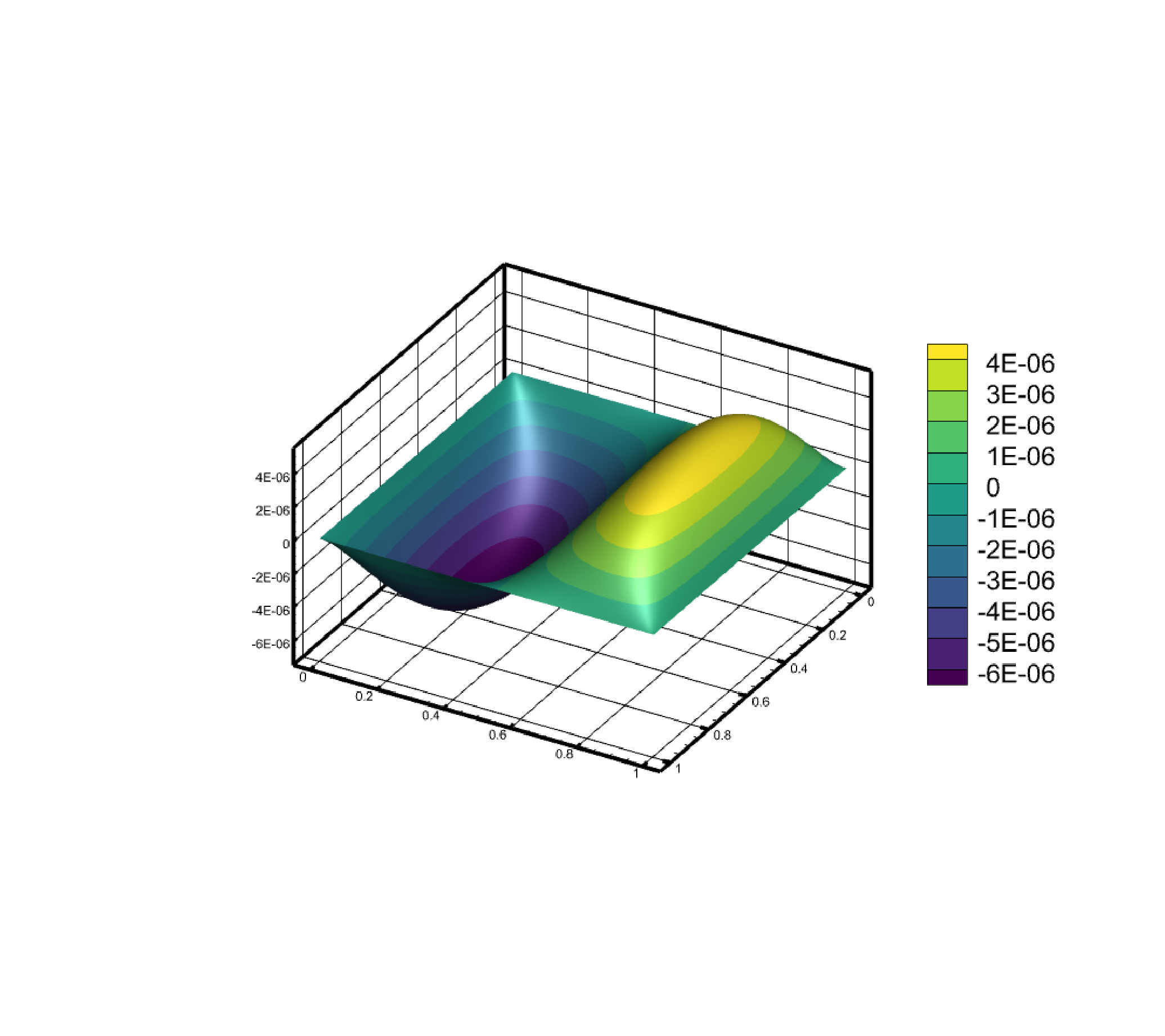} \\
		(b)
	\end{minipage}
	\begin{minipage}[b]{0.32\textwidth}
		\centering
		\includegraphics[
		width=4cm,
		trim=2.5cm 6.5cm 0cm 6.5cm, % 从左、下、右、上各裁剪1厘米
		clip % 应用裁剪
		]{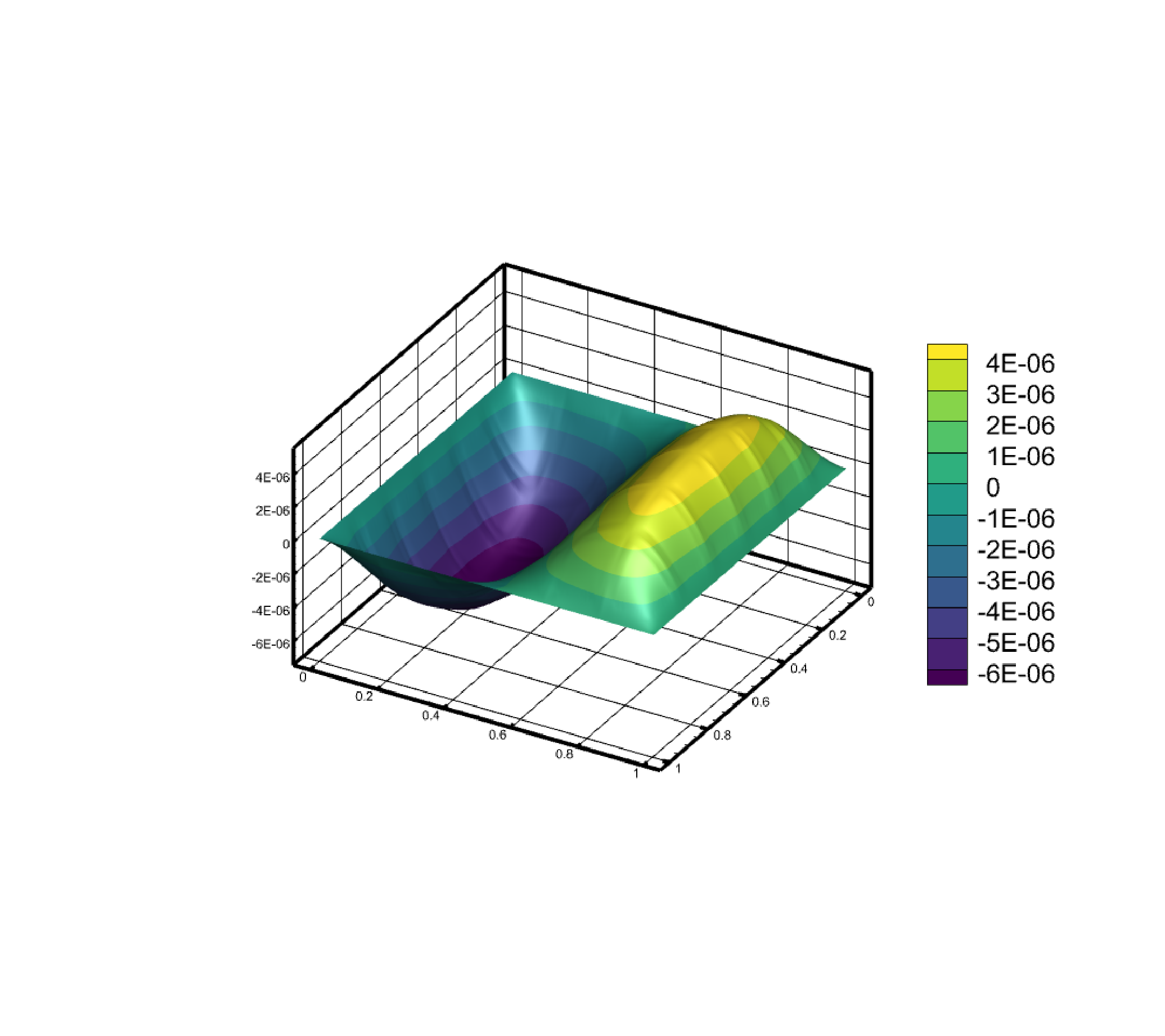} \\
		(c)
	\end{minipage}
	\begin{minipage}[b]{0.32\textwidth}
		\centering
		\includegraphics[
		width=4cm,
		trim=2.5cm 6.5cm 0cm 6.5cm, % 从左、下、右、上各裁剪1厘米
		clip % 应用裁剪
		]{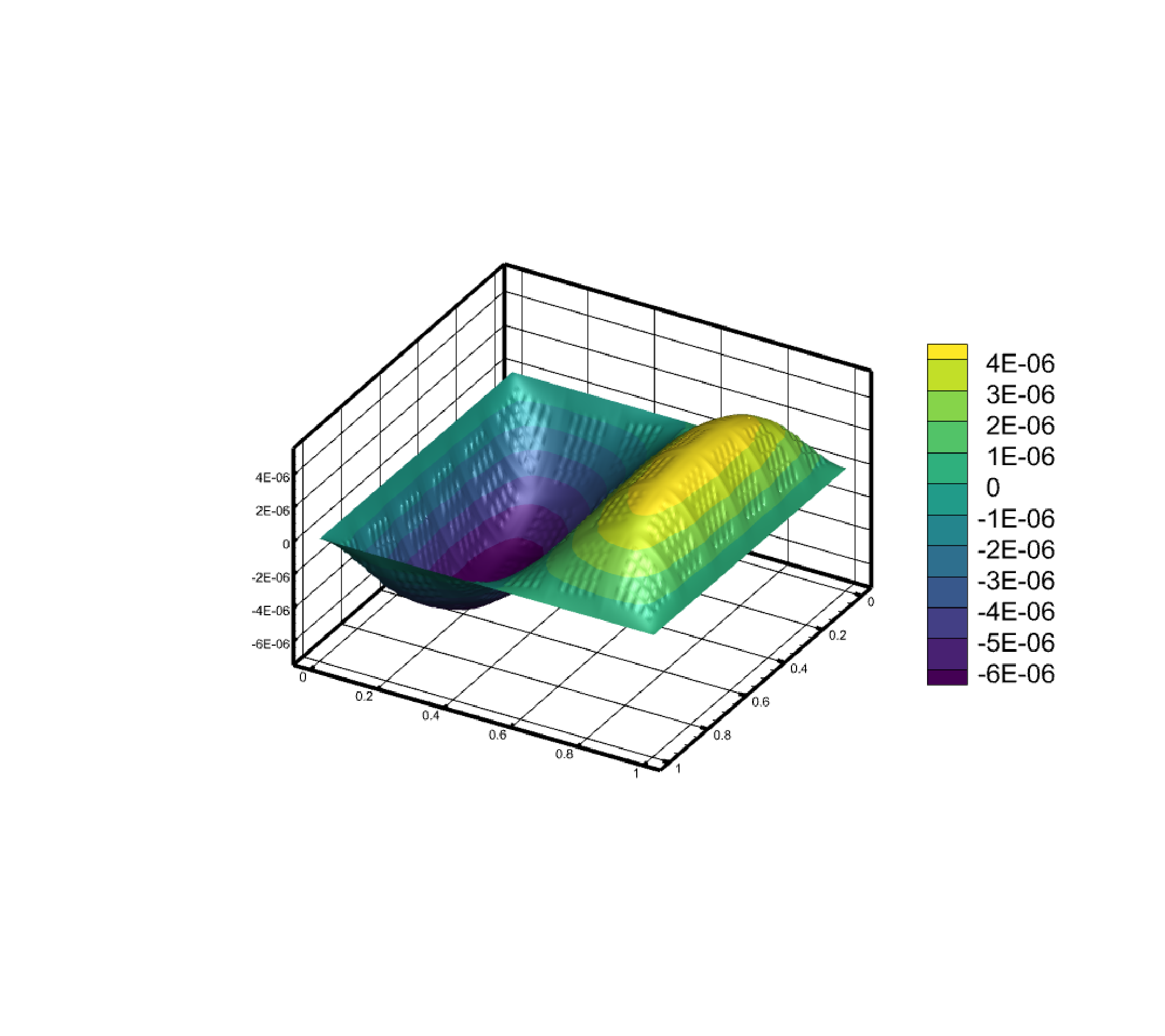} \\
		(d)
	\end{minipage}
	\begin{minipage}[b]{0.32\textwidth}
		\centering
		\includegraphics[
		width=4cm,
		trim=2.5cm 6.5cm 0cm 6.5cm, % 从左、下、右、上各裁剪1厘米
		clip % 应用裁剪
		]{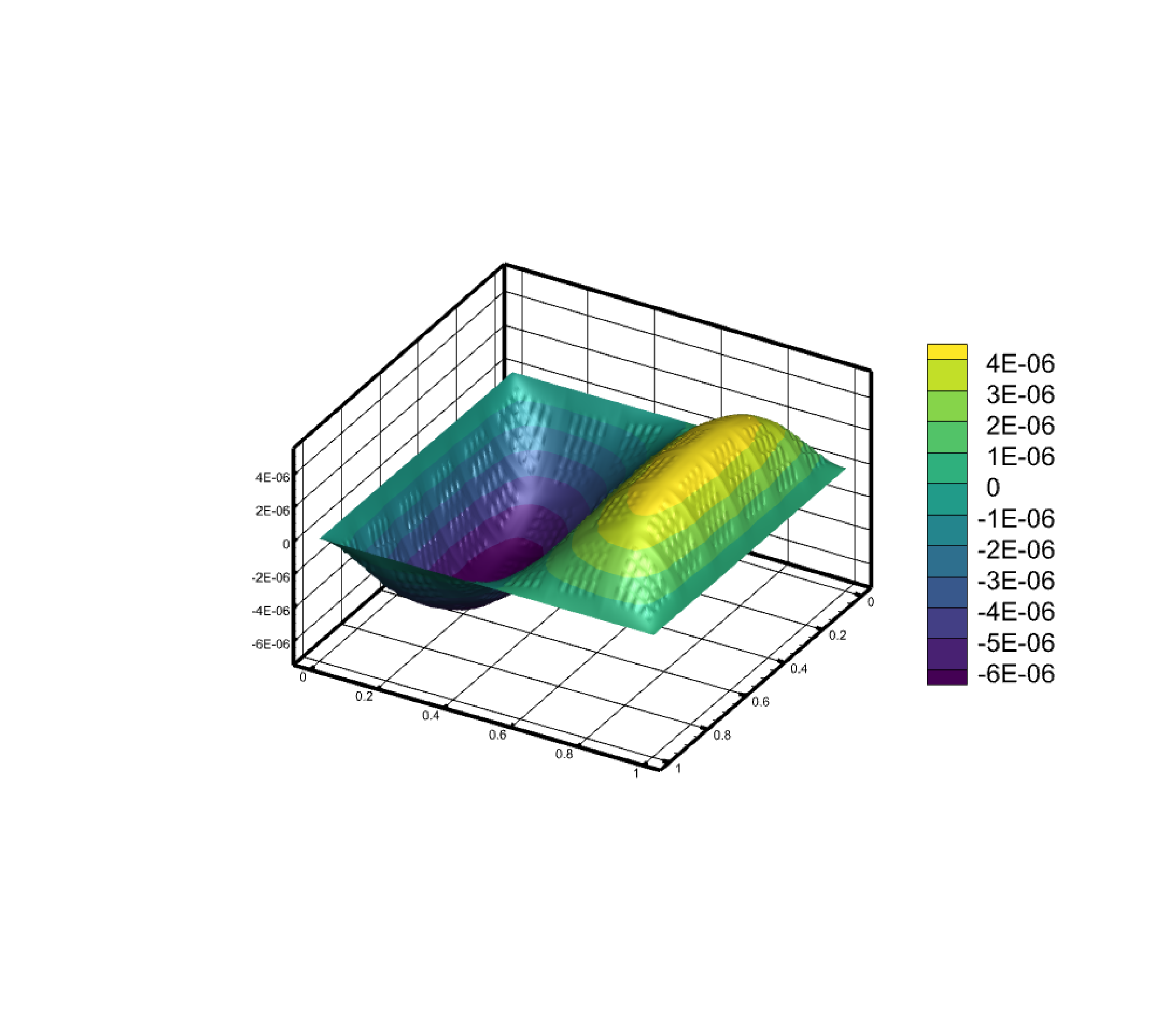} \\
		(e)
	\end{minipage}
	\caption{The second component for the displacement field: (a) ${u _{2Ref}^{{\zeta  _1}{\zeta  _2}}}$; (b) ${{u_{2}^{(0)}}}$;
		(c) ${u_{2D}^{{\zeta  _1}}}$; (d) ${u_{2t}^{{\zeta  _1}{\zeta  _2}}}$; (e) ${u_{2T}^{{\zeta  _1}{\zeta  _2}}}$.}\label{E2f4}
\end{figure}

\begin{figure}[!htb]
	\centering
	\begin{minipage}[b]{0.4\textwidth}
		\centering
		\includegraphics[width=4cm]{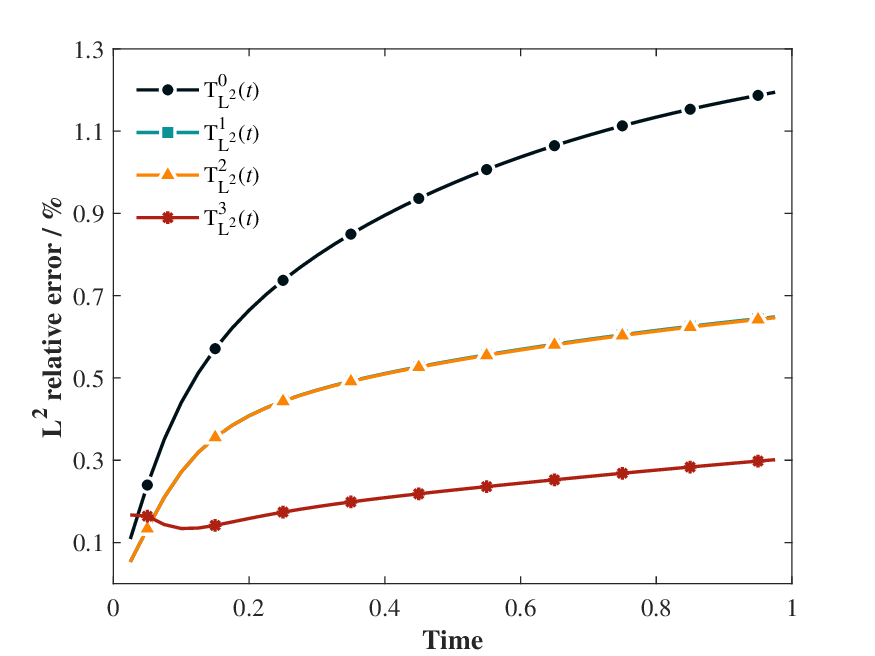} \\
		(a)
	\end{minipage}
	\begin{minipage}[b]{0.4\textwidth}
		\centering
		\includegraphics[width=4cm]{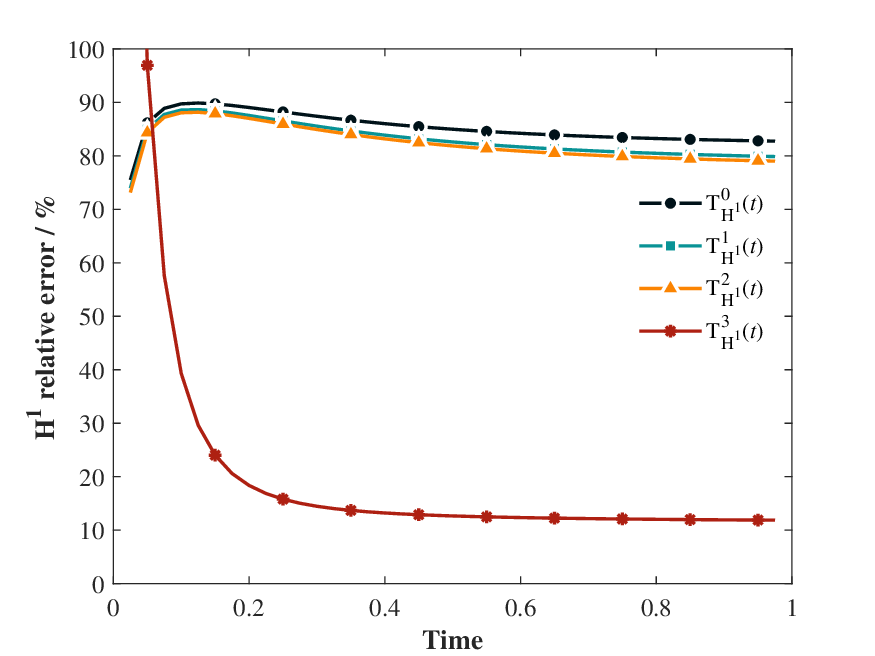} \\
		(b)
	\end{minipage}
	\begin{minipage}[b]{0.4\textwidth}
		\centering
		\includegraphics[width=4cm]{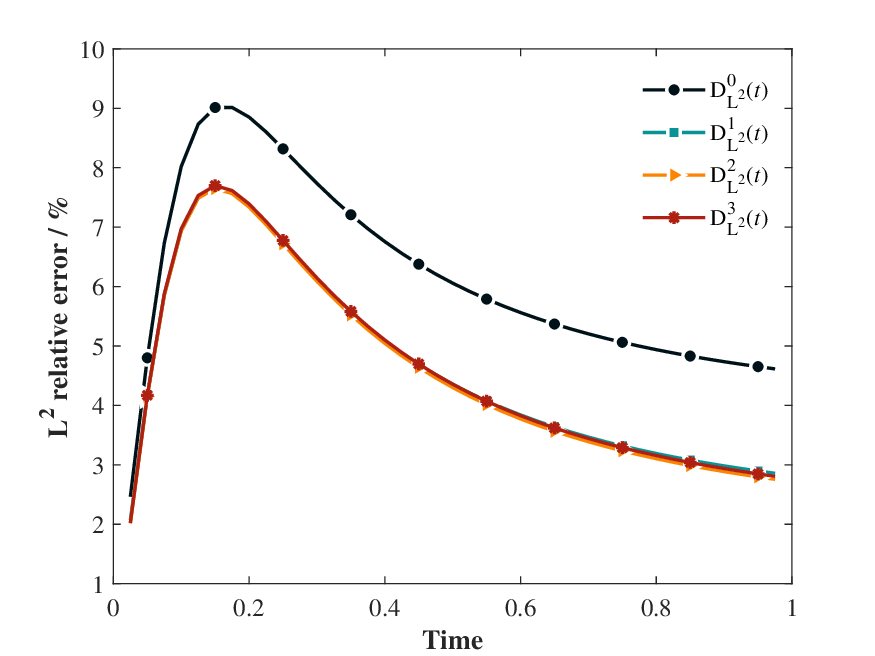} \\
		(c)
	\end{minipage}
	\begin{minipage}[b]{0.4\textwidth}
		\centering
		\includegraphics[width=4cm]{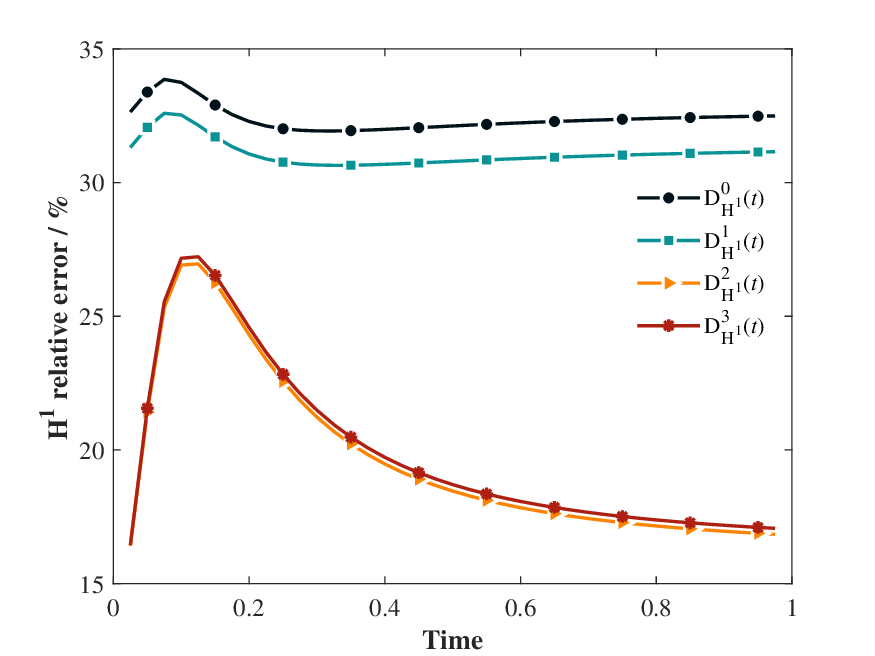} \\
		(d)
	\end{minipage}
	\caption{The evolutive relative errors of temperature and displacement fields: (a) ${\rm{T_{{L^2}}}}(t)$; (b) ${\rm{T_{{H^1}}}}(t)$; (c) ${\rm{D_{{L^2}}}}(t)$; (d) ${\rm{D_{{H^1}}}}(t)$.}\label{E2f5}
\end{figure}

Then, after numerical calculation for the novel 2D heterogeneous structure II, the simulation results of distinct multi-scale solutions at time $t=1.0s$ for physical fields are displayed in Figs.\hspace{1mm}\ref{E2f2}-\ref{E2f4}. Furthermore, Fig.\hspace{1mm}\ref{E2f5} visualizes the evolutive relative errors of physical fields with distinct multi-scale solutions.

To sum up, as shown in Table \ref{E2t1}, our HOTS method demonstrates a substantial reduction in both computational time and storage cost compared to the reference FEM, thereby accelerating the simulation process and saving the computer storage resource. From Figs.\hspace{1mm}\ref{E2f2}-\ref{E2f4}, we can conclude that our HOTS approach can precisely acquire the highly microscopic oscillating information, while the homogenized solutions, the two-scale solutions and the lower-order three-scale solutions can capture the macroscopic information without oscillations, the mesoscopic information behavior without microscopic oscillations and the low-frequency microscopic oscillating information, respectively. Moreover, as plainly depicted in Fig.\hspace{1mm}\ref{E2f5}, compared to the other three multi-scale solutions, the proposed HOTS solutions possess the optimal numerical accuracy, and associated two-stage numerical algorithm is reliable even if it runs for a long time.
\subsection{Example 3: 3D composite block structure with multiple spatial scales}
\begin{figure}[!htb]
	\centering
	\begin{minipage}[c]{0.63\textwidth}
		\centering
		\includegraphics[
		width=0.9\linewidth,
		trim=2cm 3.7cm 2cm 3.2cm, % 从左、下、右、上各裁剪1厘米
		clip % 应用裁剪
		]{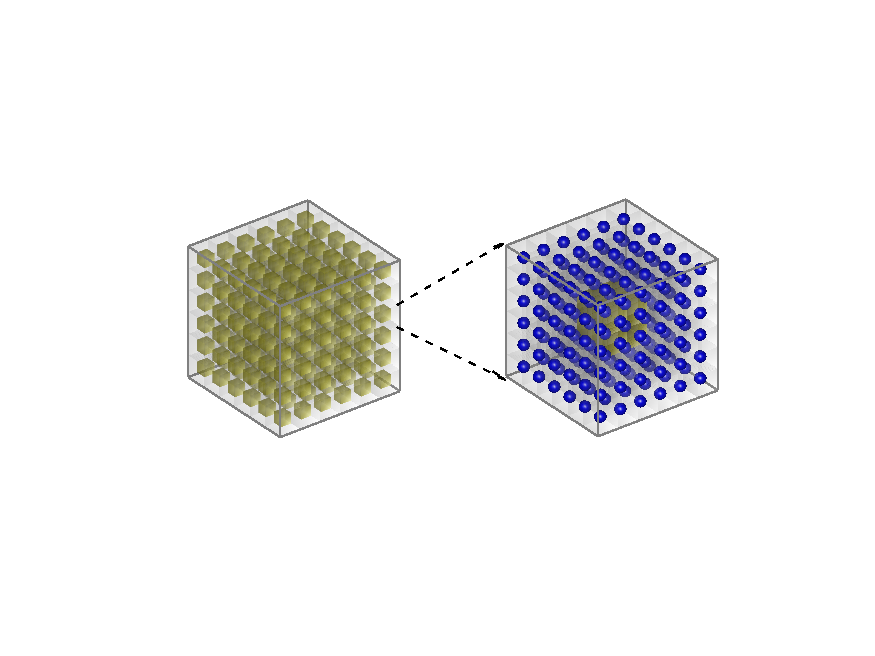} \\
		(a)
	\end{minipage}
	\begin{minipage}[c]{0.3\textwidth}
		\centering
		\includegraphics[
		width=0.9\linewidth,
		trim=2cm 1.3cm 1.8cm 1.3cm, % 从左、下、右、上各裁剪1厘米
		clip % 应用裁剪
		]{E1macro.eps} \\
		(b)
	\end{minipage}
		\begin{minipage}[c]{0.3\textwidth}
		\centering
		\includegraphics[
		width=0.9\linewidth,
		trim=2cm 1.1cm 2cm 0.5cm, % 从左、下、右、上各裁剪1厘米
		clip % 应用裁剪
		]{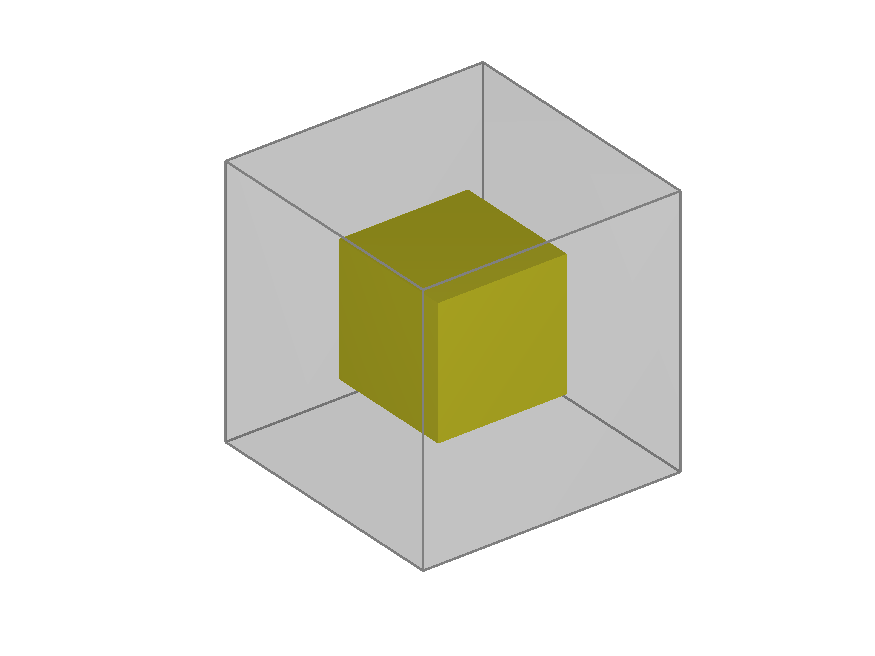} \\
		(c)
	\end{minipage}
		\begin{minipage}[c]{0.3\textwidth}
		\centering
		\includegraphics[
		width=0.9\linewidth,
		trim=2cm 1.1cm 2cm 0.5cm, % 从左、下、右、上各裁剪1厘米
		clip % 应用裁剪
		]{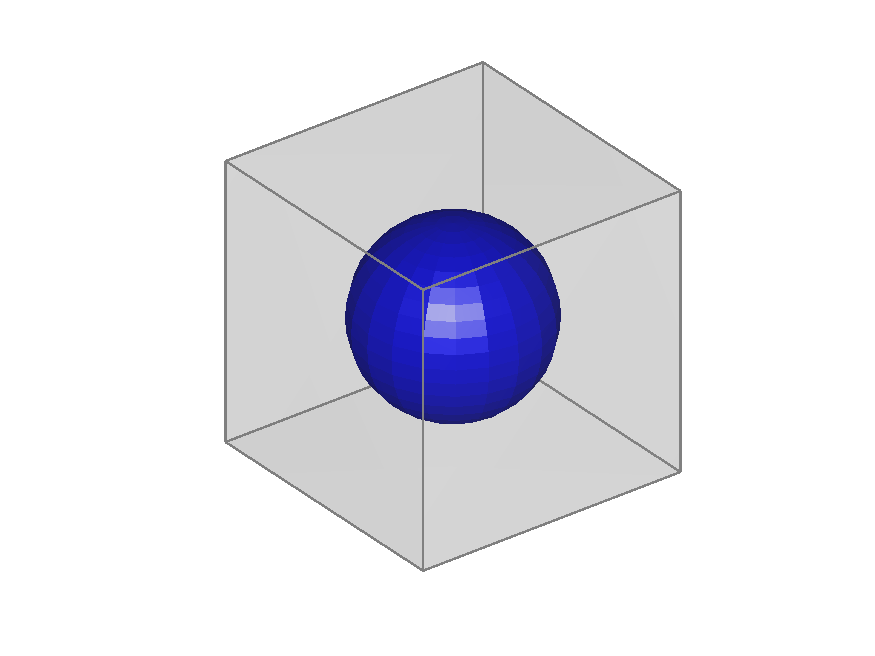} \\
		(d)
	\end{minipage}
	\caption{(a) The heterogeneous block domain $\Omega$; (b) The cross section $x_3=0.097$cm of composite block;(c) Mesoscopic unit cell $Y$; (d) Microscopic unit cell $ Z$.}\label{E3f1}
\end{figure}

This experiment studies the nonlinear thermo-mechanical behaviors of a 3D heterogeneous block structure. The detailed three-scale domain $\Omega$, mesoscopic UC $Y$ and microscopic UC $Z$ of the 3D composite block structure are shown in Fig.\hspace{1mm}\ref{E3f1}, where small periodic parameters $\zeta _1=1/6$ and $\zeta _2=1/36$, and $\Omega=(x_1,x_2,x_3)=[0,1]\times [0,1]\times[0,1]\mathrm{cm^3}$.

In this experiment, the composite block is clamped on its six surfaces, and then the boundary temperature is kept at 373.15K on its six surfaces. Furthermore, source items and initial conditions in nonlinear three-scale problem (1) are given as below.
\begin{equation}
	\begin{aligned}
		&h(\bm{x},t)=15000.0\mathrm{J/(cm^3\cdot s)},\\
        &f_i(\bm{x},t)=(0,0,-15000.0\mathrm{N/cm^3}),\;\;\text{in}\;\;\Omega \times (0,T),\\
		&\tilde{\theta}=373.15\mathrm{K},\;{\bm{u}}^{0}=0,\;{\bm{u}}^{1}(\bm{x})=0,\;\;\text{in}\;\;\Omega.
	\end{aligned}
\end{equation}
In addition, the thermal conductivity $k_{ij}^{{\zeta  _1}{\zeta  _2}}$ of three component material is defined as $1000+0.1\theta+10^{-5}\theta^2 / 10+ 10^{-3}\theta +10^{-7}\theta^2 / 200+0.02\theta+2 \times 10^{-7}\theta^2 \rm{(W/(m\cdot K))}$ and other material property parameters are the same as those of Example 1.

As for 3D heterogeneous block structure, there is a massive amount of microscopic unit cells, which is a reason that direct FE simulation cannot obtain the reference FE solutions, as shown in Table \ref{E3t1}. Next, Table \ref{E3t1} displayed the computational meshes of this heterogeneous block and it is undeniable that the presented HOTS approach can significantly reduce computer memory usage while preserving high numerical accuracy.
\begin{table}[!htb]
	\centering
	\caption{Summary of computational cost ($\Delta t = 0.01, t \in [0,1]$s).}\label{E3t1}
	\begin{tabular}{lcccc}
		\hline
		& \multirow{2}{*}{Reference FEM} & \multicolumn{3}{c}{HOTS method} \\
		\cmidrule(l){3-5}
		& & $Z$ & $Y$ & $\Omega$ \\
		\hline
		Number of nodes & $\approx$464,553,792(estimated) & 9,957 & 14,416 & 7,558,272 \\
		Number of elements & $\approx$86,593,536(estimated) & 1,856 & 2,465 & 1,295,029 \\
		\hline
	\end{tabular}
\end{table}

\begin{figure}[!htb]
	\centering
	\begin{minipage}[b]{0.32\textwidth}
		\centering
		\includegraphics[
		width=5cm,
		trim=4cm 5.9cm 0cm 5.5cm, % 从左、下、右、上各裁剪1厘米
		clip % 应用裁剪
		]{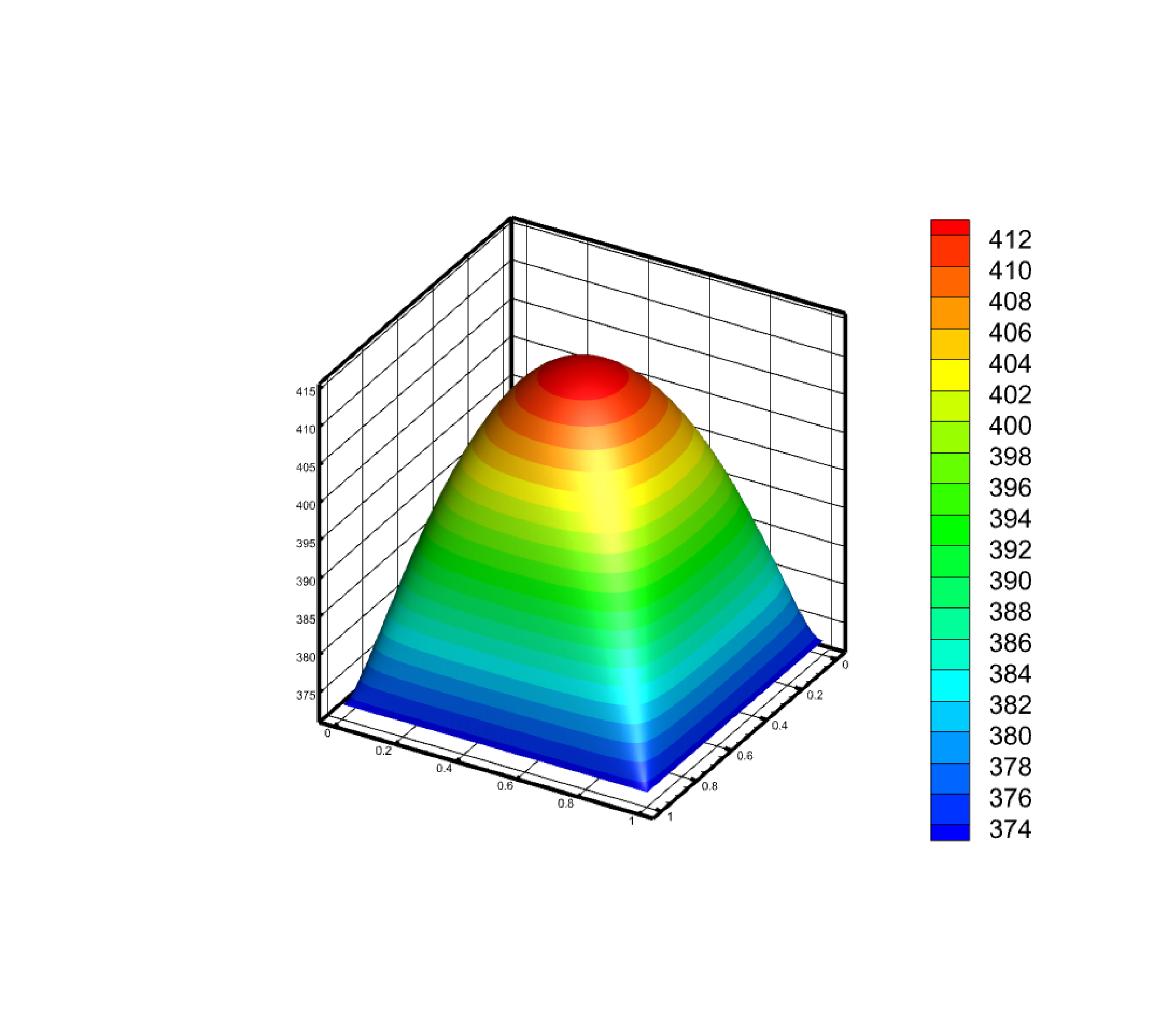} \\
		(a)
	\end{minipage}
	\begin{minipage}[b]{0.32\textwidth}
		\centering
		\includegraphics[
		width=5cm,
		trim=4cm 5.9cm 0cm 5.5cm, % 从左、下、右、上各裁剪1厘米
		clip % 应用裁剪
		]{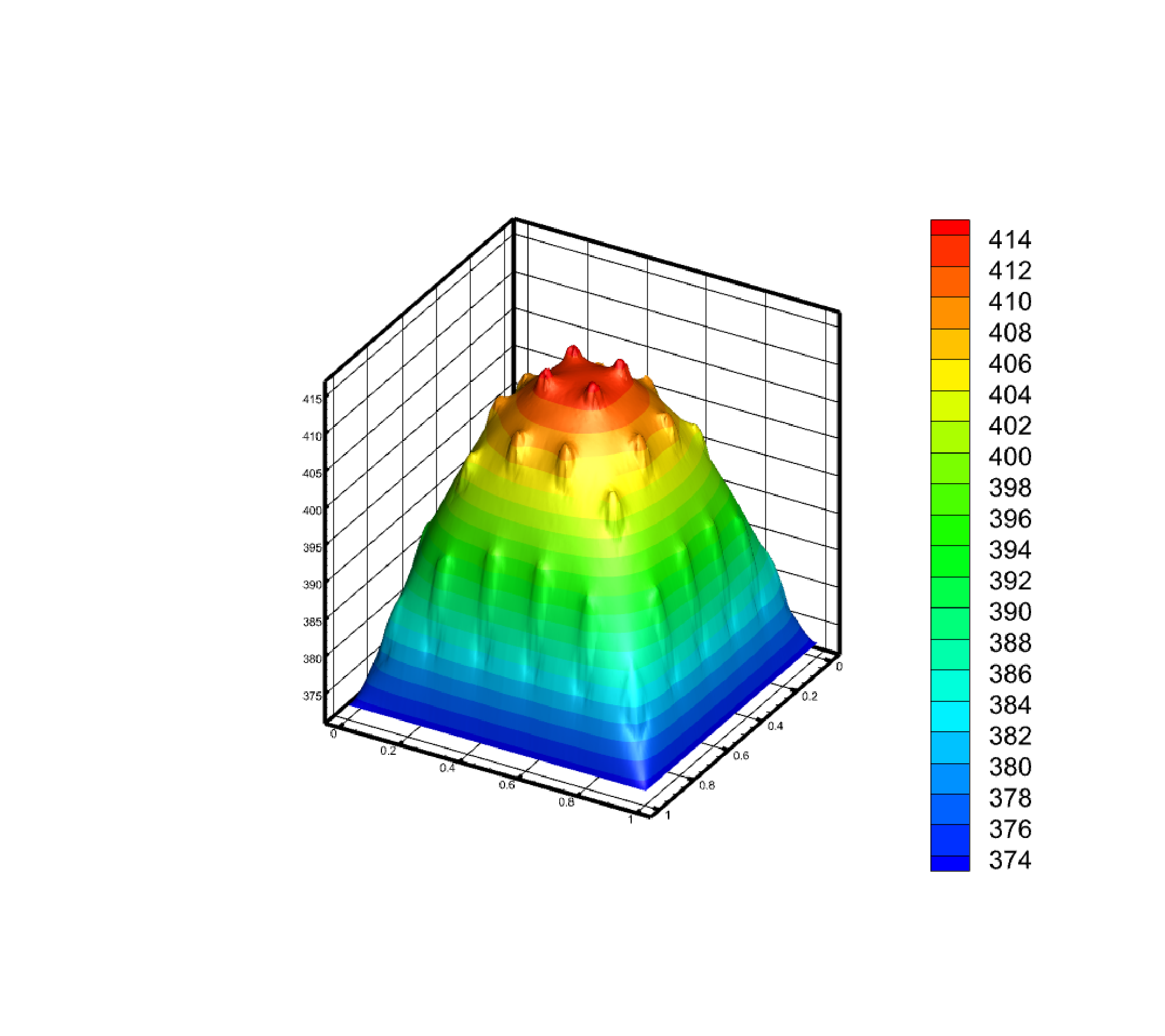} \\
		(b)
	\end{minipage}
	\begin{minipage}[b]{0.32\textwidth}
		\centering
		\includegraphics[
		width=5cm,
		trim=4cm 5.9cm 0cm 5.5cm, % 从左、下、右、上各裁剪1厘米
		clip % 应用裁剪
		]{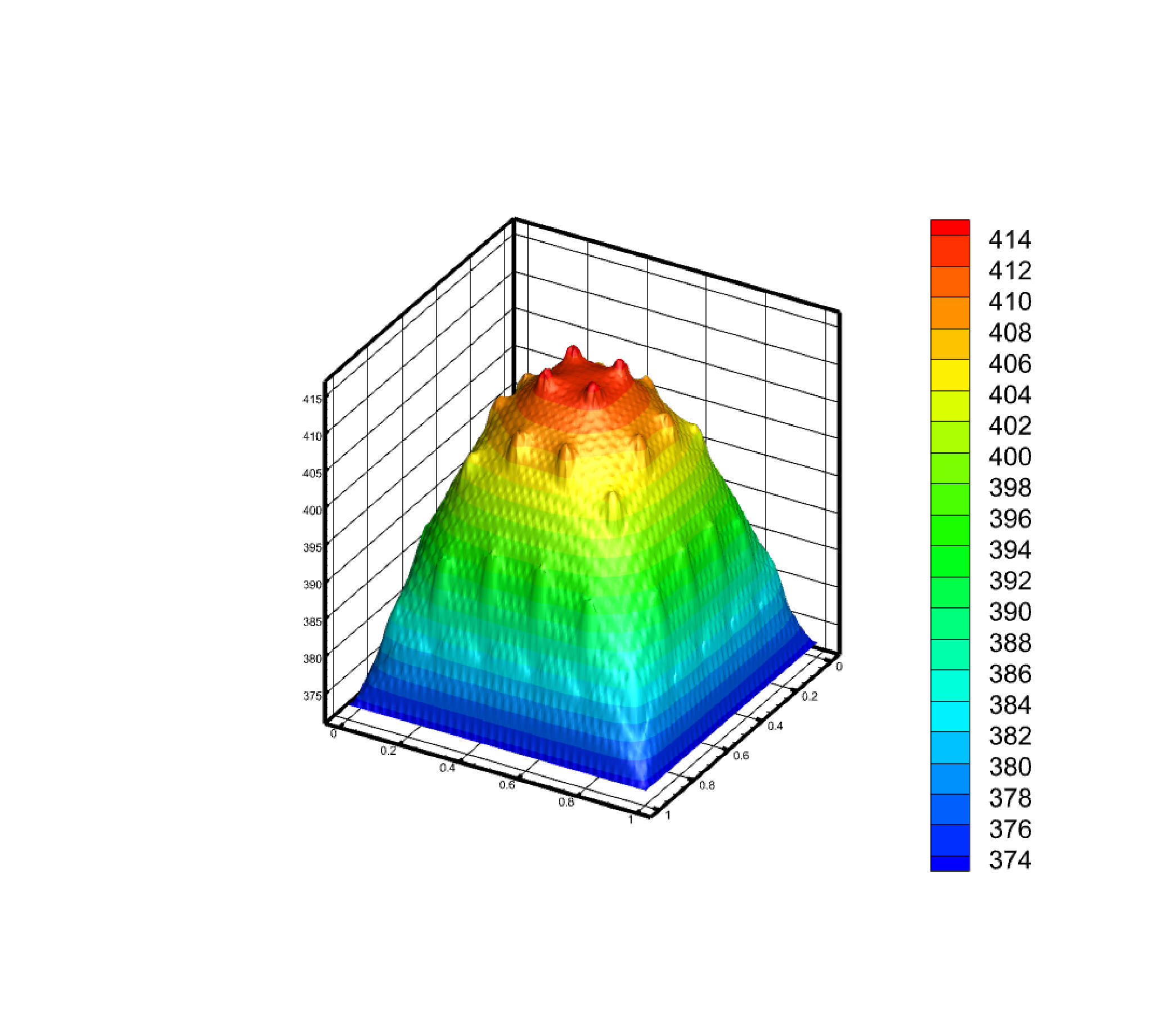} \\
		(c)
	\end{minipage}
	\begin{minipage}[b]{0.48\textwidth}
		\centering
		\includegraphics[
		width=5cm,
		trim=4cm 5.9cm 0cm 5.5cm, % 从左、下、右、上各裁剪1厘米
		clip % 应用裁剪
		]{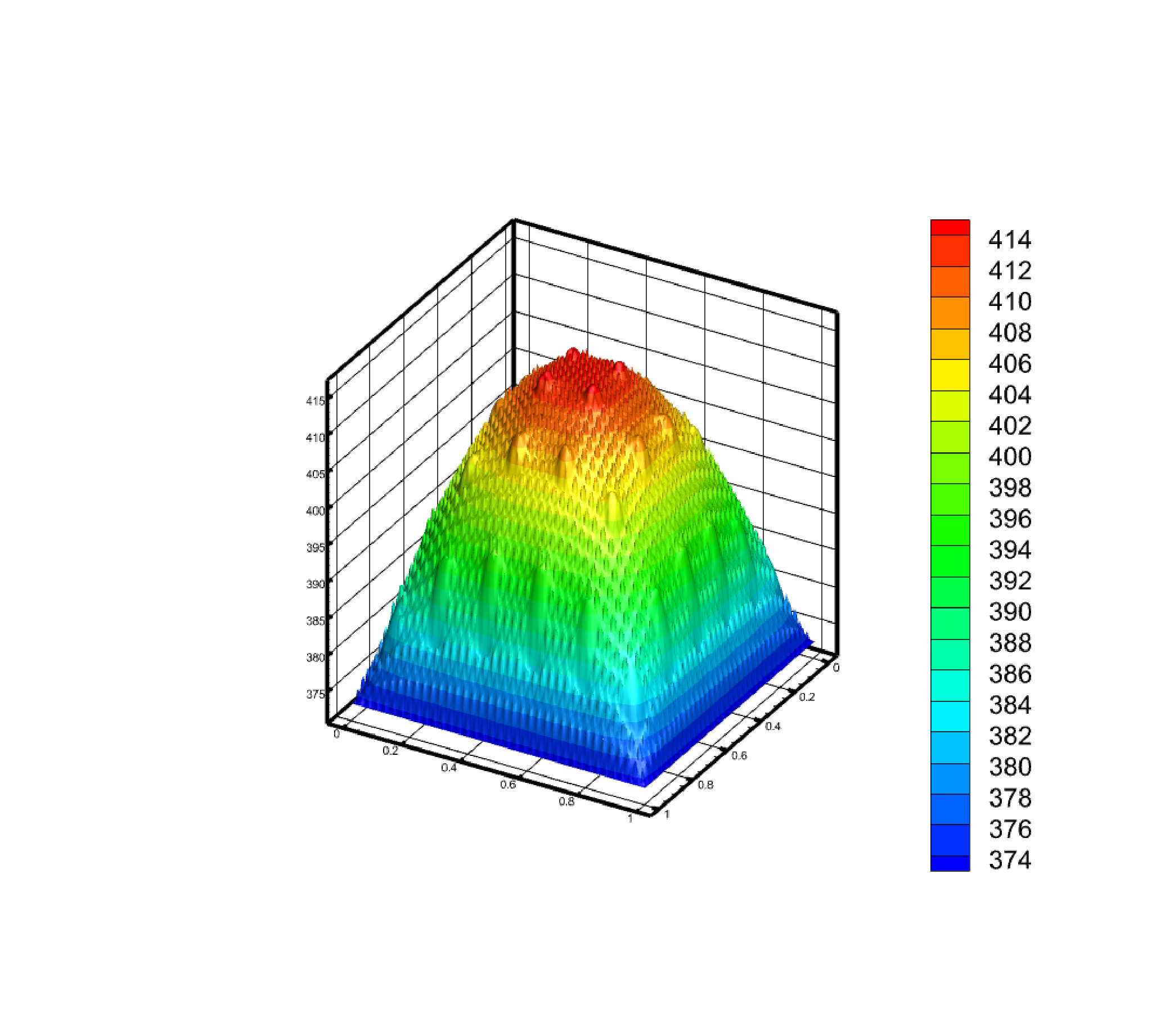} \\
		(d)
	\end{minipage}
    \begin{minipage}[b]{0.48\textwidth}
		\centering
		\includegraphics[
		width=7cm,
		trim=0cm 0cm 0cm 0cm, % 从左、下、右、上各裁剪1厘米
		clip % 应用裁剪
		]{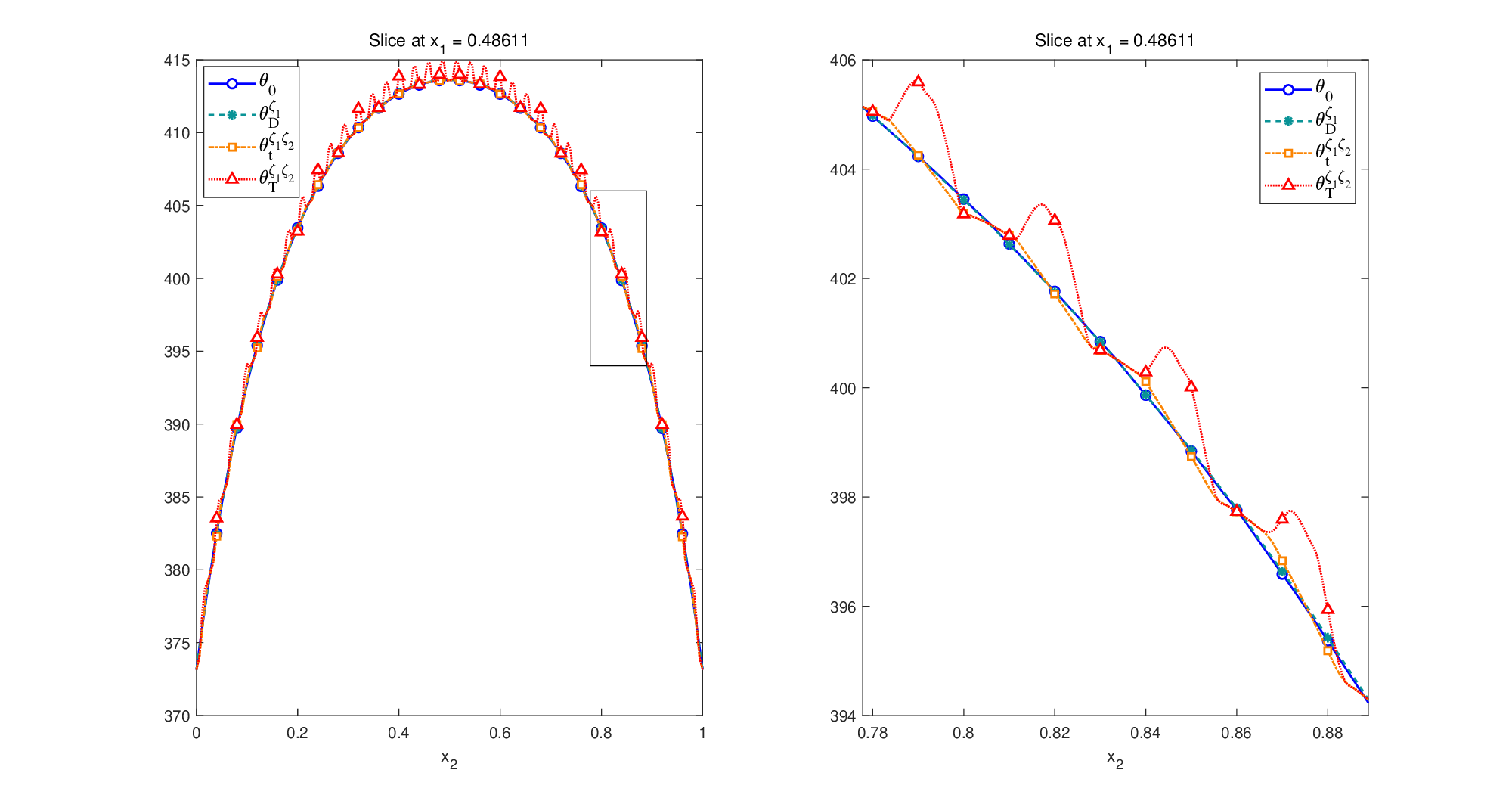} \\
		(e)
	\end{minipage}
	\caption{The temperature field in $x_3 = 0.097$cm: (a) $\theta_{0}$; (b) ${\theta _D^{{\zeta  _1}}}$; (c) ${\theta _t^{{\zeta  _1}{\zeta  _2}}}$; (d) ${\theta _T^{{\zeta  _1}{\zeta  _2}}}$; (e) Computational results on line $x_1 = 0.48611$cm and $x_3 = 0.097$cm.}\label{E3f2}
\end{figure}

\begin{figure}[!htb]
	\centering
	\begin{minipage}[b]{0.4\textwidth}
		\centering
		\includegraphics[
		width=5cm,
		trim=2.5cm 6.3cm 0cm 5.5cm, % 从左、下、右、上各裁剪1厘米
		clip % 应用裁剪
		]{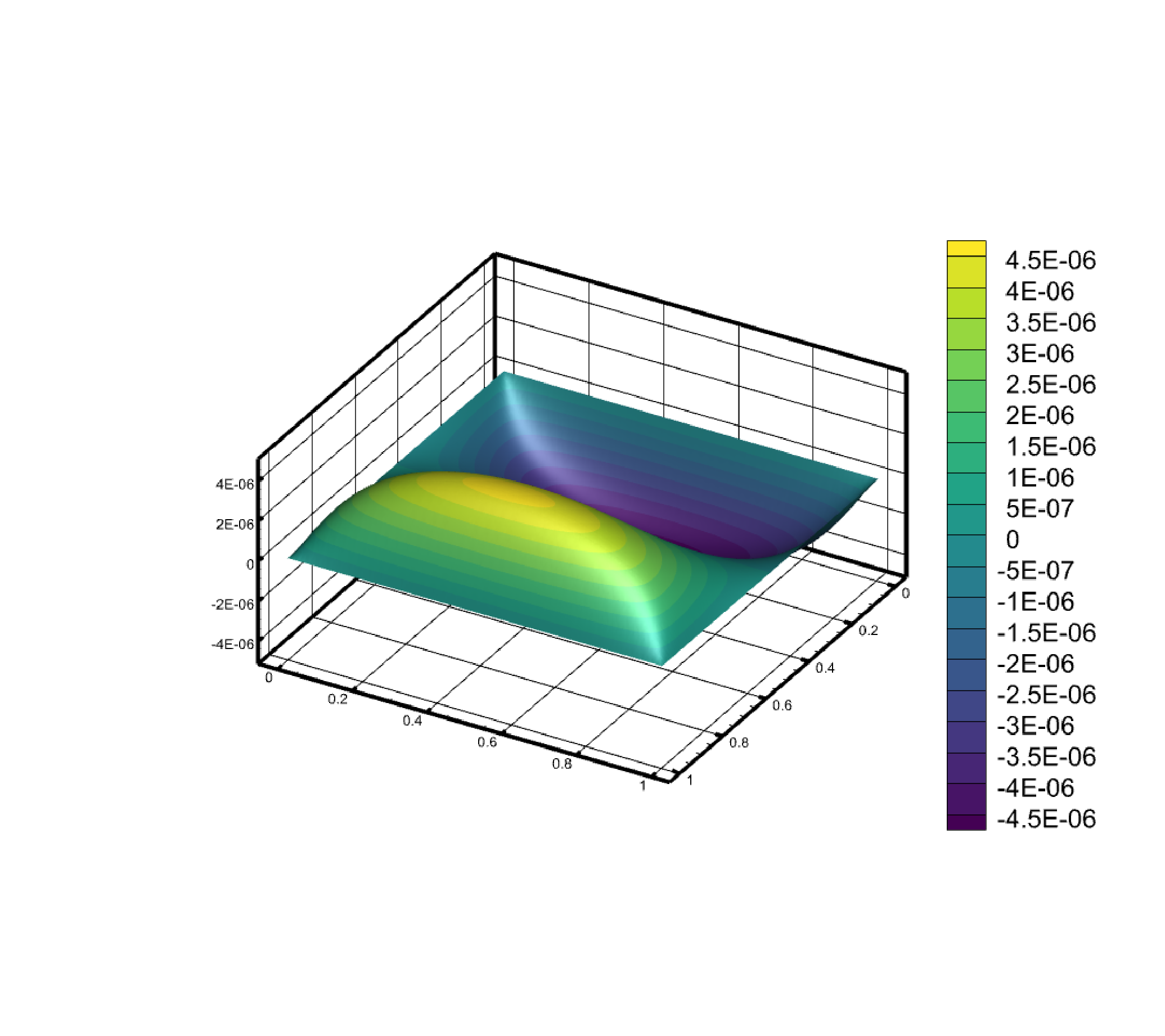} \\
		(a)
	\end{minipage}
	\begin{minipage}[b]{0.4\textwidth}
		\centering
		\includegraphics[
		width=5cm,
		trim=2.5cm 6.3cm 0cm 5.5cm, % 从左、下、右、上各裁剪1厘米
		clip % 应用裁剪
		]{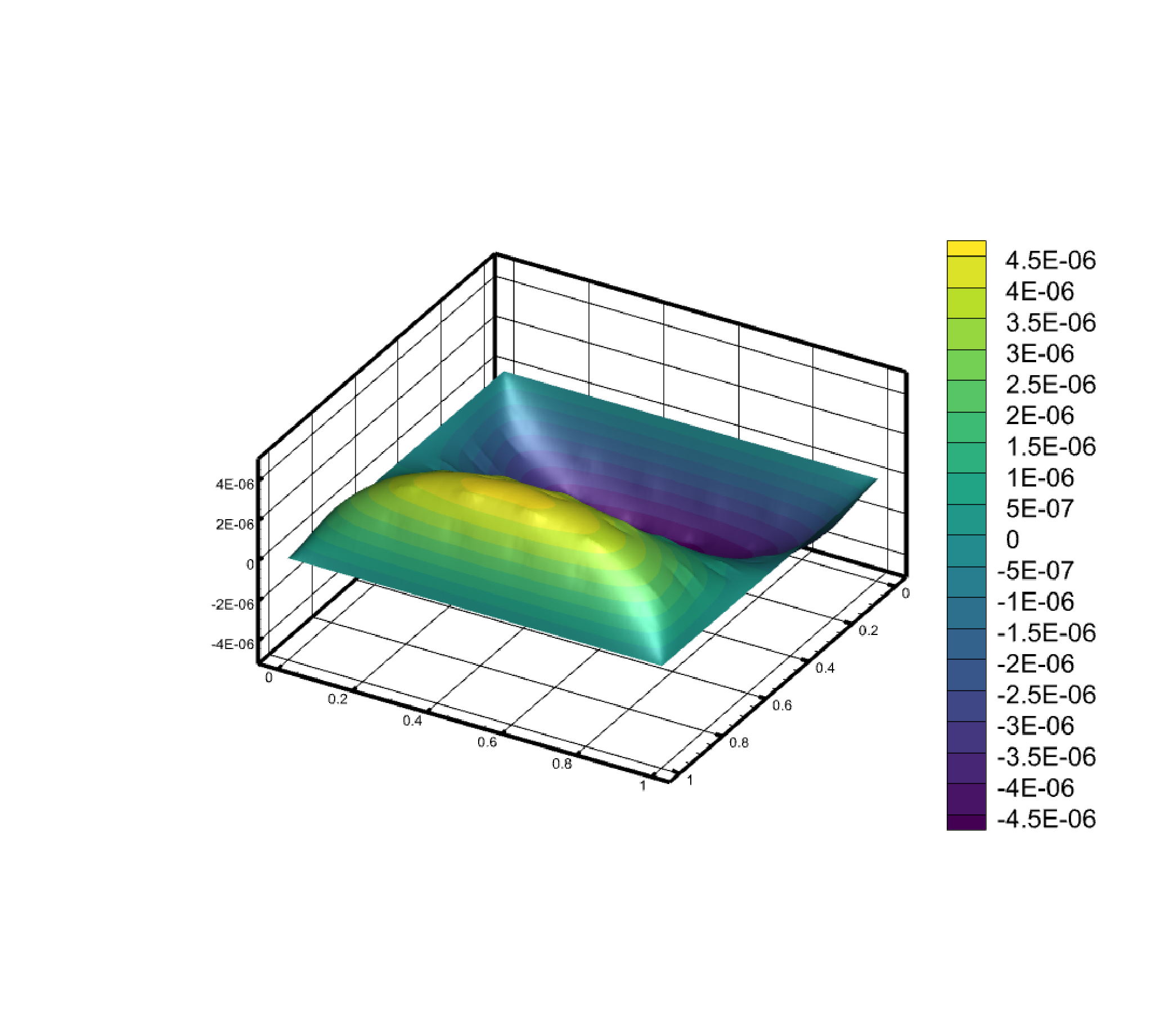} \\
		(b)
	\end{minipage}
	\begin{minipage}[b]{0.4\textwidth}
		\centering
		\includegraphics[
		width=5cm,
		trim=2.5cm 6.3cm 0cm 5.5cm, % 从左、下、右、上各裁剪1厘米
		clip % 应用裁剪
		]{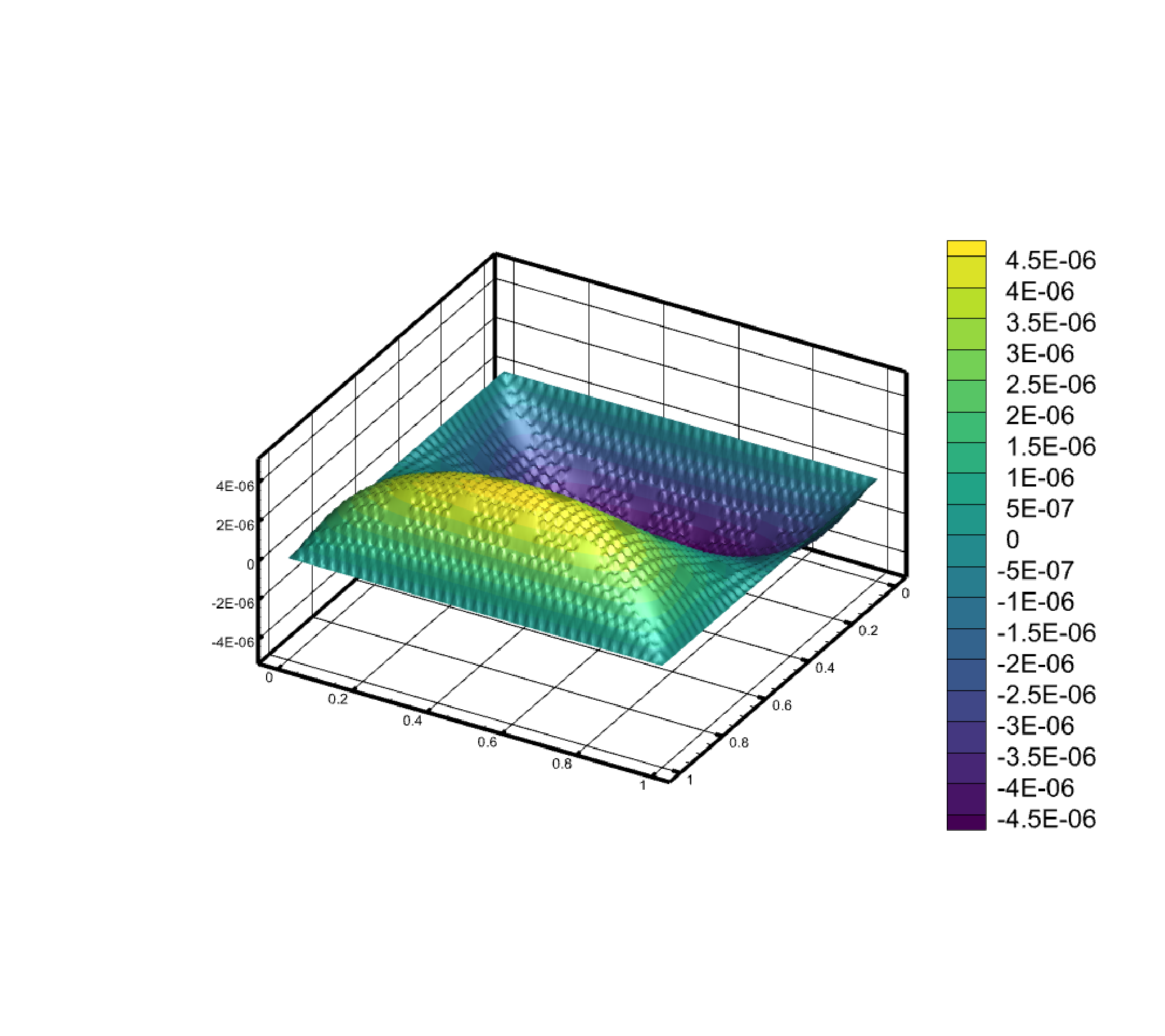} \\
		(c)
	\end{minipage}
	\begin{minipage}[b]{0.4\textwidth}
		\centering
		\includegraphics[
		width=5cm,
		trim=2.5cm 6.3cm 0cm 5.5cm, % 从左、下、右、上各裁剪1厘米
		clip % 应用裁剪
		]{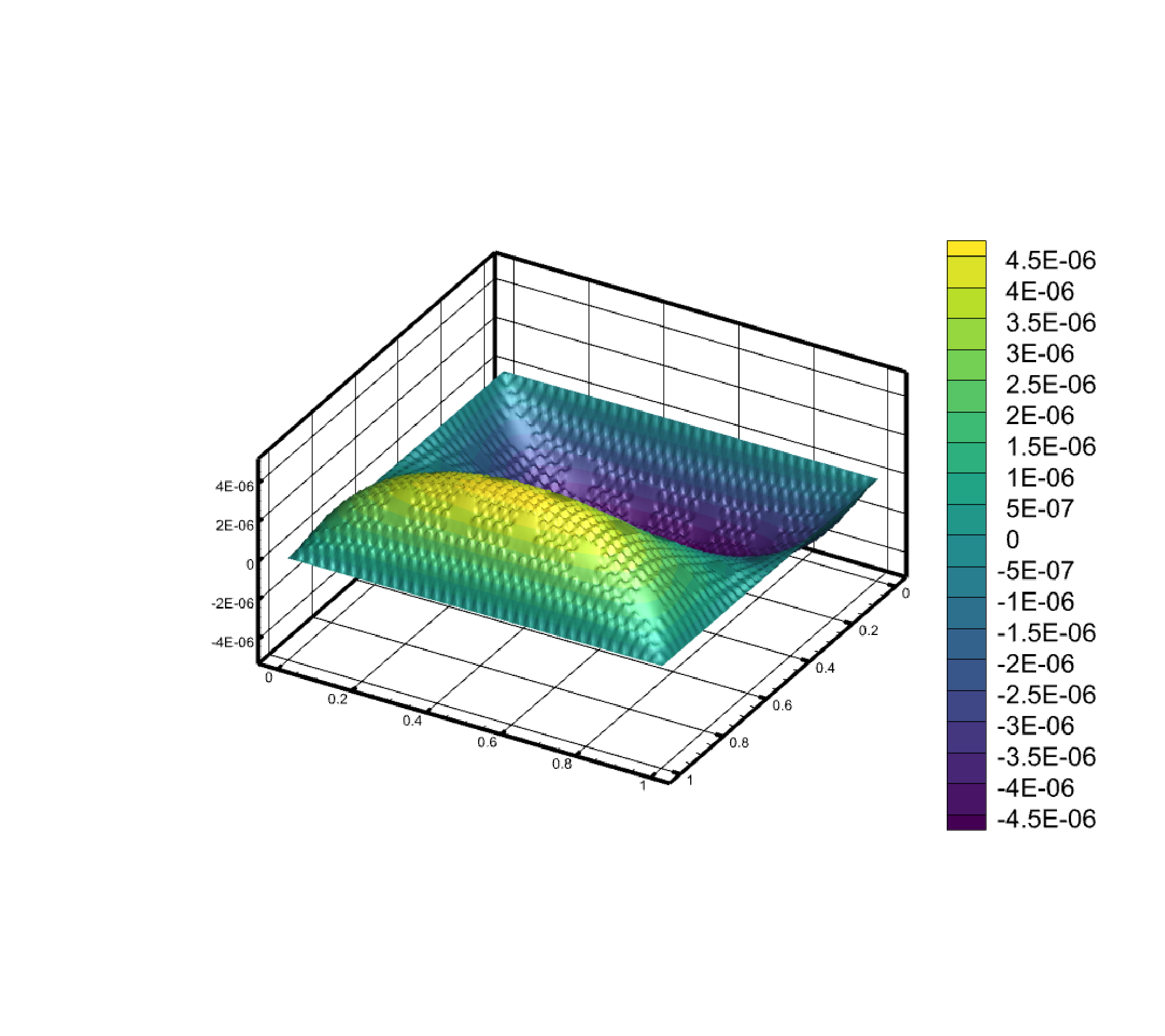} \\
		(d)
	\end{minipage}
	\caption{The first component for the displacement field in $x_3 = 0.097$cm: (a) ${{u_{1}^{(0)}}}$; (b) ${u_{1D}^{{\zeta  _1}}}$; (c) ${u_{1t}^{{\zeta  _1}{\zeta  _2}}}$; (d) ${u_{1T}^{{\zeta  _1}{\zeta  _2}}}$. }\label{E3f3}
\end{figure}

\begin{figure}[!htb]
	\centering
	\begin{minipage}[b]{0.4\textwidth}
		\centering
		\includegraphics[
		width=5cm,
		trim=2.5cm 6.3cm 0cm 5.5cm, % 从左、下、右、上各裁剪1厘米
		clip % 应用裁剪
		]{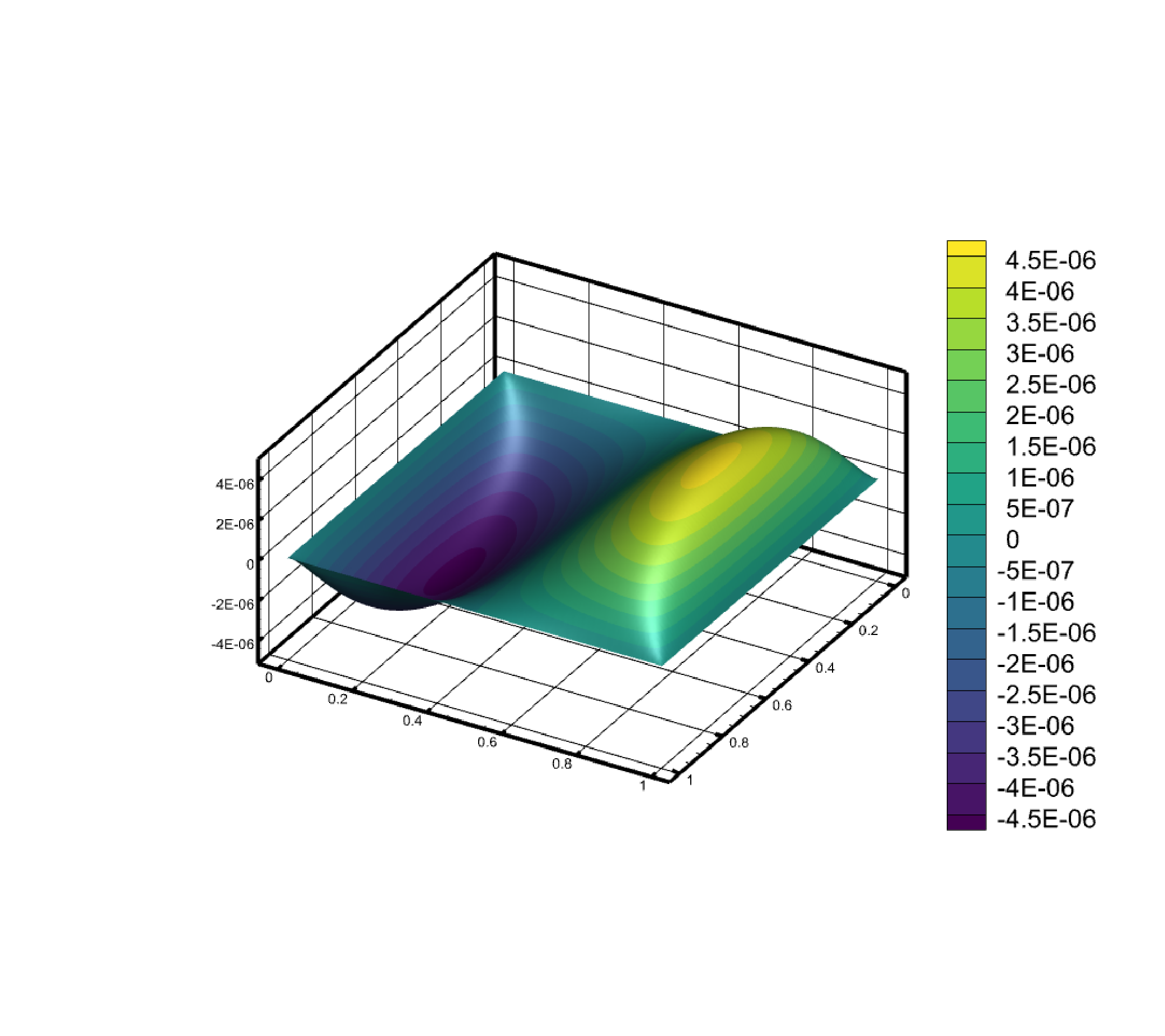} \\
		(a)
	\end{minipage}
	\begin{minipage}[b]{0.4\textwidth}
		\centering
		\includegraphics[
		width=5cm,
		trim=2.5cm 6.3cm 0cm 5.5cm, % 从左、下、右、上各裁剪1厘米
		clip % 应用裁剪
		]{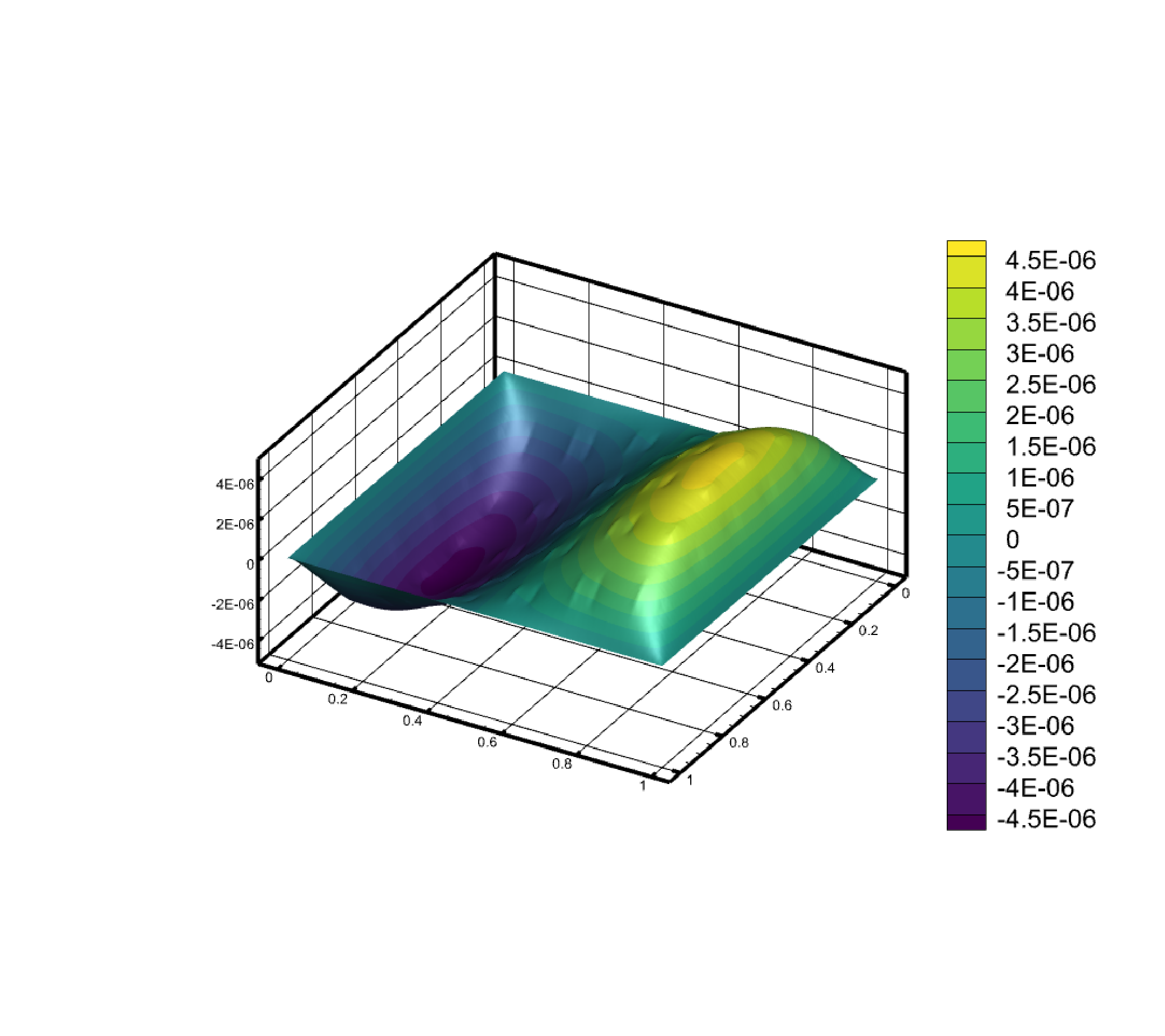} \\
		(b)
	\end{minipage}
	\begin{minipage}[b]{0.4\textwidth}
		\centering
		\includegraphics[
		width=5cm,
		trim=2.5cm 6.3cm 0cm 5.5cm, % 从左、下、右、上各裁剪1厘米
		clip % 应用裁剪
		]{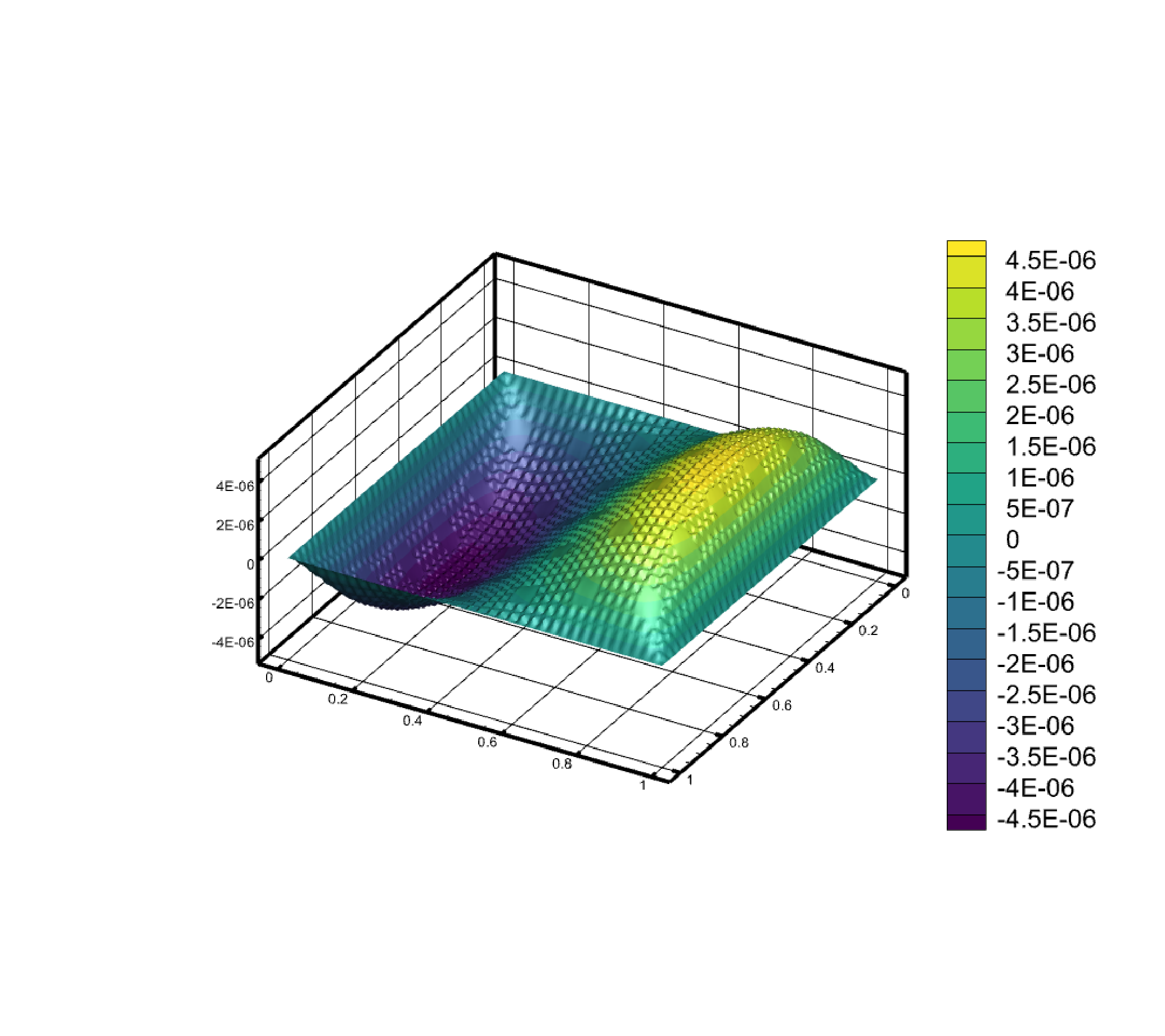} \\
		(c)
	\end{minipage}
	\begin{minipage}[b]{0.4\textwidth}
		\centering
		\includegraphics[
		width=5cm,
		trim=2.5cm 6.3cm 0cm 5.5cm, % 从左、下、右、上各裁剪1厘米
		clip % 应用裁剪
		]{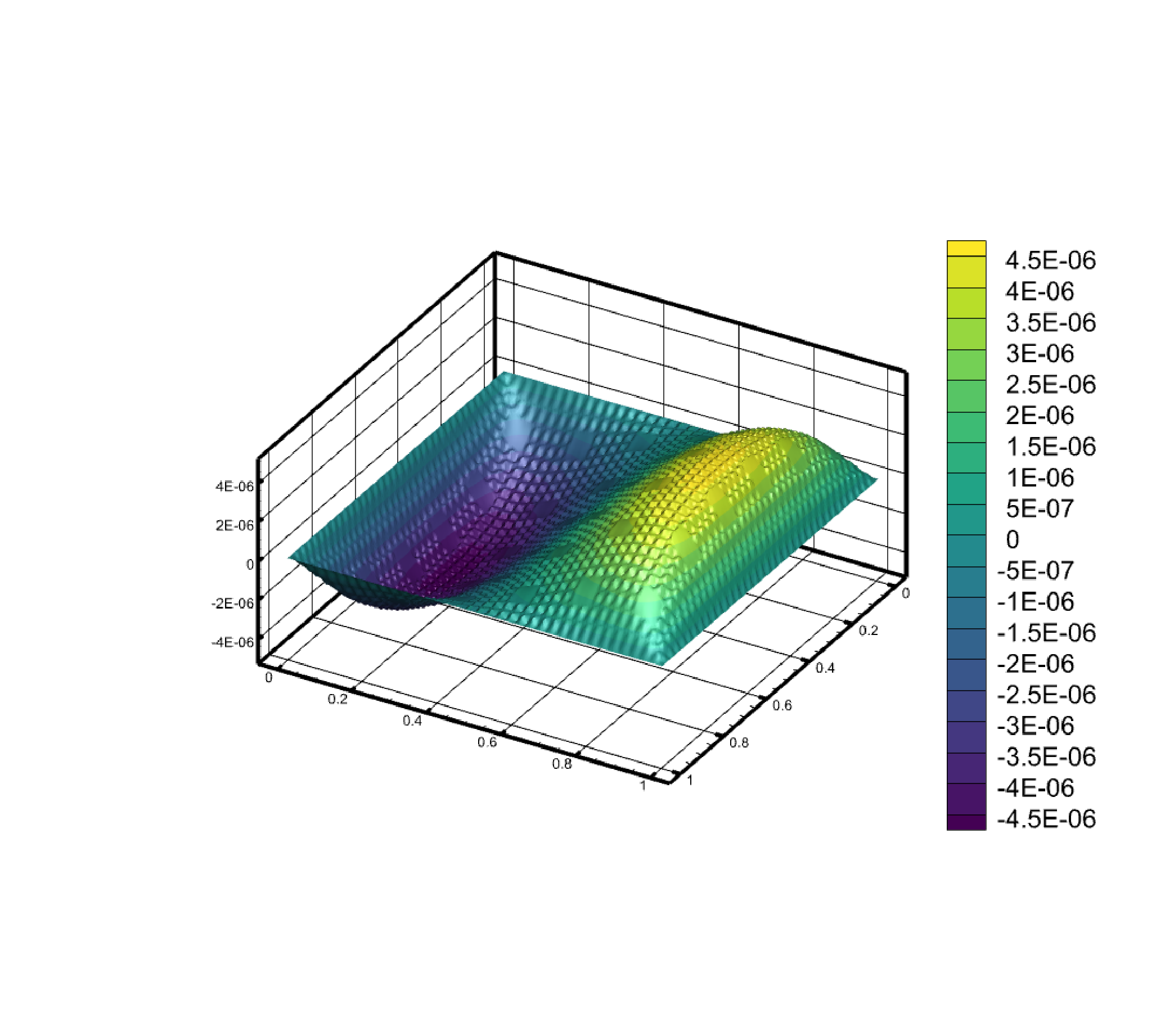} \\
		(d)
	\end{minipage}
	\caption{The second component for the displacement field in $x_3 = 0.097$cm: (a) ${{u_{2}^{(0)}}}$; (b) ${u_{2D}^{{\zeta  _1}}}$; (c) ${u_{2t}^{{\zeta  _1}{\zeta  _2}}}$; (d) ${u_{2T}^{{\zeta  _1}{\zeta  _2}}}$.}\label{E3f4}
\end{figure}

\begin{figure}[!htb]
	\centering
	\begin{minipage}[b]{0.34\textwidth}
		\centering
		\includegraphics[
		width=5cm,
		trim=2.0cm 6.5cm 0cm 5.5cm, % 从左、下、右、上各裁剪1厘米
		clip % 应用裁剪
		]{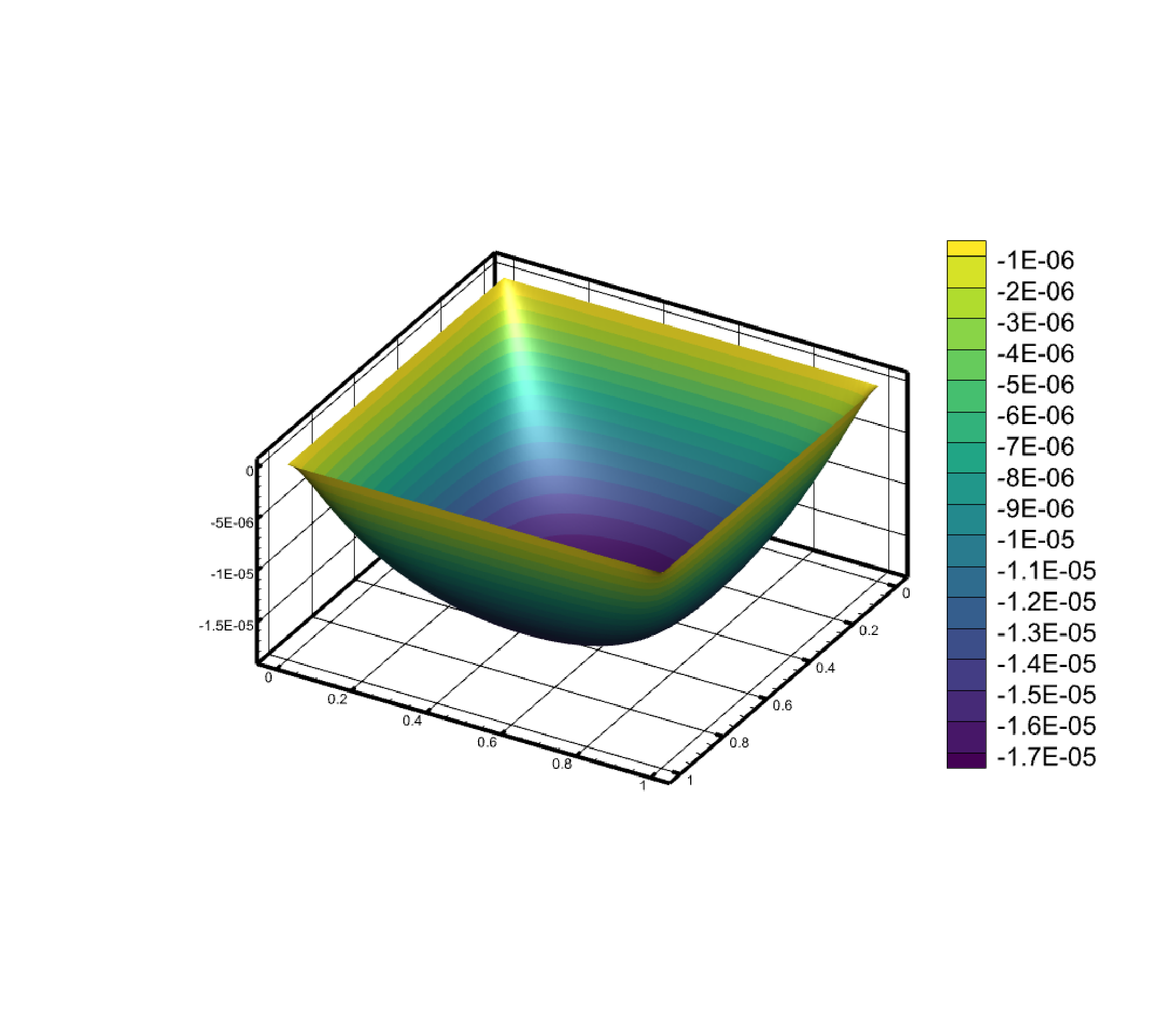} \\
		(a)
	\end{minipage}
	\begin{minipage}[b]{0.34\textwidth}
		\centering
		\includegraphics[
		width=5cm,
		trim=2.0cm 6.5cm 0cm 5.5cm, % 从左、下、右、上各裁剪1厘米
		clip % 应用裁剪
		]{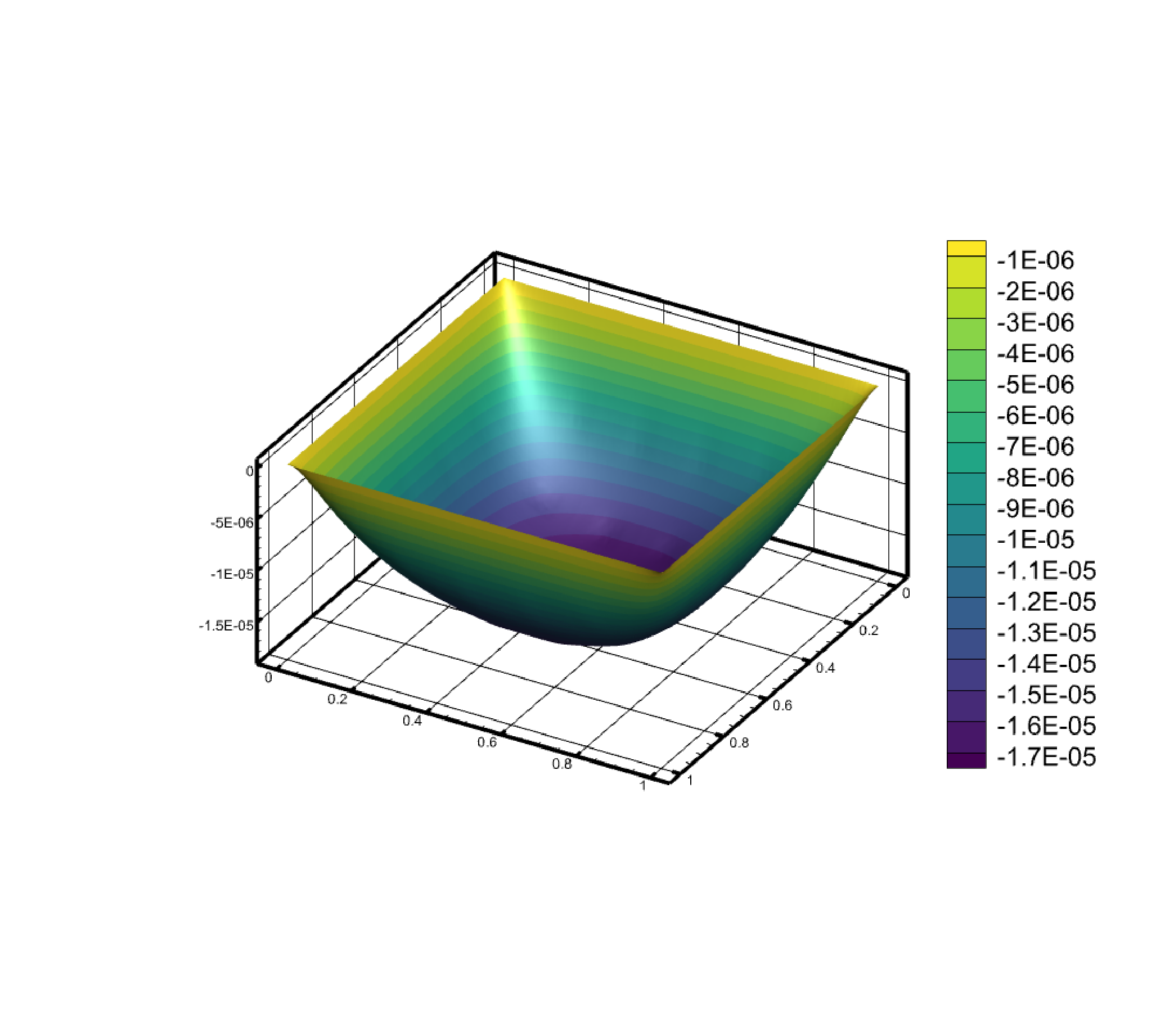} \\
		(b)
	\end{minipage}
	\begin{minipage}[b]{0.30\textwidth}
		\centering
		\includegraphics[
		width=5cm,
		trim=2.0cm 6.5cm 0cm 5.5cm, % 从左、下、右、上各裁剪1厘米
		clip % 应用裁剪
		]{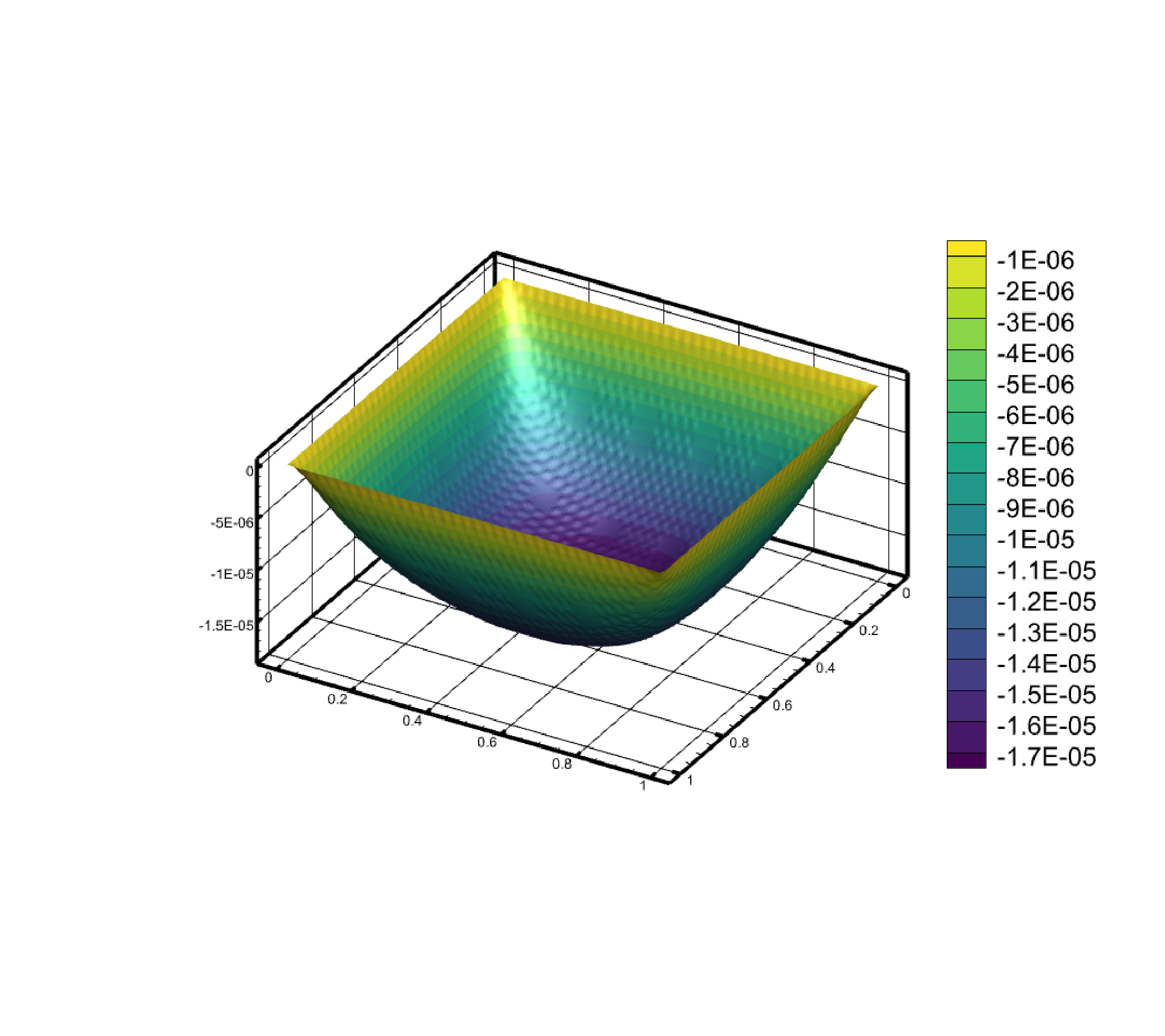} \\
		(c)
	\end{minipage}
	\begin{minipage}[b]{0.48\textwidth}
		\centering
		\includegraphics[
		width=5cm,
		trim=2.5cm 6.5cm 0cm 5.5cm, % 从左、下、右、上各裁剪1厘米
		clip % 应用裁剪
		]{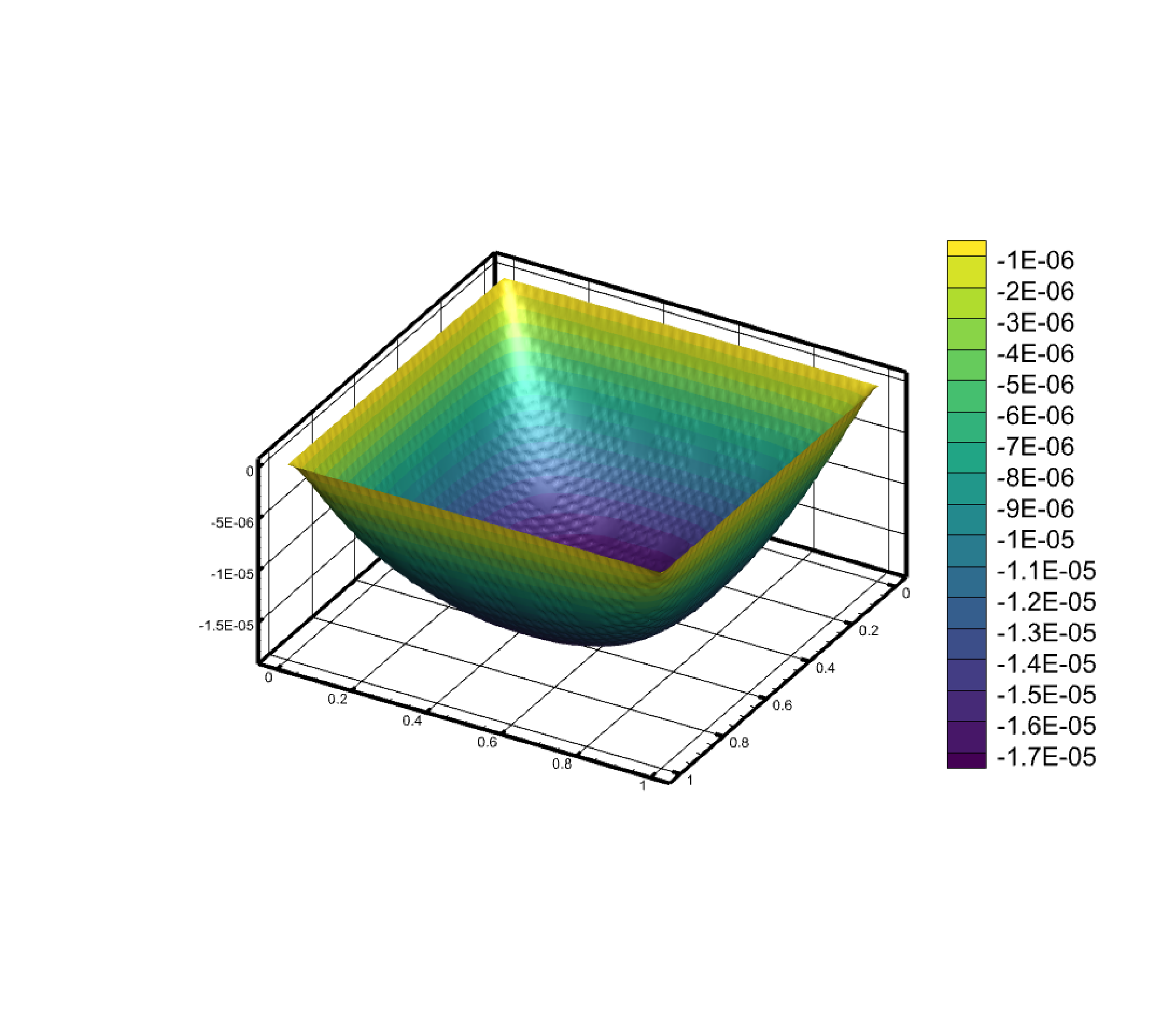} \\
		(d)
	\end{minipage}
    \begin{minipage}[b]{0.48\textwidth}
		\centering
		\includegraphics[
		width=7cm,
		trim=0cm 0cm 0cm 0cm, % 从左、下、右、上各裁剪1厘米
		clip % 应用裁剪
		]{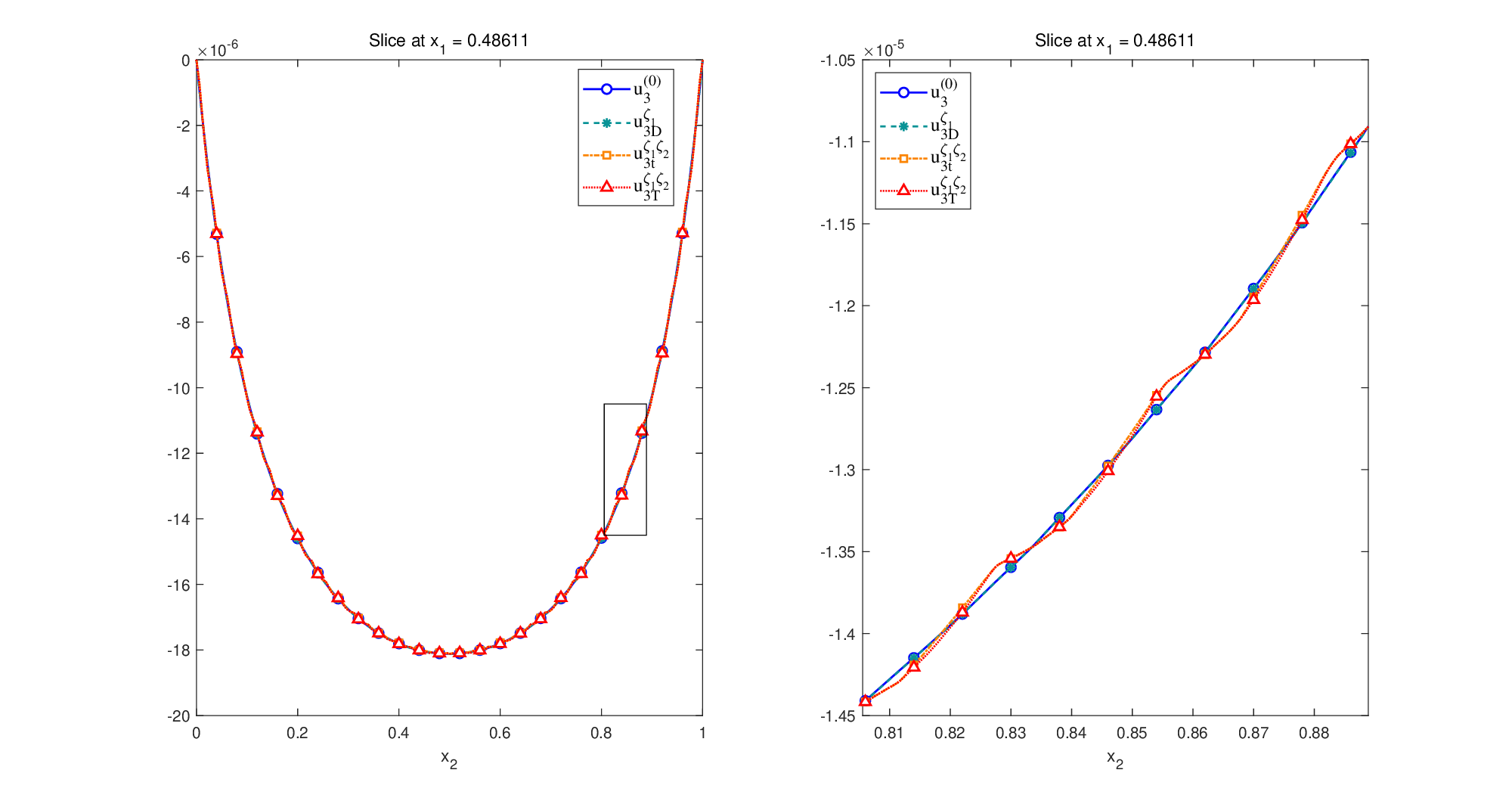} \\
		(e)
	\end{minipage}
	\caption{The third component for the displacement field in $x_3 = 0.097$cm: (a) ${{u_{3}^{(0)}}}$; (b) ${u_{3D}^{{\zeta  _1}}}$; (c) ${u_{3t}^{{\zeta  _1}{\zeta  _2}}}$; (d) ${u_{3T}^{{\zeta  _1}{\zeta  _2}}}$; (e) Computational results on line $x_1 = 0.48611$cm and $x_3 = 0.097$cm.}\label{E3f5}
\end{figure}

Now, we exhibited the distinct multi-scale solutions for temperature and displacement fields at $x_3 =0.097$cm and $t=1.0$s in Figs.\hspace{1mm}\ref{E3f2}-\ref{E3f5}.

As clearly depicted in Figs.\hspace{1mm}\ref{E3f2}-\ref{E3f5}, one can derive that our HOTS method can capture the microscopic oscillatory material behaviors of this 3D heterogeneous block, especially for temperature field. However, the homogenized approach catches only macroscopic oscillatory material behaviors, the two-scale approach characterizes mesoscopic details, and the lower-order three-scale approach provides insufficient fluctuation at the smallest scale. Moreover, it is worth mentioning that the HOTS method exhibits excellent scalability: the computational cost is independent of the number of microscopic UC within the 3D heterogeneous block, which is an ideal alternative for practical large-scale engineering computation.

\subsection{Example 4: 3D composite plate structure with multiple spatial scales}
\begin{figure}[!htb]
	\centering
	\begin{minipage}[c]{0.6\textwidth}
		\centering
		\includegraphics[
		width=0.9\linewidth,
		trim=2cm 3.8cm 2cm 3cm, % 从左、下、右、上各裁剪1厘米
		clip % 应用裁剪
		]{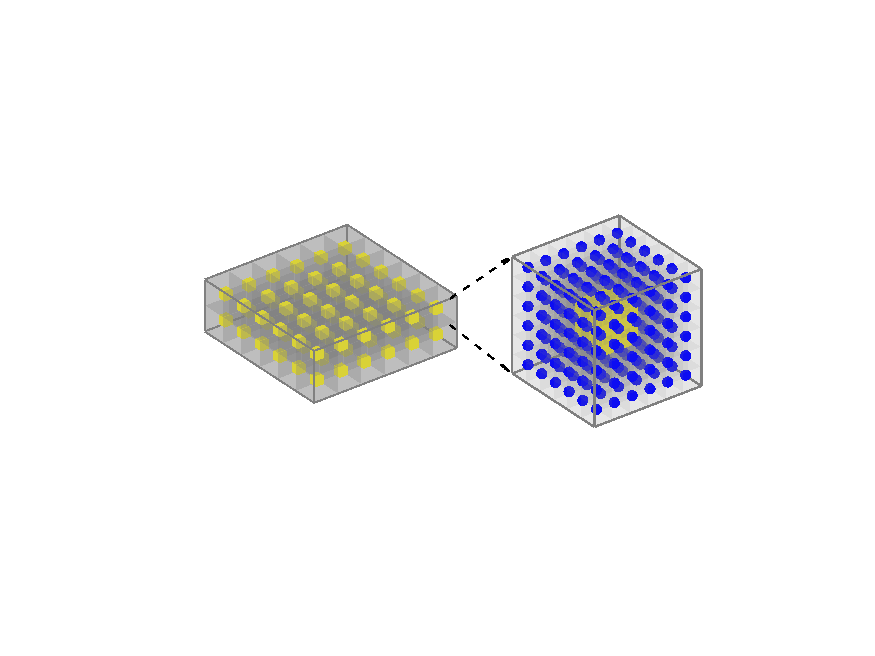} \\
		(a)
	\end{minipage}
	\begin{minipage}[c]{0.3\textwidth}
		\centering
		\includegraphics[
		width=0.9\linewidth,
		trim=2cm 1.3cm 1.8cm 1.3cm, % 从左、下、右、上各裁剪1厘米
		clip % 应用裁剪
		]{E1macro.eps} \\
		(b)
	\end{minipage}
	\begin{minipage}[c]{0.3\textwidth}
		\centering
		\includegraphics[
		width=0.9\linewidth,
		trim=2cm 1.1cm 2cm 0.5cm, % 从左、下、右、上各裁剪1厘米
		clip % 应用裁剪
		]{E3meso.eps} \\
		(c)
	\end{minipage}
	\begin{minipage}[c]{0.3\textwidth}
		\centering
		\includegraphics[
		width=0.9\linewidth,
		trim=2cm 1.1cm 2cm 0.5cm, % 从左、下、右、上各裁剪1厘米
		clip % 应用裁剪
		]{E3micro.eps} \\
		(d)
	\end{minipage}
	\caption{(a) The composite plate domain $\Omega$; (b) The cross section $x_3=0.069$cm of composite plate;(c) Mesoscopic unit cell $Y$;(d) Microscopic unit cell $ Z$.}\label{E4f1}
\end{figure}

The last experiment investigates a 3D heterogeneous plate with multiple spatial scales and spatially periodic parameters $\zeta _1=1/6$ and $\zeta _2=1/36$. The considered composite plate domain $\Omega$, mesoscopic UC $Y$ and microscopic UC $Z$ are shown in Fig.\hspace{1mm}\hspace{1mm}\ref{E4f1}, and $\Omega=(x_1,x_2,x_3)=[0,1] \times [0,1]\times[0,1/3]\mathrm{cm^3}$.

In this example, the boundary temperature is prescribed as 373.15K and 1000K at the bottom and top surface separately, its six surfaces are clamped, and data of governing equations is identical to those in Example 3. Additionally, the material parameters of 3D heterogeneous plate are equal to those in Example 1. After the computational mesh generation, we obtained the computational cost of reference FEM and our HOTS approach in Table \ref{E4t1}, and thus one can deduce that the presented HOTS method significantly reduces computer memory usage compared to the precise FEM.
\begin{table}[!htb]
	\centering
	\caption{Summary of computational cost ($\Delta t = 0.005, t \in [0,1]$s).}\label{E4t1}
	\begin{tabular}{lcccc}
		\hline
		& \multirow{2}{*}{Reference FEM} & \multicolumn{3}{c}{HOTS method} \\
		\cmidrule(l){3-5}
		& & $Z$ & $Y$ & $\Omega$ \\
		\hline
		Number of nodes & $\approx$154,851,264(estimated) & 9,957 & 14,416 & 2,519,424 \\
		Number of elements & $\approx$28,864,512(estimated)  & 1,856 & 2,465 & 439,597 \\
		\hline
	\end{tabular}
\end{table}

\begin{figure}[!htb]
	\centering
	\begin{minipage}[b]{0.32\textwidth}
		\centering
		\includegraphics[
		width=5cm,
		trim=4cm 6.4cm 0cm 5.5cm, % 从左、下、右、上各裁剪1厘米
		clip % 应用裁剪
		]{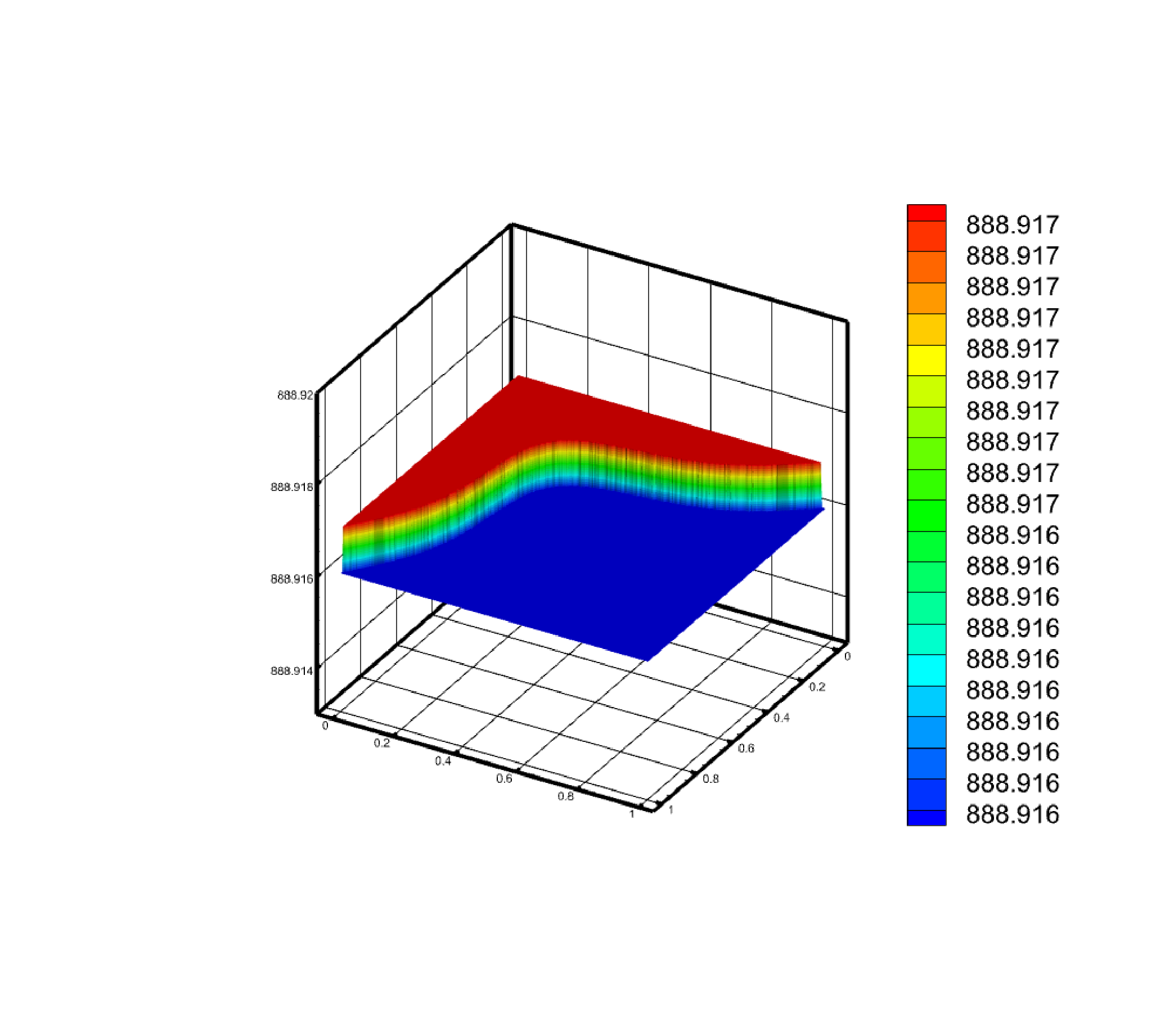} \\
		(a)
	\end{minipage}
	\begin{minipage}[b]{0.32\textwidth}
		\centering
		\includegraphics[
		width=5cm,
		trim=4cm 6.5cm 0cm 5.5cm, % 从左、下、右、上各裁剪1厘米
		clip % 应用裁剪
		]{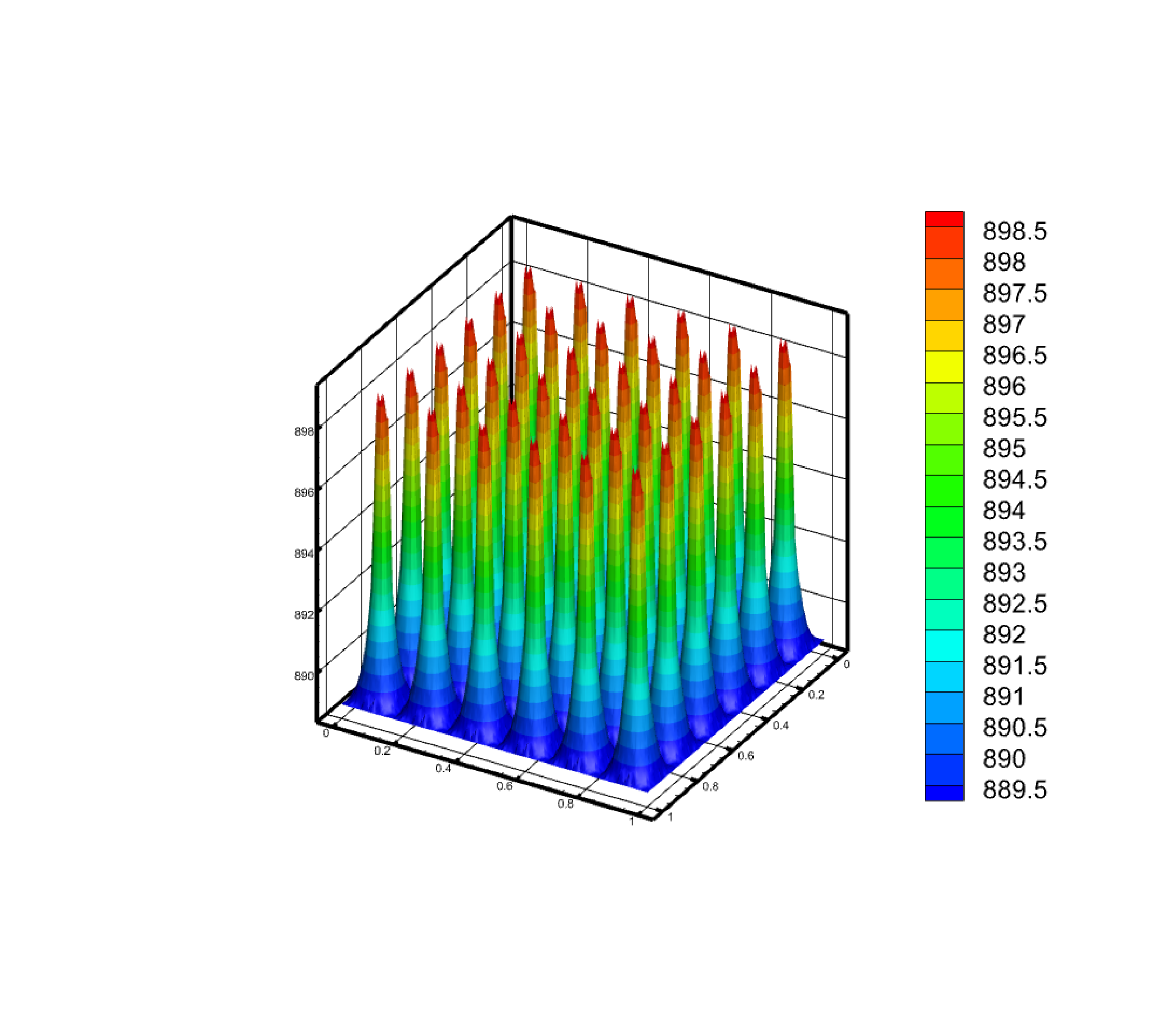} \\
		(b)
	\end{minipage}
	\begin{minipage}[b]{0.32\textwidth}
		\centering
		\includegraphics[
		width=5cm,
		trim=4cm 6.5cm 0cm 5.5cm, % 从左、下、右、上各裁剪1厘米
		clip % 应用裁剪
		]{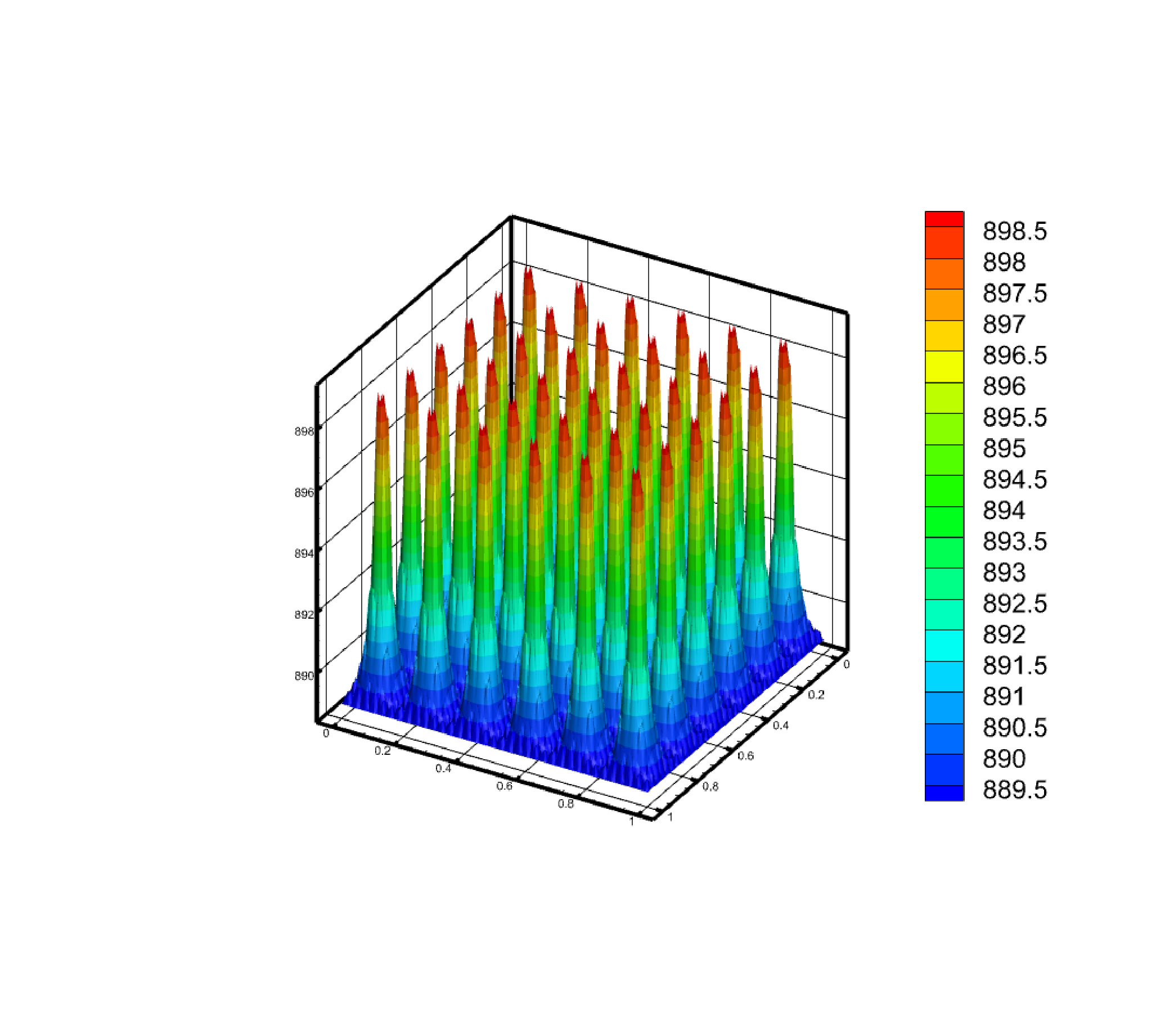} \\
		(c)
	\end{minipage}
	\begin{minipage}[b]{0.48\textwidth}
		\centering
		\includegraphics[
		width=5cm,
		trim=4cm 6.5cm 0cm 5.5cm, % 从左、下、右、上各裁剪1厘米
		clip % 应用裁剪
		]{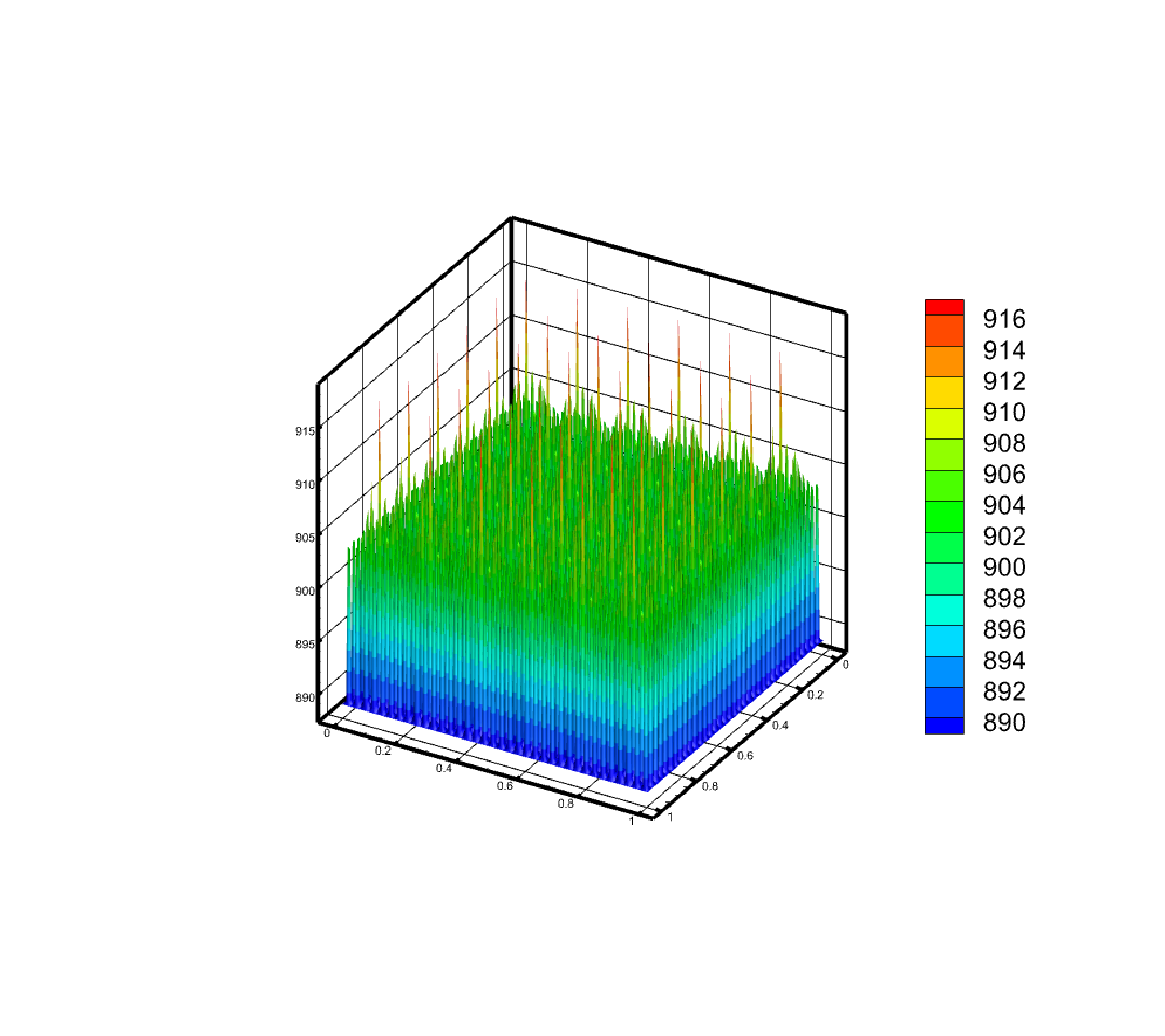} \\
		(d)
	\end{minipage}
    \begin{minipage}[b]{0.48\textwidth}
		\centering
		\includegraphics[
		width=7cm,
		trim=0cm 0cm 0cm 0cm, % 从左、下、右、上各裁剪1厘米
		clip % 应用裁剪
		]{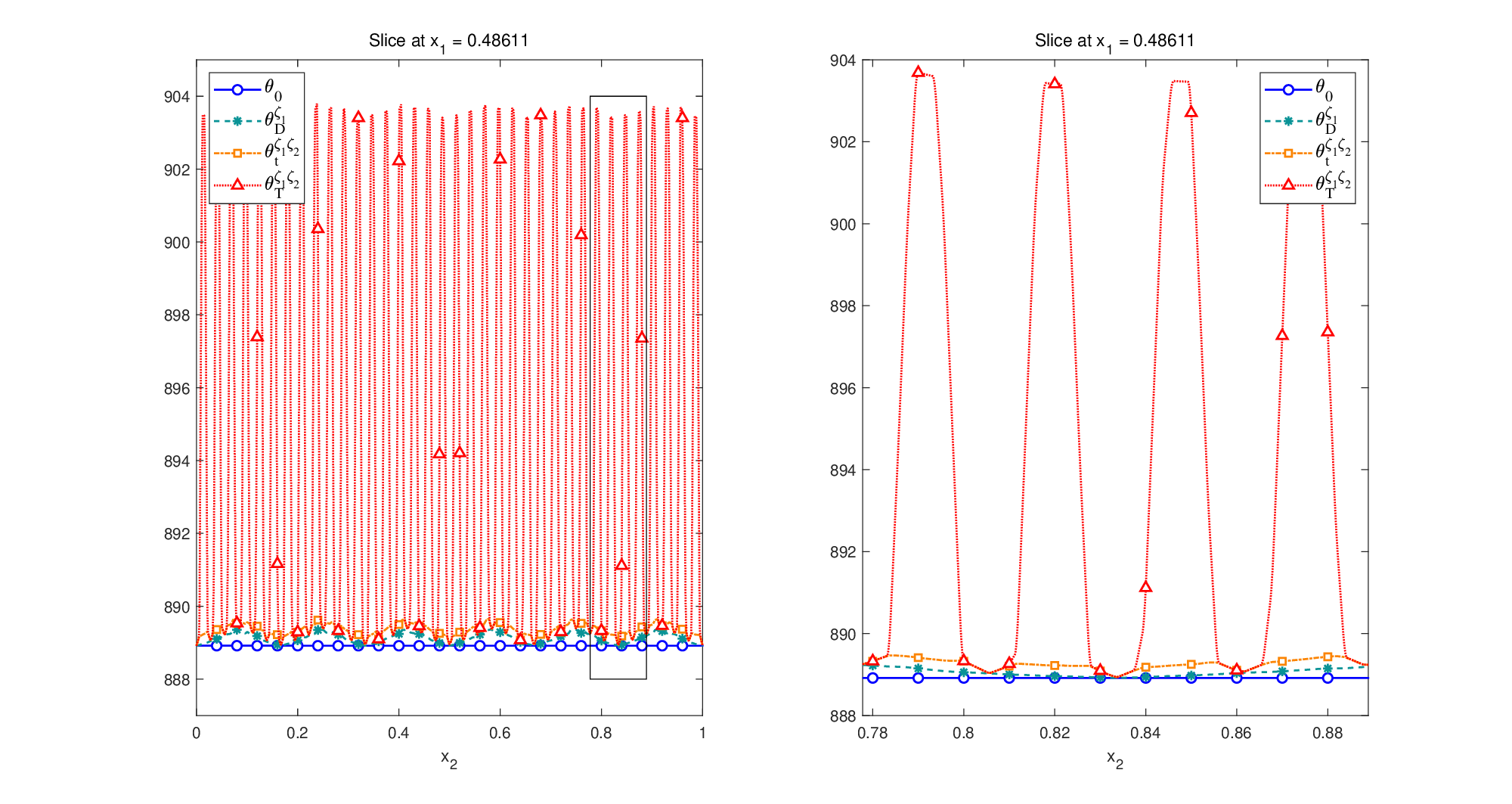} \\
		(e)
	\end{minipage}
	\caption{The temperature field in $x_3 = 0.069$cm: (a) $\theta_{0}$; (b) ${\theta _D^{{\zeta  _1}}}$; (c) ${\theta _t^{{\zeta  _1}{\zeta  _2}}}$; (d) ${\theta _T^{{\zeta  _1}{\zeta  _2}}}$; (e) Computational results on line $x_1 = 0.48611$cm and $x_3 = 0.069$cm.}\label{E4f2}
\end{figure}

\begin{figure}[!htb]
	\centering
	\begin{minipage}[b]{0.4\textwidth}
		\centering
		\includegraphics[
		width=5cm,
		trim=2.5cm 6.5cm 0cm 5.5cm, % 从左、下、右、上各裁剪1厘米
		clip % 应用裁剪
		]{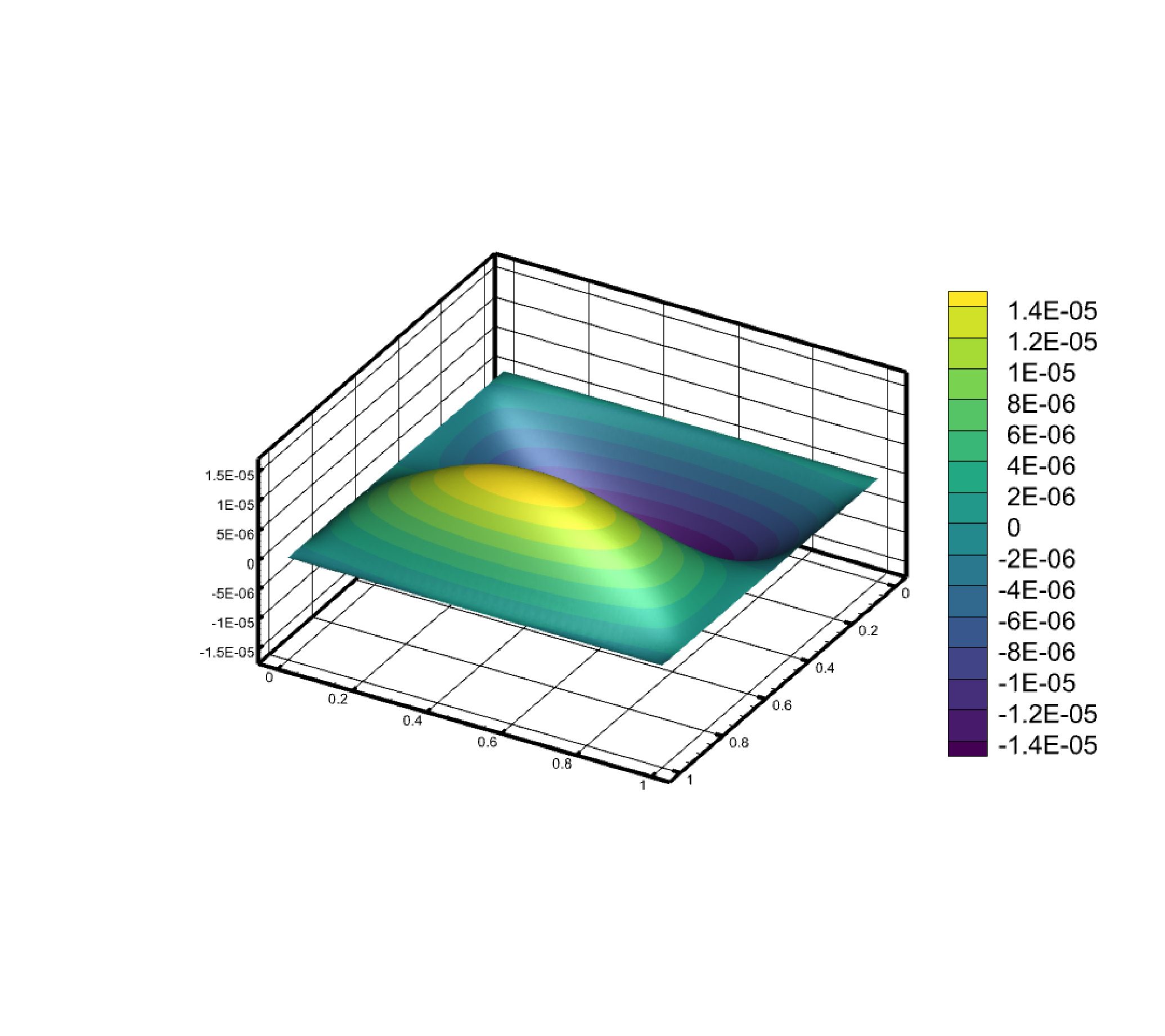} \\
		(a)
	\end{minipage}
	\begin{minipage}[b]{0.4\textwidth}
		\centering
		\includegraphics[
		width=5cm,
		trim=2.5cm 6.5cm 0cm 5.5cm, % 从左、下、右、上各裁剪1厘米
		clip % 应用裁剪
		]{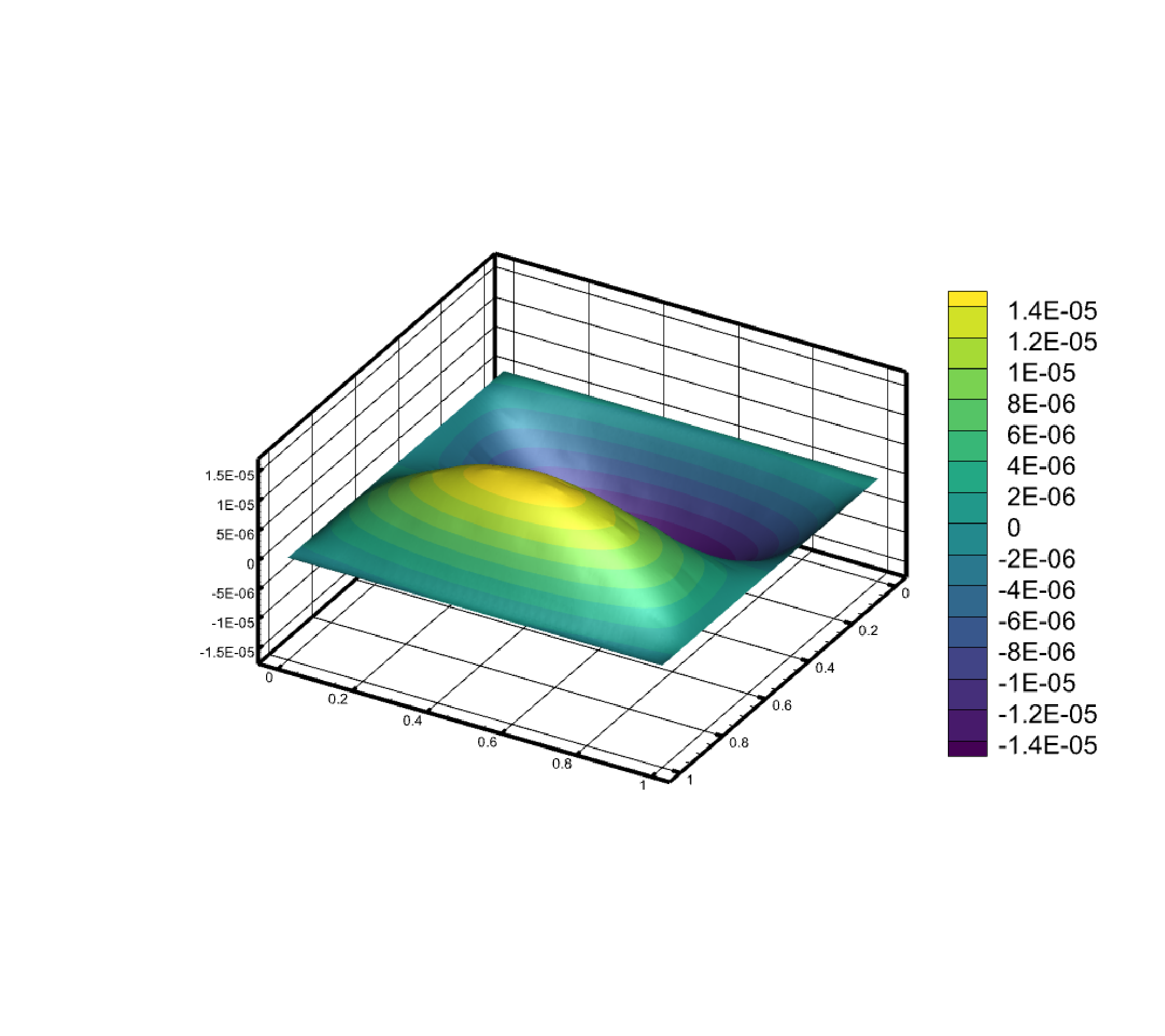} \\
		(b)
	\end{minipage}
	\begin{minipage}[b]{0.4\textwidth}
		\centering
		\includegraphics[
		width=5cm,
		trim=2.5cm 6.5cm 0cm 5.5cm, % 从左、下、右、上各裁剪1厘米
		clip % 应用裁剪
		]{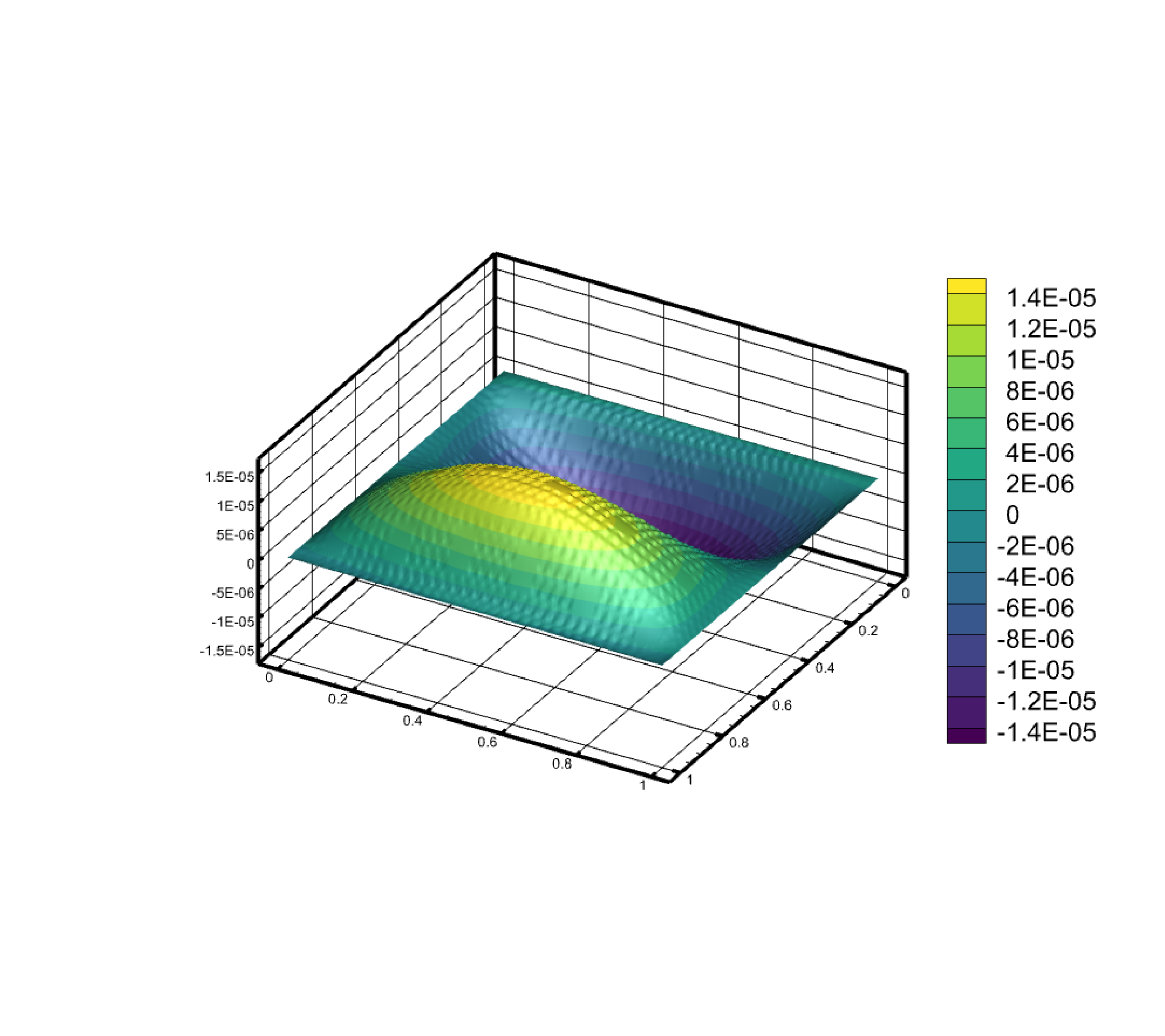} \\
		(c)
	\end{minipage}
	\begin{minipage}[b]{0.4\textwidth}
		\centering
		\includegraphics[
		width=5cm,
		trim=2.5cm 6.5cm 0cm 5.5cm, % 从左、下、右、上各裁剪1厘米
		clip % 应用裁剪
		]{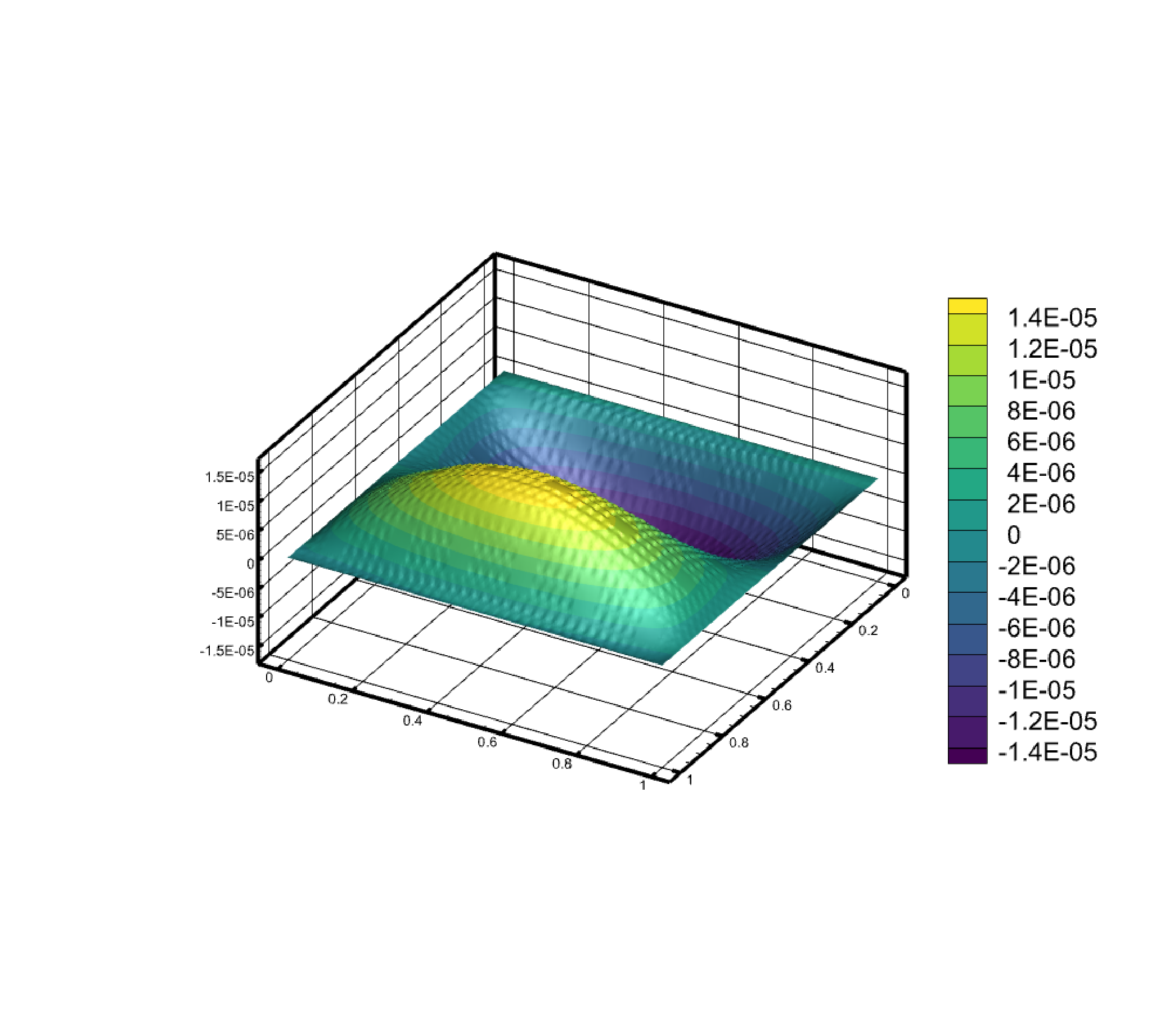} \\
		(d)
	\end{minipage}
	\caption{The first component for the displacement field in $x_3 = 0.069$cm: (a) ${{u_{1}^{(0)}}}$; (b) ${u_{1D}^{{\zeta  _1}}}$; (c) ${u_{1t}^{{\zeta  _1}{\zeta  _2}}}$; (d) ${u_{1T}^{{\zeta  _1}{\zeta  _2}}}$.}\label{E4f3}
\end{figure}

\begin{figure}[!htb]
	\centering
	\begin{minipage}[b]{0.4\textwidth}
		\centering
		\includegraphics[
		width=5cm,
		trim=2.5cm 6.5cm 0cm 5.5cm, % 从左、下、右、上各裁剪1厘米
		clip % 应用裁剪
		]{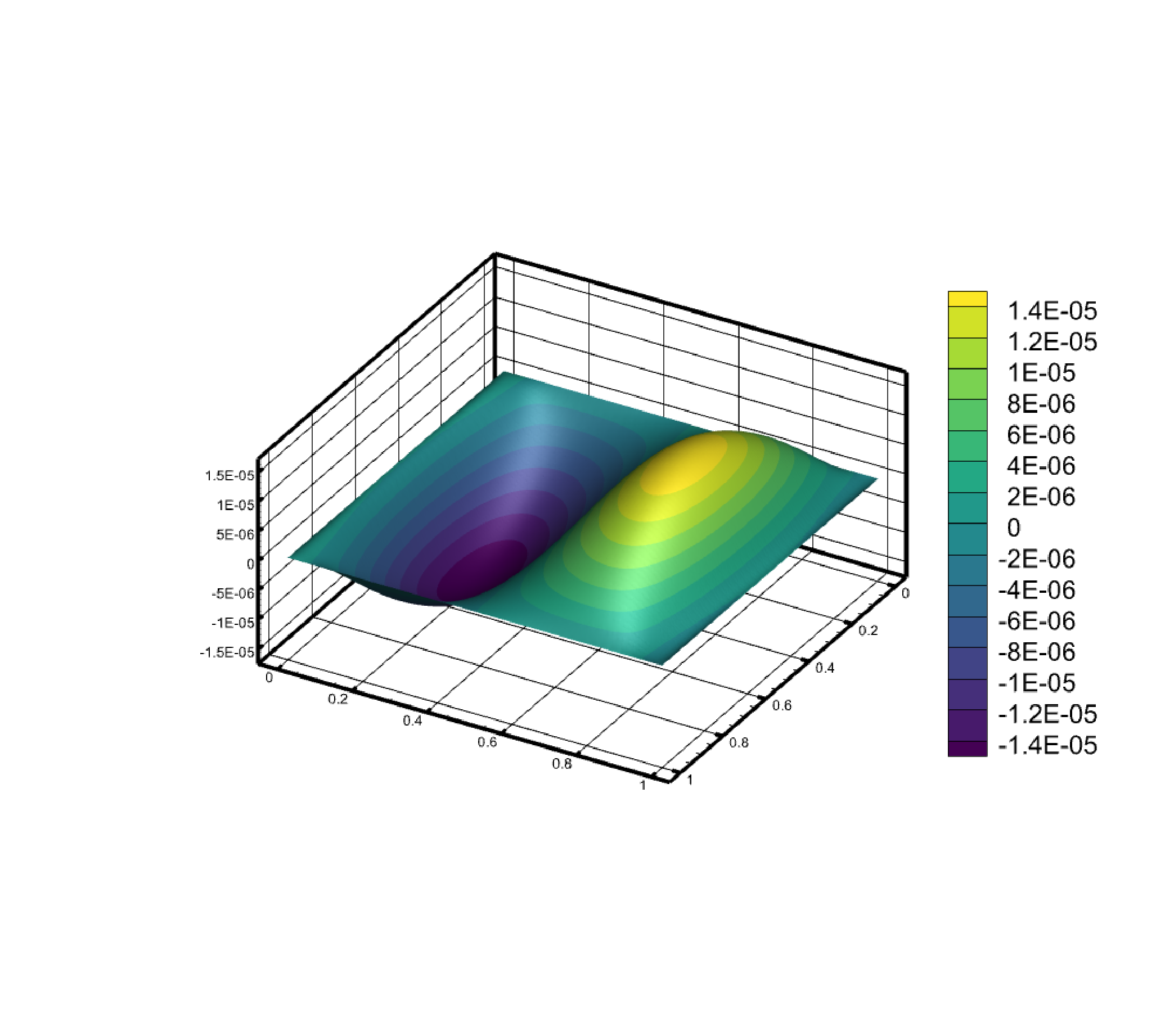} \\
		(a)
	\end{minipage}
	\begin{minipage}[b]{0.4\textwidth}
		\centering
		\includegraphics[
		width=5cm,
		trim=2.5cm 6.5cm 0cm 5.5cm, % 从左、下、右、上各裁剪1厘米
		clip % 应用裁剪
		]{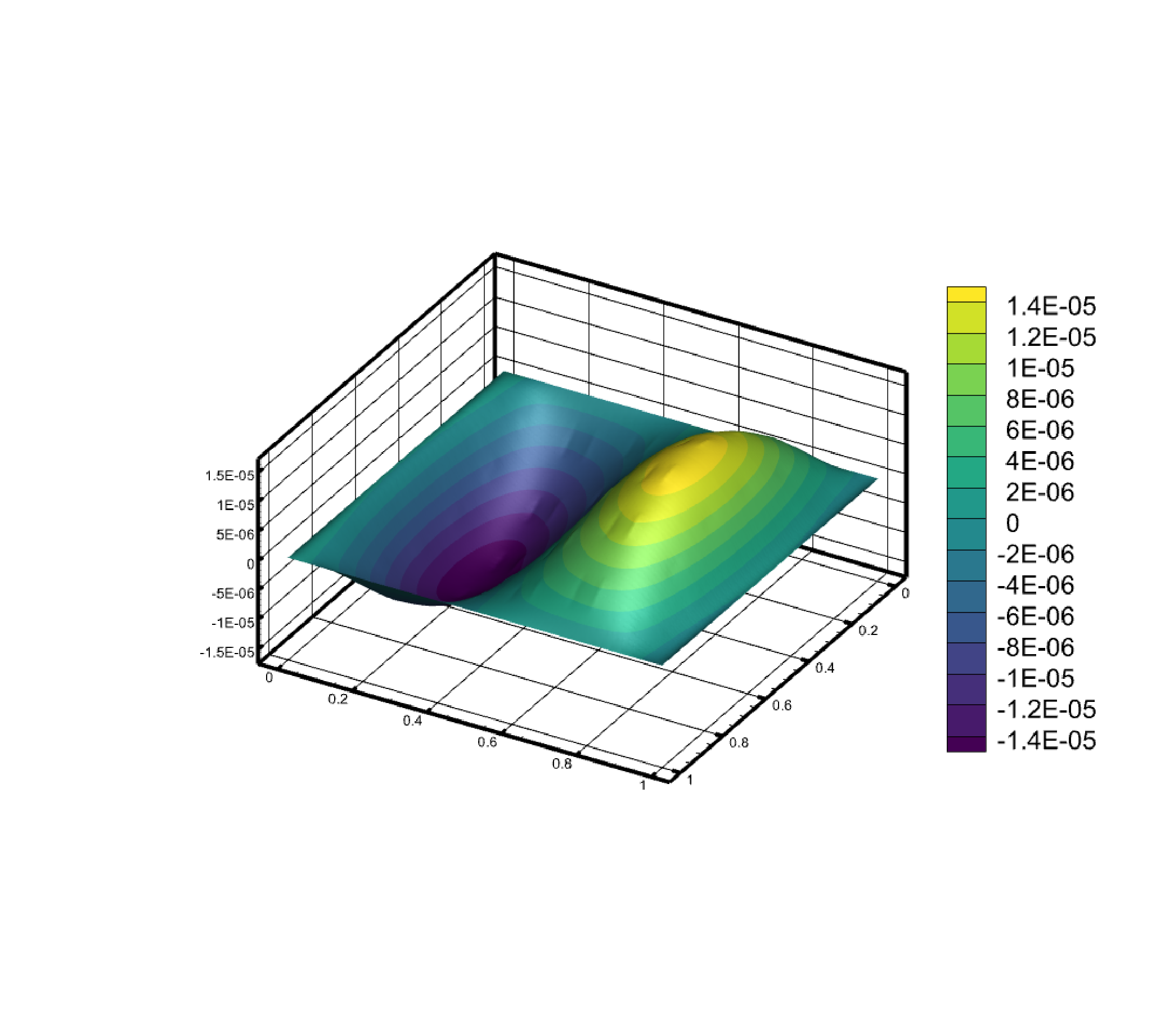} \\
		(b)
	\end{minipage}
	\begin{minipage}[b]{0.4\textwidth}
		\centering
		\includegraphics[
		width=5cm,
		trim=2.5cm 6.5cm 0cm 5.5cm, % 从左、下、右、上各裁剪1厘米
		clip % 应用裁剪
		]{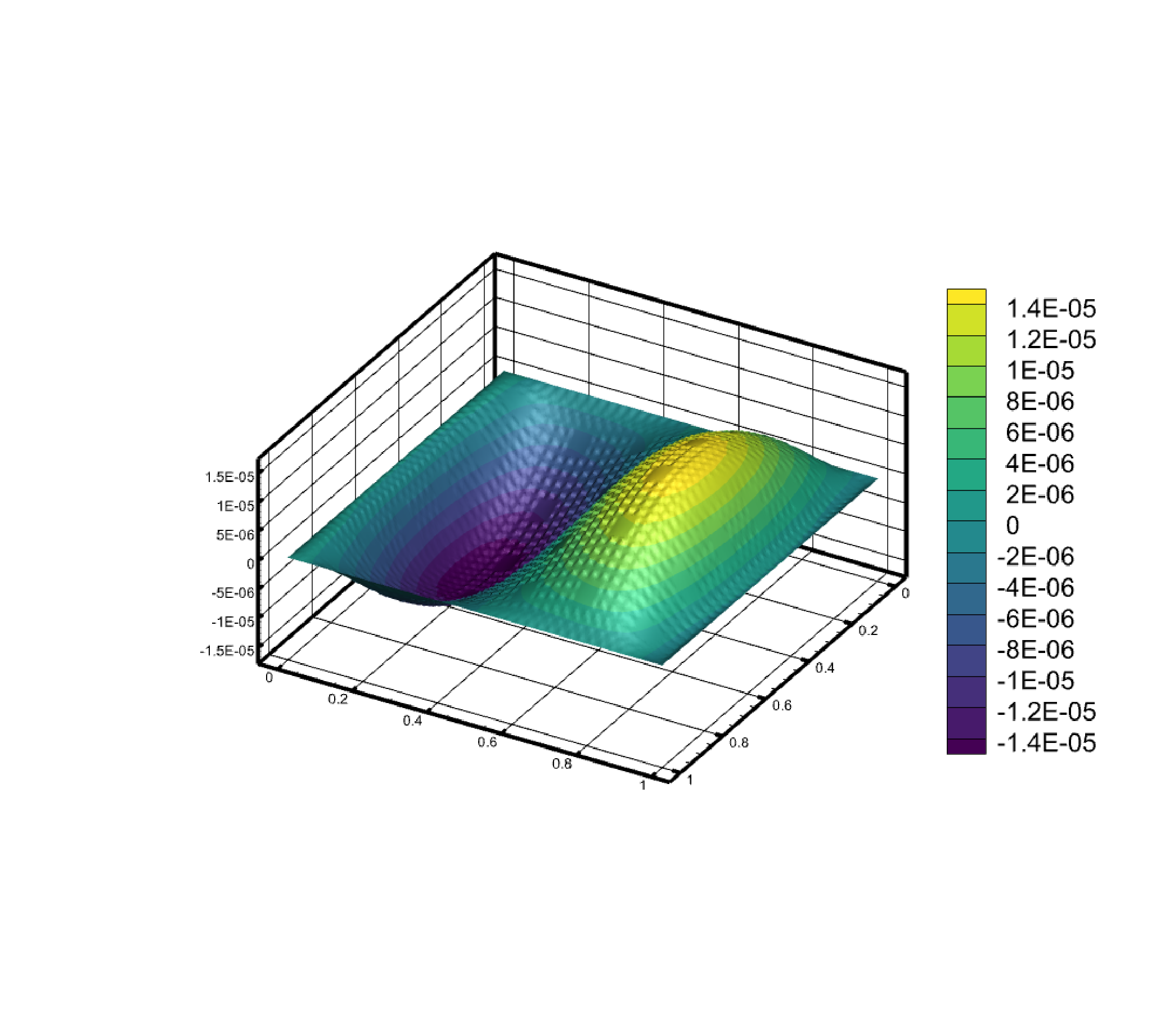} \\
		(c)
	\end{minipage}
	\begin{minipage}[b]{0.4\textwidth}
		\centering
		\includegraphics[
		width=5cm,
		trim=2.5cm 6.5cm 0cm 5.5cm, % 从左、下、右、上各裁剪1厘米
		clip % 应用裁剪
		]{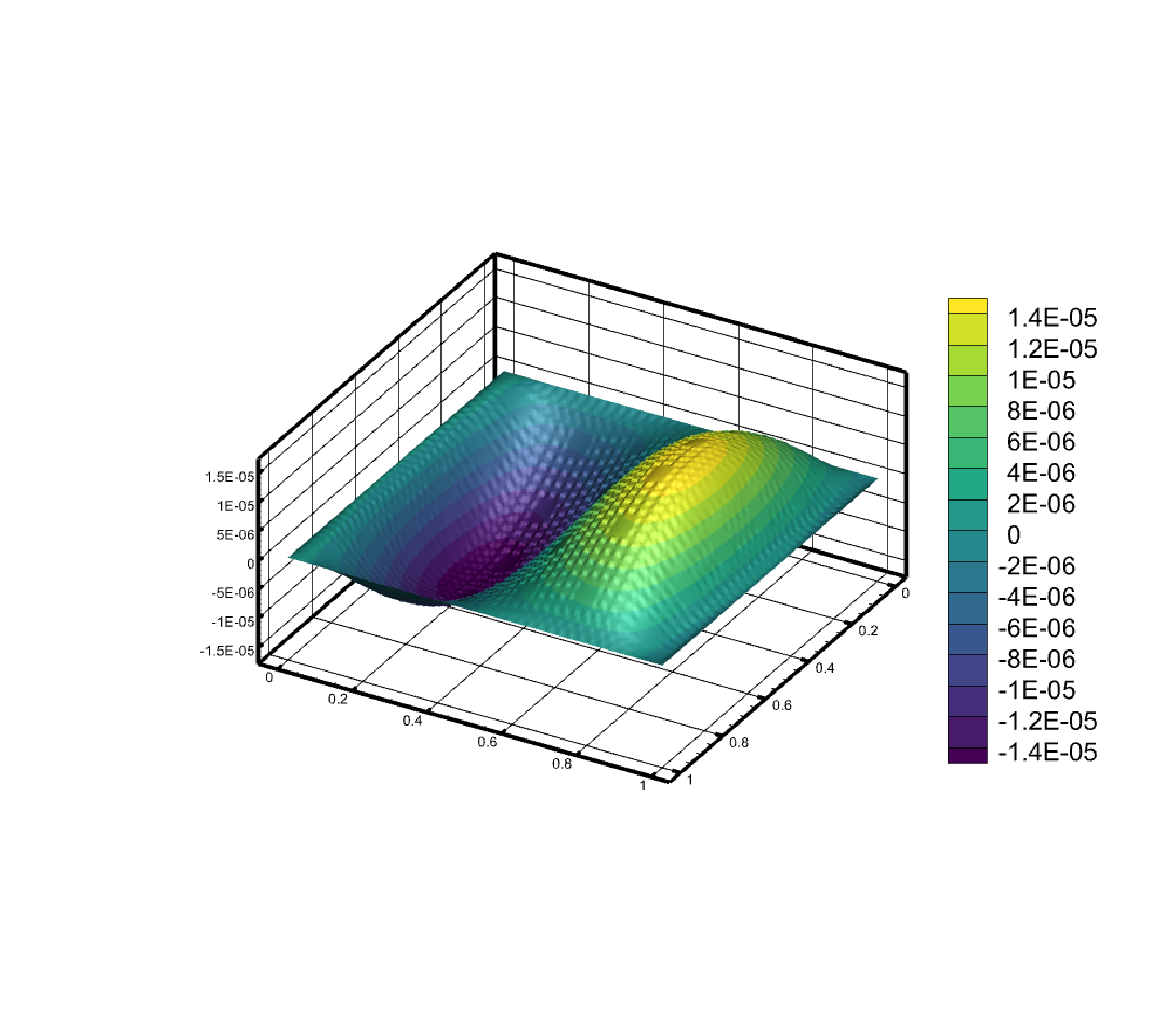} \\
		(d)
	\end{minipage}
	\caption{The second component for the displacement field in $x_3 = 0.069$cm: (a) ${{u_{2}^{(0)}}}$; (b) ${u_{2D}^{{\zeta  _1}}}$; (c) ${u_{2t}^{{\zeta  _1}{\zeta  _2}}}$; (d) ${u_{2T}^{{\zeta  _1}{\zeta  _2}}}$.}\label{E4f4}
\end{figure}

\begin{figure}[!htb]
	\centering
	\begin{minipage}[b]{0.34\textwidth}
		\centering
		\includegraphics[
		width=5cm,
		trim=2.5cm 6.5cm 0cm 5.5cm, % 从左、下、右、上各裁剪1厘米
		clip % 应用裁剪
		]{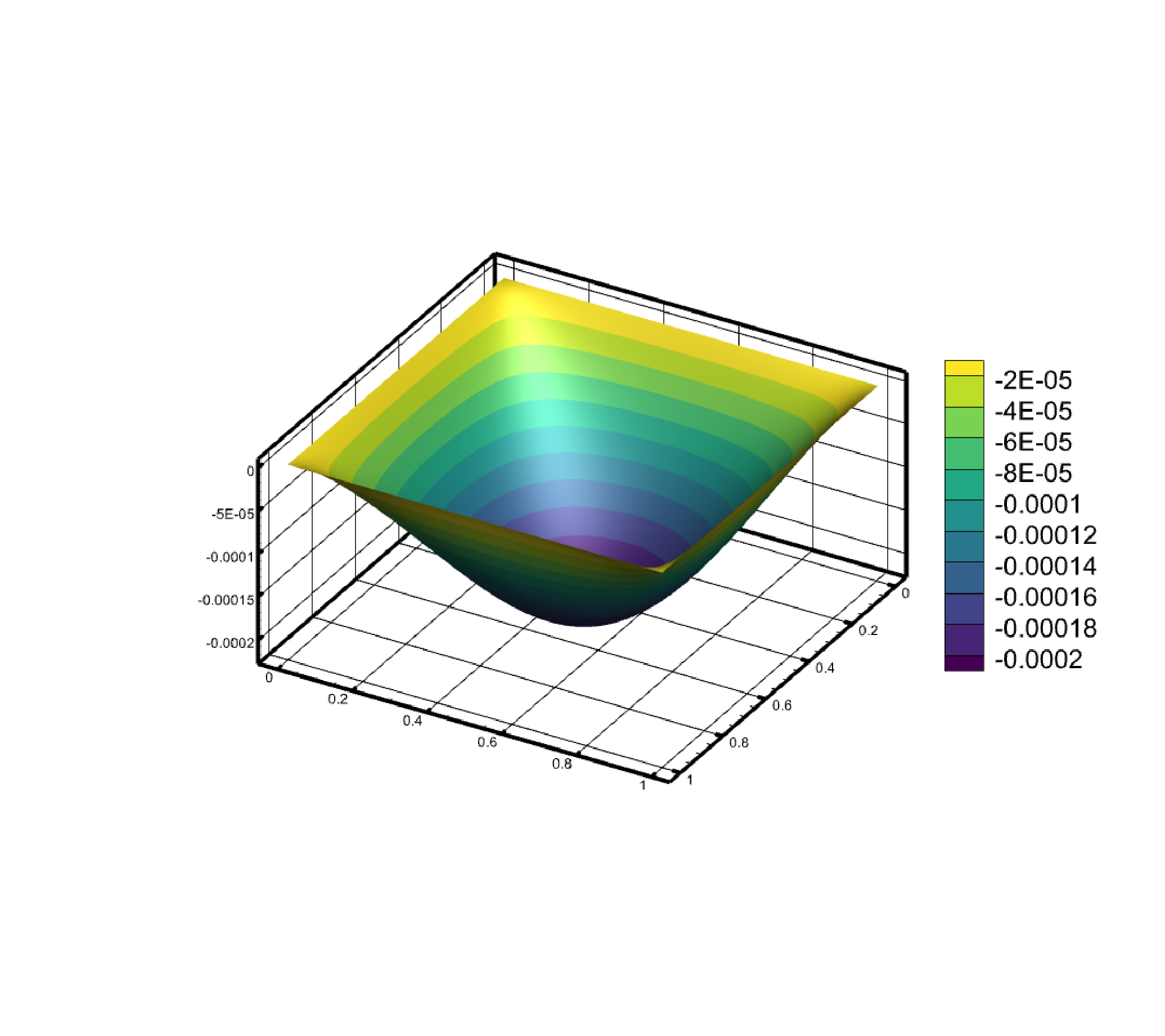} \\
		(a)
	\end{minipage}
	\begin{minipage}[b]{0.34\textwidth}
		\centering
		\includegraphics[
		width=5cm,
		trim=2.5cm 6.5cm 0cm 5.5cm, % 从左、下、右、上各裁剪1厘米
		clip % 应用裁剪
		]{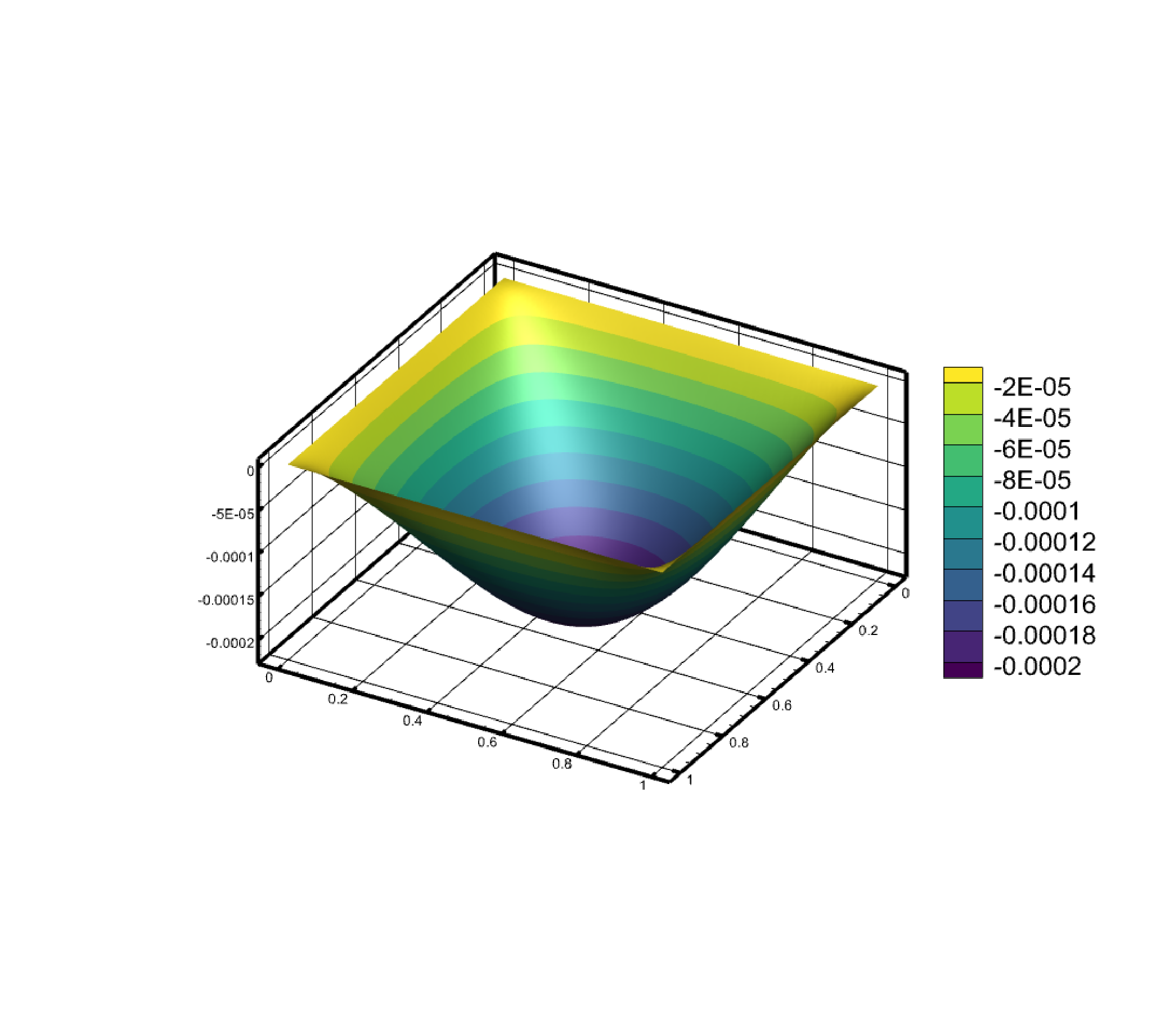} \\
		(b)
	\end{minipage}
	\begin{minipage}[b]{0.30\textwidth}
		\centering
		\includegraphics[
		width=5cm,
		trim=2.5cm 6.5cm 0cm 5.5cm, % 从左、下、右、上各裁剪1厘米
		clip % 应用裁剪
		]{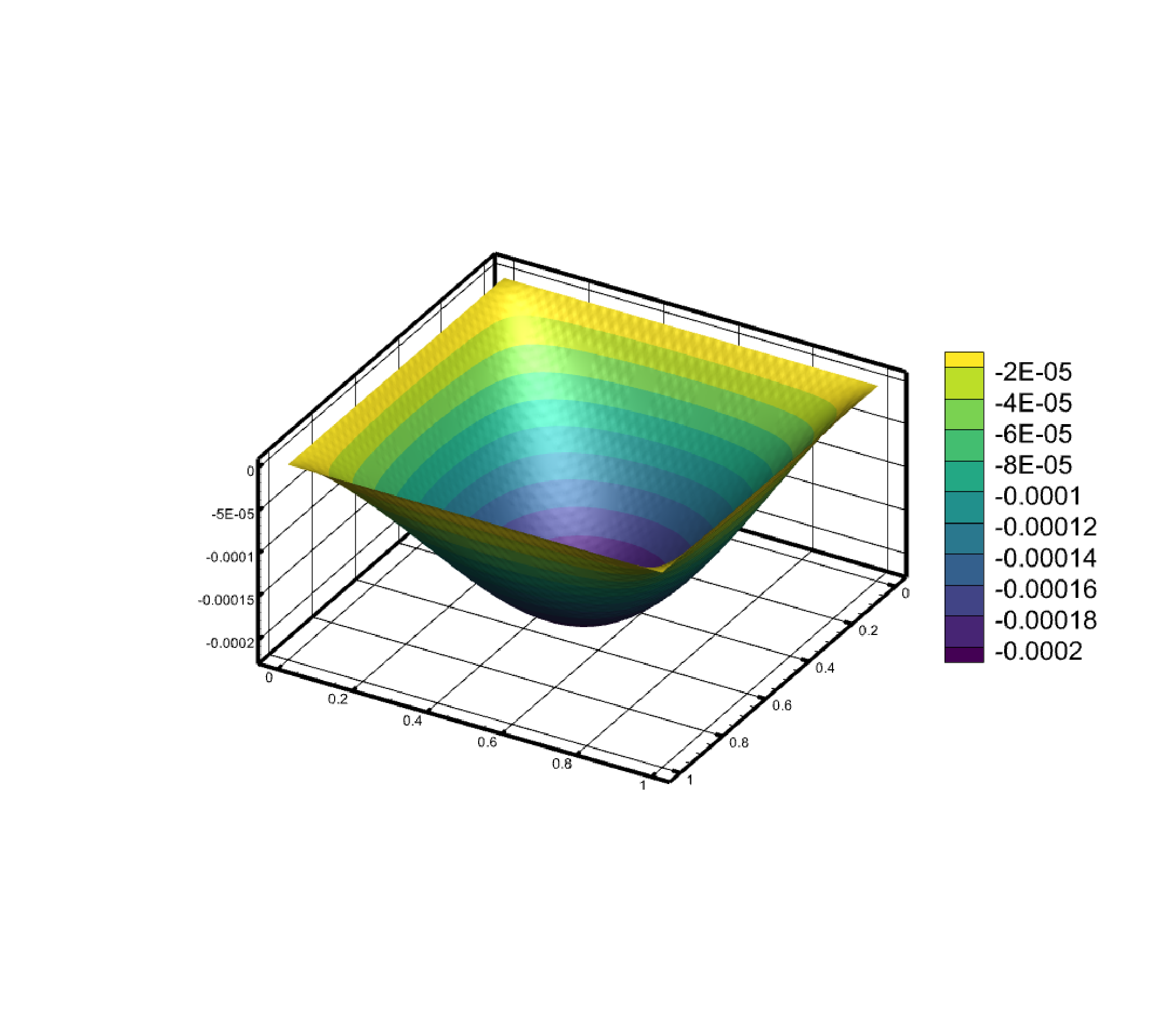} \\
		(c)
	\end{minipage}
	\begin{minipage}[b]{0.48\textwidth}
		\centering
		\includegraphics[
		width=5cm,
		trim=2.5cm 6.5cm 0cm 5.5cm, % 从左、下、右、上各裁剪1厘米
		clip % 应用裁剪
		]{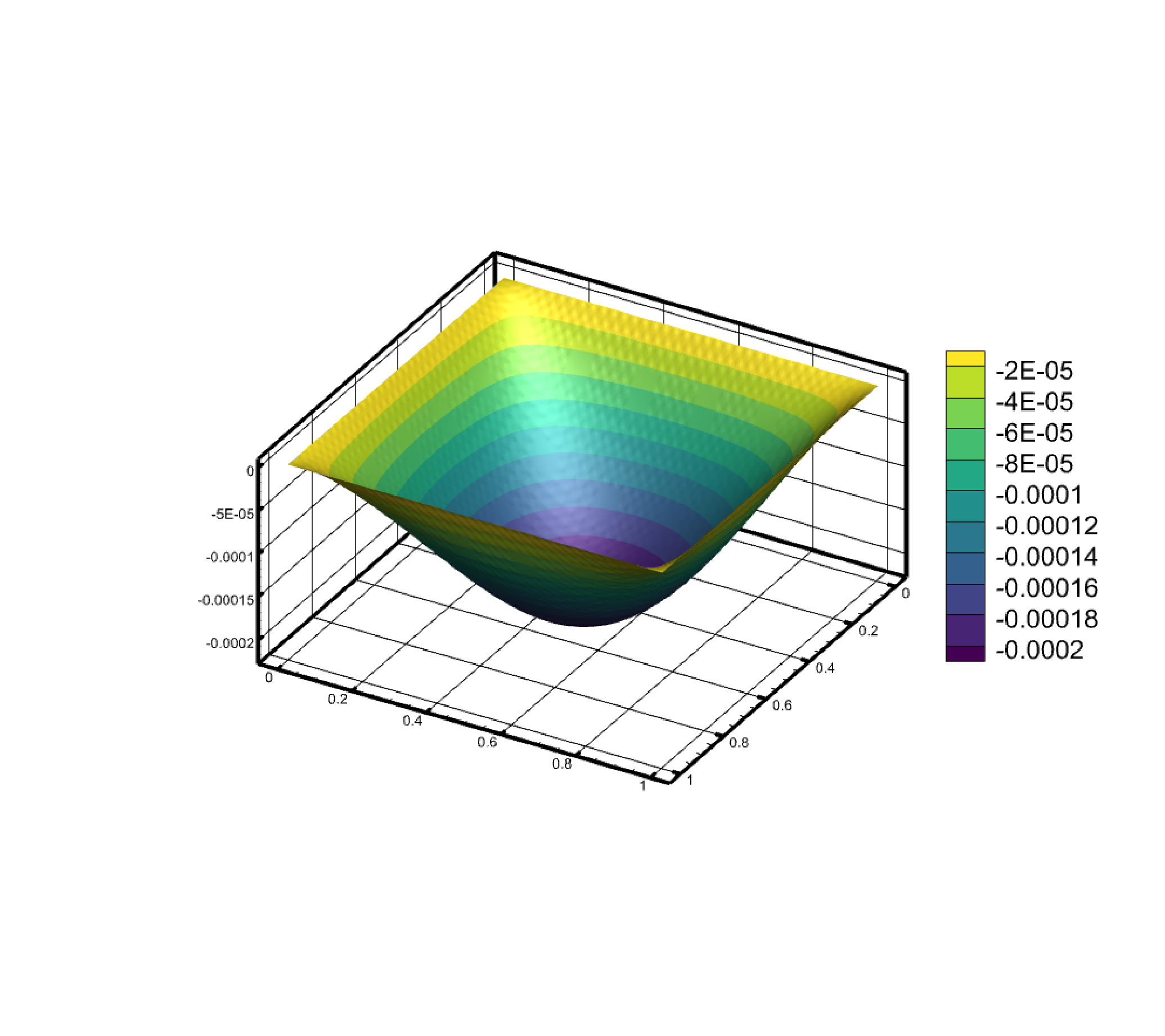} \\
		(d)
	\end{minipage}
    \begin{minipage}[b]{0.48\textwidth}
		\centering
		\includegraphics[
		width=7cm,
		trim=0cm 0cm 0cm 0cm, % 从左、下、右、上各裁剪1厘米
		clip % 应用裁剪
		]{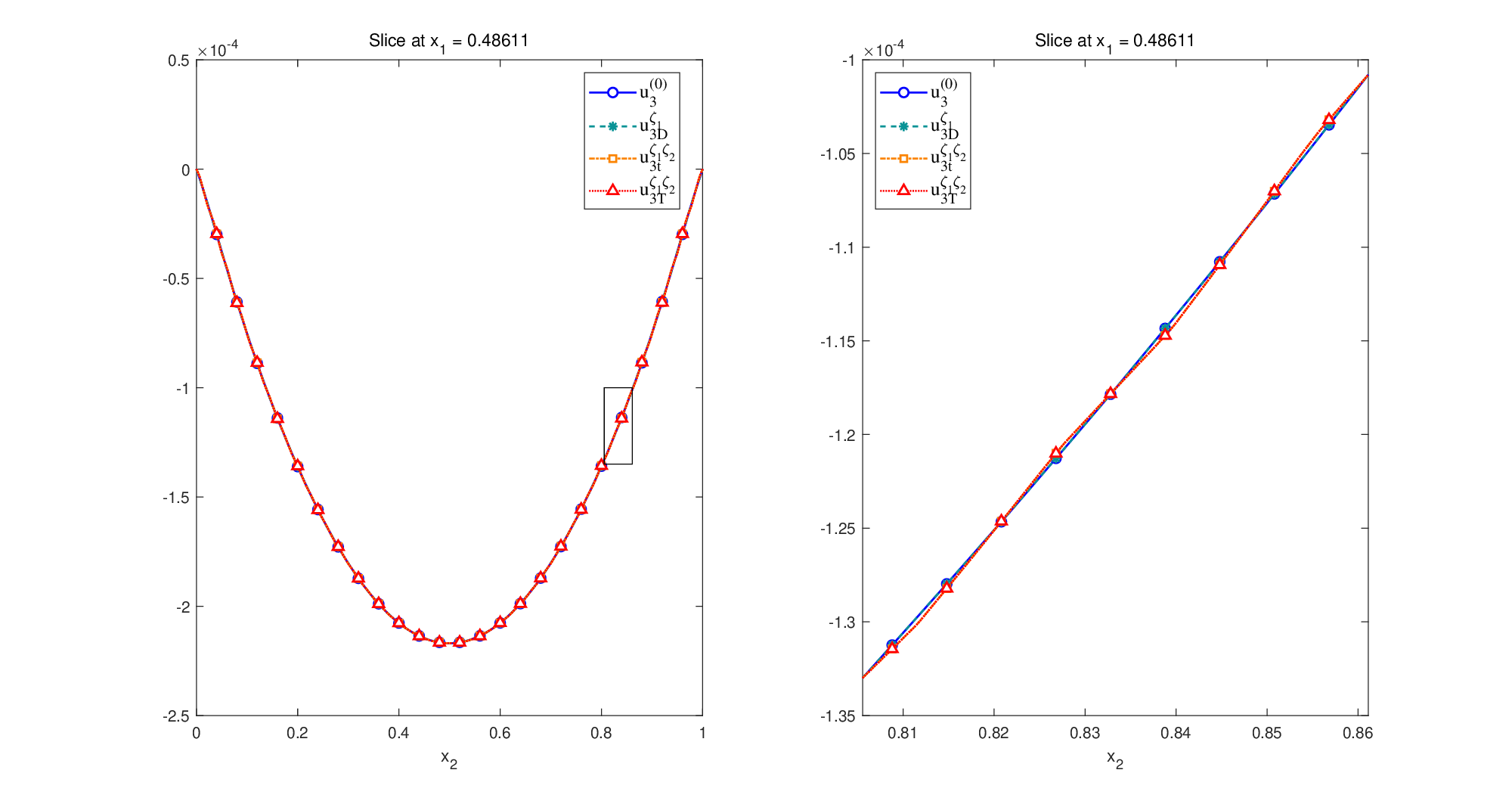} \\
		(e)
	\end{minipage}
	\caption{The third component for the displacement field in $x_3 = 0.069$cm: (a) ${{u_{3}^{(0)}}}$; (b) ${u_{3D}^{{\zeta  _1}}}$; (c) ${u_{3t}^{{\zeta  _1}{\zeta  _2}}}$; (d) ${u_{3T}^{{\zeta  _1}{\zeta  _2}}}$; (e) Computational results on line $x_1 = 0.48611$cm and $x_3 = 0.069$cm.}\label{E4f5}
\end{figure}

For the three-scale plate, Figs.\hspace{1mm}\ref{E4f2}-\ref{E4f5} displayed distinct multi-scale solutions for temperature and displacement fields in $x_3 =0.069$cm and at time $t = 1.0$s.

As illustrated in Figs.\hspace{1mm}\ref{E4f2}-\ref{E4f5}, the homogenized method, the two-scale method and the lower-order three-scale method can only acquire the macroscopic material responses, the mesoscopic material responses and the inadequate microscopic oscillation responses, respectively. However, the HOTS solutions correspond accurately to the reference numerical solutions, notably for the temperature field. Thus, the HOTS solutions would be utilized to catch the microscopic coupling behaviors of this heterogeneous plate with massive unit cells. In real-world engineering applications, it is virtually impossible to compute the precise FEM solutions for large-scale composite structures. Now, our HOTS method can effectively simulate nonlinear multi-scale multiphysics problems of large-scale heterogeneous structures with minimal computational resource consumption.

\section{Conclusion and future work}
This work proposed an innovative higher-order three-scale (HOTS) computational method  for effectively simulating nonlinear thermo-mechanical coupling problems of heterogeneous structures with three-level spatial hierarchy. The main contributions for this work are summarized as follows. Firstly, we established a high-accuracy HOTS computational model by recursively two-scale analysis between meso-micro and macro-meso scales. Secondly, a local error analysis of the proposed HOTS model is rigorously derived in the point-wise sense in order to validate its well-balanced property. Thirdly, a two-stage numerical algorithm with off-line and on-line stages is presented to effectively compute the nonlinear physical responses of heterogeneous structures with complicated three-scale configurations.

The computational efficiency and accuracy  of the presented HOTS approach are substantiated through two 2D and two 3D heterogeneous structures. On one hand, the numerical examples demonstrated that the HOTS approach could precisely catch the nonlinear thermo-mechanical coupling phenomena in three-scale heterogeneous structures, and the corresponding algorithm was stable and robust for prolonged multi-scale nonlinear multi-field coupling simulation. On the other hand, compared to direct FE calculation, the proposed HOTS method could not only significantly reduce computer storage resources but also save computational time. The high numerical accuracy, low computational cost, and effective computational performance of the proposed HOTS approach have considerable practical value for large-scale engineering computation. Additionally, the proposed higher-order three-scale framework has scalable generality for nonlinear multiphysics problems of heterogeneous structures with three-level structural configuration.

Future study will concentrate on the expansion of intricate nonlinear problems and random composite materials including thermal radiation and convection effects under high-temperature environment. Furthermore, intrinsic parallel advantage of the HOTS method in off-line stage will be exploited to further enhance its efficiency.

\section*{Acknowledgments}
This research was supported by the National Natural Science Foundation of China (No.\hspace{1mm}12471387), Xidian University Specially Funded Project for Interdisciplinary Exploration (No.\hspace{1mm}TZJH2024008), Fundamental Research Funds for the Central Universities (No.\hspace{1mm}QTZX25082), Innovation Capability Support Program of Shaanxi Province (No.\hspace{1mm}2024RS-CXTD-88), and also supported by the high performance computing center of Xidian University.
\appendix
\section{Formulation and analysis on two types of macroscopic homogenized coefficients}
\label{app1}
According to the SOTS approach and the HOTS approach in section 2, we can obtain two fundamentally distinct macroscopic homogenization coefficients. While the formula (21) has given the one type of macroscopic homogenized coefficients, we simulate the other type of macroscopic homogenized coefficients by means of the following reiteration homogenization technique, namely
\begin{equation}
	\begin{aligned}
	&{{\hat S}^*}({\theta _0}) = \frac{1}{{|Y|}}\int_Y {(\bar S(\bm{y},\theta_0) - {{\bar \vartheta }_{ij}}(\bm{y},\theta_0)\frac{{\partial {P_i}}}{{\partial {y_j}}})dY} \\
	&= \frac{1}{{|Y|}}\int_Y {\Bigl[\frac{1}{{|Z|}}\int_Z {({\rho ^{(0)}}{c^{(0)}} - \vartheta _{ij}^{(0)}\frac{{\partial {O_i}}}{{\partial {z_j}}})dZ}  - \frac{1}{{|Z|}}\int_Z {(\vartheta _{ij}^{(0)} + \vartheta _{pm}^{(0)}\frac{{\partial T_{pi}^j}}{{\partial {z_m}}})} dZ\frac{{\partial {P_i}}}{{\partial {y_j}}}\Bigl]dY} \\
	&= \frac{1}{{|Y||Z|}}\int_Y {\Bigl[\int_Z {({\rho ^{(0)}}{c^{(0)}} - \vartheta _{ij}^{(0)}\frac{{\partial {O_i}}}{{\partial {z_j}}})dZ - \int_Z {(\vartheta _{ij}^{(0)} + \vartheta _{pm}^{(0)}\frac{{\partial T_{pi}^j}}{{\partial {z_m}}})dZ} \frac{{\partial {P_i}}}{{\partial {y_j}}}} } \Bigl]dY.
	\end{aligned}
\end{equation}
\begin{equation}
	\begin{aligned}
	&{{\hat k}^*}_{ij}({\theta _0}) = \frac{1}{{|Y|}}\int_Y {({{\bar k}_{ij}}({\bm{y}},{\theta _0}) + {{\bar k}_{ik}}({\bm{y}},{\theta _0})\frac{{\partial {M_j}}}{{\partial {y_k}}})dY} \\
	&= \frac{1}{{|Y|}}\int_Y {\Bigl[\frac{1}{{|Z|}}\int_Z {(k_{ij}^{(0)} + k_{ik}^{(0)}\frac{{\partial {R_j}}}{{\partial {z_k}}})dZ}  + \frac{1}{{|Z|}}\int_Z {(k_{ik}^{(0)} + k_{im}^{(0)}\frac{{\partial {R_k}}}{{\partial {z_m}}})} dZ\frac{{\partial {M_j}}}{{\partial {y_k}}}\Bigl]dY} \\
	&= \frac{1}{{|Y||Z|}}\int_Y {\Bigl[\int_Z {(k_{ij}^{(0)} + k_{ik}^{(0)}\frac{{\partial {R_j}}}{{\partial {z_k}}})dZ + \int_Z {(k_{ik}^{(0)} + k_{im}^{(0)}\frac{{\partial {R_k}}}{{\partial {z_m}}})dZ} \frac{{\partial {P_i}}}{{\partial {y_j}}}} } \Bigl]dY.
	\end{aligned}
\end{equation}
\begin{equation}
	\begin{aligned}
	&{{\hat \vartheta }^*}_{ij}({\theta _0}) = \frac{1}{{|Y|}}\int_Y {({{\bar \vartheta }_{ij}}({\bm{y}},{\theta _0}) + {{\bar \vartheta }_{mn}}({\bm{y}},{\theta _0})\frac{{\partial N_{mi}^j}}{{\partial {y_n}}})dY} \\
	&= \frac{1}{{|Y|}}\int_Y {\Bigl[\frac{1}{{|Z|}}\int_Z {(\vartheta _{ij}^{(0)} + \vartheta _{pq}^{(0)}\frac{{\partial T_{pi}^j}}{{\partial {z_q}}})dZ}  + \frac{1}{{|Z|}}\int_Z {(\vartheta _{mn}^{(0)} + \vartheta _{pq}^{(0)}\frac{{\partial T_{pm}^n}}{{\partial {z_q}}})} dZ\frac{{\partial N_{mi}^j}}{{\partial {y_n}}}\Bigl]dY} \\
	&= \frac{1}{{|Y||Z|}}\int_Y {\Bigl[\int_Z {(\vartheta _{ij}^{(0)} + \vartheta _{pq}^{(0)}\frac{{\partial T_{pi}^j}}{{\partial {z_q}}})dZ + \int_Z {(\vartheta _{mn}^{(0)} + \vartheta _{pq}^{(0)}\frac{{\partial T_{pm}^n}}{{\partial {z_q}}})dZ} \frac{{\partial N_{mi}^j}}{{\partial {y_n}}}} }\Bigl]dY.
	\end{aligned}
\end{equation}
\begin{equation}
	\begin{aligned}
	&\hat \rho ({\theta _0}) = \frac{1}{{|Y|}}\int_Y {\bar \rho ({\bm{y}},{\theta _0})dY} \\
	&= \frac{1}{{|Y|}}\int_Y {\Bigl[\frac{1}{{|Z|}}\int_Z {{\rho ^{(0)}}({\bm{y}},{\bm{z}},{\theta _0})dZ}\Bigl]dY} \\
	&= \frac{1}{{|Y||Z|}}\int_Y {\Bigl[\int_Z {{\rho ^{(0)}}({\bm{y}},{\bm{z}},{\theta _0})dZ} }\Bigl ]dY.
	\end{aligned}
\end{equation}
\begin{equation}
	\begin{aligned}
	&{{\hat C}^*}_{ijkl}({\theta _0}) = \frac{1}{{|Y|}}\int_Y {({{\bar C}_{ijkl}}({\bm{y}},{\theta _0}) + {{\bar C}_{ijmn}}({\bm{y}},{\theta _0})\frac{{\partial N_{mk}^l}}{{\partial {y_n}}})dY} \\
	&= \frac{1}{{|Y|}}\int_Y {\Bigl[\frac{1}{{|Z|}}\int_Z {(C_{ijkl}^{(0)} + C_{ijpq}^{(0)}\frac{{\partial T_{pk}^l}}{{\partial {z_q}}})dZ}  + \frac{1}{{|Z|}}\int_Z {(C_{ijmn}^{(0)} + C_{ijpq}^{(0)}\frac{{\partial T_{pm}^n}}{{\partial {z_q}}})} dZ\frac{{\partial N_{mk}^l}}{{\partial {y_n}}}\Bigl]dY} \\
	&= \frac{1}{{|Y||Z|}}\int_Y {\Bigl[\int_Z {(C_{ijkl}^{(0)} + C_{ijpq}^{(0)}\frac{{\partial T_{pk}^l}}{{\partial {z_q}}})dZ + \int_Z {(C_{ijmn}^{(0)} + C_{ijpq}^{(0)}\frac{{\partial T_{pm}^n}}{{\partial {z_q}}})dZ} \frac{{\partial N_{mk}^l}}{{\partial {y_n}}}} } \Bigl]dY.
	\end{aligned}
\end{equation}
\begin{equation}
	\begin{aligned}
	&{{\hat \beta }^*}_{ij}({\theta _0}) = \frac{1}{{|Y|}}\int_Y {({{\bar \beta }_{ij}}({\bm{y}},{\theta _0}) + {{\bar C}_{ijkl}}({\bm{y}},{\theta _0})\frac{{\partial {P_k}}}{{\partial {y_l}}})dY} \\
	&= \frac{1}{{|Y|}}\int_Y {\Bigl[\frac{1}{{|Z|}}\int_Z {(\beta _{ij}^{(0)} + C_{ijkl}^{(0)}\frac{{\partial {O_k}}}{{\partial {z_l}}})dZ}  + \frac{1}{{|Z|}}\int_Z {(C_{ijkl}^{(0)} + C_{ijpn}^{(0)}\frac{{\partial T_{pk}^l}}{{\partial {z_n}}})} dZ\frac{{\partial {P_k}}}{{\partial {y_l}}}\Bigl]dY} \\
	&= \frac{1}{{|Y||Z|}}\int_Y {\Bigl[\int_Z {(\beta _{ij}^{(0)} + C_{ijkl}^{(0)}\frac{{\partial {O_k}}}{{\partial {z_l}}})dZ + \int_Z {(C_{ijkl}^{(0)} + C_{ijpn}^{(0)}\frac{{\partial T_{pk}^l}}{{\partial {z_n}}})dZ} \frac{{\partial {P_k}}}{{\partial {y_l}}}} } \Bigl]dY.
	\end{aligned}
\end{equation}

Then, we acquire the following approximate relationships between the two fundamentally distinct macroscopic homogenization coefficients.
\begin{equation}
\begin{aligned}
	\hat S({\theta _0}) &= \frac{1}{{|{Y^{{\zeta  _2}}}|}}\int_{{Y^{{\zeta  _2}}}} {\Bigl[{\rho ^{(0)}}({\bm{y}},{\theta _0}){c^{(0)}}({\bm{y}},{\theta _0}) - \vartheta _{ij}^{(0)}({\bm{y}},{\theta _0})\frac{{\partial {P_i}}}{{\partial {y_j}}}\Bigl]d{Y^{{\zeta  _2}}}} \\
	&= \frac{1}{{|{Y^{{\zeta  _2}}}|}}\int_{{Y^{{\zeta  _2}}}} {\Bigl[{\rho ^{(0)}}{c^{(0)}} - \vartheta _{ij}^{(0)}(\frac{\partial }{{\partial {y_j}}} + \frac{{{\zeta  _1}}}{{{\zeta  _2}}}\frac{\partial }{{\partial {z_j}}})(P_i^0 + \frac{{{\zeta  _2}}}{{{\zeta  _1}}}(T_{im}^{{\alpha _1}}\frac{{\partial P_m^0}}{{\partial {y_{{\alpha _1}}}}} + {O_i}) + {\rm{O}}(\frac{{\zeta  _2^2}}{{\zeta  _1^2}}))\Bigl]d{Y^{{\zeta  _2}}}} \\
	&= \frac{1}{{|{Y^{{\zeta  _2}}}|}}\int_{{Y^{{\zeta  _2}}}} {\Bigl[{\rho ^{(0)}}{c^{(0)}} - \vartheta _{ij}^{(0)}\frac{{\partial {O_i}}}{{\partial {z_j}}} - (\vartheta _{m{\alpha _1}}^{(0)} + \vartheta _{ij}^{(0)}\frac{{\partial T_{im}^{{\alpha _1}}}}{{\partial {z_j}}})\frac{{\partial P_m^0}}{{\partial {y_{{\alpha _1}}}}} + {\rm{O}}(\frac{{{\zeta  _2}}}{{{\zeta  _1}}})\Bigl]d{Y^{{\zeta  _2}}}} \\
	&= {{\hat S}^*}({\theta _0}) + {\rm{O}}(\zeta  _1^{ - 1}{\zeta  _2}).
\end{aligned}
\end{equation}
\begin{equation}
	\begin{aligned}
	{{\hat k}_{ij}}({\theta _0}) &= \frac{1}{{|{Y^{{\zeta  _2}}}|}}\int_{{Y^{{\zeta  _2}}}} {\Bigl[k_{ij}^{(0)}({\bm{y}},{\theta _0}) + k_{ik}^{(0)}({\bm{y}},{\theta _0})\frac{{\partial {M_j}}}{{\partial {y_k}}}\Bigl]d{Y^{{\zeta  _2}}}} \\
	&= \frac{1}{{|{Y^{{\zeta  _2}}}|}}\int_{{Y^{{\zeta  _2}}}} {\Bigl[k_{ij}^{(0)} + k_{ik}^{(0)}(\frac{\partial }{{\partial {y_k}}} + \frac{{{\zeta  _1}}}{{{\zeta  _2}}}\frac{\partial }{{\partial {z_k}}})(M_j^0 + \frac{{{\zeta  _2}}}{{{\zeta  _1}}}({R_k}\frac{{\partial M_j^0}}{{\partial {y_k}}} + {R_{{\alpha _1}}}) + {\rm{O}}(\frac{{\zeta  _2^2}}{{\zeta  _1^2}}))\Bigl]d{Y^{{\zeta  _2}}}} \\
	&= \frac{1}{{|{Y^{{\zeta  _2}}}|}}\int_{{Y^{{\zeta  _2}}}} {\Bigl[k_{ij}^{(0)} + k_{ik}^{(0)}\frac{{\partial {R_{{\alpha _1}}}}}{{\partial {z_k}}} + (k_{ik}^{(0)} + k_{ik}^{(0)}\frac{{\partial {R_k}}}{{\partial {z_k}}})\frac{{\partial M_j^0}}{{\partial {y_k}}} + {\rm{O}}(\frac{{{\zeta  _2}}}{{{\zeta  _1}}})\Bigl]d{Y^{{\zeta  _2}}}} \\
	&= {{\hat k}^*}_{ij}({\theta _0}) + {\rm{O}}(\zeta  _1^{ - 1}{\zeta  _2}).
	\end{aligned}
\end{equation}
\begin{equation}
\begin{aligned}
	{{\hat \vartheta }_{ij}}({\theta _0}) &= \frac{1}{{|{Y^{{\zeta  _2}}}|}}\int_{{Y^{{\zeta  _2}}}} {\Bigl[\vartheta _{ij}^{(0)}(} {\bm{y}},{\theta _0}) + \vartheta _{mn}^{(0)}({\bm{y}},{\theta _0})\frac{{\partial N_{mi}^j}}{{\partial {y_n}}}\Bigl]d{Y^{{\zeta  _2}}}\\
	&= \frac{1}{{|{Y^{{\zeta  _2}}}|}}\int_{{Y^{{\zeta  _2}}}} {\Bigl[\vartheta _{ij}^{(0)} + \vartheta _{mn}^{(0)}(\frac{\partial }{{\partial {y_n}}} + \frac{{{\zeta  _1}}}{{{\zeta  _2}}}\frac{\partial }{{\partial {z_n}}})(N_{mi}^{0j} + \frac{{{\zeta  _2}}}{{{\zeta  _1}}}[T_{mp}^n\frac{{\partial N_{pi}^{0j}}}{{\partial {y_n}}} + T_{mi}^j] + {\rm{O}}(\frac{{\zeta  _2^2}}{{\zeta  _1^2}}))\Bigl]d{Y^{{\zeta  _2}}}} \\
	&= \frac{1}{{|{Y^{{\zeta  _2}}}|}}\int_{{Y^{{\zeta  _2}}}} {\Bigl[\vartheta _{ij}^{(0)} + \vartheta _{mn}^{(0)}\frac{{\partial T_{mi}^j}}{{\partial {z_n}}} + (\vartheta _{pn}^{(0)} + \vartheta _{mn}^{(0)}\frac{{\partial T_{mp}^n}}{{\partial {z_n}}})\frac{{\partial N_{pi}^{0j}}}{{\partial {y_n}}} + {\rm{O}}(\frac{{{\zeta  _2}}}{{{\zeta  _1}}})\Bigl]d{Y^{{\zeta  _2}}}} \\
	&= {{\hat \vartheta }^*}_{ij}({\theta _0}) + {\rm{O}}(\zeta  _1^{ - 1}{\zeta  _2}).
\end{aligned}
\end{equation}
\begin{equation}
	\hat \rho ({\theta _0}) = \frac{1}{{|{Y^{{\zeta  _2}}}|}}\int_{{Y^{{\zeta  _2}}}} {{\rho ^{(0)}}} ({\bm{y}},{\theta _0})d{Y^{{\zeta  _2}}} = {\hat \rho ^*}({\theta _0}).
\end{equation}
\begin{equation}
\begin{aligned}
	{{\hat C}_{ijkl}}({\theta _0}) &= \frac{1}{{|{Y^{{\zeta  _2}}}|}}\int_{{Y^{{\zeta  _2}}}} {\Bigl[C_{ijkl}^{(0)}({\bm{y}},{\theta _0}) + C_{ijmn}^{(0)}({\bm{y}},{\theta _0})\frac{{\partial N_{mk}^l}}{{\partial {y_n}}}\Bigl]d{Y^{{\zeta  _2}}}} \\
	&= \frac{1}{{|{Y^{{\zeta  _2}}}|}}\int_{{Y^{{\zeta  _2}}}} {\Bigl[C_{ijkl}^{(0)} + C_{ijmn}^{(0)}(\frac{\partial }{{\partial {y_n}}} + \frac{{{\zeta  _1}}}{{{\zeta  _2}}}\frac{\partial }{{\partial {z_n}}})(N_{mk}^{0l} + \frac{{{\zeta  _2}}}{{{\zeta  _1}}}(T_{mp}^n\frac{{\partial N_{pk}^{0l}}}{{\partial {y_n}}} + T_{mk}^l)\! +\! {\rm{O}}(\frac{{\zeta  _2^2}}{{\zeta  _1^2}}))\Bigl]d{Y^{{\zeta  _2}}}} \\
	&= \frac{1}{{|{Y^{{\zeta  _2}}}|}}\int_{{Y^{{\zeta  _2}}}} {\Bigl[C_{ijkl}^{(0)} + C_{ijmn}^{(0)}\frac{{\partial T_{mk}^l}}{{\partial {z_n}}} + (C_{ijpn}^{(0)} + C_{ijmn}^{(0)}\frac{{\partial T_{mp}^n}}{{\partial {z_n}}})\frac{{\partial N_{pk}^{0l}}}{{\partial {y_n}}} + {\rm{O}}(\frac{{{\zeta  _2}}}{{{\zeta  _1}}})\Bigl]d{Y^{{\zeta  _2}}}} \\
	&= {{\hat C}^*}_{ijkl}({\theta _0}) + {\rm{O}}(\zeta  _1^{ - 1}{\zeta  _2}).
\end{aligned}
\end{equation}
\begin{equation}
\begin{aligned}
	{{\hat \beta }_{ij}}({\theta _0}) &= \frac{1}{{|{Y^{{\zeta  _2}}}|}}\int_{{Y^{{\zeta  _2}}}} {\Bigl[\beta _{ij}^{(0)}({\bm{y}},{\theta _0}) + C_{ijkl}^{(0)}({\bm{y}},{\theta _0})\frac{{\partial {P_k}}}{{\partial {y_l}}}\Bigl]d{Y^{{\zeta  _2}}}} \\
	&= \frac{1}{{|{Y^{{\zeta  _2}}}|}}\int_{{Y^{{\zeta  _2}}}} {\Bigl[\beta _{ij}^{(0)} + C_{ijkl}^{(0)}(\frac{\partial }{{\partial {y_l}}} + \frac{{{\zeta  _1}}}{{{\zeta  _2}}}\frac{\partial }{{\partial {z_l}}})(P_k^0 + \frac{{{\zeta  _2}}}{{{\zeta  _1}}}(T_{km}^{{\alpha _1}}\frac{{\partial P_m^0}}{{\partial {y_{{\alpha _1}}}}} + {O_k}) + {\rm{O}}(\frac{{\zeta  _2^2}}{{\zeta  _1^2}}))\Bigl]d{Y^{{\zeta  _2}}}} \\
	&= \frac{1}{{|{Y^{{\zeta  _2}}}|}}\int_{{Y^{{\zeta  _2}}}} {\Bigl[\beta _{ij}^{(0)} + C_{ijkl}^{(0)}\frac{{\partial {O_k}}}{{\partial {z_l}}} + (C_{ijm{\alpha _1}}^{(0)} + C_{ijkl}^{(0)}\frac{{\partial T_{km}^{{\alpha _1}}}}{{\partial {z_l}}})\frac{{\partial P_m^0}}{{\partial {y_{{\alpha _1}}}}} + {\rm{O}}(\frac{{{\zeta  _2}}}{{{\zeta  _1}}})\Bigl]d{Y^{{\zeta  _2}}}} \\
	&= {{\hat \beta }^*}_{ij}({\theta _0}) + {\rm{O}}(\zeta  _1^{ - 1}{\zeta  _2}).
\end{aligned}
\end{equation}

On the basis of mathematical relationships, tractable $\hat{S}^*(\theta_0)$, $\hat{k}_{ij}^*(\theta_0)$, ${{\hat \vartheta }^*}_{ij}({\theta _0})$, $ {\hat \rho ^*}({\theta _0})$, $ {{\hat C}^*}_{ijkl}({\theta _0})$ and ${{\hat \beta }^*}_{ij}({\theta _0})$ can be applied to replace $\hat{S}(\theta_0)$, $\hat{k}_{ij}(\theta_0)$, ${{\hat \vartheta }}_{ij}({\theta _0})$, $ {\hat \rho }({\theta _0})$, $ {{\hat C}}_{ijkl}({\theta _0})$ and ${{\hat \beta }}_{ij}({\theta _0})$ in the presented HOTS method and associated three-scale calculation.

\section{Definitions of several auxiliary cell functions and mesoscopic homogenized material parameters}
\label{app2}
Considering the SOTS analysis for mesoscopic cell function $N_{im}^{{\alpha _1}}$ in subsection 2.3, we can derive the following SOTS forms (53)-(66) for the remaining mesoscopic cell functions ${M_{\alpha_1}}$, $P_i$, $A$, $M_{\alpha_1\alpha_2}$, $C_{\alpha_1}$, ${B_{\alpha_1\alpha_2}}$, ${E_{\alpha_1\alpha_2}}$, $F_i^{\alpha _1}$, $N_{im}^{{\alpha _1}{\alpha _2}}$, $ Z_{im}^{{\alpha _1}}$, ${Q_i}$, $H_i^{{\alpha _1}}$, $W_i^{{\alpha _1}}$ and $J_{im}^{{\alpha _1}{\alpha _2}}$ in two-scale solutions (35) and (36). The detailed definitions of several auxiliary cell functions are given as follows.
\begin{equation}
	\left\{
	\begin{aligned}
		&\frac{\partial }{{\partial {z_i}}}[k_{ij}^{(0)}(\bm{y,z},{\theta _0})\frac{{\partial {R_{{\alpha _1}}}}}{{\partial {z_j}}}] = - \frac{{\partial k_{i{\alpha _1}}^{(0)}(\bm{y,z},{\theta _0})}}{{\partial {z_i}}},&\;\;\;&\bm{z}\in Z,\\
		&{R_{{\alpha _1}}}(\bm{y,z},\theta_0) = 0,&\;\;\;&\bm{z}\in\partial Z.
	\end{aligned} \right.
\end{equation}
\begin{equation}
	\left\{
	\begin{aligned}
		&\frac{\partial }{{\partial {z_j}}}[C_{ijkl}^{(0)}({\bm{y}},{\bm{z}},{\theta _0})\frac{{\partial {O_k}}}{{\partial {z_l}}}] =  - \frac{{\partial \beta _{ij}^{(0)}({\bm{y}},{\bm{z}},{\theta _0})}}{{\partial {z_j}}},&\;\;\;&\bm{z}\in Z,\\
		&{O_k}({\bm{y}},{\bm{z}},{\theta _0}) = 0,&\;\;\;&\bm{z}\in\partial Z.
	\end{aligned}  \right.
\end{equation}
\begin{equation}
	\begin{aligned}
\bar{k}_{ij}(\bm{y},{\theta_0}) = \frac{1}{{\left| Z \right|}}\int_Z {\Big( {k_{ij}^{(0)}(\bm{y,z},{\theta_0}) + k_{ik}^{(0)}(\bm{y,z},{\theta_0})\frac{{\partial {R_j}}}{{\partial {z_{k}}}}} \Big)} dZ.
	\end{aligned}
\end{equation}
\begin{equation}
	\begin{aligned}
{{\bar \beta }_{ij}}({\bm{y}},{\theta _0}) = \frac{1}{{|Z|}}\int_Z \Bigl( \beta _{ij}^{(0)}({\bm{y}},{\bm{z}},{\theta _0}) + C_{ijkl}^{(0)}({\bm{y}},{\bm{z}},{\theta _0})\frac{{\partial {O_k}}}{{\partial {z_l}}}\Bigl)dZ.
	\end{aligned}
\end{equation}
\begin{equation}
	\bar S({\bm{y}},{\theta _0}) = \frac{1}{{|Z|}}\int_Z ( {\rho ^{(0)}}({\bm{y}},{\bm{z}},{\theta _0}){c^{(0)}}({\bm{y}},{\bm{z}},{\theta _0}) - \vartheta _{ij}^{(0)}({\bm{y}},{\bm{z}},{\theta _0})\frac{{\partial {O_i}}}{{\partial {z_j}}})dZ.
\end{equation}
\begin{equation}
	\bar \rho ({\bm{y}},{\theta _0}) = \frac{1}{{|Z|}}\int_Z {{\rho ^{(0)}}({\bm{y}},{\bm{z}},{\theta _0})} dZ.
\end{equation}
\begin{equation}
	{\bar \vartheta _{{\alpha _1}{\alpha _2}}}({\bm{y}},{\theta _0}) = \frac{1}{{|Z|}}\int_Z ( \vartheta _{{\alpha _1}{\alpha _2}}^{(0)}({\bm{y}},{\bm{z}},{\theta _0}) + \vartheta _{ij}^{(0)}({\bm{y}},{\bm{z}},{\theta _0})\frac{{\partial T_{i{\alpha _1}}^{{\alpha _2}}}}{{\partial {y_j}}})dZ
\end{equation}
\begin{equation}
	\left\{
	\begin{aligned}
		&\frac{\partial }{{\partial {z_i}}}[k_{ij}^{(0)}({\bm{y}},{\bm{z}},{\theta _0})\frac{{\partial {R_{{\alpha _1}{\alpha _2}}}}}{{\partial {z_j}}}] = {{\bar k}_{{\alpha _1}{\alpha _2}}} - k_{{\alpha _1}{\alpha _2}}^{(0)} - k_{{\alpha _1}j}^{(0)}\frac{{\partial {R_{{\alpha _2}}}}}{{\partial {z_j}}} - \frac{\partial }{{\partial {z_i}}}(k_{i{\alpha _1}}^{(0)}{R_{{\alpha _2}}}),&\;\;\;&\bm{z}\in Z,\\
		&{R_{{\alpha _1}{\alpha _2}}}({\bm{y}},{\bm{z}},{\theta _0}) = 0,&\;\;\;&\bm{z}\in\partial Z.
	\end{aligned}\right.
\end{equation}
\begin{equation}
	\left\{
	\begin{aligned}
		&\frac{\partial }{{\partial {z_i}}}[k_{ij}^{(0)}({\bm{y}},{\bm{z}},{\theta _0})\frac{{\partial {D_{{\alpha _1}}}}}{{\partial {z_j}}}] = \frac{{\partial {{\bar k}_{i{\alpha _1}}}}}{{\partial {y_i}}} - \frac{{\partial k_{i{\alpha _1}}^{(0)}}}{{\partial {y_i}}} - \frac{\partial }{{\partial {y_i}}}(k_{ij}^{(0)}\frac{{\partial {R_{{\alpha _1}}}}}{{\partial {z_j}}}) - \frac{\partial }{{\partial {z_i}}}(k_{ij}^{(0)}\frac{{\partial {R_{{\alpha _1}}}}}{{\partial {y_j}}}),&\;\;\;&\bm{z}\in Z,\\
		&{D_{{\alpha _1}}}({\bm{y}},{\bm{z}},{\theta _0}) = 0,&\;\;\;&\bm{z}\in\partial Z.
	\end{aligned} \right.
\end{equation}
\begin{equation}
	\left\{
	\begin{aligned}
		&\frac{\partial }{{\partial {z_j}}}[C_{ijkl}^{(0)}({\bm{y}},{\bm{z}},{\theta _0})\frac{{\partial P_k^ * }}{{\partial {z_l}}}] = \frac{{\partial {{\bar \beta }_{ij}}}}{{\partial {y_j}}} - \frac{{\partial \beta _{ij}^{(0)}}}{{\partial {y_j}}} - \frac{\partial }{{\partial {y_j}}}(C_{ijkl}^{(0)}\frac{{\partial {O_k}}}{{\partial {z_l}}}) - \frac{\partial }{{\partial {z_j}}}(C_{ijkl}^{(0)}\frac{{\partial {O_k}}}{{\partial {y_l}}}),&\;\;\;&\bm{z}\in Z,\\
		&P_k^ * ({\bm{y}},{\bm{z}},{\theta _0}) = 0,&\;\;\;&\bm{z}\in\partial Z.
	\end{aligned} \right.
\end{equation}
\begin{equation}
	\left\{
	\begin{aligned}
		&\frac{\partial }{{\partial {z_i}}}(k_{ij}^{(0)}({\bm{y}},{\bm{z}},{\theta _0})\frac{{\partial {A^ * }}}{{\partial {z_j}}}) = {\rho ^{(0)}}{c^{(0)}} - \bar S - \vartheta _{ij}^{(0)}\frac{{\partial {O_i}}}{{\partial {z_j}}} - \vartheta _{ij}^{(0)}\frac{{\partial P_i^0}}{{\partial {y_j}}}\\
		&\quad\quad\quad\quad\quad\quad\quad - \vartheta _{ij}^{(0)}\frac{{\partial T_{im}^{{\alpha _1}}}}{{\partial {z_j}}}\frac{{\partial P_m^0}}{{\partial {y_{{\alpha _1}}}}} + {{\bar \vartheta }_{m{\alpha _1}}}\frac{{\partial P_m^0}}{{\partial {y_{{\alpha _1}}}}},&\;\;\;&\bm{z}\in Z,\\
		&{A^ * }({\bm{y}},{\bm{z}},{\theta _0}) = 0,&\;\;\;&\bm{z}\in\partial Z.
	\end{aligned}\right.
\end{equation}
\begin{equation}
	\left\{
\begin{aligned}
	&\frac{\partial }{{\partial {z_i}}}(k_{ij}^{(0)}({\bm{y}},{\bm{z}},{\theta _0})\frac{{\partial D_{{\alpha _1}}^ * }}{{\partial {z_j}}}) = \frac{{\partial {{\bar k}_{i{\alpha _1}}}}}{{\partial {x_i}}} - \frac{{\partial k_{i{\alpha _1}}^{(0)}}}{{\partial {x_i}}} - \frac{\partial }{{\partial {x_i}}}(k_{ij}^{(0)}\frac{{\partial {R_{{\alpha _1}}}}}{{\partial {z_j}}}) - \frac{\partial }{{\partial {z_i}}}(k_{ij}^{(0)}\frac{{\partial {R_{{\alpha _1}}}}}{{\partial {x_j}}}),&\;\;\;&\bm{z}\in Z,\\
	&D_{{\alpha _1}}^ * ({\bm{y}},{\bm{z}},{\theta _0}) = 0,&\;\;\;&\bm{z}\in\partial Z.
\end{aligned}\right.	
\end{equation}
\begin{equation}
	\left\{
\begin{aligned}
	&\frac{\partial }{{\partial {z_i}}}(k_{ij}^{(0)}({\bm{y}},{\bm{z}},{\theta _0})\frac{{\partial R_{{\alpha _1}}^ * }}{{\partial {z_j}}}) = \frac{{\partial {{\bf{D}}^{(0,0,1)}}k_{i{\alpha _1}}^{(0)}({\bm{y}},{\bm{z}},{\theta _0})}}{{\partial {z_i}}},&\;\;\;&\bm{z}\in Z,\\
	&R_{{\alpha _1}}^ * ({\bm{y}},{\bm{z}},{\theta _0}) = 0,&\;\;\;&\bm{z}\in\partial Z.
\end{aligned}\right.	
\end{equation}
\begin{equation}
	\left\{
	\begin{aligned}
		&\frac{\partial }{{\partial {z_i}}}(k_{ij}^{(0)}({\bm{y}},{\bm{z}},{\theta _0})\frac{{\partial E_{{\alpha _1}{\alpha _2}}^ * }}{{\partial {z_j}}}) = {{\bar \vartheta }_{{\alpha _1}{\alpha _2}}} - \vartheta _{{\alpha _1}{\alpha _2}}^{(0)} - \vartheta _{ij}^{(0)}\frac{{\partial T_{i{\alpha _1}}^{{\alpha _2}}}}{{\partial {z_j}}},&\;&\bm{z}\in Z,\\
		&E_{{\alpha _1}{\alpha _2}}^ * ({\bm{y}},{\bm{z}},{\theta _0}) = 0,&\;&\bm{z}\in\partial Z.
	\end{aligned}\right.	
\end{equation}
\begin{equation}
	\left\{
\begin{aligned}
	&\frac{\partial }{{\partial {z_j}}}(C_{ijkl}^{(0)}({\bm{y}},{\bm{z}},{\theta _0})\frac{{\partial {F_k^*}}}{{\partial {z_l}}}) = {\rho ^{(0)}} - \bar \rho ,&\;\;\;&\bm{z}\in Z,\\
	&{F_k^*}({\bm{y}},{\bm{z}},{\theta _0}) = 0,&\;\;\;&\bm{z}\in\partial Z.
\end{aligned}\right.	
\end{equation}
\begin{equation}
	\left\{
	\begin{aligned}
		&\frac{\partial }{{\partial {z_j}}}(C_{ijkl}^{(0)}({\bm{y}},{\bm{z}},{\theta _0})\frac{{\partial U_{km}^{{\alpha _1}}}}{{\partial {z_l}}}) = \frac{{\partial {{\bar C}_{ijm{\alpha _1}}}}}{{\partial {x_j}}} - \frac{{\partial C_{ijm{\alpha _1}}^{(0)}}}{{\partial {x_j}}}\\
& \quad\quad- \frac{\partial }{{\partial {x_j}}}(C_{ijkl}^{(0)}\frac{{\partial T_{km}^{{\alpha _1}}}}{{\partial {z_l}}}) - \frac{\partial }{{\partial {z_j}}}(C_{ijkl}^{(0)}\frac{{\partial T_{km}^{{\alpha _1}}}}{{\partial {x_l}}}),&\;\;\;&\bm{z}\in Z,\\
		&U_{km}^{{\alpha _1}}({\bm{y}},{\bm{z}},{\theta _0}) = 0,&\;\;\;&\bm{z}\in\partial Z.
	\end{aligned}\right.	
\end{equation}
\begin{equation}
	\left\{
\begin{aligned}
	&\frac{\partial }{{\partial {z_j}}}(C_{ijkl}^{(0)}({\bm{y}},{\bm{z}},{\theta _0})\frac{{\partial Q_k^ * }}{{\partial {z_l}}}) = \frac{{\partial {{\bar \beta }_{ij}}}}{{\partial {x_j}}} - \frac{{\partial \beta _{ij}^{(0)}}}{{\partial {x_j}}} \\
&\quad\quad- \frac{\partial }{{\partial {x_j}}}(C_{ijkl}^{(0)}\frac{{\partial {O_k}}}{{\partial {z_l}}}) - \frac{\partial }{{\partial {z_j}}}(C_{ijkl}^{(0)}\frac{{\partial {O_k}}}{{\partial {x_l}}}),&\;\;\;&\bm{z}\in Z,\\
	&Q_k^ * ({\bm{y}},{\bm{z}},{\theta _0}) = 0,&\;\;\;&\bm{z}\in\partial Z.
\end{aligned}\right.	
\end{equation}
\begin{equation}
	\left\{
	\begin{aligned}
		&\frac{\partial }{{\partial {z_j}}}(C_{ijkl}^{(0)}({\bm{y}},{\bm{z}},{\theta _0})\frac{{\partial V_k^{{\alpha _1}}}}{{\partial {z_l}}}) = {{\bar \beta }_{i{\alpha _1}}} - \beta _{i{\alpha _1}}^{(0)} - C_{i{\alpha _1}kj}^{(0)}\frac{{\partial {O_k}}}{{\partial {z_j}}}\\
 &\quad- \frac{\partial }{{\partial {z_j}}}(C_{ijk{\alpha _1}}^{(0)}{O_k}) - \frac{\partial }{{\partial {z_j}}}(\beta _{ij}^{(0)}{R_{{\alpha _1}}}),&\;\;&\bm{z}\in Z,\\
		&V_k^{{\alpha _1}}({\bm{y}},{\bm{z}},{\theta _0}) = 0,&\;\;&\bm{z}\in\partial Z.
	\end{aligned}\right.	
\end{equation}
\begin{equation}
	\left\{
	\begin{aligned}
		&\frac{\partial }{{\partial {z_j}}}(C_{ijkl}^{(0)}({\bm{y}},{\bm{z}},{\theta _0})\frac{{\partial T_{km}^{ * n}}}{{\partial {z_l}}}) = \frac{{\partial {{\bf{D}}^{(0,0,1)}}C_{ijmn}^{(0)}({\bm{y}},{\bm{z}},{\theta _0})}}{{\partial {z_j}}},&\;\;\;&\bm{z}\in Z,\\
		&T_{km}^{ * n}({\bm{y}},{\bm{z}},{\theta _0}) = 0,&\;\;\;&\bm{z}\in\partial Z.
	\end{aligned}\right.	
\end{equation}
\begin{equation}
	\left\{
	\begin{aligned}
		&\frac{\partial }{{\partial {z_j}}}(C_{ijkl}^{(0)}({\bm{y}},{\bm{z}},{\theta _0})\frac{{\partial O_k^ * }}{{\partial {z_l}}}) = \frac{{\partial {{\bf{D}}^{(0,0,1)}}\beta _{ij}^{(0)}({\bm{y}},{\bm{z}},{\theta _0})}}{{\partial {z_j}}},&\;\;\;&\bm{z}\in Z,\\
		&O_k^ * ({\bm{y}},{\bm{z}},{\theta _0}) = 0,&\;\;\;&\bm{z}\in\partial Z.
	\end{aligned}\right.	
\end{equation}
\begin{equation}
	\left\{
	\begin{aligned}
		&\frac{\partial}{\partial y_i}\Big(\bar{k}_{ij}(\bm{y},\theta_0)\frac{\partial M_{\alpha_1}^0}{\partial y_j}\Big)=-\frac{\partial\bar{k}_{i\alpha_1}}{\partial y_i},&\;\;\;&\bm{y}\in Y,\\
		&M_{\alpha_1}^0(\bm{y},\theta_0)=0,&\;\;\;&\bm{y}\in\partial Y.
	\end{aligned}\right.
\end{equation}
\begin{equation}
	\left\{
	\begin{aligned}
		&\frac{\partial }{{\partial {y_j}}}\Bigl({{\bar C}_{ijkl}}({\bm{y}},{\theta _0})\frac{{\partial P_k^0}}{{\partial {y_l}}}\Bigl) =  - \frac{{\partial {{\bar \beta }_{ij}}}}{{\partial {y_j}}},&\;\;\;&\bm{y}\in Y,\\
		&P_k^0({\bm{y}},{\theta _0}) = 0,&\;\;\;&\bm{y}\in\partial Y.
	\end{aligned}\right.
\end{equation}
\begin{equation}
	\left\{
	\begin{aligned}
		&\frac{\partial }{{\partial {y_i}}}[{{\bar k}_{ij}}({\bm{y}},{\theta _0})\frac{{\partial {A^0}}}{{\partial {y_j}}}]\! = \!\bar S({\bm{y}},{\theta _0})\! - \!\hat S({\theta _0}) \!- \!{{\bar \vartheta }_{m{\alpha _1}}}({\bm{y}},{\theta _0})\frac{{\partial P_m^0}}{{\partial {y_{{\alpha _1}}}}},&&\bm{y}\in Y,\\
		&{A^0}({\bm{y}},{\theta _0}) = 0,&&\bm{y}\in\partial Y.
	\end{aligned}\right.
\end{equation}
\begin{equation}
	\left\{
\begin{aligned}
	&\frac{\partial }{{\partial {y_i}}}[{{\bar k}_{ij}}({\bm{y}},{\theta _0})\frac{{\partial C_{{\alpha _1}}^0}}{{\partial {y_j}}}] = \frac{{\partial {{\hat k}_{i{\alpha _1}}}}}{{\partial {x_i}}} - \frac{{\partial {{\bar k}_{i{\alpha _1}}}}}{{\partial {x_i}}} - \frac{\partial }{{\partial {y_i}}}({{\bar k}_{ij}}\frac{{\partial M_{{\alpha _1}}^0}}{{\partial {x_j}}}) - \frac{\partial }{{\partial {x_i}}}({{\bar k}_{ij}}\frac{{\partial M_{{\alpha _1}}^0}}{{\partial {y_j}}}),&\;\;\;&\bm{y}\in Y,\\
	&C_{{\alpha _1}}^0({\bm{y}},{\theta _0}) = 0,&\;\;\;&\bm{y}\in\partial Y.
\end{aligned}\right.
\end{equation}
\begin{equation}
	\left\{
\begin{aligned}
	&\frac{\partial }{{\partial {y_i}}}({{\bar k}_{ij}}({\bm{y}},{\theta _0})\frac{{\partial B_{{\alpha _1}{\alpha _2}}^0}}{{\partial {y_j}}}) = \frac{\partial }{{\partial {y_i}}}(M_{{\alpha _1}}^0{{\bf{D}}^{(0,1)}}{{\bar k}_{i{\alpha _2}}} + M_{{\alpha _1}}^0{{\bf{D}}^{(0,1)}}{{\bar k}_{ij}}\frac{{\partial M_{{\alpha _2}}^0}}{{\partial {y_j}}}) \\
	&\quad\quad\quad\quad\quad\quad\quad\quad\quad\quad + \frac{\partial }{{\partial {y_i}}}(M_{{\alpha _1}}^0{{\bar k}_{i{\alpha _2}}} + M_{{\alpha _1}}^0{{\bar k}_{ij}}\frac{{\partial M_{{\alpha _2}}^0}}{{\partial {y_j}}}),&\;\;\;&\bm{y}\in Y,\\
	&B_{{\alpha _1}{\alpha _2}}^0({\bm{y}},{\theta _0}) = 0,&\;\;\;&\bm{y}\in\partial Y.
\end{aligned}\right.	
\end{equation}
\begin{equation}
	\left\{
	\begin{aligned}
		&\frac{\partial }{{\partial {y_i}}}[{{\bar k}_{ij}}({\bm{y}},{\theta _0})\frac{{\partial E_{{\alpha _1}{\alpha _2}}^0}}{{\partial {y_j}}}] = {{\hat \vartheta }_{{\alpha _1}{\alpha _2}}} - {{\bar \vartheta }_{{\alpha _1}{\alpha _2}}} - {{\bar \vartheta }_{ij}}\frac{{\partial N_{i{\alpha _1}}^{0{\alpha _2}}}}{{\partial {y_j}}},&\;\;\;&\bm{y}\in Y,\\
		&E_{{\alpha _1}{\alpha _2}}^0({\bm{y}},{\theta _0}) = 0,&\;\;\;&\bm{y}\in\partial Y.
	\end{aligned}\right.	
\end{equation}
\begin{equation}
	\left\{
\begin{aligned}
	&\frac{\partial }{{\partial {y_j}}}({{\bar C}_{ijkl}}({\bm{y}},{\theta _0})\frac{{\partial F_{k}^{{0\alpha _1}}}}{{\partial {y_l}}}) = \bar \rho  - \hat \rho ,&\;\;\;&\bm{y}\in Y,\\
	&F_{k}^{{0\alpha _1}}({\bm{y}},{\theta _0})  = 0,&\;\;\;&\bm{y}\in\partial Y.
\end{aligned}\right.	
\end{equation}
\begin{equation}
	\left\{
\begin{aligned}
	&\frac{\partial }{{\partial {y_j}}}({{\bar C}_{ijkl}}({\bm{y}},{\theta _0})\frac{{\partial Z_{km}^{0{\alpha _1}}}}{{\partial {y_l}}}) = \frac{{\partial {{\hat C}_{ijm{\alpha _1}}}}}{{\partial {x_j}}} - \frac{{\partial {{\bar C}_{ijm{\alpha _1}}}}}{{\partial {x_j}}} \\
&\quad\quad\quad\quad- \frac{\partial }{{\partial {y_j}}}({{\bar C}_{ijkl}}\frac{{\partial N_{km}^{0{\alpha _1}}}}{{\partial {x_l}}}) - \frac{\partial }{{\partial {x_j}}}({{\bar C}_{ijkl}}\frac{{\partial N_{km}^{0{\alpha _1}}}}{{\partial {y_l}}}),&\;\;\;&\bm{y}\in Y,\\
	&Z_{km}^{0{\alpha _1}}({\bm{y}},{\theta _0})= 0,&\;\;\;&\bm{y}\in\partial Y.
\end{aligned}\right.	
\end{equation}
\begin{equation}
	\left\{
\begin{aligned}
	&\frac{\partial }{{\partial {y_j}}}({{\bar C}_{ijkl}}({\bm{y}},{\theta _0})\frac{{\partial {Q_{k}^{0}}}}{{\partial {y_l}}}) = \frac{{\partial {{\hat \beta }_{ij}}}}{{\partial {x_j}}} - \frac{{\partial {{\bar \beta }_{ij}}}}{{\partial {x_j}}} - \frac{\partial }{{\partial {x_j}}}({{\bar C}_{ijkl}}\frac{{\partial P_k^0}}{{\partial {y_l}}}) - \frac{\partial }{{\partial {y_j}}}({{\bar C}_{ijkl}}\frac{{\partial P_k^0}}{{\partial {x_l}}}),&\;\;\;&\bm{y}\in Y,\\
	&{Q_{k}^{0}}({\bm{y}},{\theta _0}) = 0,&\;\;\;&\bm{y}\in\partial Y.
\end{aligned}\right.	
\end{equation}
\begin{equation}
		\left\{
	\begin{aligned}
		&\frac{\partial }{{\partial {y_j}}}({{\bar C}_{ijkl}}({\bm{y}},{\theta _0})\frac{{\partial H_{k}^{{0\alpha _1}}}}{{\partial {y_l}}}) = {{\hat \beta }_{i{\alpha _1}}} - {{\bar \beta }_{i{\alpha _1}}}\\
 &\quad\quad - \frac{\partial }{{\partial {y_j}}}({{\bar C}_{ijk{\alpha _1}}}P_k^0)- {{\bar C}_{i{\alpha _1}kl}}\frac{{\partial P_k^0}}{{\partial {y_l}}} - \frac{\partial }{{\partial {y_j}}}({{\bar \beta }_{ij}}M_{{\alpha _1}}^0),&\;\;\;&\bm{y}\in Y,\\
		&H_{k}^{{0\alpha _1}}({\bm{y}},{\theta _0}) = 0,&\;\;\;&\bm{y}\in\partial Y.
	\end{aligned}\right.	
\end{equation}
\begin{equation}
	\left\{
	\begin{aligned}
		&\frac{\partial }{{\partial {y_j}}}({{\bar C}_{ijkl}}({\bm{y}},{\theta _0})\frac{{\partial W_{k}^{{0\alpha _1}}}}{{\partial {y_l}}}) = \frac{\partial }{{\partial {y_j}}}[M_{{\alpha _1}}^0{{\bf{D}}^{(0,1)}}{{\bar \beta }_{ij}} + M_{{\alpha _1}}^0{{\bf{D}}^{(0,1)}}{{\bar C}_{ijkl}}\frac{{\partial P_k^0}}{{\partial {y_l}}}] \\
		&\quad\quad\quad\quad\quad\quad\quad\quad\quad\quad+ \frac{\partial }{{\partial {y_j}}}[{{\bar \beta }_{ij}}M_{{\alpha _1}}^0 + {{\bar C}_{ijkn}}M_{{\alpha _1}}^0\frac{{\partial P_k^0}}{{\partial {y_n}}}],&\;\;\;&\bm{y}\in Y,\\
		&W_{k}^{{0\alpha _1}}({\bm{y}},{\theta _0}) = 0,&\;\;\;&\bm{y}\in\partial Y.
	\end{aligned}\right.	
\end{equation}
\begin{equation}
	\left\{
	\begin{aligned}
		&\frac{\partial }{{\partial {y_j}}}({{\bar C}_{ijkl}}({\bm{y}},{\theta _0})\frac{{\partial J_{km}^{0{\alpha _1}{\alpha _2}}}}{{\partial {y_l}}}) = \frac{\partial }{{\partial {y_j}}}[M_{{\alpha _1}}^0{{\bf{D}}^{(0,1)}}{{\bar C}_{ijm{\alpha _2}}} + M_{{\alpha _1}}^0{{\bf{D}}^{(0,1)}}{{\bar C}_{ijkl}}\frac{{\partial N_{km}^{0{\alpha _2}}}}{{\partial {y_l}}}]\\
		&\quad\quad\quad\quad\quad\quad\quad\quad\quad\quad\quad + \frac{\partial }{{\partial {y_j}}}[M_{{\alpha _1}}^0{{\bar C}_{ijm{\alpha _2}}} + M_{{\alpha _1}}^0{{\bar C}_{ijkl}}\frac{{\partial N_{km}^{0{\alpha _2}}}}{{\partial {y_l}}}],&\;\;\;&\bm{y}\in Y,\\
		&J_{km}^{0{\alpha _1}{\alpha _2}}({\bm{y}},{\theta _0}) = 0,&\;\;\;&\bm{y}\in\partial Y.
	\end{aligned}\right.	
\end{equation}

\bibliographystyle{model1-num-names}
\bibliography{ref-1}

\begin{thebibliography}{50}
\expandafter\ifx\csname natexlab\endcsname\relax\def\natexlab#1{#1}\fi
\providecommand{\url}[1]{\texttt{#1}}
\providecommand{\href}[2]{#2}
\providecommand{\path}[1]{#1}
\providecommand{\DOIprefix}{doi:}
\providecommand{\ArXivprefix}{arXiv:}
\providecommand{\URLprefix}{URL: }
\providecommand{\Pubmedprefix}{pmid:}
\providecommand{\doi}[1]{\href{http://dx.doi.org/#1}{\path{#1}}}
\providecommand{\Pubmed}[1]{\href{pmid:#1}{\path{#1}}}
\providecommand{\bibinfo}[2]{#2}
\ifx\xfnm\relax \def\xfnm[#1]{\unskip,\space#1}\fi
%Type = Book
\bibitem[{Fan(2009)}]{R1}
\bibinfo{author}{X.~Fan}, \bibinfo{title}{{Thermal Structures Analysis and
  Applications of Highspeed Vehicles}}, \bibinfo{publisher}{National Defense
  Industry Press}, \bibinfo{year}{{2009}}.
%Type = Article
\bibitem[{Zhang et~al.(2020)Zhang, Wang, Wang, and Zhang}]{R2}
\bibinfo{author}{Y.~Zhang}, \bibinfo{author}{K.~Wang},
  \bibinfo{author}{B.~Wang}, \bibinfo{author}{C.~Zhang},
\newblock \bibinfo{title}{Thermal shock resistance of porous ceramic foams with
  temperature-dependent material properties},
\newblock \bibinfo{journal}{Ceramics International} \bibinfo{volume}{46}
  (\bibinfo{year}{2020}) \bibinfo{pages}{1503--1511}.
%Type = Article
\bibitem[{Abdoun and and(2020)}]{R3}
\bibinfo{author}{F.~Abdoun}, \bibinfo{author}{L.~A. and},
\newblock \bibinfo{title}{Thermal buckling and vibration of laminated composite
  plates with temperature dependent properties by an asymptotic numerical
  method},
\newblock \bibinfo{journal}{International Journal for Computational Methods in
  Engineering Science and Mechanics} \bibinfo{volume}{21}
  (\bibinfo{year}{2020}) \bibinfo{pages}{43--57}.
%Type = Article
\bibitem[{Dai and Huang(2021)}]{R4}
\bibinfo{author}{G.~Dai}, \bibinfo{author}{J.~Huang},
\newblock \bibinfo{title}{Nonlinear thermotics: designing thermal metamaterials
  with temperature response},
\newblock \bibinfo{journal}{Journal of Nantong University (Natural Science
  Edition)} \bibinfo{volume}{20} (\bibinfo{year}{2021}) \bibinfo{pages}{1--18}.
%Type = Article
\bibitem[{Nasirov et~al.(2020)Nasirov, Gupta, Hasanov, and Fidan}]{R5}
\bibinfo{author}{A.~Nasirov}, \bibinfo{author}{A.~Gupta},
  \bibinfo{author}{S.~Hasanov}, \bibinfo{author}{I.~Fidan},
\newblock \bibinfo{title}{Three-scale asymptotic homogenization of short fiber
  reinforced additively manufactured polymer composites},
\newblock \bibinfo{journal}{Composites Part B: Engineering}
  \bibinfo{volume}{202} (\bibinfo{year}{2020}) \bibinfo{pages}{108269}.
%Type = Article
\bibitem[{Song et~al.(2023)Song, Li, Zhang, Liu, Tang, Xie, and Liao}]{R6}
\bibinfo{author}{K.~Song}, \bibinfo{author}{D.~Li}, \bibinfo{author}{C.~Zhang},
  \bibinfo{author}{T.~Liu}, \bibinfo{author}{Y.~Tang}, \bibinfo{author}{Y.~M.
  Xie}, \bibinfo{author}{W.~Liao},
\newblock \bibinfo{title}{Bio-inspired hierarchical honeycomb metastructures
  with superior mechanical properties},
\newblock \bibinfo{journal}{Composite Structures} \bibinfo{volume}{304}
  (\bibinfo{year}{2023}) \bibinfo{pages}{116452}.
%Type = Article
\bibitem[{Liu et~al.(2024)Liu, Hou, Sapanathan, Nie, Meng, and Xu}]{R7}
\bibinfo{author}{Y.~Liu}, \bibinfo{author}{Y.~Hou},
  \bibinfo{author}{T.~Sapanathan}, \bibinfo{author}{R.~Nie},
  \bibinfo{author}{L.~Meng}, \bibinfo{author}{Y.~Xu},
\newblock \bibinfo{title}{A multiscale strategy for exploring the mechanical
  behavior of 3d braided composite thin-walled cylinders},
\newblock \bibinfo{journal}{Thin-Walled Structures} \bibinfo{volume}{198}
  (\bibinfo{year}{2024}) \bibinfo{pages}{111705}.
%Type = Article
\bibitem[{Li et~al.(2024)Li, Zuo, Li, and Liu}]{R8}
\bibinfo{author}{X.-N. Li}, \bibinfo{author}{X.-B. Zuo},
  \bibinfo{author}{L.~Li}, \bibinfo{author}{J.-H. Liu},
\newblock \bibinfo{title}{Multiscale modeling and simulation on mechanical
  behavior of fiber reinforced concrete},
\newblock \bibinfo{journal}{International Journal of Solids and Structures}
  \bibinfo{volume}{286} (\bibinfo{year}{2024}) \bibinfo{pages}{112569}.
%Type = Article
\bibitem[{Yang et~al.(2024)Yang, Liu, Xia, Fan, Taylor, and Han}]{R9}
\bibinfo{author}{H.~Yang}, \bibinfo{author}{Z.~Liu}, \bibinfo{author}{Y.~Xia},
  \bibinfo{author}{W.~Fan}, \bibinfo{author}{A.~C. Taylor},
  \bibinfo{author}{X.~Han},
\newblock \bibinfo{title}{Mechanical properties of hierarchical lattice via
  strain gradient homogenization approach},
\newblock \bibinfo{journal}{Composites Part B: Engineering}
  \bibinfo{volume}{271} (\bibinfo{year}{2024}) \bibinfo{pages}{111153}.
%Type = Article
\bibitem[{Zou et~al.(2015)Zou, Ainslie, Fujishiro, Bhagurkar, Naito, Babu,
  Fagnard, Vanderbemden, and Yamamoto}]{R10}
\bibinfo{author}{J.~Zou}, \bibinfo{author}{M.~Ainslie},
  \bibinfo{author}{H.~Fujishiro}, \bibinfo{author}{A.~Bhagurkar},
  \bibinfo{author}{T.~Naito}, \bibinfo{author}{N.~H. Babu},
  \bibinfo{author}{J.-F. Fagnard}, \bibinfo{author}{P.~Vanderbemden},
  \bibinfo{author}{A.~Yamamoto},
\newblock \bibinfo{title}{Numerical modelling and comparison of mgb2 bulks
  fabricated by hip and infiltration growth},
\newblock \bibinfo{journal}{Superconductor Science and Technology}
  \bibinfo{volume}{28} (\bibinfo{year}{2015}) \bibinfo{pages}{075009}.
%Type = Article
\bibitem[{Li et~al.(2014)Li, Wang, and Gao}]{R11}
\bibinfo{author}{Y.~Li}, \bibinfo{author}{X.~Wang}, \bibinfo{author}{Y.~Gao},
\newblock \bibinfo{title}{Computational method for elastic--plastic and
  anisotropic superconducting cable under simple load},
\newblock \bibinfo{journal}{International Journal of Computational Methods}
  \bibinfo{volume}{11} (\bibinfo{year}{2014}) \bibinfo{pages}{1344006}.
%Type = Book
\bibitem[{Bensoussan et~al.(2011)Bensoussan, Lions, and Papanicolaou}]{R12}
\bibinfo{author}{A.~Bensoussan}, \bibinfo{author}{J.-L. Lions},
  \bibinfo{author}{G.~Papanicolaou}, \bibinfo{title}{Asymptotic analysis for
  periodic structures}, volume \bibinfo{volume}{374},
  \bibinfo{publisher}{American Mathematical Soc.}, \bibinfo{year}{2011}.
%Type = Article
\bibitem[{Liu et~al.(2021)Liu, Chung, and Zhang}]{R13}
\bibinfo{author}{X.~Liu}, \bibinfo{author}{E.~Chung},
  \bibinfo{author}{L.~Zhang},
\newblock \bibinfo{title}{Iterated numerical homogenization for multiscale
  elliptic equations with monotone nonlinearity},
\newblock \bibinfo{journal}{Multiscale Modeling \& Simulation}
  \bibinfo{volume}{19} (\bibinfo{year}{2021}) \bibinfo{pages}{1601--1632}.
%Type = Article
\bibitem[{Efendiev et~al.(2000)Efendiev, Hou, and Wu}]{R14}
\bibinfo{author}{Y.~R. Efendiev}, \bibinfo{author}{T.~Y. Hou},
  \bibinfo{author}{X.-H. Wu},
\newblock \bibinfo{title}{Convergence of a nonconforming multiscale finite
  element method},
\newblock \bibinfo{journal}{SIAM Journal on Numerical Analysis}
  \bibinfo{volume}{37} (\bibinfo{year}{2000}) \bibinfo{pages}{888--910}.
%Type = Article
\bibitem[{Efendiev and Wu(2002)}]{R15}
\bibinfo{author}{Y.~R. Efendiev}, \bibinfo{author}{X.-H. Wu},
\newblock \bibinfo{title}{Multiscale finite element for problems with highly
  oscillatory coefficients},
\newblock \bibinfo{journal}{Numerische Mathematik} \bibinfo{volume}{90}
  (\bibinfo{year}{2002}) \bibinfo{pages}{459--486}.
%Type = Article
\bibitem[{Guan et~al.(2024)Guan, Jiang, and Wang}]{R16}
\bibinfo{author}{X.~Guan}, \bibinfo{author}{L.~Jiang},
  \bibinfo{author}{Y.~Wang},
\newblock \bibinfo{title}{Regularized coupling multiscale method for
  thermomechanical coupled problems},
\newblock \bibinfo{journal}{Journal of Computational Physics}
  \bibinfo{volume}{499} (\bibinfo{year}{2024}) \bibinfo{pages}{112737}.
%Type = Article
\bibitem[{Abdulle et~al.(2012)Abdulle, Weinan, Engquist, and
  Vanden-Eijnden}]{R17}
\bibinfo{author}{A.~Abdulle}, \bibinfo{author}{E.~Weinan},
  \bibinfo{author}{B.~Engquist}, \bibinfo{author}{E.~Vanden-Eijnden},
\newblock \bibinfo{title}{The heterogeneous multiscale method},
\newblock \bibinfo{journal}{Acta Numerica} \bibinfo{volume}{21}
  (\bibinfo{year}{2012}) \bibinfo{pages}{1--87}.
%Type = Article
\bibitem[{Zheng et~al.(2020)Zheng, Yifeng, Xiao, and Peng}]{R18}
\bibinfo{author}{S.~Zheng}, \bibinfo{author}{Z.~Yifeng},
  \bibinfo{author}{P.~Xiao}, \bibinfo{author}{W.~Peng},
\newblock \bibinfo{title}{Dimensional reduction analyzing the thermoelastic
  behavior of wind turbine blades based on the variational asymptotic
  multiscale method},
\newblock \bibinfo{journal}{Composite Structures} \bibinfo{volume}{254}
  (\bibinfo{year}{2020}) \bibinfo{pages}{112835}.
%Type = Article
\bibitem[{M{\aa}lqvist and Peterseim(2014)}]{R19}
\bibinfo{author}{A.~M{\aa}lqvist}, \bibinfo{author}{D.~Peterseim},
\newblock \bibinfo{title}{Localization of elliptic multiscale problems},
\newblock \bibinfo{journal}{Mathematics of Computation} \bibinfo{volume}{83}
  (\bibinfo{year}{2014}) \bibinfo{pages}{2583--2603}.
%Type = Article
\bibitem[{Xing et~al.(2010)Xing, Yang, and Wang}]{R20}
\bibinfo{author}{Y.~Xing}, \bibinfo{author}{Y.~Yang},
  \bibinfo{author}{X.~Wang},
\newblock \bibinfo{title}{A multiscale eigenelement method and its application
  to periodical composite structures},
\newblock \bibinfo{journal}{Composite Structures} \bibinfo{volume}{92}
  (\bibinfo{year}{2010}) \bibinfo{pages}{2265--2275}.
%Type = Article
\bibitem[{Geers et~al.(2010)Geers, Kouznetsova, and Brekelmans}]{R21}
\bibinfo{author}{M.~G. Geers}, \bibinfo{author}{V.~G. Kouznetsova},
  \bibinfo{author}{W.~Brekelmans},
\newblock \bibinfo{title}{Multi-scale computational homogenization: Trends and
  challenges},
\newblock \bibinfo{journal}{Journal of computational and applied mathematics}
  \bibinfo{volume}{234} (\bibinfo{year}{2010}) \bibinfo{pages}{2175--2182}.
%Type = Article
\bibitem[{Jin(2022)}]{R22}
\bibinfo{author}{S.~Jin},
\newblock \bibinfo{title}{Asymptotic-preserving schemes for multiscale physical
  problems},
\newblock \bibinfo{journal}{Acta Numerica} \bibinfo{volume}{31}
  (\bibinfo{year}{2022}) \bibinfo{pages}{415--489}.
%Type = Article
\bibitem[{Chen et~al.(2020)Chen, Ma, and Zhang}]{R23}
\bibinfo{author}{J.~Chen}, \bibinfo{author}{D.~Ma}, \bibinfo{author}{Z.~Zhang},
\newblock \bibinfo{title}{A multiscale reduced basis method for the
  schrödinger equation with multiscale and random potentials},
\newblock \bibinfo{journal}{Multiscale Modeling \& Simulation}
  \bibinfo{volume}{18} (\bibinfo{year}{2020}) \bibinfo{pages}{1409--1434}.
%Type = Article
\bibitem[{Dong et~al.(2019)Dong, Zheng, Cui, Nie, Yang, and Yang}]{R24}
\bibinfo{author}{H.~Dong}, \bibinfo{author}{X.~Zheng},
  \bibinfo{author}{J.~Cui}, \bibinfo{author}{Y.~Nie},
  \bibinfo{author}{Z.~Yang}, \bibinfo{author}{Z.~Yang},
\newblock \bibinfo{title}{Multiscale computational method for dynamic
  thermo-mechanical problems of composite structures with diverse periodic
  configurations in different subdomains},
\newblock \bibinfo{journal}{Journal of Scientific Computing}
  \bibinfo{volume}{79} (\bibinfo{year}{2019}) \bibinfo{pages}{1630--1666}.
%Type = Article
\bibitem[{Dong et~al.(2023)Dong, Cui, Nie, Ma, Jin, and Huang}]{R25}
\bibinfo{author}{H.~Dong}, \bibinfo{author}{J.~Cui}, \bibinfo{author}{Y.~Nie},
  \bibinfo{author}{R.~Ma}, \bibinfo{author}{K.~Jin},
  \bibinfo{author}{D.~Huang},
\newblock \bibinfo{title}{Multi-scale computational method for nonlinear
  dynamic thermo-mechanical problems of composite materials with
  temperature-dependent properties},
\newblock \bibinfo{journal}{Communications in Nonlinear Science and Numerical
  Simulation} \bibinfo{volume}{118} (\bibinfo{year}{2023})
  \bibinfo{pages}{107000}.
%Type = Article
\bibitem[{Yang et~al.(2016)Yang, Cui, Sun, Liang, and Yang}]{R26}
\bibinfo{author}{Z.~Yang}, \bibinfo{author}{J.~Cui}, \bibinfo{author}{Y.~Sun},
  \bibinfo{author}{J.~Liang}, \bibinfo{author}{Z.~Yang},
\newblock \bibinfo{title}{Multiscale analysis method for thermo-mechanical
  performance of periodic porous materials with interior surface radiation},
\newblock \bibinfo{journal}{International Journal for Numerical Methods in
  Engineering} \bibinfo{volume}{105} (\bibinfo{year}{2016})
  \bibinfo{pages}{323--350}.
%Type = Article
\bibitem[{Wang et~al.(2015)Wang, Cao, and Wong}]{R27}
\bibinfo{author}{X.~Wang}, \bibinfo{author}{L.~Cao}, \bibinfo{author}{Y.~Wong},
\newblock \bibinfo{title}{Multiscale computation and convergence for coupled
  thermoelastic system in composite materials},
\newblock \bibinfo{journal}{Multiscale Modeling \& Simulation}
  \bibinfo{volume}{13} (\bibinfo{year}{2015}) \bibinfo{pages}{661--690}.
%Type = Article
\bibitem[{Gaka et~al.(2001)Gaka, Telega, and Tokarzewski}]{R28}
\bibinfo{author}{A.~Gaka}, \bibinfo{author}{J.~J. Telega},
  \bibinfo{author}{S.~Tokarzewski},
\newblock \bibinfo{title}{Heat equation with temperature-dependent conductivity
  coefficients and macroscopic properties of microheterogeneous media},
\newblock \bibinfo{journal}{Mathematical and Computer Modelling}
  \bibinfo{volume}{33} (\bibinfo{year}{2001}) \bibinfo{pages}{927--942}.
%Type = Article
\bibitem[{Chung et~al.(2001)Chung, Tamma, and Namburu}]{R29}
\bibinfo{author}{P.~W. Chung}, \bibinfo{author}{K.~K. Tamma},
  \bibinfo{author}{R.~R. Namburu},
\newblock \bibinfo{title}{Homogenization of temperature-dependent thermal
  conductivity in composite materials},
\newblock \bibinfo{journal}{Journal of Thermophysics and Heat Transfer}
  \bibinfo{volume}{15} (\bibinfo{year}{2001}) \bibinfo{pages}{10--17}.
%Type = Article
\bibitem[{Abdulle et~al.(2015)Abdulle, Bai, and Vilmart}]{R30}
\bibinfo{author}{A.~Abdulle}, \bibinfo{author}{Y.~Bai},
  \bibinfo{author}{G.~Vilmart},
\newblock \bibinfo{title}{Reduced basis finite element heterogeneous multiscale
  method for quasilinear elliptic homogenization problems},
\newblock \bibinfo{journal}{Discrete and Continuous Dynamical Systems-Series S}
  \bibinfo{volume}{8} (\bibinfo{year}{2015}) \bibinfo{pages}{91--118}.
%Type = Article
\bibitem[{Abdulle and Huber(2016)}]{R31}
\bibinfo{author}{A.~Abdulle}, \bibinfo{author}{M.~E. Huber},
\newblock \bibinfo{title}{Finite element heterogeneous multiscale method for
  nonlinear monotone parabolic homogenization problems},
\newblock \bibinfo{journal}{ESAIM: Mathematical Modelling and Numerical
  Analysis} \bibinfo{volume}{50} (\bibinfo{year}{2016})
  \bibinfo{pages}{1659--1697}.
%Type = Article
\bibitem[{Sengupta et~al.(2012)Sengupta, Papadopoulos, and Taylor}]{R32}
\bibinfo{author}{A.~Sengupta}, \bibinfo{author}{P.~Papadopoulos},
  \bibinfo{author}{R.~L. Taylor},
\newblock \bibinfo{title}{A multiscale finite element method for modeling fully
  coupled thermomechanical problems in solids},
\newblock \bibinfo{journal}{International Journal for Numerical Methods in
  Engineering} \bibinfo{volume}{91} (\bibinfo{year}{2012})
  \bibinfo{pages}{1386--1405}.
%Type = Article
\bibitem[{Verfürth(2021)}]{R33}
\bibinfo{author}{B.~Verfürth},
\newblock \bibinfo{title}{Numerical homogenization for nonlinear strongly
  monotone problems},
\newblock \bibinfo{journal}{IMA Journal of Numerical Analysis}
  \bibinfo{volume}{42} (\bibinfo{year}{2021}) \bibinfo{pages}{1313--1338}.
%Type = Article
\bibitem[{Chen et~al.(2022)Chen, Li, Lu, and Wright}]{R34}
\bibinfo{author}{S.~Chen}, \bibinfo{author}{Q.~Li}, \bibinfo{author}{J.~Lu},
  \bibinfo{author}{S.~J. Wright},
\newblock \bibinfo{title}{Manifold learning and nonlinear homogenization},
\newblock \bibinfo{journal}{Multiscale Modeling \& Simulation}
  \bibinfo{volume}{20} (\bibinfo{year}{2022}) \bibinfo{pages}{1093--1126}.
%Type = Article
\bibitem[{Dong et~al.(2025)Dong, Guan, and Nie}]{R35}
\bibinfo{author}{H.~Dong}, \bibinfo{author}{X.~Guan}, \bibinfo{author}{Y.~Nie},
\newblock \bibinfo{title}{Multiscale method and convergence analysis for
  coupled nonlinear thermomechanical problems in heterogeneous shells},
\newblock \bibinfo{journal}{SIAM Journal on Scientific Computing}
  \bibinfo{volume}{47} (\bibinfo{year}{2025}) \bibinfo{pages}{B190--B219}.
%Type = Article
\bibitem[{Dong et~al.(2017)Dong, Cui, Nie, and Yang}]{R36}
\bibinfo{author}{H.~Dong}, \bibinfo{author}{J.~Cui}, \bibinfo{author}{Y.~Nie},
  \bibinfo{author}{Z.~Yang},
\newblock \bibinfo{title}{Second-order two-scale computational method for
  nonlinear dynamic thermo-mechanical problems of composites with cylindrical
  periodicity},
\newblock \bibinfo{journal}{Communications in Computational Physics}
  \bibinfo{volume}{21} (\bibinfo{year}{2017}) \bibinfo{pages}{1173--1206}.
%Type = Article
\bibitem[{Liu et~al.(2013)Liu, Liu, Guan, He, and Yuan}]{R37}
\bibinfo{author}{S.~Liu}, \bibinfo{author}{X.~Liu}, \bibinfo{author}{X.~Guan},
  \bibinfo{author}{P.~He}, \bibinfo{author}{Y.~Yuan},
\newblock \bibinfo{title}{A stochastic multi-scale model for predicting the
  thermal expansion coefficient of early-age concrete},
\newblock \bibinfo{journal}{CMES Comput. Model. Eng. Sci} \bibinfo{volume}{92}
  (\bibinfo{year}{2013}) \bibinfo{pages}{173--191}.
%Type = Article
\bibitem[{Rodr{\'\i}guez et~al.(2016)Rodr{\'\i}guez, Cruz, and
  Bravo-Castillero}]{R38}
\bibinfo{author}{E.~I. Rodr{\'\i}guez}, \bibinfo{author}{M.~E. Cruz},
  \bibinfo{author}{J.~Bravo-Castillero},
\newblock \bibinfo{title}{Reiterated homogenization applied to heat conduction
  in heterogeneous media with multiple spatial scales and perfect thermal
  contact between the phases},
\newblock \bibinfo{journal}{Journal of the Brazilian Society of Mechanical
  Sciences and Engineering} \bibinfo{volume}{38} (\bibinfo{year}{2016})
  \bibinfo{pages}{1333--1343}.
%Type = Article
\bibitem[{Iglesias-Rodr{\'\i}guez et~al.(2020)Iglesias-Rodr{\'\i}guez,
  Bravo-Castillero, Cruz, P{\'e}rez-Fern{\'a}ndez, and Sabina}]{R39}
\bibinfo{author}{E.~Iglesias-Rodr{\'\i}guez},
  \bibinfo{author}{J.~Bravo-Castillero}, \bibinfo{author}{M.~E. Cruz},
  \bibinfo{author}{L.~D. P{\'e}rez-Fern{\'a}ndez}, \bibinfo{author}{F.~J.
  Sabina},
\newblock \bibinfo{title}{Reiterated homogenization applied to nanofluids with
  an interfacial thermal resistance},
\newblock \bibinfo{journal}{International Journal for Multiscale Computational
  Engineering} \bibinfo{volume}{18} (\bibinfo{year}{2020}).
%Type = Article
\bibitem[{Nascimento et~al.(2017)Nascimento, Cruz, and Bravo-Castillero}]{R40}
\bibinfo{author}{E.~S. Nascimento}, \bibinfo{author}{M.~E. Cruz},
  \bibinfo{author}{J.~Bravo-Castillero},
\newblock \bibinfo{title}{Calculation of the effective thermal conductivity of
  multiscale ordered arrays based on reiterated homogenization theory and
  analytical formulae},
\newblock \bibinfo{journal}{International Journal of Engineering Science}
  \bibinfo{volume}{119} (\bibinfo{year}{2017}) \bibinfo{pages}{205--216}.
%Type = Article
\bibitem[{Yang et~al.(2017)Yang, Zhang, Dong, Cui, Guan, and Yang}]{R41}
\bibinfo{author}{Z.~Yang}, \bibinfo{author}{Y.~Zhang},
  \bibinfo{author}{H.~Dong}, \bibinfo{author}{J.~Cui},
  \bibinfo{author}{X.~Guan}, \bibinfo{author}{Z.~Yang},
\newblock \bibinfo{title}{High-order three-scale method for mechanical behavior
  analysis of composite structures with multiple periodic configurations},
\newblock \bibinfo{journal}{Composites Science and Technology}
  \bibinfo{volume}{152} (\bibinfo{year}{2017}) \bibinfo{pages}{198--210}.
%Type = Article
\bibitem[{Dong et~al.(2021)Dong, Cui, Nie, Jin, Guan, and Yang}]{R42}
\bibinfo{author}{H.~Dong}, \bibinfo{author}{J.~Cui}, \bibinfo{author}{Y.~Nie},
  \bibinfo{author}{K.~Jin}, \bibinfo{author}{X.~Guan},
  \bibinfo{author}{Z.~Yang},
\newblock \bibinfo{title}{High-order three-scale computational method for
  elastic behavior analysis and strength prediction of axisymmetric composite
  structures with multiple spatial scales},
\newblock \bibinfo{journal}{Mathematics and Mechanics of Solids}
  \bibinfo{volume}{26} (\bibinfo{year}{2021}) \bibinfo{pages}{905--936}.
%Type = Article
\bibitem[{Yang et~al.(2018)Yang, Sun, Cui, Yang, and Guan}]{R43}
\bibinfo{author}{Z.~Yang}, \bibinfo{author}{Y.~Sun}, \bibinfo{author}{J.~Cui},
  \bibinfo{author}{Z.~Yang}, \bibinfo{author}{T.~Guan},
\newblock \bibinfo{title}{A three-scale homogenization algorithm for coupled
  conduction-radiation problems in porous materials with multiple
  configurations},
\newblock \bibinfo{journal}{International Journal of Heat and Mass Transfer}
  \bibinfo{volume}{125} (\bibinfo{year}{2018}) \bibinfo{pages}{1196--1211}.
%Type = Article
\bibitem[{Dong et~al.(2019)Dong, Zheng, Cui, Nie, Yang, and Yang}]{R44}
\bibinfo{author}{H.~Dong}, \bibinfo{author}{X.~Zheng},
  \bibinfo{author}{J.~Cui}, \bibinfo{author}{Y.~Nie},
  \bibinfo{author}{Z.~Yang}, \bibinfo{author}{Z.~Yang},
\newblock \bibinfo{title}{High-order three-scale computational method for
  dynamic thermo-mechanical problems of composite structures with multiple
  spatial scales},
\newblock \bibinfo{journal}{International Journal of Solids and Structures}
  \bibinfo{volume}{169} (\bibinfo{year}{2019}) \bibinfo{pages}{95--121}.
%Type = Article
\bibitem[{Zuo et~al.(2024)Zuo, Yang, Cui, Deng, Li, and Guo}]{R45}
\bibinfo{author}{H.~Zuo}, \bibinfo{author}{Z.~Yang}, \bibinfo{author}{J.~Cui},
  \bibinfo{author}{S.~Deng}, \bibinfo{author}{H.~Li}, \bibinfo{author}{Z.~Guo},
\newblock \bibinfo{title}{High-order models for hydro-mechanical coupling
  problems in multiscale porous media},
\newblock \bibinfo{journal}{International Journal for Numerical Methods in
  Engineering} \bibinfo{volume}{125} (\bibinfo{year}{2024})
  \bibinfo{pages}{e7456}.
%Type = Article
\bibitem[{Ye et~al.(2024)Ye, Ma, Tang, Cui, and Li}]{R46}
\bibinfo{author}{S.~Ye}, \bibinfo{author}{Q.~Ma}, \bibinfo{author}{Q.~Tang},
  \bibinfo{author}{J.~Cui}, \bibinfo{author}{Z.~Li},
\newblock \bibinfo{title}{Second-order three-scale asymptotic analysis and
  algorithms for steklov eigenvalue problems in composite domain with
  hierarchical cavities},
\newblock \bibinfo{journal}{Journal of Scientific Computing}
  \bibinfo{volume}{98} (\bibinfo{year}{2024}) \bibinfo{pages}{61}.
%Type = Article
\bibitem[{Dong et~al.(2025)Dong, Wang, Ye, Nie, and Gao}]{R47}
\bibinfo{author}{H.~Dong}, \bibinfo{author}{Y.~Wang}, \bibinfo{author}{C.~Ye},
  \bibinfo{author}{Y.~Nie}, \bibinfo{author}{P.~Gao},
\newblock \bibinfo{title}{Higher-order three-scale asymptotic model and
  efficient two-stage numerical algorithm for transient nonlinear thermal
  conduction problems of composite structures},
\newblock \bibinfo{journal}{Computers \& Mathematics with Applications}
  \bibinfo{volume}{192} (\bibinfo{year}{2025}) \bibinfo{pages}{72--103}.
%Type = Article
\bibitem[{Cao(2006)}]{R48}
\bibinfo{author}{L.-Q. Cao},
\newblock \bibinfo{title}{Multiscale asymptotic expansion and finite element
  methods for the mixed boundary value problems of second order elliptic
  equation in perforated domains},
\newblock \bibinfo{journal}{Numerische Mathematik} \bibinfo{volume}{103}
  (\bibinfo{year}{2006}) \bibinfo{pages}{11--45}.
%Type = Article
\bibitem[{Dong et~al.(2018)Dong, Cao, Wang et~al.}]{R49}
\bibinfo{author}{Q.-l. Dong}, \bibinfo{author}{L.-q. Cao},
  \bibinfo{author}{X.~Wang}, et~al.,
\newblock \bibinfo{title}{Multiscale numerical algorithms for elastic wave
  equations with rapidly oscillating coefficients},
\newblock \bibinfo{journal}{Applied Mathematics and Computation}
  \bibinfo{volume}{336} (\bibinfo{year}{2018}) \bibinfo{pages}{16--35}.
%Type = Inproceedings
\bibitem[{Cui(2001)}]{R50}
\bibinfo{author}{J.~Cui},
\newblock \bibinfo{title}{{Multiscale computational method for unified design
  of structure, components and their materials}},
\newblock in: \bibinfo{booktitle}{Proceedings on Computational Mechanics in
  Science and Engineering, CCCM-2001, Guangzhou, 5-8 December},
  \bibinfo{publisher}{Peking University Press}, \bibinfo{year}{2001}, pp.
  \bibinfo{pages}{33--43}.

\end{thebibliography}

\end{document}